\documentclass[12pt,a4paper]{article}
\usepackage{amssymb,amsmath,latexsym,times,graphicx,standard}
\usepackage{pdfsync}

\textheight 
            675pt
\textwidth 
           460pt

\newcommand{\dk}{\, \mathrm{d}k}

\newcommand{\dx}{\, \mathrm{d}x}

\newcommand{\ds}{\, \mathrm{d}s}
\newcommand{\dy}{\, \mathrm{d}y}

\newcommand{\sm}{\!\setminus\!}

\newcommand{\R}{\mathbb R}
\newcommand{\N}{\mathbb N}

\DeclareMathOperator{\esssup}{ess\, sup}
\DeclareMathOperator{\sgn}{sgn}
\DeclareMathOperator{\sech}{sech}
\DeclareMathOperator{\supp}{supp}
\DeclareMathOperator{\dist}{dist}

\newcommand{\EE}{\mathcal E}
\newcommand{\FF}{\mathcal F}
\newcommand{\GG}{\mathcal G}
\newcommand{\HH}{\mathcal H}
\newcommand{\II}{\mathcal I}
\newcommand{\JJ}{\mathcal J}
\newcommand{\KK}{\mathcal K}
\newcommand{\LL}{\mathcal L}
\newcommand{\MM}{\mathcal M}

\newcommand{\PP}{\mathcal P}
\newcommand{\RR}{\mathcal R}
\renewcommand{\SS}{\mathcal S}
\newcommand{\TT}{\mathcal T}

\newcommand{\e}{\mathrm{e}}
\renewcommand{\i}{\mathrm{i}}

\newcommand{\qed}{\hfill$\Box$\bigskip}

\newcommand{\eqn}[1]{(\ref{#1})}

\newcommand{\nn}{|{\mskip-2mu}|{\mskip-2mu}|}

\newtheorem{theorem}{Theorem}[section]
\newtheorem{lemma}[theorem]{Lemma}
\newtheorem{definition}[theorem]{Definition}
\newtheorem{proposition}[theorem]{Proposition}
\newtheorem{corollary}[theorem]{Corollary}
\newtheorem{remark}[theorem]{Remark}

\newcommand\blfootnote[1]{%
  \begingroup
  \renewcommand\thefootnote{}\footnote{#1}%
  \addtocounter{footnote}{-1}%
  \endgroup
}

\begin{document}
\newcounter{count}
\newcounter{spice}

\title{Existence and conditional energetic stability of solitary gravity-capillary water waves with constant vorticity}

\author{M. D. Groves\thanks{FR 6.1 - Mathematik, Universit\"{a}t des Saarlandes, Postfach 151150, 66041 Saarbr\"{u}cken, Germany;
Department of Mathematical Sciences, Loughborough University, Loughborough, Leics, LE11 3TU, UK},
\quad E. Wahl\'{e}n\thanks{Centre for Mathematical Sciences, Lund University, 22100 Lund, Sweden}}
\date{}
\maketitle{}

\begin{abstract}
We present an existence and stability theory for
gravity-capillary solitary waves with constant vorticity on the surface of a body of water
of finite depth. Exploiting a rotational version of the classical variational principle,
we prove the existence of a minimiser of the wave energy $\HH$
subject to the constraint $\II=2\mu$, where $\II$ is the
wave momentum and $0< \mu \ll 1$. Since $\HH$ and
$\II$ are both conserved quantities a standard argument asserts
the stability of the set $D_\mu$ of minimisers: solutions starting
near $D_\mu$ remain close to $D_\mu$ in a suitably defined
energy space over their interval of existence.

In the applied mathematics literature solitary water waves of the
present kind are described by solutions of a
Korteweg-deVries equation (for strong surface tension) or
a nonlinear Schr\"{o}dinger equation (for weak surface tension).
We show that the waves detected by
our variational method converge (after an appropriate rescaling)
to solutions of the appropriate model equation as $\mu \downarrow 0$.
\blfootnote{To appear in Proceedings of the Royal Society
of Edinburgh: Section A. 
\copyright The Royal Society of Edinburgh.}
\end{abstract}

%
%
\newpage

\section{Introduction}

\subsection{Variational formulation of the hydrodynamic problem} \label{Variational formulation}

\subsubsection{The water-wave problem}

In this paper we consider a two-dimensional perfect fluid bounded below by 
a flat rigid bottom $\{y=0\}$ and above by a  free surface $\{y=d+\eta(x,t)\}$.
The fluid has unit density and flows under
the influence of gravity and surface tension with constant vorticity
$\omega$, so that the velocity field $(u(x,y,t),v(x,y,t))$ in the fluid domain
$\Sigma_\eta=\{0 < y < d + \eta(x,t)\}$ satisfies $v_x-u_y = \omega$. We study
waves which are perturbations of underlying shear flows given by $\eta=0$ and $(u,v)=(\omega (d-y), 0)$
(which may be a good description of tidal currents (see Constantin \cite[Chapter 2.3.2]{Constantin})) 
and are evanescent as $x \rightarrow \pm \infty$.
In terms of a generalised velocity potential $\phi$ such that
$(u,v)=(\phi_x + \omega (d-y), \phi_y)$ and stream function $\psi$ such that $(u,v)=(\psi_y,-\psi_x)$,
the governing equations are
\begin{eqnarray*}
\Delta \phi & = & \parbox{90mm}{$0$,} 0 < y < d+\eta, \\
\phi_{y} & = & \parbox{90mm}{$0$,} y=0, \\
\eta_t & = & \parbox{90mm}{$\phi_{y} - \eta_x\phi_x+\omega \eta \eta_x$,}  y=d+\eta, \\
\phi_t & = & \parbox{90mm}{$\displaystyle-\frac{1}{2}|\nabla \psi|^2 -\omega \psi-g\eta+ \beta\left[\frac{\eta_x}{\sqrt{1+\eta_x^2}}\right]_x$,} y=d+\eta,
\end{eqnarray*}
with $\eta(x,t)$, $\phi(x,y,t), \psi(x,y,t) +\frac{1}{2}\omega (d-y)^2 \rightarrow 0$ as $x \rightarrow \pm \infty$,
where $g$ and $\beta$ are respectively the acceleration due to
gravity and the (positive) coefficient of surface tension (see Constantin, Ivanov \& Prodanov 
\cite{ConstantinIvanovProdanov08}). 

At this point it is convenient to introduce dimensionless variables
$$(x^\prime,y^\prime)=\frac{1}{d}(x,y),\qquad t^\prime=\left(\frac{g}{d}\right)^{\!\frac{1}{2}}t,$$
$$\eta^\prime(x^\prime,t^\prime) = \frac{1}{d}\eta(x,t),\quad
\phi^\prime(x^\prime,t^\prime) = \frac{1}{(gd^3)^\frac{1}{2}}\phi(x,t),\quad
\psi^\prime(x^\prime,t^\prime) = \frac{1}{(gd^3)^\frac{1}{2}}\psi(x,t)
$$
and parameters $\omega^\prime=\omega(d/g)^\frac{1}{2}$, $\beta^\prime=\beta/gd^2$; one obtains the equations
\begin{eqnarray}
\Delta \phi & = & \parbox{90mm}{$0$,} 0 < y < 1+\eta, \label{Laplace} \\
\phi_{y} & = & \parbox{90mm}{$0$,} y=0, \label{Impermeable} \\
\eta_t & = & \parbox{90mm}{$\phi_{y} - \eta_x\phi_x+\omega \eta \eta_x$,}  y=1+\eta, \label{Kinematic BC}\\
\phi_t & = & \parbox{90mm}{$\displaystyle-\frac{1}{2}|\nabla \psi|^2 -\omega \psi-\eta+ \beta\left[\frac{\eta_x}{\sqrt{1+\eta_x^2}}\right]_x$,} y=1+\eta, \label{Dynamic BC}
\end{eqnarray}
in which the primes have been dropped for notational simplicity. In particular we seek \emph{solitary-wave solutions}
of \eqn{Laplace}--\eqn{Dynamic BC}, that is
waves of permanent form which propagate from right to left with constant (dimensionless) speed $\nu$, so that
$\eta(x,t)=\eta(x+\nu t)$ (and of course $\eta(x+\nu t) \rightarrow 0$ as
$x+\nu t \rightarrow \pm \infty$).

\subsubsection{Formulation as a Hamiltonian system}

We proceed by reducing the hydrodynamic problem to a pair of nonlocal, coupled evolutionary equations for
the variables $\eta$ and $\xi=\phi|_{y=1+\eta}$. For fixed $\eta$ and $\xi$, let $\phi$ denote the unique solution to
the boundary-value problem
\begin{eqnarray*}
& & \parbox{5cm}{$\Delta \phi=0,$}0<y<1+\eta, \\
& & \parbox{5cm}{$ \phi = \xi,$}y=1+\eta, \\
& & \parbox{5cm}{$\phi_y =0,$}y=0
\end{eqnarray*}
and denote the harmonic conjugate of $\phi$ by $\tilde{\psi}$. We define the \emph{Hilbert transform} $H(\eta)$
and \emph{Dirichlet-Neumann operator} $G(\eta)$ for this boundary-value problem by
$$H(\eta)\xi = \tilde{\psi}|_{y=1+\eta}, \qquad G(\eta)\xi=(\phi_y-\eta_x \phi_x)|_{y=1+\eta},$$
so that $G(\eta)=-\partial_x H(\eta)$, and note that the boundary conditions \eqn{Kinematic BC},
\eqn{Dynamic BC} can be written as
\begin{eqnarray*}
\eta_t & = & \parbox{90mm}{$G(\eta)\xi+\omega \eta \eta_x$,} \\
\xi_t & = & 
-\frac{1}{2(1+\eta_x^2)}\big(\xi_x^2 - (G(\eta)\xi)^2 - 2 \eta_x\xi_x G(\eta)\xi\big) \\
& &
\parbox{90mm}{$\displaystyle\qquad\mbox{} + \omega \eta \xi_x - \omega H(\eta) \xi 
-\eta+ \beta\left[\frac{\eta_x}{\sqrt{1+\eta_x^2}}\right]_x$.}
\end{eqnarray*}

Wahl\'{e}n \cite{Wahlen07} observed that the above equations can be formulated as the Hamiltonian system
\begin{equation}
\left(\begin{array}{c} \eta_t\\ \xi_t \end{array}\right)=\left(\begin{array}{cc}0 & 1\\ -1 & \omega \partial_x^{-1} \end{array}\right)
\!\!\left(\begin{array}{c}\delta_\eta \mathcal H\\ 
\delta_\xi \mathcal H\end{array}\right), \label{Hamiltonian system}
\end{equation}
in which
\begin{equation}
\HH(\eta, \xi) =\int_{-\infty}^\infty \left(\frac{\xi G(\eta)\xi}{2}+\omega
  \xi \eta \eta_x + \frac{\omega^2}{6}\eta^3+\frac{\eta^2}{2}+\beta(\sqrt{1+\eta_x^2}-1)\right )\dx,
\label{Definition of HH}
\end{equation}
(note that the well-known formulation of the water-wave problem by Zakharov \cite{Zakharov68} is
recovered in the irrotational case $\omega=0$). This Hamiltonian system
has the conserved quantities $\HH(\eta,\xi)$ (total energy) and
\begin{equation}
\II(\eta,\xi)  =\int_{-\infty}^\infty \left(\xi \eta_x+\frac{\omega}{2}\eta^2\right)\dx,
\label{Definition of II}
\end{equation}
(total horizontal momentum), which satisfies the equation
\begin{equation}
\left(\begin{array}{c} \eta_x \\ \xi_x \end{array}\right)=\left(\begin{array}{cc}0 & 1\\ -1 & \omega \partial_x^{-1} \end{array}\right)
\!\!\left(\begin{array}{c}\delta_\eta \mathcal I\\ 
\delta_\xi \mathcal I\end{array}\right); \label{Noether for I}
\end{equation}
these quantities are associated with its independence of
respectively $t$ and $x$. According to \eqn{Hamiltonian system} and
 \eqn{Noether for I}, a solution of the form $\eta(x,t)=\eta(x+\nu t)$, $\xi(x,t)=\xi(x+\nu t)$
is characterised as a critical point of the total energy
subject to the constraint of fixed momentum (cf.\ Benjamin \cite{Benjamin84}).
It is therefore a critical point of the functional $\HH- \nu \II$, where the
the speed of the wave is given by the Lagrange multiplier $\nu$. This
functional depends on the single independent variable $x+\nu t$, which we
now abbreviate to $x$.

A similar variational principle for waves of permanent form with a general distribution of
vorticity has been used by Groves \& Wahl\'{e}n \cite{GrovesWahlen07} in an existence theory
for solitary waves. Groves \& Wahl\'{e}n
interpreted their variational functional as an action functional and derived a formulation
of the hydrodynamic problem as an infinite-dimensional spatial Hamiltonian system; a rich solution
set is found using a centre-manifold reduction technique to convert it into a
Hamiltonian system with a finite number of degrees of freedom. 

 In this paper we present a direct existence theory
for minimisers of $\HH$ subject to the constraint
$\II=2\mu$ for $0 < \mu < \mu_0$, where $\mu_0$ is a fixed positive constant chosen
small enough for for the validity of our calculations. We seek constrained minimisers in a two-step
approach.

\begin{enumerate}
\item \emph{Fix $\eta \neq 0$ and minimise $\HH(\eta,\cdot)$ over
$T_\mu=\{\xi: \II(\eta,\xi)=2\mu\}$.}
This problem (of minimising a quadratic functional over a linear manifold)
admits a unique global minimiser $\xi_\eta$.
\item \emph{Minimise $\JJ_\mu(\eta):=\HH(\eta,\xi_\eta)$ over $\eta \in U\sm\{0\}$,} where $U$ is a fixed
ball centred upon the origin in a suitable function space. 
Because $\xi_\eta$ minimises
$\HH(\eta,\cdot)$ over $T_\mu$ there exists a Lagrange multiplier $\nu_\eta$ such that
$$G(\eta)\xi_\eta+\omega \eta \eta^\prime=\nu_\eta \eta^\prime,$$
and straightforward calculations show that
$$
\xi_\eta=\nu_\eta G(\eta)^{-1}\eta^\prime-\frac{\omega}{2} G(\eta)^{-1}(\eta^2)_x,
$$
\[
\nu_\eta= \left(\frac{1}{2}\int_{-\infty}^\infty \eta^\prime
G(\eta)^{-1}\eta^\prime \dx\right)^{\!\!-1}\!\!\left(\mu -\frac{\omega}{4}\int_{-\infty}^\infty \eta^2 \dx+\frac{\omega}{4}\int_{-\infty}^\infty
(\eta^2)_x G(\eta)^{-1}\eta^\prime \dx\right),
\]
so that
\begin{equation}
\JJ_\mu(\eta) = \KK(\eta)+\frac{(\mu+\GG(\eta))^2}{\LL(\eta)}, \label{Definition of J}
\end{equation}
where
\begin{eqnarray}
\GG(\eta) & = & \frac{\omega}{4}\int_{-\infty}^\infty \eta^2 K(\eta)\eta \dx -\frac{\omega}{4}\int_{-\infty}^\infty \eta^2\dx,
\label{Definition of GG} \\
\KK(\eta) & = & \int_{-\infty}^\infty \left( \frac{1}{2}\eta^2 + \beta[\sqrt{1+\eta^{\prime 2}}-1] \right)\dx \\
& & \qquad\quad\mbox{}
-\frac{\omega^2}{2} \int_{-\infty}^\infty \frac{\eta^2}{2} K(\eta)\frac{\eta^2}{2} \dx
+ \frac{\omega^2}{6}\int_{-\infty}^\infty \eta^3 \dx,\quad \label{Definition of KK} \\
\LL(\eta) & = & \frac{1}{2} \int_{-\infty}^\infty \eta K(\eta)\eta \dx \label{Definition of LL}
\end{eqnarray}
and $K(\eta)=-\partial_x G(\eta)^{-1} \partial_x$.
This computation also shows that the dimensionless speed of a solitary wave corresponding to
a constrained minimiser $\eta$ of $\HH$ is
$$\nu=\frac{\mu+\GG(\eta)}{\LL(\eta)}.$$
\end{enumerate}

This two-step approach to the constrained minimisation problem
was introduced  in a corresponding theory for irrotational solitary waves by Buffoni \cite{Buffoni04a},
who used a conformal mapping due to Babenko \cite{Babenko87a,Babenko87b} to transform
$\JJ_\mu$ into another functional $\tilde{\JJ}_\mu$ depending only upon $H(0)$ and hence simplify
the necessary variational analysis.
Buffoni established the existence of a (non-zero) minimiser of $\tilde{\JJ}_\mu$ for strong surface tension
(Buffoni \cite{Buffoni04a}) and obtained partial results in this direction for weak surface tension
(Buffoni \cite{Buffoni05,Buffoni09}). A method for completing his results for weak surface tension
was sketched in a short note by Groves \& Wahl\'{e}n \cite{GrovesWahlen10}; in the present paper
we give complete details, including non-zero vorticity in our treatment and working directly with the
original physical variables.  Although versions of the Babenko transformation for non-zero constant
vorticity have been published (Constantin \& Varvaruca \cite{ConstantinVarvaruca11}, Martin \cite{Martin13}),
finding minimisers of $\JJ_\mu$ over $U \sm \{0\}$ has the advantage of immediately yielding 
precise information on solutions to the original water-wave equations
\eqn{Laplace}--\eqn{Dynamic BC}.

\subsubsection{Functional-analytic framework} \label{FAF}
An appropriate functional-analytic framework for the above variational problem is  introduced in
Section \ref{FA setting}. We work with the function spaces
$$H^r(\R)=\overline{(\SS(\R), \|\cdot\|_r)}, \qquad
\|\eta\|_r^2:=\int_{-\infty}^\infty (1+k^2)^r |\hat{\eta}|^2 \dk$$
for $r \in \R$ (the standard Sobolev spaces), and
$$H_\star^{1/2}(\R)=\overline{(\SS(\R), \|\cdot\|_{H_\star^{1/2}(\R)})}, \qquad
\|\eta\|_{H_\star^{1/2}(\R)}^2:= \int_{-\infty}^\infty (1+k^2)^{-\frac{1}{2}}k^2|\hat{\eta}|^2\dk,$$
$$H_\star^{-1/2}(\R)=\overline{(\overline\SS(\R), \|\cdot\|_{H_\star^{-1/2}(\R)})}, \qquad
\|\eta\|_{H_\star^{-1/2}(\R)}^2:= \int_{-\infty}^\infty (1+k^2)^{\frac{1}{2}}k^{-2}|\hat{\eta}|^2\dk;$$
here $\overline{(\SS(\R),\|\cdot\|)}$ denotes the completion of the inner product space
constructed by equipping the Schwartz class $\SS(\R)$ (or the subclass $\bar{\SS}(\R)$ of
Schwartz-class functions with zero mean) with the norm $\|\cdot\|$ and $\hat{\eta}=\FF[\eta]$
is the Fourier transform of $\eta$. 

The mathematical analysis of $G(\eta)$ and $K(\eta)$ is complicated by the fact that they are
defined in terms of boundary-value problems in the variable domain $\Sigma_\eta$. Lannes
\cite[Chapters 2 and 3]{Lannes} has presented a comprehensive theory for handling such 
such boundary-value problems by transforming them into 
serviceable nonlinear elliptic problems in the fixed domain $\Sigma_0$, and here we adapt
Lannes's methods to our specific requirements. Our main results are 
stated in the following theorem, according to which
equations \eqn{Definition of GG}--\eqn{Definition of LL} define analytic functionals
$\GG,\KK,\LL: W^{s+3/2} \rightarrow \R$ for $s>0$.
In accordance with this theorem we take $U=B_M(0) \subseteq H^2(\R)$,
where $M>0$ is chosen small enough
so that $\overline{B}_M(0) \subseteq H^2(\R)$ lies in $W^{s+3/2}$
and for for the validity of our calculations.

\begin{theorem} \label{Main analyticity result}
Choose $h_0 \in (0,1)$ and define $W=\{\eta \in W^{1,\infty}(\R) \colon 1+\inf \eta>h_0\}$
and $W^r = H^{r} \cap W$ for $r \geq 0$.
\begin{list}{(\roman{count})}{\usecounter{count}}
\item
The Dirichlet-Neumann operator $G(\eta)$ is an isomorphism
$H_\star^{1/2}(\R) \to H_\star^{-1/2}(\R)$ for each $\eta \in W$.
\item
The Dirichlet-Neumann operator $G(\cdot): W \rightarrow
\LL(H_\star^{1/2}(\R), H_\star^{-1/2}(\R))$
and Neumann-Dirichlet operator
$G(\cdot)^{-1}: W \rightarrow
\LL(H_\star^{-1/2}(\R), H_\star^{1/2}(\R))$ are analytic.
\item
The operator $K(\cdot): W^{s+3/2} \rightarrow
\LL(H^{s+3/2}(\R), H^{s+1/2}(\R))$ is analytic for each $s>0$.
\end{list}
\end{theorem}

\subsection{Heuristics}

\begin{figure}[h]
\centering
\includegraphics[width=7cm]{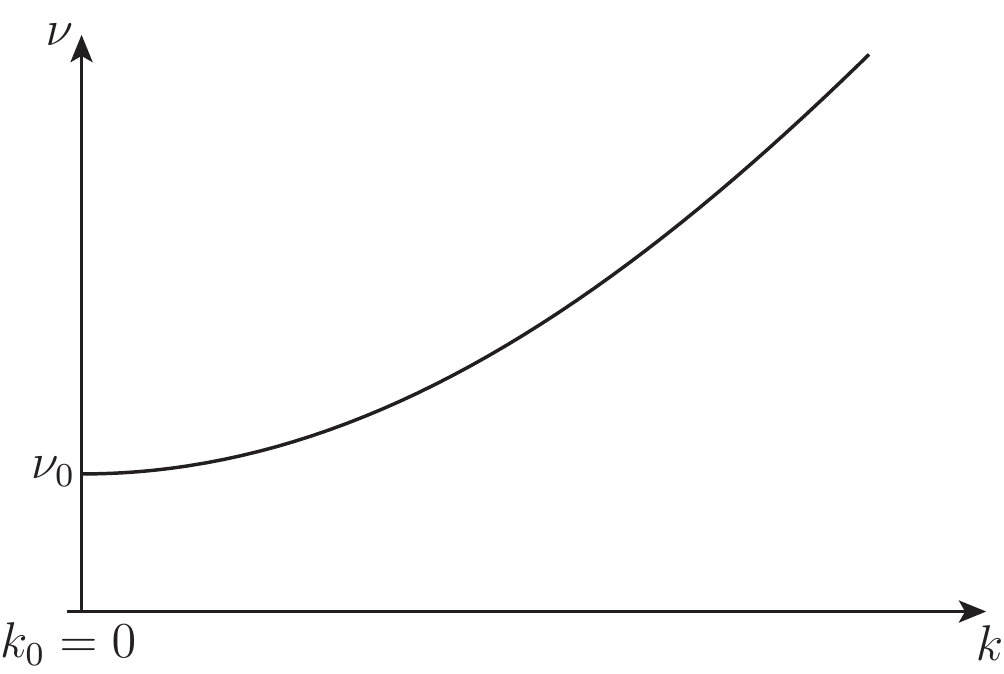}
\hspace{1cm}
\includegraphics[width=7cm]{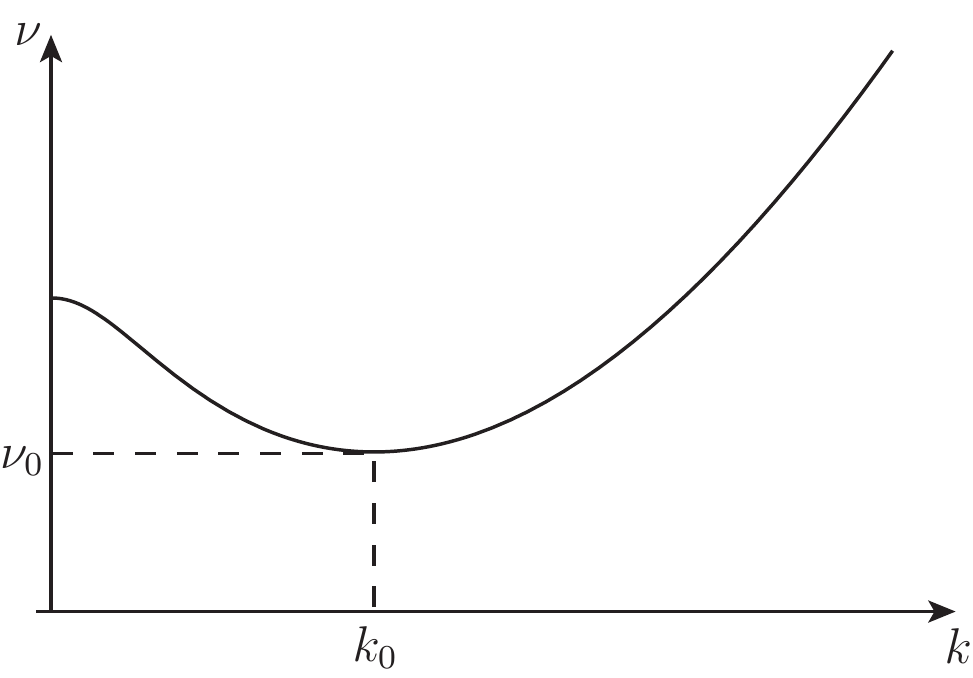}

$\beta>\beta_\mathrm{c}\hspace{6.5cm}\beta<\beta_\mathrm{c}$

{\it
\caption{Dispersion relation for linear water waves}
\label{Dispersion curves}}
\end{figure}

The existence of small-amplitude solitary waves is predicted by studying the dispersion relation for the
linearised version of \eqn{Laplace}--\eqn{Dynamic BC}. Linear waves of the form $\eta(x,t) = \cos k (x +\nu t)$ exist whenever
$$1+\beta k^2 - \omega \nu - \nu^2 f(k) =0, \qquad f(k)=|k|\coth |k|,$$
that is whenever
$$\nu = -\frac{\omega}{2f(k)} + \frac{1}{2}\left( \frac{\omega^2}{f(k)^2} + \frac{4(1+\beta k^2)}{f(k)}\right)^{\!\!\frac{1}{2}}.$$
The function $k \mapsto \nu(k)$, $k \geq 0$ has a unique global minimum $\nu_0=\nu(k_0)$, and one finds that
$k_0>0$ for $\beta<\beta_\mathrm{c}$ and $k_0=0$ (with $\nu_0 =\nu(0) = \frac{1}{2}(-\omega+\sqrt{\omega^2+4})$) for $\beta>\beta_\mathrm{c}$, where
$$\beta_\mathrm{c} = {\textstyle\frac{1}{6}}(\omega^2 +2 - \omega\sqrt{\omega^2+4})$$
(see Figure~\ref{Dispersion curves}).
For later use let us also note that
$$g(k):=1+\beta k^2 - \omega \nu_0 - \nu_0^2f(k) \geq 0, \qquad k \in \R,$$
with equality precisely when $k=\pm k_0$.

Bifurcations of nonlinear solitary waves are are expected whenever the
linear group and phase speeds are equal, so that $\nu^\prime(k)=0$ (see Dias \& Kharif \cite[\S 3]{DiasKharif99}).
We therefore expect the existence of small-amplitude solitary waves with speed near $\nu_0$; the waves bifurcate from
laminar flow when $\beta>\beta_\mathrm{c}$ and from a linear periodic wave train with frequency $k_0 \nu(k_0)$ when
$\beta<\beta_\mathrm{c}$. Model equations for both types of solution have been derived by Johnson \cite[\S\S 4--5]{Johnson12}.\\

\begin{figure}[h]
\begin{center}
\includegraphics{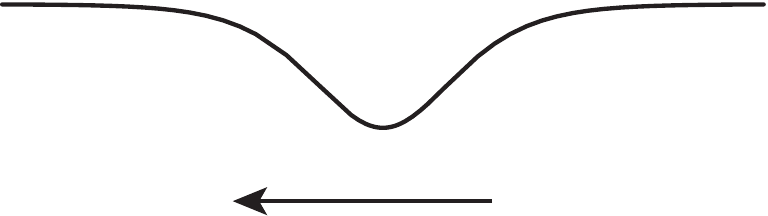}
\end{center}
{\it
\caption{Korteweg-deVries theory predicts the existence of small-amplitude solitary waves of depression for strong surface tension} \label{Sketch of KdV wave}}
\end{figure}

\noindent
\underline{$\beta>\beta_\mathrm{c}$}:
The appropriate model equation is the Korteweg-deVries equation
\begin{equation}
-2u_T-\left(\beta-\frac{\nu_0^2}{3}\right)u_{XXX}+(\omega^2+3)uu_X=0, \label{T-KdV}
\end{equation}
in which
$$\eta=\mu^\frac{2}{3}u(X,T)+O(\mu^\frac{4}{3}), \qquad X=\mu^\frac{1}{3}(x+\nu_0 t),\quad
 T=2(\omega^2+4)^{-\frac{1}{2}}\mu^\frac{2}{3}t.$$
At this level of approximation a solution
to \eqn{T-KdV} of the form $u(X,T)=\phi(X+\nu_\mathrm{KdV}T)$ with $\phi(X) \rightarrow 0$ as $X \rightarrow \pm \infty$
corresponds to a solitary water wave with speed
$$\nu = \nu_0 + 2(\omega^2+4)^{-\frac{1}{2}}\mu^\frac{2}{3}\nu_\mathrm{KdV}
=-{\textstyle\frac{1}{2}}\omega + {\textstyle\frac{1}{2}}(\omega^2+4)^{1/2} + 2(\omega^2+4)^{-\frac{1}{2}}\mu^\frac{2}{3}\nu_\mathrm{KdV}.$$
The following lemma gives a variational description of the set of such solutions; the corresponding solitary waves are sketched in Figure \ref{Sketch of KdV wave}.

\begin{lemma} \label{Variational KdV}
\quad
\begin{list}{(\roman{count})}{\usecounter{count}}
\item
The set of solutions to the ordinary differential equation
$$
-\left(\beta-\frac{\nu_0^2}{3}\right)\phi^{\prime\prime}-2\nu_\mathrm{KdV}\phi+\frac{3}{2}\left(
\frac{\omega^2}{3}+1\right)\phi^2=0
$$
satisfying $\phi(X) \rightarrow 0$ as $X \rightarrow \infty$ is
$D_\mathrm{KdV} = \{\phi_\mathrm{KdV}(\cdot + y)\colon y \in \R\},$
where
\begin{eqnarray*}
\nu_\mathrm{KdV} & =& -\frac{\displaystyle 2\left(\frac{3}{16}\right)^{\!\!\frac{2}{3}}\left(\frac{\omega^2}{3}+1\right)^{\!\!\frac{4}{3}}}
{\displaystyle\left(\beta-\frac{\nu_0^2}{3}\right)^{\!\!\frac{1}{3}}(\omega^2+4)^\frac{1}{3}},\\
\phi_\mathrm{KdV}(x) & = & -\frac{\displaystyle\sqrt{3}\left(\frac{3}{16}\right)^{\!\!\frac{1}{6}}\left(\frac{\omega^2}{3}+1\right)^{\!\!\frac{1}{3}}}
{\displaystyle\left(\beta-\frac{\nu_0^2}{3}\right)^{\!\!\frac{1}{3}}(\omega^2+4)^\frac{1}{3}}
\sech^2 \left(\frac{\displaystyle\left(\frac{3}{16}\right)^{\!\!\frac{1}{3}}\left(\frac{\omega^2}{3}+1\right)^{\!\!\frac{2}{3}}x}
{\displaystyle\left(\beta-\frac{\nu_0^2}{3}\right)^{\!\!\frac{2}{3}}(\omega^2+4)^\frac{1}{6}}\right).
\end{eqnarray*}
These functions are precisely the minimisers of the functional
$\EE_\mathrm{KdV}:H^1(\R) \rightarrow \R$ given by
$$
\EE_\mathrm{KdV}(\phi)=\frac{1}{2}\int_{-\infty}^\infty \left(\left(\beta-\frac{\nu_0^2}{3}\right)(\phi^\prime)^2+\left(\frac{\omega^2}{3}+1\right)\phi^3
\right) \dx
$$
over the set $N_\mathrm{KdV} = \{\phi \in H^1(\R): \|\phi\|_0^2=2\alpha_\mathrm{KdV}\}$;
the constant $2\nu_\mathrm{KdV}$ is the Lagrange multiplier in this constrained variational
principle and
$$
c_\mathrm{KdV}:=\inf\left\{\EE_\mathrm{KdV}(\phi)\colon \phi\in N_\mathrm{KdV}\right\}
=-\frac{\displaystyle\frac{9}{5}\left(\frac{2}{3}\right)^{\!\!\frac{1}{3}}\left(\frac{\omega^2}{3}+1\right)^{\!\!\frac{4}{3}}}
{\displaystyle\left(\beta-\frac{\nu_0^2}{3}\right)^{\!\!\frac{1}{3}}(\omega^2+4)^\frac{5}{6}}.
$$
Here the numerical value $\alpha_\mathrm{KdV} = 2(\omega^2+4)^{-\frac{1}{2}}$ is chosen for
compatibility with an estimate (Proposition \ref{Strong ST convergence in L2})
in the following water-wave theory.
\item
Suppose that $\{\phi_m\}\subset N_\mathrm{KdV}$ is a minimising sequence for $\EE_\mathrm{KdV}$.
There exists a sequence $\{x_m\}$ of real numbers with the property that a subsequence of $\{\phi_m(\cdot +x_m)\}$
converges in $H^1(\R)$ to an element of $D_\mathrm{KdV}$.
\end{list}
\end{lemma}

\begin{figure}[h]
\begin{center}
\includegraphics{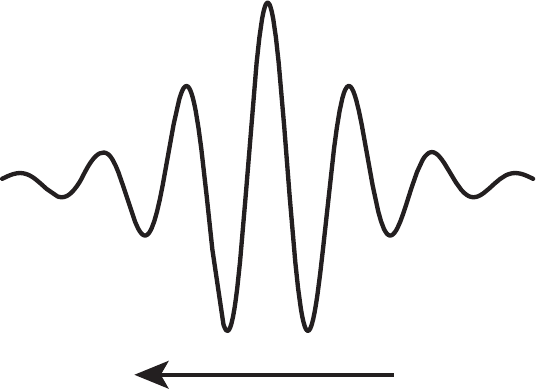}
\end{center}
{\it
\caption{Nonlinear Schr\"{o}dinger theory predicts the existence of small-amplitude envelope solitary waves for weak surface tension} \label{Sketch of NLS wave}}
\end{figure}

\noindent
\underline{$\beta<\beta_\mathrm{c}$}:
The appropriate model equation is the cubic nonlinear Schr\"{o}dinger equation
\begin{equation}
2\mathrm{i}A_T -\frac{1}{4}g^{\prime\prime}(k_0) A_{XX} + \frac{3}{2}\left(\frac{A_3}{2}+A_4\right)|A|^2 A =0, \label{T-NLS}
\end{equation}
in which
$$\eta=\frac{1}{2}\mu (A(X,T)\mathrm{e}^{\mathrm{i}k_0(x+\nu_0 t)}+\mathrm{c.c.})+O(\mu^2),$$
$$X=\mu(x+\nu_0 t),\qquad T=4k_0(\omega+2\nu_0 f(k_0))^{-1} \mu^2 t$$
and $A_3$, $A_4$ are functions of $\beta$ and $\omega$ which are given in Corollary \ref{lot in fours}
and Proposition \ref{lot in threes} below; the abbreviation `$\mathrm{c.c.}$' denotes the complex conjugate of the
preceding quantity. (It is demonstrated in Appendix B that $A_3+2A_4$ is negative.) At this level of approximation a
solution to \eqn{T-NLS} of the form
$A(X,T)=\mathrm{e}^{\mathrm{i}\nu_\mathrm{NLS} T}\phi(X)$
with $\phi(X) \rightarrow 0$ as $X \rightarrow \pm\infty$ corresponds
to a solitary water wave with speed
$$\nu=\nu_0 + 4(\omega+2\nu_0 f(k_0))^{-1} \mu^2\nu_\mathrm{NLS}.$$
The following lemma gives a variational description of the set of such solutions (see Cazenave \cite[\S 8]{Cazenave});
the corresponding solitary waves are sketched in Figure \ref{Sketch of NLS wave}.

\begin{lemma} \label{Variational NLS}
\quad
\begin{list}{(\roman{count})}{\usecounter{count}}
\item
The set of complex-valued solutions to the ordinary differential equation
$$
-\frac{1}{4}g^{\prime\prime}(k_0)\phi^{\prime\prime}-2\nu_\mathrm{NLS}\phi+\frac{3}{2}\left(
\frac{A_3}{2}+A_4\right)|\phi|^2\phi=0
$$
satisfying $\phi(X) \rightarrow 0$ as $X \rightarrow \infty$ is
$D_\mathrm{NLS} = \{\mathrm{e}^{\mathrm{i}\omega}\phi_\mathrm{NLS}(\cdot + y)\colon \omega \in [0,2\pi), y \in \R\},$
where
\begin{eqnarray*}
\nu_\mathrm{NLS} & =& -\frac{9\alpha_\mathrm{NLS}^2}{8g^{\prime\prime}(k_0)}\left(\frac{A_3}{2}+A_4\right)^{\!\!2}, \\
\phi_\mathrm{NLS}(x) & = & \alpha_\mathrm{NLS}\left(-\frac{3}{g^{\prime\prime}(k_0)}\left(\frac{A_3}{2}+A_4\right)\right)^{\!\!\frac{1}{2}}
\sech\left(-\frac{3\alpha_\mathrm{NLS}}{g^{\prime\prime}(k_0)}\left(\frac{A_3}{2}+A_4\right)x\right)
\end{eqnarray*}
These functions are precisely the minimisers of the functional
$\EE_\mathrm{NLS}:H^1(\R) \rightarrow \R$ given by
$$
\EE_\mathrm{NLS}(\phi)=\int_{-\infty}^\infty \left(\frac{1}{8}g^{\prime\prime}(k_0)|\phi^\prime|^2
+\frac{3}{8}\left(\frac{A_3}{2}+A_4\right)|\phi|^4\right) \dx
$$
over the set $N_\mathrm{NLS} = \{\phi \in H^1(\R): \|\phi\|_0^2=2\alpha_\mathrm{NLS}\}$;
the constant $2\nu_\mathrm{NLS}$ is the Lagrange multiplier in this constrained variational
principle and
$$
c_\mathrm{NLS}:=\inf\left\{\EE_\mathrm{NLS}(\phi)\colon \phi\in N_\mathrm{NLS}\right\}
=-\frac{3\alpha_\mathrm{NLS}^3}{4g^{\prime\prime}(k_0)}\left(\frac{A_3}{2}+A_4\right)^{\!\!2}.
$$
Here the numerical value $\alpha_\mathrm{NLS}=\frac{1}{2}\left(\frac{1}{4}\nu_0f(k_0)+\frac{\omega}{8}\right)^{-1}$
is chosen for
compatibility with an estimate (Proposition \ref{Weak ST convergence in L2})
in the following water-wave theory.
\item
Suppose that $\{\phi_n\}\subset N_\mathrm{NLS}$ is a minimising sequence for $\EE_\mathrm{NLS}$.
There exists a sequence $\{x_m\}$ of real numbers with the property that a subsequence of $\{\phi_m(\cdot +x_m)\}$
converges in $H^1(\R)$ to an element of $D_\mathrm{NLS}$.
\end{list}
\end{lemma}

\subsection{The main results}

In this paper we establish the existence of minimisers of the functional $\JJ_\mu$ over $U \sm \{0\}$
and confirm that the corresponding solitary water waves
are approximated by suitable scalings of the functions $\phi_\mathrm{KdV}$ (for $\beta>\beta_\mathrm{c}$)
and $\phi_\mathrm{NLS}$ (for $\beta<\beta_\mathrm{c}$). The following theorem states these results more precisely.

\begin{theorem} \label{Main results}
\quad
\begin{list}{(\roman{count})}{\usecounter{count}}
\item
The set $B_\mu$ of minimisers of $\JJ_\mu$ over $U \sm \{0\}$ is non-empty.
\item
Suppose that $\{\eta_m\}$ is a minimising sequence for $\JJ_\mu$ on $U\sm\{0\}$ which satisfies
$$
\sup_{m\in{\mathbb N}} \|\eta_m\|_2 < M.
$$
There exists a sequence $\{x_m\} \subset \R$ with the property that
a subsequence of $\{\eta_m(x_m+\cdot)\}$ converges
in $H^r(\R)$, $r \in [0,2)$, to a function $\eta \in B_\mu$.
\item
Suppose that $\beta>\beta_\mathrm{c}$.
The set $B_\mu$ of minimisers of $\JJ_\mu$ over $U\sm\{0\}$ satisfies
$$
\sup_{\eta \in B_\mu} \inf_{x \in \R} \|\phi_\eta - \phi_\mathrm{KdV}(\cdot+x)\|_1 \rightarrow 0
$$
as $\mu \downarrow 0$, where we write
$$\eta_1(x)=\mu^\frac{2}{3}\phi_\eta(\mu^\frac{1}{3} x)$$
and $\eta_1$ is obtained from $\eta$ by multiplying its Fourier transform by the characteristic function of the interval
$[-\delta_0,\delta_0]$ with $\delta_0>0$. Furthermore, the speed $\nu_\mu$ of the corresponding solitary water waves
satisfies
$$\nu_\mu = \nu_0 + 2(\omega^2+4)^{-\frac{1}{2}}\nu_\mathrm{KdV}\mu^\frac{2}{3} + o(\mu^\frac{2}{3})$$
uniformly over $\eta \in B_\mu$.
\item
Suppose that $\beta<\beta_\mathrm{c}$. The set $B_\mu$ of minimisers of $\JJ_\mu$ over $U\sm\{0\}$ satisfies
$$\sup_{\eta \in B_\mu} \inf_{\substack{\omega \in [0,2\pi],\\ x \in \R}} \|\phi_\eta-e^{\mathrm{i}\omega}\phi_\mathrm{NLS}(\cdot+x)\|_1 \rightarrow 0$$
as $\mu \downarrow 0$, where we write
$$\eta_1^+(x) = \frac{1}{2}\mu \phi_\eta(\mu x)\e^{\i k_0 x},$$
and $\eta_1^+$ is obtained from $\eta$ by multiplying its Fourier transform by the characteristic function of the interval
$[k_0-\delta_0,k_0+\delta_0]$ with $\delta_0 \in (0,k_0/3)$. Furthermore, the speed $\nu_\mu$ of the corresponding solitary water waves
satisfies
$$\nu_\mu = \nu_0 + 4(\omega+2\nu_0f(k_0))^{-1}\nu_\mathrm{NLS}\mu^2 + o(\mu^2)$$
uniformly over $\eta \in B_\mu$.
\end{list}
\end{theorem}

The first part of Theorem \ref{Main results} is proved by reducing it to a special case of the second.
We proceed by introducing the coercive penalised functional
$\JJ_{\rho,\mu}: H^2(\R) \rightarrow \R \cup \{\infty\}$ defined by
$$\JJ_{\rho,\mu}(\eta) = \left\{\begin{array}{lll}\displaystyle \KK(\eta)+\frac{(\mu+\GG(\eta))^2}{\LL(\eta)}
+ \rho(\|\eta\|_2^2), & & \eta \in U\sm\{0\}, \\
\\
\infty, & & \eta \not\in U\sm\{0\},
\end{array}\right.$$
where $\rho: [0,M^2) \rightarrow {\mathbb R}$ is a smooth, increasing `penalisation' function which
explodes to infinity as $t \uparrow M^2$ and vanishes for $0 \leq t \leq \tilde{M}^2$; the number
$\tilde{M}$ is chosen very close to $M$.
Minimising sequences $\{\eta_m\}$ for $\JJ_{\rho,\mu}$,
which clearly satisfy $\sup_{m \in {\mathbb N}} \|\eta_m\|_2 < M$, are studied in detail in Section \ref{Minimising sequences} with the help
of the concentration-compactness principle (Lions \cite{Lions84a,Lions84b}). The main difficulty
here lies in discussing the consequences of `dichotomy'.

On the one hand the functionals
$\GG$, $\KK$ and $\LL$ are nonlocal and therefore do not act linearly when applied
to the sum of two functions with disjoint supports.
They are however `pseudolocal' in the sense that
$$\begin{Bmatrix} \GG \\ \KK \\ \LL \end{Bmatrix}
(\eta_m^{(1)}+\eta_m^{(2)})
-\begin{Bmatrix} \GG \\ \KK \\ \LL \end{Bmatrix}
(\eta_m^{(1)})
-\begin{Bmatrix} \GG \\ \KK \\ \LL \end{Bmatrix}
(\eta_m^{(2)}) \rightarrow 0$$
as $m \rightarrow \infty$,
where $\{\eta_m^{(1)}\}$, $\{\eta_m^{(2)}\}$ have the properties that
$\supp \eta_m^{(1)} \subset [-R_m,R_m]$ and $\supp \eta_m^{(2)} \subset \R \sm (-S_m,S_m)$
for sequences $\{R_m\}$, $\{S_m\}$ of positive real
numbers with $R_m$, $S_m \rightarrow \infty$, $R_m/S_m \rightarrow 0$ as $m \rightarrow \infty$
(Lemma \ref{Splitting properties}(iii)).
This result is established in Section \ref{Pseudo-local properties} by a new method which involves
studying the weak formulation of the boundary-value problems defining the terms in the power-series
expansion of $K$ about $\eta_0 \in W^{s+3/2}$.
On the other hand no \emph{a priori} estimate is available to rule out `dichotomy'
at this stage; proceeding iteratively we find that minimising sequences can theoretically have profiles
with infinitely many `bumps'. In particular we show that $\{\eta_m\}$
asymptotically lies in the region unaffected by the penalisation
and construct a special minimising sequence $\{\tilde{\eta}_m\}$ for $\JJ_{\rho,\mu}$
which lies in a neighbourhood of the origin with radius $O(\mu^\frac{1}{2})$ in $H^2(\R)$
and satisfies $\|\JJ^\prime_\mu(\tilde{\eta}_m)\|_0 \rightarrow 0$ as $n \rightarrow \infty$. The fact that
the construction is independent of the choice of $\tilde{M}$ allows us to conclude that $\{\tilde{\eta}_m\}$
is also a minimising sequence for $\JJ_\mu$ over $U\sm\{0\}$.

The special minimising sequence $\{\tilde{\eta}_m\}$ is used in Section \ref{SSA section} to establish the
\emph{strict sub-additivity} of the infimum $c_\mu$ of $\JJ_\mu$ over $U \sm \{0\}$, that is the inequality
$$
c_{\mu_1+\mu_2} < c_{\mu_1} + c_{\mu_2}, \qquad 0<\mu_1, \mu_2, \mu_1+\mu_2 < \mu_0.
$$
The strict sub-additivity of $c_\mu$ follows from the fact that the function
\begin{equation}
a \mapsto a^{-q}\MM_{a^2\mu}(a\tilde{\eta}_m), \qquad a \in [1,a_0], \label{Decreasing function of a - Intro}
\end{equation}
is decreasing and strictly negative for some $q>2$ and $a_0 \in (1,2]$, where
$$\MM_\mu(\eta) := \JJ_\mu(\eta) - \KK_2(\eta) - \frac{(\mu+\GG_2(\eta))^2}{\LL_2(\eta)}$$
is the `nonlinear' part of $\JJ_\mu(\eta)$ (see Section \ref{SSA subsection}). We proceed by
approximating $\MM_\mu(\eta_m)$ with its dominant term and showing that this term has the required
property.

The heuristic arguments given above suggest firstly that the spectrum of
minimisers of $\JJ_\mu$ over $U\sm\{0\}$ (that is, the support of their Fourier transform)
is concentrated near wavenumbers $k=\pm k_0$, and secondly that they have the KdV or nonlinear
Schr\"{o}dinger length scales; the same should be true of the functions $\tilde{\eta}_m$,
which approximate minimisers. We therefore
decompose $\tilde{\eta}_m$ into the sum
of a function $\tilde{\eta}_{m,1}$ whose spectrum is compactly supported near $k=\pm k_0$
and a function $\tilde{\eta}_{m,2}$ whose spectrum is bounded away from these points,
and study $\tilde{\eta}_{m,1}$ using the weighted norm
$$
\nn \eta \nn_\alpha^2 := 
\int _{-\infty}^\infty (1+\mu^{-4\alpha}(|k|-k_0)^4) |\hat{\eta}(k)|^2\dk.
$$
A careful analysis of the equation $\JJ^\prime_\mu(\tilde{\eta}_m)=O(\mu^N)$ in $L^2(\R)$ shows that
$\nn \tilde{\eta}_{m,1} \nn_\alpha^2 = O(\mu)$ and $\| \tilde{\eta}_{m,2} \|_2 = O(\mu^{2+\alpha})$
for $\alpha<\frac{1}{3}$ when $\beta>\beta_\mathrm{c}$ and for $\alpha<1$ when $\beta<\beta_\mathrm{c}$.
Using these estimates on the size of $\tilde{\eta}_n$, we find that
$$
\MM_\mu(\tilde{\eta}_m) = \left\{\begin{array}{ll} 
\displaystyle c\!\! \int_{-\infty}^\infty \tilde{\eta}_{m,1}^3 \dx + o(\mu^\frac{5}{3}), & \beta>\beta_\mathrm{c}, \\
\\
\displaystyle -c\!\! \int_{-\infty}^\infty \tilde{\eta}_{m,1}^4 \dx + o(\mu^3),\quad\mbox{} & \beta<\beta_\mathrm{c}. \end{array}\right.
$$
That the function \eqn{Decreasing function of a - Intro} is decreasing and strictly negative follows
from the above estimate and the fact that $\MM_\mu(\eta_m)$ is negative for any minimising sequence
$\{\eta_m\}$ for $\JJ_\mu$ over $U\sm\{0\}$.

Knowledge of the strict sub-additivity property of $c_\mu$ (and general estimates for general minimising sequences)
reduces the proof of part (ii) of Theorem \ref{Main results} to a straightforward application of the concentration-compactness
principle (see Section \ref{Minimisation}). Parts (iii) and (iv) are derived from Lemmata \ref{Variational KdV}(ii) and \ref{Variational NLS}(ii)
by means of a scaling and contradiction argument from the estimates
$$
\|\phi_\eta\|_0^2
= 2\begin{Bmatrix} \alpha_\mathrm{KdV} \\ \alpha_\mathrm{NLS} \end{Bmatrix}  + o(1), \qquad
\begin{Bmatrix}\EE_\mathrm{KdV} \\ \EE_\mathrm{NLS} \end{Bmatrix}(\phi_\eta)
= \begin{Bmatrix} c_\mathrm{KdV} \\ c_\mathrm{NLS} \end{Bmatrix} + o(1), \qquad
\eta \in B_\mu,
$$
which emerge as part of the proof of Theorem \ref{Main results}(i) (see Section \ref{Convergence}).

Some of the techniques used in the present paper were developed by Buffoni \emph{et al.}
\cite{BuffoniGrovesSunWahlen13} in an existence theory for three-dimensional irrotational solitary waves.
While we make reference to relevant parts of that paper, many aspects of our construction differ
significantly from theirs. In particular, our treatment of nonlocal analytic operators is more comprehensive.
Their version of Theorem \ref{Main analyticity result} (see Lemmata 1.1 and 1.4 in that reference) is obtained using a less
sophisticated `flattening' transformation and shows only that the operators are
analytic at the origin. Correspondingly, `pseudo-localness' in the sense described above is established
there only for constant-coefficient boundary-value problems (using an explcit representation of the solution by
means of Green's functions). Our treatment of the consequences of `dichotomy' in the concentration-compactness
principle (Section \ref{Minimising sequences}) is on the other hand similar to that given by
Buffoni \emph{et al.} \cite{BuffoniGrovesSunWahlen13}, and we omit proofs which are straightforward
modifications of theirs; the main difference here is that negative
values of the parameter $\mu$  emerge in our iterative construction of the special minimising
sequence (see the remarks below Lemma \ref{Vanishing and concentration}).

\subsection{Conditional energetic stability}

Our original problem of finding minimisers of $\HH(\eta,\xi)$ subject to the constraint $\II(\eta,\xi) = 2\mu$
is also solved as a corollary to Theorem \ref{Main results}(ii); one follows the two-step minimisation
procedure described in Section \ref{Variational formulation} (see Section \ref{Minimisation}).

\begin{theorem} \label{Intro - main result}
\quad
\begin{list}{(\roman{count})}{\usecounter{count}}
\item
The set $D_\mu$ of minimisers of $\HH$ on the set
$$S_\mu = \{(\eta,\xi) \in U \times H_\star^{1/2}(\R): \II(\eta,\xi) = 2\mu\}$$
is non-empty.
\item
Suppose that $\{(\eta_m,\xi_m)\} \subset S_\mu$ is a minimising sequence for $\HH$ with the
property that $\sup_{m \in \N} \|\eta_m\|_2 < M$. There exists a sequence $\{x_m\} \subset \R$
with the property that a subsequence of $\{(\eta_m(x_m + \cdot), \xi_m(x_m+\cdot)\}$ converges
in $H^r(\R) \times H_\star^{1/2}(\R)$, $r \in [0,2)$, to a function in $D_\mu$.
\end{list}
\end{theorem}

It is a general principle that the solution set of a constrained minimisation
problem constitutes a stable set of solutions of the corresponding
initial-value problem (e.g.\ see Cazenave \& Lions \cite{CazenaveLions82}).
The usual informal interpretation of the statement that a set $X$ of solutions to an initial-value problem
is `stable' is that a solution which begins close to a solution in $X$
remains close to a solution in $X$ at all subsequent times. Implicit in this statement
is the assumption that the initial-value problem is globally well-posed, that is every pair
$(\eta_0,\Phi_0)$ in an appropriately chosen set is indeed the initial datum
of a unique solution $t \mapsto (\eta(t),\Phi(t))$, $t \in [0,\infty)$.
At present there is no global well-posedness theory for gravity-capillary water
waves with constant vorticity (although there is a large and growing body of literature
concerning well-posedness issues for water-wave problems in general).
Assuming the existence of solutions, we obtain the following stability result as a corollary of
Theorem \ref{Intro - main result} using the argument given by Buffoni \emph{et al.}
\cite[Theorem 5.5]{BuffoniGrovesSunWahlen13}.
(The only property of a solution $(\eta,\xi)$ to the initial-value problem which is
relevant to stability theory is that $\HH(\eta(t),\xi(t))$ and $\II(\eta(t),\xi(t))$ are constant; we
therefore adopt this property as the definition of a solution.)

\begin{theorem} \label{CES - Intro}
Suppose that $(\eta,\xi): [0,T] \rightarrow U \times H_\star^{1/2}(\R)$
has the properties that
$$\HH(\eta(t),\xi(t)) = \HH(\eta(0),\xi(0)),\ \II(\eta(t),\xi(t))=\II(\eta(0),\xi(0)), \qquad t \in [0,T]$$
and 
$$\sup_{t \in [0,T]} \|\eta(t)\|_2 < M.$$
Choose $r \in [0,2)$, and let `$\dist$' denote the distance in $H^r(\R) \times H_\star^{1/2}(\R)$.
For each
$\varepsilon>0$ there exists $\delta>0$ such that
$$\dist((\eta(0),\xi(0)), D_\mu) < \delta \quad \Rightarrow \quad
\dist((\eta(t),\xi(t)), D_\mu)<\varepsilon$$
for $t\in[0,T]$.
\end{theorem}

This result is a statement of
the \emph{conditional, energetic stability of the set $D_\mu$}. Here
\emph{energetic} refers to the fact that the
distance in the statement of stability is measured in the `energy space'
$H^r(\R) \times H_\star^{1/2}(\R)$, while \emph{conditional}
alludes to the well-posedness issue. Note that the solution $t \mapsto (\eta(t),\xi(t))$
may exist in a smaller space over the interval $[0,T]$, at each instant of which it remains close
(in energy space) to a solution in $D_\mu$. Furthermore, Theorem \ref{CES - Intro}
is a statement of the stability of the \emph{set} of constrained minimisers $D_\mu$;
establishing the uniqueness of the constrained minimiser would imply
that $D_\mu$ consists of translations of a single solution, so that the statement
that $D_\mu$ is stable is equivalent to classical orbital stability of this unique
solution (Benjamin \cite{Benjamin74}).
The phrase `conditional, energetic stability' was introduced by
Mielke \cite{Mielke02} in his study of the stability of irrotational
solitary water waves with strong surface tension using dynamical-systems
methods.

\section{The functional-analytic setting} \label{FA setting}

\subsection{Nonlocal operators}
The goal of this section is to introduce rigorous definitions of the Dirichlet-Neumann operator $G(\eta)$,
its inverse $N(\eta)$ and the operator $K(\eta):=-\partial_x (N(n) \partial_x)$.

\subsubsection{Function spaces} \label{Function spaces}

Choose $h_0 \in (0,1)$. We consider the class
$$W=\{\eta \in W^{1,\infty}(\R) \colon 1+\inf \eta>h_0\}$$
of surface profiles and  denote the fluid domain by
\[
\Sigma_\eta=\{(x,y)\in \mathbb R^2 \colon 0<y<1+\eta(x)\}, \qquad \eta \in W.
\]
The observation that velocity potentials are unique only up to additive constants leads us to introduce the completion
$H_\star^1(\Sigma_\eta)$ of
$$\SS(\Sigma_\eta)=\{\phi \in C^\infty(\overline{\Sigma}_\eta):
|x|^m|\partial_x^{\alpha_1}\partial_y^{\alpha_2} \phi|\mbox{ is bounded for all }m,\alpha_1,\alpha_2 \in {\mathbb N}_0 \}$$
with respect to the Dirichlet norm
as an appropriate function space for $\phi$. The corresponding space for the trace $\phi|_{y=1+\eta}$ is the space
$H_\star^{1/2}(\R)$ defined in Section \ref{FAF}.

\begin{proposition}
\label{prop:trace}
Fix $\eta \in W$. 
The trace map $\phi\mapsto \phi|_{y=1+\eta}$ 
defines a continuous operator $H^1_\star(\Sigma_\eta)\rightarrow H_\star^{1/2}(\R)$ with a continuous right inverse $H_\star^{1/2}(\R)\rightarrow H_\star^1(\Sigma_\eta)$.
\end{proposition}

We also use anisotropic function spaces for functions defined in the strip $\Sigma_0 = \R \times (0,1)$.

\begin{definition}
Suppose that $r\in \R$ and $n \in \N_0$. 
\begin{list}{(\roman{count})}{\usecounter{count}}
\item
The Banach space $(L^\infty H^r,\|\cdot\|_{r,\infty})$ is defined by
$$L^\infty H^r=L^\infty((0,1), H^r(\R)), \qquad
\|u\|_{r,\infty} 
=\mathop{\esssup}_{y\in (0,1)} \|u(\cdot,y)\|_{H^r(\R)}.
$$
\item
The Banach space $(H^{r,m}, \|\cdot\|_{r,m})$ is defined by
\[
H^{r,m}=\bigcap_{j=0}^n H^j((0,1), H^{r-j}(\R)), \qquad \|u\|_{r,m}=\sum_{j=0}^n \|\Lambda^{r-j}\partial_y^j u\|_{L^2(\Sigma)},
\]
where $\Lambda f=\FF^{-1}[(1+k^2)^\frac{1}{2}\hat f(k)]$.
\end{list}
\end{definition}

The following propositions state some properties of these function spaces which are used in the subsequent analysis;
they are deduced from results for standard Sobolev spaces (see H\"{o}rmander \cite[Theorem 8.3.1]{Hoermander} for
Proposition \ref{Hprop}).

\begin{proposition}
\quad
\begin{list}{(\roman{count})}{\usecounter{count}}
\item
The space $C_0^\infty(\overline \Sigma)$ is dense in $H^{r,1}$ for each $r \in \R$.
\item
For each $r \in \R$ the mapping $u\mapsto u|_{y=1}$, $u\in C_0^\infty(\overline \Sigma)$, extends continuously to an operator $H^{r+1,1} \rightarrow H^{r+1/2}(\R)$.
\item
The space $H^{r+1,1}$ is continuously embedded in $L^\infty H^{r+1/2}$ for each $r \in \R$.
\item
The space $H^{r+1,1}$ is a Banach algebra for each $r > 0$.
\end{list}
\end{proposition}
\begin{proposition} \label{Hprop}
Suppose that $r_0$, $r_1$, $r_2$ satisfy $r_0 \leq r_1$, $r_0 \leq r_2$, $r_1+r_2 \geq 0$ and $r_0 < r_1+r_2 - \frac{1}{2}$. The
product $u_1u_2$ of each $u_1 \in L^\infty H^{r_1}$ and $u_2 \in H^{r_2,0}$ lies in $H^{r_0,0}$ and satisfies
\[
\|u_1u_2\|_{r_0,0}\le c\|u_1\|_{r_1,\infty} \|u_2\|_{r_2,0}.
\]
\end{proposition}
\begin{proposition} \label{Bprop}
For each bounded linear function $L: L^2(\R) \rightarrow L^\infty H^0$ the formula $(\eta,w) \mapsto L(\eta)w$ defines a bounded bilinear
function $L^2(\R) \times H^1(\Sigma) \rightarrow L^2(\Sigma)$ which satisfies the estimate
$$\|L(\eta)w\|_0 \leq c\|L\|\|w\|_0^\frac{1}{2}\|w\|_1^\frac{1}{2}\|\eta\|_0.$$
The assertion remains valid when $\Sigma$ is replaced by $\{|x| < M\}$ or $\{|x| > M\}$ and the estimate
holds uniformly over all values of $M$ greater than unity.
\end{proposition}

\subsubsection{The Dirichlet-Neumann operator}

The Dirichlet-Neumann operator $G(\eta)$ for the boundary-value problem 
\begin{eqnarray}
& & \parbox{6cm}{$\Delta \phi=0,$}0<y<1+\eta, \label{BC for DNO 1} \\
& & \parbox{6cm}{$ \phi = \xi,$}y=1+\eta, \label{BC for DNO 2}\\
& & \parbox{6cm}{$\phi_y =0,$}y=0 \label{BC for DNO 3}
\end{eqnarray}
is defined formally as follows: fix $\xi=\xi(x)$, solve \eqn{BC for DNO 1}--\eqn{BC for DNO 3} and set
$$G(\eta)\xi=(\phi_y-\eta^\prime\phi_x)|_{y=1+\eta}.$$
Our rigorous definition of $G(\eta)$ is given in terms of weak solutions to \eqn{BC for DNO 1}--\eqn{BC for DNO 3}
(see Lannes \cite[Proposition 2.9]{Lannes} for the proof of Lemma \ref{thm:Dirichlet weak solution}).
\begin{definition}
Suppose that $\xi\in H_\star^{1/2}(\R)$ and $\eta \in W$. A 
\underline{weak solution} of \eqref{BC for DNO 1}--\eqref{BC for DNO 3} is a function 
$\phi \in H^1_\star(\Sigma_\eta)$ with $\phi|_{y=1+\eta}=\xi$ which satisfies
$$
\int_{\Sigma_\eta} \nabla \phi \cdot \nabla \psi \dx \dy=0
$$
for all $\psi \in H^1_\star(\Sigma_\eta)$ with $\psi|_{y=1+\eta}=0$.
\end{definition}
\begin{lemma}
\label{thm:Dirichlet weak solution}
For each $\xi\in H_\star^{1/2}(\R)$ and $\eta \in W$ there
exists a unique weak solution $\phi$ of \eqref{BC for DNO 1}--\eqref{BC for DNO 3}.
The solution satisfies the estimate
\[
\|\phi\|_{H_\star^1(\Sigma_\eta)}\le C\|\xi\|_{H_\star^{1/2}(\R)},
\]
where $C=C(\|\eta\|_{1,\infty})$.
\end{lemma}
\begin{definition}
\label{def:Dirichlet-Neumann}
Suppose that $\eta\in W$ and $\xi \in H_\star^{1/2}(\R)$.
The \underline{Dirichlet-Neumann operator} is the bounded linear operator
$G(\eta)\colon H_\star^{1/2}(\R)\rightarrow H_\star^{-1/2}(\R)$ 
defined by
$$
\int_{-\infty}^\infty (G(\eta)\xi_1)\, \xi_2 \dx=\int_{\Sigma_\eta} \nabla \phi_1\cdot \nabla \phi_2 \dx \dy,
$$
where $\phi_j  \in H^1_\star(\Sigma_\eta)$ is the unique weak solution of \eqref{BC for DNO 1}--\eqref{BC for DNO 3} with $\xi=\xi_j$, $j=1,2$.
\end{definition}

\subsubsection{The Neumann-Dirichlet operator}

The Neumann-Dirichlet operator $N(\eta)$ for the the boundary-value problem 
\begin{eqnarray}
& & \parbox{6cm}{$\Delta \phi=0,$}0<y<1+\eta, \label{BC for NDO 1} \\
& & \parbox{6cm}{$\phi_y - \eta^\prime\phi_x = \xi,$}y=1+\eta, \label{BC for NDO 2}\\
& & \parbox{6cm}{$\phi_y =0,$}y=0 \label{BC for NDO 3}
\end{eqnarray}
is defined formally as follows: fix $\xi=\xi(x)$, solve \eqn{BC for NDO 1}--\eqn{BC for NDO 3}
and set
$$N(\eta)\xi=\phi|_{y=1+\eta}.$$
Our rigorous definition of $N(\eta)$ is also given in terms of weak solutions; Lemma \ref{thm:Neumann weak solution}
is proved in the same fashion as Lemma \ref{thm:Dirichlet weak solution}.

\begin{definition}
Suppose that $\xi\in H_\star^{-1/2}(\R)$ and $\eta \in W$. A 
\underline{weak solution} of \eqref{BC for NDO 1}--\eqref{BC for NDO 3} is a function 
$\phi \in H^1_\star(\Sigma_\eta)$ which satisfies
$$
\int_{\Sigma_\eta} \nabla \phi \cdot \nabla \psi \dx \dy=\int_{-\infty}^\infty \xi \psi|_{y=1+\eta} \dx
$$
for all $\psi\in H^1_\star(\Sigma_\eta)$.
\end{definition}
\begin{lemma}
\label{thm:Neumann weak solution}
For each $\xi\in H_\star^{-1/2}(\R)$ and $\eta \in W$ there
exists a unique weak solution $\phi$ of \eqref{BC for NDO 1}--\eqref{BC for NDO 3}. The solution satisfies the estimate
\[
\|\phi\|_{H_\star^1(\Sigma_\eta)}\le C\|\xi\|_{H_\star^{-1/2}(\R)},
\]
where $C=C(\|\eta\|_{1,\infty})$.
\end{lemma}
\begin{definition}
\label{def:Neumann-Dirichlet}
Suppose that $\eta\in W$ and $\xi \in H_\star^{-1/2}(\R)$.
The \underline{Neumann-Dirichlet operator} is the bounded linear operator
$N(\eta)\colon H_\star^{-1/2}(\R)\rightarrow H_\star^{1/2}(\R)$  defined by
$$
N(\eta)\xi=\phi|_{y=1+\eta},
$$
where $\phi \in H^1_\star(\Sigma_\eta)$ is the unique weak solution of \eqref{BC for NDO 1}--\eqref{BC for NDO 3}.
\end{definition}

The relationship between $G(\eta)$ and $N(\eta)$ is clarified by the following result, which follows from the
definitions of these operators.
\begin{lemma} \label{DNO is inverse of NDO}
Suppose that $\eta\in W$. The operator $G(\eta)\in \LL(H_\star^{1/2}(\R), H_\star^{-1/2}(\R))$ is 
invertible with $G(\eta)^{-1}=N(\eta)$.
\end{lemma}

\subsubsection{Analyticity of the operators}

Let us begin by recalling the definition
of analyticity given by Buffoni \& Toland \cite[Definition 4.3.1]{BuffoniToland} together with
a precise formulation of our result in their terminology.

\begin{definition} \label{BT defn of analyticity}
Let $X$ and $Y$ be Banach spaces,  $U$ be a non-empty, open subset of $X$ and
$\LL_\mathrm{s}^k(X,Y)$ be the space of bounded, $k$-linear
symmetric operators $X^k \rightarrow Y$ with norm
$$\nn m\nn:=\inf \{c: \|m(\{f\}^{(k)})\|_Y \leq c \|f\|_X^k \mbox{ \rm{for all} } f \in X\}.$$
A function $F:U \rightarrow Y$ is \underline{analytic at a point $x_0 \in U$} if there exist
real numbers $\delta, r>0$ and a sequence $\{m_k\}$, where $m_k \in \LL_\mathrm{s}^k(X,Y)$,
$k=0,1,2,\ldots$, with the properties that
$$F(x) = \sum_{k=0}^\infty m_k(\{x-x_0\}^{(k)}), \qquad x \in B_\delta(x_0)$$
and
$$\sup_{k \geq 0} r^k\nn m_k\nn < \infty.$$
The function is \underline{analytic} if it is analytic at each point $x_0\in U$.
\end{definition}

\begin{theorem} \label{G and N are analytic}
\quad
\begin{list}{(\roman{count})}{\usecounter{count}}
\item
The Dirichlet-Neumann operator $G(\cdot): W \rightarrow
\LL(H_\star^{1/2}(\R), H_\star^{-1/2}(\R))$ is analytic.
\item
The Neumann-Dirichlet operator $N(\cdot): W \rightarrow
\LL(H_\star^{-1/2}(\R), H_\star^{1/2}(\R))$ is analytic.
\end{list}
\end{theorem}

To prove this theorem we study the dependence of solutions to the boundary-value
problems \eqn{BC for DNO 1}--\eqn{BC for DNO 3}
and \eqn{BC for NDO 1}--\eqn{BC for NDO 3} on $\eta$
by transforming them into equivalent problems in the fixed domain $\Sigma:=\Sigma_0$.
For this purpose we define a change of variable $(x,y)=F^\delta(x,y^\prime)$ in the following way.
Choose $\delta > 0$ and an even function $\chi \in C_0^\infty(\mathbb R)$ with $\chi(k) \in [0,1]$ for $k\in \R$, 
$\supp \chi \in [-2,2]$ and $\chi(x)\equiv 1$ for $|x|\le 1$, write 
\[
\eta^\delta(x,y^\prime)=\FF^{-1}[\chi(\delta (1-y^\prime)k)\hat \eta(k)](x)
\]
and define
$$F^\delta(x,y^\prime)=(x, y^\prime(1+\eta^\delta(x,y^\prime)))=(x,y^\prime +f^\delta(x,y^\prime)),$$
in which 
$f^\delta(x,y^\prime)=y^\prime\eta^\delta(x,y^\prime)$.

\begin{lemma}
Suppose that $\eta \in W$. The mapping $F^\delta$  is a bijection $\Sigma \rightarrow \Sigma_\eta$  and $\overline{\Sigma} \rightarrow \overline{\Sigma}_\eta$ with $y \in C_\mathrm{b}^1(\Sigma)$, $y^\prime \in C_\mathrm{b}^1(\Sigma_\eta)$ and
$$\inf\limits_{(x,y^\prime) \in \overline{\Sigma}} y_{y^\prime}(x,y^\prime)
=\inf\limits_{(x,y^\prime) \in \overline{\Sigma}} (1+f_{y^\prime}^\delta(x,y))>0$$
for each $\delta \in (0,\delta_\mathrm{max})$, where $\delta_\mathrm{max}=\delta_\mathrm{max}(\|\eta^\prime\|_\infty^{-1})$.
\end{lemma}
{\bf Proof.} Writing
$$\eta^\delta(x,y^\prime) = \int_{-\infty}^\infty K(s) \eta(x-\delta(1-y^\prime)s) \ds,$$
where $K = (2\pi)^{-\frac{1}{2}} \FF^{-1}[\chi] \in \SS(\R)$,
one finds that $\eta^\delta \in C^\infty(\Sigma) \cap C_\mathrm{b}^1(\Sigma)$ with\linebreak
$\|\eta^\delta\|_\infty \leq c\|\eta\|_\infty$, $\|\eta^\delta_x\|_\infty \leq c\|\eta^\prime\|_\infty$, $\|\eta^\delta_{y^\prime}\|_\infty \leq c\delta \|\eta^\prime\|_\infty$.
It follows that $F^\delta \in C^\infty(\Sigma)$ and $y \in C_\mathrm{b}^1(\Sigma)$.
Furthermore $y(x,0)=0$, $y(x,1) = 1+ \eta(x)$ and
\begin{eqnarray*}
\partial_{y^\prime}{y} & = & 1+y^\prime \eta^\delta_{y^\prime} + \eta^\delta \\
& = & 1+y^\prime \eta^\delta_{y^\prime} + \eta - \int_{y^\prime}^1 \eta^\delta_{y^\prime} \\
& \geq & h_0 - c\delta\|\eta^\prime\|_\infty \\
& \geq & \tfrac{1}{2}h_0 \\
& > & 0
\end{eqnarray*}
for sufficiently small $\delta$ (depending only upon $\|\eta^\prime\|_\infty^{-1}$), so that $F^\delta$ is a bijection
$\Sigma \rightarrow \Sigma_\eta$  and $\overline{\Sigma} \rightarrow \overline{\Sigma}_\eta$. It follows from the inverse function theorem
that $(F^\delta)^{-1} \in C^\infty(\Sigma_\eta)$; the estimate
$$\det \mathrm{d}F^\delta[x,y^\prime]=\partial_{y^\prime} y(x,y^\prime) \geq \tfrac{1}{2}h_0$$
and the fact that $\mathrm{d}F^\delta$ is bounded on $\Sigma$ imply that
$\mathrm{d}(F^\delta)^{-1} \in C_\mathrm{b}(\Sigma_\eta)$, whereby $y^\prime \in C_\mathrm{b}^1(\Sigma_\eta)$.
\qed

The change of variable $(x,y)=F^\delta(x,y^\prime)$ transforms 
the boundary-value problem \eqn{BC for NDO 1}--\eqn{BC for NDO 3} into
\begin{eqnarray}
& & \parbox{90mm}{$\nabla\cdot((I+Q)\nabla u)=0$}0 < y <1,
\label{BC for u 1} \\
& & \parbox{90mm}{$(I+Q)\nabla u \cdot (0, 1) = \xi,$}y=1, \label{BC for u 2} \\
& & \parbox{90mm}{$(I+Q)\nabla u \cdot (0,-1)=0,$}y=0, \label{BC for u 3}
\end{eqnarray}
where
\[
Q=\begin{pmatrix}f^\delta_y & - f^\delta_x \\ 
-  f^\delta_x & \displaystyle{\frac{- f^\delta_y+( f^\delta_x)^2}{1+ f^\delta_y}}\end{pmatrix}
\]
and the primes have been dropped for notational simplicity.

\begin{lemma}\label{Q is analytic 1}
The mapping $W\rightarrow (L^\infty(\overline{\Sigma}))^{2\times 2}$ given by $\eta \mapsto Q(\eta)$ is analytic.
\end{lemma}

 It is helpful to consider the more general boundary-value problem
\begin{eqnarray}
& & \parbox{90mm}{$\nabla\cdot((I+Q)\nabla u)=\nabla\cdot G$}0 < y <1,
\label{BC for u 4} \\
& & \parbox{90mm}{$(I+Q)\nabla u \cdot (0, 1) = \xi+G \cdot (0, 1) ,$}y=1, \label{BC for u 5} \\
& & \parbox{90mm}{$(I+Q)\nabla u \cdot (0,-1)=G\cdot (0,-1),$}y=0, \label{BC for u 6}
\end{eqnarray}
where $I+Q \in(L^\infty(\overline{\Sigma}))^{2\times 2}$ is uniformly positive definite, that is, there
exists a constant $p_0>0$ such that
\[
(I+Q)(x, y)\nu\cdot \nu \ge p_0|\nu|^2,
\]
for all $(x,y)\in \overline \Sigma$ and all $\nu \in \R^2$.

\begin{definition} \label{Weak soln}
Suppose that
$\xi \in H_\star^{-1/2}(\R)$ and $G\in (L^2(\Sigma))^2$.
A \underline{weak solution} of \eqn{BC for u 4}--\eqn{BC for u 6} is a
function $u \in H_\star^1(\Sigma)$ which satisfies
$$
\int_\Sigma (I+Q)\nabla u\cdot \nabla w \dx\dy = \int_\Sigma G \cdot \nabla w\dx\dy +\int_{-\infty}^\infty \xi w|_{y=1} \dx
$$
for all $w \in H_\star^1(\Sigma)$.
\end{definition}
\begin{lemma} \label{Flattened generalised BVP has a WS}
For each $\xi \in H_\star^{-1/2}(\R)$ and $G\in (L^2(\Sigma))^2$
the boundary-value problem \eqn{BC for u 4}--\eqn{BC for u 6} has a unique weak solution
$u \in H_\star^1(\Sigma)$. The solution satisfies the estimate
\[
\|u\|_{H_\star^1(\Sigma)}\le C(\|\xi\|_{H_\star^{-1/2}(\R)}+\|G\|_{L^2(\R)}),
\]
where $C=C(p_0^{-1})$.
\end{lemma}

Lemma \ref{Flattened generalised BVP has a WS} applies in particular to \eqn{BC for u 1}--\eqn{BC for u 3}
for each fixed $\eta \in W$ (the matrix $I+Q$ is uniformly positive definite since it is uniformly bounded
above, its determinant is unity and its upper left entry is positive). The next theorem shows that its unique weak solution depends analytically upon $\eta$.

\begin{theorem} \label{Weak solns analytic 1}
The mapping $W \rightarrow
\LL(H_\star^{-1/2}(\R), H_\star^1(\Sigma))$ given by
$\eta \mapsto (\xi \mapsto u)$, where $u \in H_\star^1(\Sigma)$
is the unique weak solution of \eqn{BC for u 1}--\eqn{BC for u 3},
is analytic.
\end{theorem}
{\bf Proof.} Choose $\eta_0 \in W$ and write $\tilde{\eta}=\eta-\eta_0$ and
\[
Q(x,y)=\sum_{n=0}^\infty Q^n(x,y), \qquad Q^n = \tilde{m}_n(\tilde{\eta}^{\{(n)\}})
\]
where $\tilde{m}_n(\tilde{\eta}^{\{(n)\}}) \in \LL_s^n(W^{1,\infty}(\R), (L^\infty(\overline{\Sigma}))^{2\times 2})$ satisfies
$$\nn\tilde{m}^n\nn \leq C_2 r^{-n} \|\tilde{\eta}\|_{1,\infty}^n$$
(see Lemma \ref{Q is analytic 1}). We proceed by seeking a solution of \eqn{BC for u 1}--\eqn{BC for u 3}
of the form
\begin{equation}
u(x,y) = \sum_{n=0}^\infty u^n(x,y), \qquad u^n = m_1^n(\{\tilde{\eta}\}^{(n)}) \label{Formal expansion}
\end{equation}
where
$m_1^n \in \LL_\mathrm{s}^n(W^{1,\infty}(\R),H_\star^1(\Sigma))$ is linear in $\xi$ and satisfies
$$\nn m_1^n \nn \leq C_1B^n \|\xi\|_{H_\star^{-1/2}(\R)}$$
for some constant $B>0$.

Substituting the \emph{Ansatz} \eqn{Formal expansion} into the equations, one finds that
\begin{eqnarray} 
& & \parbox{60mm}{$\nabla \cdot ((I+Q^0)\nabla u^0) = 0,$}0 < y <1, \label{BC for u0 1}\\
& & \parbox{60mm}{$(I+Q^0) \nabla u^0\cdot (0,1) = \xi,$}y=1, \label{BC for u0 2}\\
& & \parbox{60mm}{$(I+Q^0) \nabla u^0\cdot (0,-1)=0,$}y=0 \label{BC for u0 3}
\end{eqnarray}
and
\begin{eqnarray}
& & \parbox{90mm}{$\nabla \cdot ((I+Q^0)\nabla u^n) = \nabla \cdot G^n,$}0 < y <1, \label{BC for un 1} \\
& & \parbox{90mm}{$(I+Q^0) \nabla u^n\cdot (0,1) =  G^n\cdot (0,1),$}y=1, \label{BC for un 2} \\
& & \parbox{90mm}{$(I+Q^0) \nabla u^n\cdot (0,-1)=G^n\cdot (0,-1),$}y=0 \label{BC for un 3}
\end{eqnarray}
for $n=1,2,3,\ldots$, where
$$
G^n = - \sum_{k=1}^n Q^k \nabla u^{n-k}.
$$
The estimate for $m^0$ follows directly from Lemma \ref{Flattened generalised BVP has a WS}. Proceeding
inductively, suppose the result for $m^n$ is true for all $k<n$. Estimating
\begin{eqnarray}
\|G^n\|_0 & \leq & \sum_{k=1}^n \|Q^k\|_\infty \|\nabla u^{n-k}\|_0 \label{Expansion of Gn} \\
& \leq & C_1 C_2 B^n \|\xi\|_{H_\star^{-1/2}(\R)}\|\tilde{\eta}\|_{1,\infty}^n \sum_{k=1}^n (Br)^{-k} \nonumber
\end{eqnarray}
and using Lemma \ref{Flattened generalised BVP has a WS} again, we find that
\begin{eqnarray*}
\|u^n\|_{H_\star^1(\Sigma)} & \leq & C_1 C_2 C_3 B^n \|\xi\|_{H_\star^{-1/2}(\R)}\|\tilde{\eta}\|_{1,\infty}^n \sum_{k=1}^\infty (Br)^{-k} \\
& \leq & C_1 B^n \|\xi\|_{H_\star^{-1/2}(\R)}\|\tilde{\eta}\|_{1,\infty}^n
\end{eqnarray*}
for sufficiently large values of $B$ (independently of $n$).

A straightforward supplementary argument shows that the expansion \eqn{Formal expansion}
defines a weak solution $u$ of \eqn{BC for un 1}--\eqn{BC for un 3}. \qed

Theorem \ref{G and N are analytic}(ii)
follows from the above theorem, the formula $N(\eta)\xi = u|_{y=1}$ and the continuity of the trace operator
$H_\star^1(\Sigma) \rightarrow H_\star^{1/2}(\R)$, while
Theorem \ref{G and N are analytic}(i) follows from the inverse function theorem for analytic functions.

Finally, we record another useful result.

\begin{theorem} \label{Properties of DNO}

For each $\eta \in W$ the norms
$$\xi \mapsto \left(\int_{-\infty}^\infty  \xi G(\eta) \xi \dx\right)^{\!\!\frac{1}{2}}, \qquad
\kappa \mapsto \left(\int_{-\infty}^\infty  \kappa N(\eta) \kappa \dx\right)^{\!\!\frac{1}{2}}$$
are equivalent to the usual norms for respectively $H_\star^{1/2}(\R)$ and $H_\star^{-1/2}(\R)$ .

\end{theorem}
{\bf Proof.} Let $T:H_\star^{-1/2}(\R) \mapsto H_\star^{1/2}(\R)$ be the isometric isomorphism 
$\eta \mapsto \FF^{-1}[(1+k^2)^\frac{1}{2}k^{-2}\hat \eta]$,
which has the property that
\[
\int_{-\infty}^\infty \psi\, \xi\dx = \langle T\psi, \xi\rangle_{H_\star^{1/2}(\R)}, \qquad \psi \in H_\star^{-1/2}(\R),
\xi\in H_\star^{1/2}(\R). 
\]
It follows from Definition \ref{def:Dirichlet-Neumann}, Lemma \ref{DNO is inverse of NDO}
and the calculation
\[
\langle TG(\eta)\xi, \xi\rangle_{H_\star^{1/2}(\R)}=\int_{-\infty}^\infty (G(\eta)\xi)\, \xi \dx=\int_{\Sigma_\eta} |\nabla \phi|^2 \dx \dy\ge 0,
\]
where $\phi$ is the unique weak solution of \eqref{BC for DNO 1}--\eqref{BC for DNO 3},
that $TG(\eta)$ is a self-adjoint, positive, isomorphism $H_\star^{1/2}(\R) \rightarrow H_\star^{1/2}(\R)$.
The spectral theory for bounded, self-adjoint operators shows that $\xi \mapsto \langle TG(\eta)\xi,\xi\rangle_{H_\star^{1/2}(\R)}^\frac{1}{2}$
and $\xi \mapsto \langle N(\eta)T^{-1}\xi,\xi\rangle_{H_\star^{1/2}(\R)}^\frac{1}{2}$ are equivalent to the usual norm for $H_\star^{1/2}(\R)$, so that
$\kappa \mapsto \langle N(\eta)\kappa,T\kappa\rangle_{H_\star^{1/2}(\R)}^\frac{1}{2}$ is equivalent to the usual norm for $H_\star^{-1/2}(\R)$. The assertion now
follows from the first equality in the previous equation and the calculation
$$
\langle N(\eta)\kappa, T\kappa\rangle_{H_\star^{1/2}(\R)}=\int_{-\infty}^\infty (N(\eta)\kappa)\, \kappa \dx.\eqno{\Box}
$$

\subsubsection{The operator \boldmath$K(\eta) = -\partial_x (N(\eta)\partial_x)$\unboldmath}

Our first result for this operator is obtained from the material presented above for $N$.

\begin{theorem} \label{Analyticity of K 1}
\quad
\begin{list}{(\roman{count})}{\usecounter{count}}
\item
The operator $K(\cdot): W \rightarrow
\LL(H^{1/2}(\R), H^{-1/2}(\R))$ is analytic.
\item
For each $\eta \in W$ the operator 
$K(\eta): H^{1/2}(\R) \rightarrow H^{-1/2}(\R)$
is an isomorphism and the norm
$$\zeta \mapsto \left(\int_{-\infty}^\infty  \zeta K(\eta) \zeta \dx\right)^{\!\!\frac{1}{2}}$$
is equivalent to the usual norm for $H^{1/2}(\R)$.
\end{list}
\end{theorem}
{\bf Proof.} (i) This result follows from the definition of $K$
and the continuity of the operators\linebreak
$\partial_x: H^{1/2}(\R) \rightarrow H_\star^{-1/2}(\R)$
and $\partial_x: H_\star^{1/2}(\R) \rightarrow H^{-1/2}(\R)$.

(ii) This follows by writing 
\begin{align*}
\int_{-\infty}^\infty  \zeta K(\eta) \zeta \dx&=
\int_{-\infty}^\infty  \zeta^\prime  N(\eta) \zeta^\prime \dx\\
&\ge c\|\zeta^\prime\|_{H_\star^{-1/2}(\R)}^2\\
&=c\|\zeta\|_{1/2}^2,
\end{align*}
in which Theorem \ref{Properties of DNO} has been used.
\qed

In the remainder of this section we establish the following result concerning the analyticity of
$K$ in higher-order Sobolev spaces, using the symbol $W^r$ as an abbreviation for $W \cap H^r(\R)$.

\begin{theorem} \label{Analyticity of K 2}
The operator $K(\cdot): W^{s+3/2} \rightarrow
\LL(H^{s+3/2}(\R), H^{s+1/2}(\R))$ is analytic for each $s>0$.
\end{theorem}

To prove Theorem \ref{Analyticity of K 2} it is necessary to establish additional regularity of the weak solutions
$u^n$, $n=1,2,\ldots$ of the boundary-value problems \eqn{BC for u0 1}--\eqn{BC for u0 3} and
\eqn{BC for un 1}--\eqn{BC for un 3}.
We proceed by examining the general boundary-value problem \eqn{BC for u 4}--\eqn{BC for u 6} under additional
regularity assumptions on $\zeta$ and $G$. Our result is stated in Lemma \ref{Regularity of u} below, whose proof requires
an \emph{a priori} estimate and a commutator estimate (see Lannes \cite[Proposition B.10(2)] {Lannes} for a derivation of
the latter).

\begin{lemma}
\label{regularity lemma}
Suppose that $Q \in (H^{s+1,2})^{2\times 2}$ and $G \in (H^{t,1})^2$ for some $t \in (\frac{1}{2}-s, s+1]$. The weak solution $u$ to
\eqn{BC for u 4}--\eqn{BC for u 6} satisfies the \emph{a priori} estimate
\[
\|\nabla u\|_{t,1}\le C(\|G\|_{t,1}+\|\nabla u\|_{t,0}),
\]
where $C=C(p_0^{-1},\|Q\|_{s+1,2})$.
\end{lemma}
{\bf Proof.} Note that
\begin{eqnarray*}
\| \nabla u \|_{t,1} & = & \| u_x \|_{t,1} + \|u_y \|_{t,1} \\
& = & \|u_x\|_{t,0} + \|u_{xy}\|_{t-1,0} + \|u_y\|_{t,0} + \|u_{yy}\|_{t-1,0} \\
& \leq & C(\|\nabla u\|_{t,0}+\|u_{yy}\|_{t-1,0})
\end{eqnarray*}
because $\|u_{xy}\|_{t-1,0} \leq \|u_y\|_{t,0}$, and to estimate 
$\|u_{yy}\|_{t-1,0}$ we use equation \eqn{BC for u 4}, which we write in the form
$$
(1+q_{22})u_{yy}=\nabla \cdot G-\partial_x [(1+q_{11})u_x+q_{12} u_y] -
\partial_y(q_{12}u_x)-q_{22y}u_y.
$$
Denoting the right hand side of this equation by $H$, one finds that
\begin{eqnarray*}
\|u_{yy}\|_{t-1,0} & = & \|(1+q_{22})^{-1} H\|_{t-1,0} \\
& \leq & \|H\|_{t-1,0} +\|\tilde q_{22} H\|_{t-1,0} \\
& \leq & (1+ \|\tilde q_{22}\|_{s+1/2,\infty})\|H\|_{t-1,0}\\
& \leq & C \|H\|_{t-1,0},
\end{eqnarray*}
where $\tilde q_{22}=-q_{22}(1+q_{22})^{-1}$ and we have used the interpolation estimate
$$\left\| \frac{p}{1+p} \right\|_r \leq C_1(p_0^{-1},\|p\|_\infty) \|p\|_r \leq C_2(p_0^{-1}, \|p\|_r)$$
for $p \in H^r(\R)$, $r>\frac{1}{2}$ with $1+p(x) \geq p_0$ for all $x \in \R$.

It remains to estimate $\|H\|_{t-1,0}$. Observe that $\|\nabla \cdot G\|_{t-1,0} \leq \|G\|_{t,1}$,
$\|u_{xx}\|_{t-1,0} \leq \|\nabla u\|_{t,0}$ and
\begin{eqnarray}
\|q_{ij} \nabla u_x\|_{t-1,0}
& \leq &
C\|Q\|_{s+1/2,\infty} \|\nabla u_x\|_{t-1,0} \nonumber \\
&\le &
C\|Q\|_{s+1,1}\|\nabla u\|_{t,0}. \label{H estimate 1}
\end{eqnarray}
The terms in $H$ involving derivatives of $Q$ are treated differently.

Suppose first that $t\le s+\frac{1}{2}$. Combining the estimate
\begin{eqnarray*}
\|\begin{Bmatrix}\partial_x\\ \partial_y\end{Bmatrix} q_{ij} \nabla u\|_{t-1,0}
&\le & C\|\begin{Bmatrix}\partial_x\\ \partial_y\end{Bmatrix} q_{ij}\|_{s-1/2,\infty}\|\nabla u\|_{t,0}\\
&\le &
C\|Q\|_{H^{s+1,2}}\|\nabla u\|_{t,0}
\end{eqnarray*}
(Proposition \ref{Hprop}) and the estimate \eqn{H estimate 1},
one obtains the required result
$$\|u_{yy}\|_{t-1,0} \leq  \|H\|_{t-1,0}\le C(\|G\|_{t,1}+ \|\nabla u\|_{t,0}).$$

In the case $t \in (s+\frac{1}{2},s+1]$ we instead estimate
\begin{align*}
\|\begin{Bmatrix}\partial_x\\ \partial_y\end{Bmatrix} q_{ij} \nabla u\|_{t-1,0}
&\le C\|\begin{Bmatrix}\partial_x\\ \partial_y\end{Bmatrix} q_{ij}\|_{s,0} \|\nabla u\|_{t-1/2-\varepsilon,\infty}\\
&\le 
C\|Q\|_{s+1,1}\|\nabla u\|_{t-\varepsilon,1}
\end{align*}
with 
$0<\varepsilon< \min\{\frac{1}{2},s\}$
by Proposition \ref{Hprop} to find that
\begin{align*}
\|u_{yy}\|_{t-1,0}&\le C(\|G\|_{t,1}+\|\nabla u\|_{t,0}+ \|\nabla u\|_{t-\varepsilon,1})\\
&\le C(\|G\|_{t,1}+\|\nabla u\|_{t,0}+ \|u_{yy}\|_{t-1-\varepsilon,0}).
\end{align*}
The result follows by repeating this argument a finite number of times and using the already
established result for $t=s+\frac{1}{2}$.\qed

\begin{lemma}
\label{commutator estimate}
Suppose that $r_0>\frac{1}{2}$, $\Delta \in [0,1]$ and $r \in (-\frac{1}{2},r_0+\Delta]$
and define $\Lambda^r_\varepsilon=\Lambda^r \chi(\varepsilon \Lambda)$ for $\varepsilon \in [0,\varepsilon_0)$.
The estimate
\[
\|[\Lambda^r_\varepsilon, u] v\|_0\le c\|u\|_{r_0+\Delta}\|v\|_{r-\Delta}.
\]
holds for each $u \in H^{r_0+\Delta}$ and each $v \in H^{r-\Delta}$, where the constant $c$
does not depend upon $\varepsilon$.
\end{lemma}

\begin{lemma}\label{Regularity of u}
Suppose that $Q \in (H^{s+1,2})^{2\times 2}$ and $\zeta \in H^{t+3/2}(\mathbb R)$, $G\in (H^{t+1,1})^2$ for some $t \in [0,s]$.
The weak solution
$u$ of \eqn{BC for u 4}--\eqn{BC for u 6} with $\xi=\zeta^\prime$ satisfies 
$\nabla u\in H^{t+1,1}$ with 
\[
\|\nabla u\|_{t+1,1}\le C(\|G\|_{t+1,1}+\|\zeta\|_{t+3/2}),
\]
where $C=C(p_0^{-1},\|Q\|_{s+1,2})$.
\end{lemma}
{\bf Proof.}
Choose $ r \in (0,t+1]$, $\varepsilon >0$ and note that $\Lambda^r_\varepsilon$ is well defined as an operator on $H_\star^1(\Sigma)$. 
Writing $w=(\Lambda_\varepsilon^r)^2 u$ in Definition \ref{Weak soln}, we find that
\[
\int_\Sigma \Lambda_\varepsilon^r (P\nabla u)\cdot \nabla \Lambda_\varepsilon^r u\dx \dy =
\int_\Sigma \Lambda_\varepsilon^r G\cdot \nabla \Lambda_\varepsilon^r u\dx \dy 
+\int_{-\infty}^\infty \Lambda_\varepsilon^r \xi \Lambda_\varepsilon^r u|_{y=1}\dx
\]
because
$\Lambda_\varepsilon^r$ commutes with partial derivatives and is symmetric with respect to the 
$L^2$-inner product. This equation can be rewritten as
\begin{align*}
\int_\Sigma P\nabla \Lambda_\varepsilon^r u \cdot \nabla \Lambda_\varepsilon^r u\dx \dy =&
-\int_\Sigma [\Lambda_\varepsilon^r, Q]\nabla u \cdot \nabla \Lambda_\varepsilon^r u\dx \dy 
+
\int_\Sigma \Lambda_\varepsilon^r G\cdot \nabla \Lambda_\varepsilon^r u\dx \dy \\
& 
-\int_{-\infty}^\infty \Lambda_\varepsilon^r \zeta (\Lambda_\varepsilon^r  u|_{y=1})_x\dx,
\end{align*}
and it follows from the coercivity of $P$ and the continuity of the trace map $H^1_\star(\Sigma)\rightarrow H_\star^{1/2}(\R)$ that
\begin{align*}
\|\Lambda_\varepsilon^r\nabla u\|_{L^2(\Sigma)}
&\le 
C(\| [\Lambda_\varepsilon^r, Q]\nabla u\|_{L^2(\Sigma)} 
+
\| \Lambda_\varepsilon^r G\|_{L^2(\Sigma)}
+\|\Lambda_\varepsilon^r\Lambda^\frac{1}{2} \zeta \|_{L^2(\R)})\\
&\le C(\| [\Lambda_\varepsilon^r, Q]\nabla u\|_{L^2(\Sigma)} 
+ \|G\|_{t+1,1} + \| \zeta \|_{t+3/2}).
\end{align*}

The next step is to estimate the commutator $[\Lambda_\varepsilon^r, Q]$.
For $r\le s+\frac{1}{2}$ we choose\linebreak $\tilde{\Delta} \in (0,\min(s,1))$ and estimate
\begin{eqnarray*}
 \| [\Lambda_\varepsilon^r, Q]\nabla u\|_{L^2(\Sigma)}
 & \leq & C\|Q\|_{s+1/2,\infty} \|\nabla u\|_{r-\tilde{\Delta},0} \\
 & \leq & C\|Q\|_{s+1,1} \|\nabla u\|_{r-\tilde{\Delta},0}
 \end{eqnarray*}
using Lemma \ref{commutator estimate} (with $r_0=s+\frac{1}{2}-\tilde{\Delta}$, $\Delta=\tilde{\Delta}$).
In the case $r \in (s+\frac{1}{2},s+1]$ on the other hand, we choose
$\tilde{\Delta} \in (0,\min(s,\frac{1}{2}))$ and estimate
\begin{align*}
 \| [\Lambda_\varepsilon^r, Q]\nabla u\|_{L^2(\Sigma)}
&\le
    C \|Q\|_{s+1,0} \|\nabla u\|_{r-\tilde{\Delta}-1/2,\infty}\\
&\le 
C\|Q\|_{s+1,0} \|\nabla u\|_{r-\tilde{\Delta},1}
\end{align*}
using Lemma \ref{commutator estimate} (with $r_0=s+\frac{1}{2}-\tilde{\Delta}$
and $\Delta=\tilde{\Delta}+\frac{1}{2}$) and
\[
\|\nabla u\|_{r-\tilde{\Delta},1}\le C(\|G\|_{t+1,1}
+\|\nabla u\|_{r-\tilde{\Delta},0})
\]
using Lemma \ref{regularity lemma}.

Combing the above estimates yields
$$
\|\Lambda_\varepsilon^r\nabla u\|_{L^2(\Sigma)}
\leq
C(\|\nabla u\|_{r-\tilde{\Delta},0}
+ \|G\|_{t+1,1} + \| \zeta \|_{t+3/2}),
$$
where $\tilde{\Delta} \in (0,\min(s,\frac{1}{2}))$, and letting $\varepsilon \rightarrow 0$
and using the resulting estimate iteratively, we find that
$$
\|\nabla u\|_{t+1,0}
\leq
C(\|G\|_{t+1,1}
+\| \zeta \|_{t+3/2} + \|u\|_{H_\star^1(\Sigma)}),
$$
from which the result follows by Lemma \ref{Flattened generalised BVP has a WS} and Lemma \ref{regularity lemma}.\qed

The following result shows that Lemma \ref{Regularity of u} is applicable to the boundary-value problems
\eqn{BC for u0 1}--\eqn{BC for u0 3} and \eqn{BC for un 1}--\eqn{BC for un 3}.

\begin{lemma}\label{Q is analytic 2}
The mapping $W^{s+3/2} \rightarrow(H^{s+1,2})^{2 \times 2}$ given by $\eta \mapsto Q(\eta)$ is analytic.
\end{lemma}
\begin{remark} \label{Formulae for Qx and Qy}
Observe that
\begin{eqnarray*}
Q_x(\eta) & = & S_0(\eta) + R_0(\eta)L_0^\delta \eta^{\prime\prime} + R_1(\eta)L_1^\delta \eta^{\prime\prime}, \\
Q_y(\eta) & = & T_0(\eta) + R_0(\eta)L_1^\delta \eta^{\prime\prime} + R_1(\eta)L_2^\delta \eta^{\prime\prime},
\end{eqnarray*}
where
$$
L_j^\delta(\cdot) = \FF^{-1}[(\i \delta)^j \chi^{(j)} ((1-y)\delta k) \FF[\cdot]], \qquad j=0,1,2,
$$
are bounded bilinear functions $L^2(\R) \rightarrow L^\infty H^0$
and
\begin{eqnarray*}
& & S_0:\eta \rightarrow \begin{pmatrix} \eta_x^\delta & 0 \\
0 & \displaystyle{-\frac{\eta_x^\delta}{1+f_y^\delta} - \frac{(-f_y^\delta+(f_x^\delta)^2)\eta_x^\delta}{(1+f_y^\delta)^2}}\end{pmatrix},\\
& & T_0:\eta \rightarrow \begin{pmatrix} 2L_1^\delta \eta^\prime & -\eta_x^\delta \\
-\eta_x^\delta & \displaystyle{-\frac{2L_1^\delta \eta^\prime}{1+f_y^\delta} +\frac{2f_x^\delta \eta_x^\delta}{1+f_y^\delta}
- \frac{2(-f_y^\delta+(f_x^\delta)^2)L_1^\delta \eta^\prime}{(1+f_y^\delta)^2}}\end{pmatrix},\\
& & R_0:\eta \rightarrow \begin{pmatrix} 0 & -y \\
-y & \displaystyle{\frac{2yf_x^\delta}{1+f_y^\delta}}\end{pmatrix}, \\
& & R_1:\eta \rightarrow \begin{pmatrix} y & 0 \\
0 & \displaystyle{-\frac{y}{1+f_y^\delta} - \frac{y(-f_y^\delta+(f_x^\delta)^2)}{(1+f_y^\delta)^2}}\end{pmatrix}
\end{eqnarray*}
are analytic functions $W \rightarrow (L^\infty(\overline{\Sigma}))^{2\times 2}$.
\end{remark}
The regularity assertion in Theorem \ref{Analyticity of K 2} now
follows from the next result and the continuity of the trace operator
$H^{s+1,1} \rightarrow H^{s+1/2}(\R)$. 

\begin{theorem}
The mapping $W^{s+3/2} \rightarrow
\LL(H^{s+3/2}(\R), (H^{s+1,1})^2)$ given by
$\eta \mapsto (\zeta \mapsto \nabla u)$,
where $u \in H_\star^1(\Sigma)$ is the unique weak solution of \eqn{BC for u 1}--\eqn{BC for u 3}
with $\xi=\zeta^\prime$, is analytic.
\end{theorem}
{\bf Proof.} Repeating the proof of Theorem \ref{Weak solns analytic 1}, replacing
Lemma \ref{Flattened generalised BVP has a WS} by Lemma \ref{Regularity of u},
Lemma \ref{Q is analytic 1} by Lemma \ref{Q is analytic 2} and inequality \eqn{Expansion of Gn} by
$$
\|G^n\|_{s+1,1} \leq \sum_{k=1}^n \|Q^k\|_{s+1,1} \|\nabla u^{n-k}\|_{s+1,1}
$$
($H^{s+1,1}$ is a Banach algebra), we obtain the representation
$$
\nabla u(x,y) = \sum_{n=0}^\infty \nabla u^n(x,y), \qquad \nabla u^n = m^n_2(\{\tilde{\eta}\}^{(n)})
$$
where
$m^n_2 \in \LL_\mathrm{s}^n(H^{s+3/2}(\R),(H^{s+1,1})^2)$ is linear in $\zeta$ and satisfies
$$\nn m^n_2 \nn \leq C_1B^n \|\zeta\|_{s+3/2}$$
for some constant $B>0$.\qed

We conclude this section with a useful supplementary estimate for $\|K^n(\tilde{\eta})\|$.

\begin{proposition} \label{Analyticity of K 3}
There exists a constant $B>0$ such that
$$\|K^n(\tilde{\eta})\zeta\|_0 \leq C_1 B^n (\|\tilde{\eta}\|_{1,\infty} + \|\tilde{\eta}^{\prime\prime}+k_0^2\tilde{\eta}\|_0)^n \|\zeta\|_{3/2}, \qquad n=0,1,2,\ldots.$$
\end{proposition}
{\bf Proof.} It suffices to establish the estimate
$$\|\nabla u^n\|_1 \leq C_1  B^n(\|\tilde{\eta}\|_{1,\infty} + \|\tilde{\eta}^{\prime\prime}+k_0^2\tilde{\eta}\|_0)^n \|\zeta\|_{3/2}, \qquad n=0,1,2,\ldots;$$
for $n=0$ this result follows from Lemma \ref{Regularity of u} (with $t=0$ and $s=\frac{1}{2}$).

Proceeding inductively, suppose the estimate for $\|\nabla u^k\|_1$ is true for all $k<n$, and recall from the proof of
Theorem \ref{Weak solns analytic 1} that
$$\|Q^k\|_\infty \leq C_2 r^{-k} \|\tilde{\eta}\|_{1,\infty}^k,
\qquad
\|G^n\|_0 \leq C_1C_2 B^n \|\zeta\|_{3/2} \|\tilde{\eta}\|_{1,\infty}^n \sum_{k=1}^n (B r)^{-k}.$$
Writing
\begin{eqnarray*}
Q_x^k & = & S_0^k + R_0^k L_0^\delta \eta_0^{\prime\prime}  + R_0^{k-1}L_0^\delta \tilde{\eta}^{\prime\prime}
+ R_1^k L_0^\delta \eta_0^{\prime\prime}  + R_1^{k-1}L_0^\delta \tilde{\eta}^{\prime\prime} \\
& = & S_0^k + \sum_{j=0}^1 \left(
- k_0^2 R_j^{k-1} L_j^\delta \tilde{\eta} + R_j^k  L_j^\delta \eta_0^{\prime\prime}
+ R_j^{k-1}L_j^\delta(\tilde{\eta}^{\prime\prime} + k_0^2 \tilde{\eta}) \right),
\end{eqnarray*}
where
$$\|S_0^k\|_\infty \leq C_2 r^{-k} \|\tilde{\eta}\|_{1,\infty}^k, \qquad
\|R_j^k\|_\infty \leq C_2 r^{-k} \|\tilde{\eta}\|_{1,\infty}^k, \quad j=0,1,$$
(see Remark \ref{Formulae for Qx and Qy}), we find that
\begin{eqnarray*}
G^n_x & = & \!\!\!\!- \sum_{k=1}^n \big(Q_x^k \nabla u^{n-k} + Q^k \nabla u_x^{n-k}\big) \\
& = & \!\!\!\!\sum_{k=1}^n \!\left(
\!S_0^k \nabla u^{n-k}\!+\! \sum_{j=0}^1\! \left(
- k_0^2 R_j^{k-1} L_j^\delta \tilde{\eta} + R_j^k  L_j^\delta \eta_0^{\prime\prime}
+ R_j^{k-1}L_j^\delta(\tilde{\eta}^{\prime\prime} + k_0^2 \tilde{\eta}) \right)\! \nabla u^{n-k} \!+\! Q^k \nabla u_x^{n-k}\!\right)\!\!.
\end{eqnarray*}
It follows that
\begin{eqnarray*}
\|G^n_x\|_0 & \leq & \sum_{k=1}^n \Big((\|S_0^k\|_\infty + k_0^2 (\|R_0^{k-1}\|_\infty +\|R_1^{k-1}\|_\infty)
\|\tilde{\eta}\|_{\infty})\|\nabla u^{n-k}\|_0\\[-4mm]
& & \hspace{1.5cm}\mbox{} 
+ \vphantom{\Big(}(\|R_0^k\|_\infty\|L_0^\delta\| +\|R_1^k\|_\infty\|\|L_1^\delta\|) \|\eta_0^{\prime\prime}\|_0 \|\nabla u^{n-k}\|_1 \\
& & \hspace{1.5cm}\mbox{} + (\|R_0^{k-1}\|_\infty \|L_0^\delta\|+\|R_1^{k-1}\|_\infty\|L_1^\delta\|) \|\tilde{\eta}^{\prime\prime} + k_0^2 \tilde{\eta}\|_0 \|\nabla u^{n-k} \|_1 \\
& & \hspace{1.5cm}\mbox{} 
+ \|Q^k\|_\infty \|\nabla u_x^{n-k}\|_0 \Big) \\
& \leq & C_1 C_2 B^n \big(1+2k_0^2 r +( \|L_0^\delta\|+\|L_1^\delta\|)(\|\eta_0^{\prime\prime}\|_0+r)
+ 1\big) \\
& & \qquad \mbox{}\times
\|\zeta\|_{3/2}(\|\tilde{\eta}\|_{1,\infty} + \|\tilde{\eta}^{\prime\prime}+k_0^2\tilde{\eta}\|_0)^n \sum_{k=1}^n (B r)^{-k},
\end{eqnarray*}
in which Proposition \ref{Bprop} has been used.
A similar calculation yields the same estimate for $\|G_y^n\|_0$.

Combining the estimates for $\|G^n\|_0$, $\|G_x^n\|_0$ and $\|G_y^n\|_0$ and applying Lemma \ref{Regularity of u}
(with $t=0$ and $s=\frac{1}{2}$), one finds that
\begin{eqnarray*}
\|\nabla u^n\|_1 & \leq & \sqrt{3}C_1 C_2 C_3 B^n 
\big(1+2k_0^2 r +( \|L_0^\delta\|+\|L_1^\delta\|)(\|\eta_0^{\prime\prime}\|_0+r)
+ 1\big) \\
& & \qquad \mbox{}\times
\|\zeta\|_{3/2}(\|\tilde{\eta}\|_{1,\infty} + \|\tilde{\eta}^{\prime\prime}+k_0^2\tilde{\eta}\|_0)^n \sum_{k=1}^n (B r)^{-k},
\end{eqnarray*}
so that
$$
\|\nabla u^n\|_1 \leq C_1 B^n (\|\tilde{\eta}\|_{1,\infty} + \|\tilde{\eta}^{\prime\prime}+k_0^2\tilde{\eta}\|_0)^n \|\zeta\|_{3/2}
$$
for sufficiently large values of $B$ (independently of $n$).\qed

\subsection{Variational functionals} \label{Properties of J, K, L}

In this section we study the functional
\begin{equation}
\TT(\eta) = \int_{-\infty}^\infty f_1(\eta) K(\eta) f_2(\eta) \dx, \label{Definition of T}
\end{equation}
where $f_1, f_2: \R \rightarrow \R$ are polynomials with $f_1(0)=f_2(0)=0$,
and apply our results to the functionals $\GG$, $\KK$ and $\LL$.

\subsubsection{Analyticity of the functionals}

In this section we again suppose that $s>0$.
The first result follows from Theorem \ref{Analyticity of K 1}(i).

\begin{lemma} \label{Functionals are analytic}
Equation \eqn{Definition of T} defines a functional
$\TT: W^{s+3/2} \rightarrow \R$
which is analytic and satisfies $\TT(0)=0$.
\end{lemma}

We now turn to the construction of the gradient $\TT^\prime(\eta)$ in $L^2(\R)$,
the main step of which is accomplished by the following lemma.

\begin{lemma} \label{Expression for HHprime}
Define $\HH: W^{s+3/2} \rightarrow \LL_s^2(H^{s+3/2}(\R), \R)$
by the formula
$$\HH(\eta)(\zeta_1,\zeta_2)=\langle \zeta_1, K(\eta) \zeta_2 \rangle_0.$$
The gradient $\HH^\prime(\eta)(\zeta_1,\zeta_2)$ in $L^2(\R)$ exists for
each $\eta \in W^{s+3/2}$ and $\zeta_1$, $\zeta_2 \in H^{s+3/2}(\R)$ and is given by the formula
$$
\HH^\prime(\eta)(\zeta_1,\zeta_2) = 
-u_{1x}u_{2x}+\frac{1+\eta^{\prime 2}}{(1+\eta)^2} u_{1y }u_{2y}\Bigg|_{y=1},
\label{Formula for Hprime}
$$
where $u_j$ is the weak solution of \eqn{BC for u 1}--\eqn{BC for u 3}
with $\xi=\zeta_j^\prime$, $j=1,2$. This formula defines an analytic function
$\HH^\prime: W^{s+3/2} \rightarrow \LL_\mathrm{s}^2(H^{s+3/2}(\R), H^{s+1/2}(\R))$.
\end{lemma}
{\bf Proof.} 
It follows from the formula
$$
\HH(\eta)
=\int_{\Sigma} (I+Q(\eta))\nabla u_1 \cdot \nabla u_2 \dx \dy
$$
that
\begin{eqnarray}
\lefteqn{\mathrm{d}\HH[\eta](\omega) = 
\int_{\Sigma} \mathrm{d}Q[\eta](\omega) \nabla u_1 \cdot \nabla u_2\dx \dy} \nonumber \\
& & \hspace{0.75in}\mbox{}
+\int_{\Sigma} (I+Q(\eta)) \nabla w_1 \cdot \nabla u_2\dx \dy
+\int_{\Sigma} (I+Q(\eta)) \nabla u_1 \cdot \nabla w_2\dx \dy,
\label{Basic formula for dH} 
\end{eqnarray}
where $w_j=\mathrm{d} u_j(\eta)[\omega]$, $j=1,2$.
Recall that
$$\int_\Sigma (I+Q(\eta))\nabla u_j\cdot \nabla v \dx\dy = \int_{-\infty}^\infty \zeta_j^\prime v|_{y=1}\dx, \qquad j=1,2,$$
for every $v \in H_\star^1(\Sigma)$ (Definition \ref{Weak soln} with $\xi=\zeta_j^\prime$ and $G=0$), so that
\begin{equation}
 \int_\Sigma \Bigg\{ \mathrm{d}Q[\eta](\omega) \nabla u_j\cdot \nabla v+ (I+Q(\eta))\nabla w_j\cdot \nabla v \Bigg\} \dx\dy=0
\qquad j=1,2, \label{Weak 1}
\end{equation}
for every $v \in H_\star^1(\Sigma)$.
Subtracting \eqn{Weak 1} with $j=1$, $v=u_2$ and $j=2$, $v=u_1$ from
\eqn{Basic formula for dH} yields
$$
\mathrm{d}\HH[\eta](\omega)= -\int_\Sigma \mathrm{d}Q[\eta](\omega)\nabla u_1\cdot\nabla u_2 \dx\dy.
$$

Finally, write $h^\delta(x,y)=y\omega^\delta(x,y)$, where
$\omega^\delta(x, y)=\FF^{-1}[\chi(\delta (y-1)|k|)\hat \omega(k)](x)$, so that
$h^\delta=\mathrm{d}f^\delta[\eta](\omega)$, and observe that
\begin{eqnarray*}
\lefteqn{\int_{-\infty}^\infty\Bigg\{-h^\delta\left(u_{1x} - \frac{f^\delta_x y u_{1y}}{1+f^\delta_y}\right)\!\!\!\left(u_{2x} - \frac{f^\delta_x y u_{2y}}{1+f^\delta_y}\right)
+ \frac{h^\delta u_{1y}u_{2y}}{(1+f^\delta_y)^2} \Bigg\}\Bigg|_{y=1}\dx}\qquad \\
& = & \frac{1}{2}\int_\Sigma\frac{\mathrm{d}}{\mathrm{d}y}\Bigg\{
-h^\delta\left(u_{1x} - \frac{f^\delta_x y u_{1y}}{1+f^\delta_y}\right)\!\!\!\left(u_{2x} - \frac{f^\delta_x y u_{2y}}{1+f^\delta_y}\right)
+ \frac{h^\delta u_{1y}u_{2y}}{(1+f^\delta_y)^2}
\Bigg\}\dx\dy \nonumber \\
& = & \int_\Sigma \left\{ -h^\delta_xu_{1x}u_{2x} +h^\delta_x u_{1x}u_{2y} + h^\delta_x u_{1y}u_{2x} + \frac{h^\delta_y u_{1y}u_{2y}}{(1+f^\delta_y)^2} \right.\\
& & \qquad\qquad\mbox{}
+ \left.\frac{2(f^\delta_x)^2h^\delta_y u_{1y}u_{2y}}{(1+f^\delta_y)^2} - \frac{2f^\delta_xh^\delta_x u_{1y}u_{2y}}{1+f^\delta_y}\right\}\dx\dy\\
& & + \int_\Sigma \frac{h^\delta u_{1y}}{1+f^\delta_y}\left\{((1+f^\delta_y)u_{2x}-f^\delta_xu_{2y})_x +\left (-f^\delta_xu_{2x}+\frac{1+(f^\delta_x)^2}{1+f^\delta_y}u_{2y}\right)_y
 \right\} \dx\dy\\
 & & + \int_\Sigma \frac{h^\delta u_{2y}}{1+f^\delta_y}\left\{((1+f^\delta_y)u_{1x}-f^\delta_xu_{1y})_x +\left (-f^\delta_xu_{1x}+\frac{1+(f^\delta_x)^2}{1+f^\delta_y}u_{1y}\right)_y
 \right\} \dx\dy\\
& & +\int_{-\infty}^\infty\Bigg\{\frac{h^\delta f^\delta_xu_{1y}}{1+f^\delta_y}\left(u_{2x} - \frac{f^\delta_x y u_{2y}}{1+f^\delta_y}\right)
+\frac{h^\delta f^\delta_xu_{2y}}{1+f^\delta_y}\left(u_{1x} - \frac{f^\delta_x y u_{1y}}{1+f^\delta_y}\right)
\Bigg\}\Bigg|_{y=1}\dx
\end{eqnarray*}
\begin{eqnarray*}
& = & -\int_{\Sigma} \Bigg\{\mathrm{d}Q[\eta](\omega) \nabla u_1\cdot \nabla u_2 \\
& & \qquad\qquad\mbox{}
+ \frac{h^\delta u_{1y}}{1+f^\delta_y} \underbrace{\nabla \cdot ((I+Q(\eta))\nabla u_1)}_{\displaystyle =0}
+ \frac{h^\delta u_{2y}}{1+f^\delta_y} \underbrace{\nabla \cdot ((I+Q(\eta))\nabla u_2)}_{\displaystyle =0} \Bigg\}\dx\dy \\
& & +\int_{-\infty}^\infty\Bigg\{\frac{h^\delta f^\delta_xu_{1y}}{1+f^\delta_y}\left(u_{2x} - \frac{f^\delta_x y u_{2y}}{1+f^\delta_y}\right)
+\frac{h^\delta f^\delta_xu_{2y}}{1+f^\delta_y}\left(u_{1x} - \frac{f^\delta_x y u_{1y}}{1+f^\delta_y}\right)
\Bigg\}\Bigg|_{y=1}\dx,
\end{eqnarray*}
in which the third line follows from the second by differentiating the term in braces with respect to $y$ (note that $h^\delta|_{y=0}=0$)
and
integrating by parts. One concludes that
$$
\mathrm{d}\HH[\eta](\omega) = \int_{-\infty}^\infty \left\{ -u_{1x}u_{2x} + \frac{1+(f^\delta_x)^2}{(1+f^\delta_y)^2}u_{1y}u_{2y}\right\}h^\delta\Bigg|_{y=1} \dx,
$$
and the stated formula follows from this result and the facts that $f^\delta|_{y=1}=\eta$ and $h^\delta|_{y=1}=\omega$.

The hypotheses of the lemma imply that $\nabla u_j\in H^{s+1,1}$ and $\nabla u_j|_{u=1} \in H^{s+1/2}(\R)$, $j=1,2$.
This observation ensures that the above algebraic manipulations are valid and that $\mathrm{d}\HH[\eta]$ belongs to $H^{s+1/2}(\R)$
because $H^{s+1,1}$ and $H^{s+1/2}(\R)$ are Banach algebras.\qed

\begin{corollary} \label{Gradient of T}
The gradient $\TT^\prime(\eta)$ in $L^2(\R)$ exists for each $\eta \in W^{s+3/2}$ and is given by the formula
$$\TT^\prime(\eta) = \HH^\prime(\eta)(f_1(\eta),f_2(\eta))+f_1^\prime(\eta)K(\eta)f_2(\eta)+f_2^\prime(\eta)K(\eta)f_1(\eta).$$
This formula defines an analytic function $\TT^\prime: W^{s+3/2} \rightarrow H^{s+1/2}(\R)$ which satisfies $\TT^\prime(0)=0$.
\end{corollary}

\begin{theorem} \label{Main result for functionals}
\quad
\begin{list}{(\roman{count})}{\usecounter{count}}
\item 
Equations \eqn{Definition of GG}--\eqn{Definition of LL} define analytic functionals
$\GG,\KK,\LL: W^{s+3/2} \rightarrow \R$
which satisfy $\GG(0),\KK(0),\LL(0)=0$.
\item
Equation \eqn{Definition of J} defines an analytic functional $\JJ_\mu: W^{s+3/2}\sm\{0\} \rightarrow \R$.
\item
The gradients $\GG^\prime(\eta)$ and $\LL^\prime(\eta)$ in $L^2(\R)$ exist for each $\eta \in W^{s+3/2}$
and are given by the formulae
\begin{eqnarray}
\GG^\prime(\eta) & = & \frac{\omega}{4}\HH^\prime(\eta)(\eta^2,\eta)
+\frac{\omega}{4}K(\eta)\eta^2 + \frac{\omega}{2}\eta K(\eta)\eta - \frac{\omega}{2}\eta, \label{Definition of Gprime} \\
\LL^\prime(\eta) & = & \frac{1}{2}\HH^\prime(\eta)(\eta,\eta) + K(\eta)\eta. \label{Definition of Lprime}
\end{eqnarray}
These formulae define analytic functions $\GG^\prime, \LL^\prime: W^{s+3/2} \rightarrow H^{s+1/2}(\R)$ which
satisfy $\GG^\prime(0)=0$ and $\LL^\prime(0)=0$.
\item
The gradient $\KK^\prime(\eta)$ in $L^2(\R)$ exists for each $\eta \in W^2$ and is given by the formula
\begin{equation}
\KK^\prime(\eta) = \eta - \beta\left(\frac{\eta^\prime}{\sqrt{1+\eta^{\prime 2}}}\right)^{\!\!\prime}
- \frac{\omega^2}{8}\HH^\prime(\eta)(\eta^2,\eta^2) -\frac{\omega^2}{2}\eta^2K(\eta)\eta + \frac{\omega^2}{3}\eta^2.
 \label{Definition of Kprime} \\
\end{equation}
This formula defines an analytic function $\KK^\prime: W^2 \rightarrow L^2(\R)$ which
satisfies $\KK^\prime(0)=0$.
\item
The gradient $\JJ_\mu^\prime(\eta)$ in $L^2(\R)$ exists for each $\eta \in W^2 \sm \{0\}$ and defines an analytic
function $\JJ_\mu^\prime: W^2 \sm\{0\} \rightarrow L^2(\R)$.
\end{list}
\end{theorem}

\begin{corollary}
Choose $M>0$ so that $\overline{B}_M(0) \subseteq H^2(\R)$ lies in $W^{s+3/2}$ and define $U=B_M(0)$. Equations
\eqn{Definition of GG}--\eqn{Definition of LL} define analytic functionals
$\GG, \KK, \LL: U \rightarrow \R$ while equations \eqn{Definition of Gprime}--\eqn{Definition of Kprime} define
analytic functions $\GG^\prime, \KK^\prime, \LL^\prime: U \rightarrow L^2(\R)$.
\end{corollary}

Finally, we state some further useful estimates for the operators $\GG$, $\KK$ and $\LL$. Here, and in the
remainder of this paper, the constant $M$ is chosen small enough
for the validity of our calculations.

\begin{proposition} \label{Quadratic estimates}
The estimates
$$|\GG(\eta)| \leq c\|\eta\|_{1/2}^2, \qquad
\KK(\eta) \geq c\|\eta\|_1^2, \qquad
c\|\eta\|_{1/2}^2 \leq \LL(\eta) \leq c\|\eta\|_{1/2}^2
$$
hold for each $\eta \in U$.
\end{proposition}
{\bf Proof.} The estimate for $\GG$ follows from the calculation
$$|\GG(\eta)| \ \leq\  c(\|\eta\|_0^2 \|K(\eta)\eta\|_0 + \|\eta\|_0^2)
\ \leq\  c\|\eta\|_{1/2}^2,$$
while that for $\LL$ is a direct consequence of Theorem \ref{Analyticity of K 1}(ii).
Turning to the estimate for $\KK$, observe that
$$
\KK(\eta) = \underbrace{\int_{-\infty}^\infty  \left\{\frac{\beta\eta^{\prime 2}}{1+\sqrt{1+\eta^{\prime 2}}}+\frac{\eta^2}{2}\right\}\dx}_{\displaystyle \geq c\|\eta\|_1^2} - \frac{\omega^2}{8} \int_{-\infty}^\infty \eta^2K(\eta)\eta^2 \dx
+ \frac{\omega^2}{6}\int_{-\infty}^\infty \eta^3 \dx
$$
and
$$
\left|\int_{-\infty}^\infty \eta^3 \dx\right| \leq c \|\eta\|_1^3,\qquad
\left|\int_{-\infty}^\infty \eta^2K(\eta)\eta^2 \dx\right|\ \leq\  c\|\eta^2\|_{1/2}^2 \ \leq\  c\|\eta\|_1^4,
$$
for each $\eta \in U$, so that $\KK(\eta) \geq c\|\eta\|_1^2$.\qed

\subsubsection{Pseudo-local properties of the operator ${\mathcal T}$} \label{Pseudo-local properties}

In this section we consider sequences $\{\eta_m^{(1)}\}$, $\{\eta_m^{(2)}\} \subset U$ with the properties that\linebreak
$\supp \eta_m^{(1)} \subset [-R_m,R_m]$, $\supp \eta_m^{(2)} \subset \R \sm (-S_m,S_m)$
and $\sup_{m \in {\mathbb N}}\|\eta_m^{(1)}+\eta_m^{(1)}\|_2 < M$, where $\{R_m\}$, $\{S_m\}$ are sequences of positive real
numbers with $R_m$, $S_m \rightarrow \infty$, $R_m/S_m \rightarrow 0$ as $m \rightarrow \infty$.
We establish the following `pseudo-local' property of the operator $\TT$.

\begin{theorem} \label{Main splitting theorem}
The operator ${\mathcal T}$ satisfies
$$\lim_{m \rightarrow \infty} \Big(\TT(\eta_m^{(1)}+\eta_m^{(2)}) - \TT(\eta_m^{(1)})-\TT(\eta_m^{(2)}) \Big)=0,$$
$$\lim_{m \rightarrow \infty} \|\TT^\prime(\eta_m^{(1)}+\eta_m^{(2)}) - \TT^\prime(\eta_m^{(1)})-\TT^\prime(\eta_m^{(2)}) \|_0=0,$$
$$\lim_{m \rightarrow \infty} \langle \TT^\prime(\eta_m^{(2)}),\eta_m^{(1)} \rangle_0 =0.$$
In particular, this result applies to $\GG$, $\KK$ and $\LL$.
\end{theorem}

We begin the proof of Theorem \ref{Main splitting theorem} by re-examining the general boundary-value problem
\eqn{BC for u 4}--\eqn{BC for u 6}.

\begin{lemma} \label{Key splitting lemma}
Suppose that $\{R_m\}$, $\{S_m\}$ and $\{U_m\}$ are sequences of positive real numbers and
$\{Q_m\} \subseteq (L^\infty(\overline{\Sigma}))^{2 \times 2}$, $\{G_m\} \subseteq L^2(\Sigma)$,
$\{\zeta_m^{(1)}\}$, $\{\zeta_m^{(2)}\} \subseteq H^{1/2}(\R)$ are bounded sequences with the properties that
\begin{list}{(\roman{count})}{\usecounter{count}}
\item
$S_m-U_m$, $U_m - R_m \rightarrow \infty$ as $m \rightarrow \infty$;
\item
$\supp \zeta_m^{(1)}  \subset [-R_m,R_m]$ and $\supp \zeta_m^{(2)}\subset \R \sm (-S_m,S_m)$;
\item
$\|G_m^{(1)}\|_{L^2(|x| > R_m)},\ \|G_m^{(2)}\|_{L^2(|x| <S_m)} \rightarrow 0$ as $m \rightarrow \infty$;
\item
there exists a constant $p_0>0$ such that
\[
(I+Q_m)(x, y)\nu\cdot \nu \ge p_0|\nu|^2
\]
for all $(x,y)\in \overline \Sigma$, all $m \in {\mathbb N}$ and all $\nu \in \R^2$.
\end{list}
The unique weak solutions $u_m^{(j)} \in H^1_\star(\Sigma)$ of the boundary-value problems
\begin{eqnarray}
\label{eq:split 1}
& & \parbox{90mm}{$\nabla\cdot((I+Q_m)\nabla u_m^{(j)})=\nabla\cdot G_m^{(j)},$}0 < y <1,\\
\label{eq:split 2}
& & \parbox{90mm}{$(I+Q_m)\nabla u_m^{(j)} \cdot (0, 1) = \zeta_{m,x}^{(j)}+G_m^{(j)}\cdot(0,1),$}y=1, \\
\label{eq:split 3}
& & \parbox{90mm}{$(I+Q_m)\nabla u_m^{(j)} \cdot (0,-1)=G_m^{(j)}\cdot(0,-1),$}y=0,
\end{eqnarray}
$j=1$, $2$, satisfy the estimates
$$
\lim_{m \rightarrow \infty} \|\nabla u_{m}^{(1)}\|_{L^2(|x| >U_m)} = 0, \qquad
\lim_{m \rightarrow \infty} \|\nabla u_{m}^{(2)}\|_{L^2(|x| <U_m)} = 0.
$$
\end{lemma}
{\bf Proof.} Write $\zeta_m^{(2)}=\zeta_{m,+}^{(2)}+\zeta_{m,-}^{(2)}$, where
$\supp \zeta_{m,+}^{(2)}\subseteq [S_m,\infty)$
and $\supp \zeta_{m,-}^{(2)}\subseteq (-\infty, -S_m]$, and let $u_{m,+}^{(2)}$, $u_{m,-}^{(2)}$ be the weak solutions
of the boundary-value problem \eqn{eq:split 1}--\eqn{eq:split 2} with\linebreak $\zeta_m^{(2)}$, $G_m^{(2)}$ replaced by
respectively $\zeta_{m,+}^{(2)}$, $G_{m,+}^{(2)}:=G_m^{(2)}\chi_{\{x>0\}}$ and $\zeta_{m,-}^{(2)}$, $G_{m,-}^{(2)}:=G_m^{(2)}\chi_{\{x<0\}}$,
so that $u_m^{(2)}=u_{m,+}^{(2)}+u_{m,-}^{(2)}$.

Choose $T>0$ and take $m$ large enough so that $T+1<S_m$. Define  $\phi \in C^\infty(\R)$ by the formula
$$\phi_T(x) =\left\{\begin{array}{lll}1, & \qquad x \leq T, \\[2mm]
\chi(2(x-T)), & \qquad x > T
\end{array}\right.$$
and set
\[
w_m(x,y)=\phi_T^2(x) (u_{m,+}^{(2)}(x,y)-M_T),
\]
where 
\[
M_T=\int_{T\le x\le T+1} u_{m,+}^{(2)} (x,y)\dx \dy,
\]
so that $\supp w_m \subseteq (-\infty, T+1] \times [0,1]$ and the mean value of $u_{m,+}^{(2)}(x,y)-M_T$ over\linebreak
 $(T,T+1)\times(0,1)$ is zero.
Using Definition \ref{Weak soln}, we find that
$$
\int_{\Sigma} (I+Q_m) \nabla u_{m,+}^{(2)}\cdot \nabla w_m\dx\dy = \int_\Sigma G_{m,+}^{(2)} \cdot \nabla w_m\dx\dy
+\underbrace{\int_{-\infty}^\infty \partial_x \zeta_{m,+}^{(2)} w_m|_{y=1} \dx}_{\displaystyle =0},
$$
from which it follows that
\begin{eqnarray*}
\lefteqn{\int_{\Sigma} (I+Q_m) \phi_T^2|\nabla u_{m,+}^{(2)}|^2 \dx\dy} \\
& & \leq c \left( \left(\int_\Sigma \phi_T^2|\nabla u_{m,+}^{(2)}|^2 \dx\dy\right)^{\!\!\frac{1}{2}}\!\!
\left(\int_{T \leq x \leq T+1} |u_{m,+}^{(2)}-M_T|^2 \dx\dy\right)^{\!\!\frac{1}{2}}\right. \\
& & \qquad\quad\mbox{}+ \left(\int_{x \leq T+1} |G_{m,+}^{(2)}|^2 \dx\dy\right)^{\!\!\frac{1}{2}}\!\!
\left(\int_{T \leq x \leq T+1} |u_{m,+}^{(2)}-M_T|^2 \dx\dy\right)^{\!\!\frac{1}{2}} \\
& & \qquad\quad\mbox{}+ \left.\left(\int_{x \leq T+1} |G_{m,+}^{(2)}|^2 \dx\dy\right)^{\!\!\frac{1}{2}}\!\!
\left(\int_\Sigma \phi_T^2|\nabla u_{m,+}^{(2)}|^2 \dx\dy\right)^{\!\!\frac{1}{2}}\right)
\end{eqnarray*}
and hence that
$$
\int_{\Sigma} \phi_T^2|\nabla u_{m,+}^{(2)}|^2 \dx\dy
\leq c \left( \int_{T \leq x \leq T+1} |\nabla u_{m,+}^{(2)}|^2 \dx\dy
+\int_{x \leq T+1} |G_{m,+}^{(2)}|^2 \dx\dy\right),
$$
where the Poincar\'{e} inequality
$$
\int_{T \leq x \leq T+1} |u_{m,+}^{(2)}-M_T|^2 \dx\dy  \leq c \int_{T \leq x \leq T+1} |\nabla u_{m,+}^{(2)}|^2 \dx\dy$$
has been used.

The above inequality implies that
$$\Phi(T) \leq c_\star \big(\Phi(T+1) - \Phi(T) + \Psi(T+1)\big),$$
for some $c_\star>0$, where
$$\Phi(T) = \int_{x \leq T} |\nabla u_{m,+}^{(2)}|^2 \dx\dy, \qquad \Psi(T) = \int_{x \leq T} |G_{m,+}^{(2)}|^2 \dx\dy,$$
so that
$$\Phi(T) \leq d_\star \big(\Phi(T+1) + \Psi(T+1)\big),$$
where $d_\star = c_\star/(c_\star+1) \in (0,1)$, and using this inequality recursively, one finds that
$$\Phi(T) \leq d_\star^{[r]} \Phi(T+r) + \frac{d_\star}{1-d_\star}\Psi(T+r), \qquad r \geq 1.$$
In particular, this result asserts that
$$\Phi(U_m) \leq d_\star^{S_m-U_m-1} \Phi(S_m) + \frac{d_\star}{1-d_\star}\Psi(S_m),$$
and because
$$\Phi(S_m)\ =\ \int_{x < S_m} |\nabla u_{m,+}^{(2)}|^2 \dx\dy\ \leq\ \int_\Sigma |\nabla u_{m,+}^{(2)}|^2 \dx\dy
\ \leq\ \|\zeta_m^{(2)}\|_{1/2}\ =\ O(1)$$
and
$$\Psi(S_m)\ =\ \int_{x < S_m} |G_{m,+}^{(2)}|^2 \dx\dy\ \leq\ \int_{|x|<S_m} |G_m^{(2)}|^2 \dx\dy\ =\ o(1)$$
as $m \rightarrow \infty$, we conclude that
$$\Phi(U_m)\ =\ \int_{x < U_m} |\nabla u_{m,+}^{(2)}|^2 \dx\dy\ =\ o(1)$$
as $m \rightarrow \infty$.

A similar argument shows that
$$\int_{x > -U_m} |\nabla u_{m,-}^{(2)}|^2 \dx\dy = o(1)$$
as $m \rightarrow \infty$, so that
\begin{eqnarray*}
\int_{|x|< U_m} |\nabla u_m^{(2)}|^2 & \leq & \int_{|x| < U_m} |\nabla u_{m,+}^{(2)}|^2 \dx\dy
+ \int_{|x| < U_m} |\nabla u_{m,-}^{(2)}|^2 \dx\dy \\
& \leq & \int_{x < U_m} |\nabla u_{m,+}^{(2)}|^2 \dx\dy
+ \int_{x > -U_m} |\nabla u_{m,-}^{(2)}|^2 \dx\dy \\
& \rightarrow & 0
\end{eqnarray*}
as $m \rightarrow \infty$.

The complementary estimate
$$\int_{|x| > U_m} |\nabla u_m^{(1)}|^2 \rightarrow 0$$
as $m \rightarrow \infty$ is obtained in a similar fashion.\qed

The next step is to apply Lemma \ref{Key splitting lemma} to the boundary-value problem \eqn{BC for u 1}--\eqn{BC for u 3}.

\begin{lemma} \label{First splitting result}
Let $u(\eta)$ be the solution to \eqn{BC for u 1}--\eqn{BC for u 3} with $\xi = \partial_x f(\eta)$, $\eta \in U$,
where $f$ is a real polynomial. The estimates
$$
\lim_{m \rightarrow \infty} \|\nabla u(\eta_m^{(1)})\|_{H^1(|x|>T_m)}=0,\qquad \lim_{m \rightarrow \infty} \|\nabla u(\eta_m^{(2)})\|_{H^1(|x|<T_m)} = 0
$$
hold for each sequence $\{T_m\}$ of positive real numbers with $S_m-T_m$, $T_m - R_m \rightarrow \infty$ as $m \rightarrow \infty$.
\end{lemma}
{\bf Proof.} Choose sequences $\{\tilde{R}_m\}$, $\{\tilde{S}_m\}$ of positive real numbers with
$S_m - \tilde{S}_m$, $\tilde{S}_m - T_m \rightarrow \infty$ and $T_m-\tilde{R}_m$, $\tilde{R}_m - R_m \rightarrow \infty$ as
$m \rightarrow \infty$.
The quantities $u_m^{(j)}=u(\eta_m^{(j)})$, $j=1,2$, satisfy the boundary-value problems
\begin{eqnarray*}
& & \parbox{90mm}{$\nabla\cdot((I+Q_m^{(j)})\nabla u_m^{(j)})=0$}0 < y <1,\\
& & \parbox{90mm}{$(I+Q_m^{(j)})\nabla u_m^{(j)} \cdot (0, 1) = f(\eta_m^{(j)})_x,$}y=1, \\
& & \parbox{90mm}{$(I+Q_m^{(j)})\nabla u_m^{(j)} \cdot (0,-1)=0,$}y=0,
\end{eqnarray*}
where
$Q_m^{(j)}=Q(\eta_m^{(j)})$,
and Lemma \ref{Key splitting lemma} asserts that
$$
\lim_{m \rightarrow \infty} \|\nabla u_m^{(1)}\|_{L^2(|x| >\tilde{R}_m)} = 0, \qquad \lim_{m \rightarrow \infty} \|\nabla u_m^{(2)}\|_{L^2(|x| <\tilde{S}_m)} = 0.
$$

The derivatives $u_{mx}^{(j)}$, $j=1,2$ are weak solutions of the boundary-value problems
\begin{eqnarray*}
& & \parbox{95mm}{$\nabla\cdot((I+Q_m^{(j)})\nabla u_{mx}^{(j)})=\nabla\cdot G_m^{(j)}$,}0 < y <1,\\
& & \parbox{95mm}{$(I+Q_m^{(j)})\nabla u_{mx}^{(j)} \cdot (0, 1) = f(\eta_m^{(j)})_{xx}+G_m^{(j)}\cdot(0,1),$}y=1, \\
& & \parbox{95mm}{$(I+Q_m^{(j)})\nabla u_{mx}^{(j)} \cdot (0,-1)=G_m^{(j)}\cdot(0,-1),$}y=0,
\end{eqnarray*}
where $G_m^{(j)} = -Q_{mx}^{(j)} \nabla u_m^{(j)}$. Using Remark \ref{Formulae for Qx and Qy}
and writing $S_{0m}^{(j)} = S_0(\eta_m^{(j)})$, $R_{0m}^{(j)} = R_0(\eta_m^{(j)})$, $R_{1m}^{(j)} = R_1(\eta_m^{(j)})$,
one finds that
\begin{eqnarray}
\lefteqn{\|Q_{mx}^{(1)} \nabla u_m^{(1)}\|_{L^2(|x| > \tilde{R}_m)}} \nonumber \\
& \leq & \|S_{0m}^{(1)}\|_\infty \|\nabla u_m^{(1)}\|_{L^2(|x| > \tilde{R}_m)} \nonumber \\
& & \qquad\mbox{}
+ c( \|R_{0m}^{(1)}\|_\infty \|L_0^\delta\| + \|R_{1m}^{(1)}\|_\infty \|L_1^\delta\|)\|(\eta_m^{(1)})^{\prime\prime}\|_0
\|\nabla u_m^{(1)}\|_{L^2(|x| > \tilde{R}_m)}^\frac{1}{2}\|\nabla u_m^{(1)}\|_{H^1(|x| > \tilde{R}_m)}^\frac{1}{2} \nonumber \\
& = & o(1) \label{Estimate Qxnabla}
\end{eqnarray}
as $m \rightarrow \infty$. (Lemma \ref{Regularity of u} asserts that $\{\nabla u_m^{(j)}\} \subseteq H^{3/2,1}$
and hence $\{\nabla u_m^{(j)}\} \subseteq H^1(\Sigma)$ is bounded; it follows that
$\|\nabla u_{m}^{(1)}\|_{H^1(|x|>\tilde{R}_m)}=O(1)$ as $m \rightarrow \infty$.)
A similar calculation shows that $\|Q_{mx}^{(2)} \nabla u_m^{(2)}\|_{L^2(|x| < \tilde{S}_m)} = o(1)$
as $m \rightarrow \infty$, and Lemma  \ref{Key splitting lemma} yields the estimates
$$
\lim_{m \rightarrow \infty} \|\nabla u_{mx}^{(1)}\|_{L^2(|x| >T_m)} = 0, \qquad \lim_{m \rightarrow \infty} \|\nabla u_{mx}^{(2)}\|_{L^2(|x| <T_m)} = 0.
$$

The calculation
$$u_{myy}^{(j)} = -\frac{1}{1+q_{m22}^{(j)}} \left( \partial_x [(1+q_{m11}^{(j)})u_{mx}^{(j)}+q_{m12}^{(j)} u_{my}^{(j)}] +
\partial_y(q_{m12}^{(j)}u_{mx}^{(j)})-q_{m22y}^{(j)}u_{my}^{(j)}\right)$$
(see equation \eqn{BC for u 4}) and estimates
$$\|q_{mij}^{(1)} \nabla u_{mx}^{(1)} \|_{L^2(|x| > T_m)} \leq \|q_{mij}^{(1)}\|_\infty \| \nabla u_{mx}^{(1)} \|_{L^2(|x| > T_m)} = o(1),$$
$$\|\begin{Bmatrix}\partial_x\\ \partial_y\end{Bmatrix} q_{mij}^{(1)} \nabla u_m^{(1)}\|_{L^2(|x|>T_m)} = o(1) $$
as $m \rightarrow \infty$ (cf.\ \eqn{Estimate Qxnabla})
show that
$$\lim_{m \rightarrow \infty} \|u_{myy}^{(1)} \|_{L^2(|x|>T_m)} = 0$$
(recall that $\|(1+q_{m22}^{(j)})^{-1}\|_\infty$ is bounded); the complementary limit
$$\lim_{m \rightarrow \infty} \|u_{myy}^{(2)} \|_{L^2(|x|<T_m)} = 0$$
is obtained in a similar fashion.\qed

Lemma \ref{Second splitting result} below states another useful application of Lemma
\ref{Key splitting lemma} to the boundary-value problem \eqn{BC for u 1}--\eqn{BC for u 3}; the following proposition is used in its proof.

\begin{proposition} \label{Q splitting prop}
Choose $N \in {\mathbb N}$. The estimates
\[
|(Q(\eta_m^{(1)}+\eta_m^{(2)})-Q(\eta_m^{(2)}))(x,y)|\le c\, \dist(x, [-R_m, R_m])^{-N}
\]
and
\[
|(Q(\eta_m^{(1)}+\eta_m^{(2)})-Q(\eta_m^{(1)}))(x,y)|\le c\, \dist(x, \R\sm (-S_m, S_m))^{-N}
\]
hold for all $(x,y) \in \overline{\Sigma}$, where $|\cdot|$ denotes the $2\times2$ matrix maximum norm,
and remain valid when $Q$ is replaced by $Q_x$ or $Q_y$.
\end{proposition}
{\bf Proof.} Observe that
$$\eta^\delta(x,y^\prime) = \frac{1}{1-y}\int_{\supp \eta} K\left(\frac{x-s}{1-y}\right)\eta(s) \mathrm{d}s,$$
where $K = (2\pi)^{-1/2} \delta^{-1}\FF^{-1}[\chi] \in \SS(\R)$. The above formula shows that
$\eta^\delta \in C^\infty(\overline{\Sigma}\setminus \supp \eta \times \{1\})$ with
\[
|\partial_x^j\partial_y^k\eta^\delta(x,y)|\le c\, \dist(x, \supp \eta)^{-N}
\|\eta\|_\infty
\]
for each $N \in {\mathbb N}$.

Note that
\begin{eqnarray*}
|(Q(\eta_1+\eta_2)-Q(\eta_2))(x,y)| 
& = & 
\left|\begin{pmatrix}f_{1y}^\delta & - f_{1x}^\delta \\ 
-  f_{1x}^\delta & \displaystyle \frac{- f_{3y}^\delta +(f_{3x}^\delta)^2}{1+ f_{3y}^\delta}
-\frac{- f_{2y}^\delta +(f_{2x}^\delta)^2}{1+ f_{2y}^\delta}\end{pmatrix}\!\!(x,y)\right| \\
& \leq & c |(f_{1x}^\delta, f_{1y}^\delta)(x,y)|
\end{eqnarray*}
for all $\eta_1$, $\eta_2$ and $\eta_3:=\eta_1+\eta_2 \in U$. It follows that
\begin{eqnarray*}
|(Q(\eta_m^{(1)}+\eta_m^{(2)})-Q(\eta_m^{(2)}))(x,y)|
& \leq & |((\eta_m^{(1)})^\delta(x,y), (\eta_{mx}^{(1)})^\delta(x,y),(\eta_{my}^{(1)})^\delta(x,y))| \\
& \leq & c\, \dist(x, [-R_m, R_m])^{-N}.
\end{eqnarray*}
The same argument yields the estimate for $Q(\eta_m^{(1)}+\eta_m^{(2)})-Q(\eta_m^{(1)})$ and the
corresponding results for $Q_x$ and $Q_y$.\qed

\begin{lemma} \label{Second splitting result}
Let $u(\eta)$ be the solution to \eqn{BC for u 1}--\eqn{BC for u 3} with $\xi = \partial_x f(\eta)$, $\eta \in U$,
where $f$ is a real polynomial. The estimates
\begin{eqnarray*}
& & \lim_{m \rightarrow \infty} \|\nabla u(\eta_m^{(1)}+\eta_m^{(2)})-\nabla u(\eta_m^{(1)})\|_{H^1(|x| < T_m)} = 0, \\
& & \lim_{m \rightarrow \infty} \|\nabla u(\eta_m^{(1)}+\eta_m^{(2)})-\nabla u(\eta_m^{(2)})\|_{H^1(|x| > T_m)} = 0
\end{eqnarray*}
hold for each sequence $\{T_m\}$ of positive real numbers with $S_m-T_m$, $T_m - R_m \rightarrow \infty$ as $m \rightarrow \infty$.
\end{lemma}
{\bf Proof.} Choose sequences $\{\tilde{R}_m\}$, $\{\tilde{S}_m\}$ of positive real numbers with
$S_m - \tilde{S}_m$, $\tilde{S}_m - T_m \rightarrow \infty$ and $T_m-\tilde{R}_m$, $\tilde{R}_m - R_m \rightarrow \infty$ as
$m \rightarrow \infty$. The quantities $w_m^{(1)}=u(\eta_m^{(1)}+\eta_m^{(2)})-u(\eta_m^{(2)})$ and $w_m^{(2)}=u(\eta_m^{(1)}+\eta_m^{(2)})-u(\eta_m^{(1)})$
satisfy the boundary-value problems
\begin{eqnarray*}
& & \parbox{90mm}{$\nabla\cdot((I+Q_m)\nabla w_m^{(j)})=
\nabla\cdot G_m^{(j)},$}0 < y <1,\\
& & \parbox{90mm}{$(I+Q_m)\nabla w_m^{(j)} \cdot (0, 1)= f(\eta_{m}^{(j)})_x+G_m^{(j)}\cdot (0,1) ,$}y=1, \\
& & \parbox{90mm}{$(I+Q_m)\nabla w_m^{(j)} \cdot (0,-1)=G_m^{(j)}\cdot(0,-1),$}y=0, 
\end{eqnarray*}
where $Q_m=Q(\eta_m^{(1)}+\eta_m^{(2)})$ and
$$
G_m^{(1)} = (Q_m^{(2)}-Q_m)\nabla u_m^{(2)}, \qquad
G_m^{(2)} = (Q_m^{(1)}-Q_m)\nabla u_m^{(1)}.
$$
Using the estimate
\[
|(Q_m^{(2)}-Q_m)(x,y)|\le c\,\dist(x, [-R_m, R_m])^{-N}
\]
(Proposition \ref{Q splitting prop}), one finds that
\[
\|G_m^{(1)}\|_{L^2(|x|>\tilde{R}_m)}^2
\ \leq\ c(\tilde{R}_m-R_m)^{-N} \|\nabla u_m^{(2)}\|_0^2
\ \leq\ c(\tilde{R}_m-R_m)^{-N} \|f(\eta_{m}^{(2)})\|_{1/2}^2
\ =\ o(1)
\]
as $m \rightarrow \infty$ and a similar argument shows that
$\|G_m^{(2)}\|_{L^2(|x|<\tilde{S}_m)}^2 = o(1)$
as $m \rightarrow \infty$. It follows from Lemma \ref{Key splitting lemma} that
$$
\lim_{m \rightarrow \infty} \|w_m^{(1)}\|_{L^2(|x| >T_m)} = 0, \qquad \lim_{m \rightarrow \infty} \|w_m^{(2)}\|_{L^2(|x| <T_m)} = 0.
$$

The derivatives $w_{mx}^{(j)}$, $j=1,2$ are weak solutions of the boundary-value problems
\begin{eqnarray*}
& & \parbox{95mm}{$\nabla\cdot((I+Q_m^{(j)})\nabla w_{mx}^{(j)})=\nabla\cdot H_m^{(j)}$,}0 < y <1,\\
& & \parbox{95mm}{$(I+Q_m^{(j)})\nabla w_{mx}^{(j)} \cdot (0, 1) = \partial_x^2 f(\eta_m^{(j)})+H_m^{(j)}\cdot(0,1),$}y=1, \\
& & \parbox{95mm}{$(I+Q_m^{(j)})\nabla w_{mx}^{(j)} \cdot (0,-1)=H_m^{(j)}\cdot(0,-1),$}y=0,
\end{eqnarray*}
where
\begin{eqnarray*}
H_m^{(1)} & = & -Q_{mx} \nabla w_m^{(1)} + (Q_m^{(2)}-Q_m) \nabla u_{mx}^{(2)} + (Q_{mx}^{(2)}-Q_{mx}) \nabla u_m^{(2)}, \\
H_m^{(2)} & = & -Q_{mx} \nabla w_m^{(2)} + (Q_m^{(1)}-Q_m) \nabla u_{mx}^{(1)} + (Q_{mx}^{(1)}-Q_{mx}) \nabla u_m^{(1)}.
\end{eqnarray*}
Treating $\|Q_{mx} \nabla w_m^{(1)}\|_{L^2(|x|>\tilde{R}_m)}$ using the method given in the proof
of Lemma \ref{First splitting result} (estimate \eqn{Estimate Qxnabla}) and
$\|(Q_m^{(2)}-Q_m) \nabla u_{mx}^{(2)}\|_{L^2(|x|>\tilde{R}_m)}$, $\|(Q_{mx}^{(2)}-Q_{mx}) \nabla u_m^{(2)}\|_{L^2(|x|>\tilde{R}_m)}$
using the method given above, one finds that
 $\|H_m^{(1)}\|_{L^2(|x|>\tilde{R}_m)} = o(1)$ as $m \rightarrow \infty$. A similar argument yields
$\|H_m^{(2)}\|_{L^2(|x|<\tilde{S}_m)} = o(1)$ as $m \rightarrow \infty$, and it follows from Lemma  \ref{Key splitting lemma}
that
$$
\lim_{m \rightarrow \infty} \|\nabla u_{mx}^{(1)}\|_{L^2(|x| >T_m)} = 0, \qquad \lim_{m \rightarrow \infty} \|\nabla u_{mx}^{(2)}\|_{L^2(|x| <T_m)} = 0.
$$

Finally, observe that
\begin{eqnarray*}
& & w_{myy}^{(1)} = -\frac{1}{1+q_{m22}^{(1)}} \Big( \partial_x [(1+q_{m11}^{(1)})w_{mx}^{(1)}+q_{m12}^{(1)} w_{my}^{(1)}] +
\partial_y(q_{m12}^{(1)}w_{mx}^{(1)})-q_{m22y}^{(1)}w_{my}^{(1)} \\
& & \hspace{4cm} \mbox{}+\nabla(Q_m^{(2)}-Q_m)\cdot\nabla u_m^{(1)} + (Q_m^{(2)}-Q_m)\Delta u_m^{(1)} \Big).
\end{eqnarray*}
The argument given in the proof of Lemma \ref{First splitting result} shows that
$$\|\partial_x [(1+q_{m11}^{(1)})w_{mx}^{(1)}+q_{m12}^{(1)} w_{my}^{(1)}] +
\partial_y(q_{m12}^{(1)}w_{mx}^{(1)})-q_{m22y}^{(1)}w_{my}^{(1)}\|_{L^2(|x|>T_m)} = o(1),$$
and the method given above shows that
$$\|\nabla(Q_m^{(2)}-Q_m)\cdot\nabla u_m^{(1)})\|_{L^2(|x|>T_m)},\ 
\|(Q_m^{(2)}-Q_m)\Delta u_m^{(1)}\|_{L^2(|x|>T_m)} = o(1)$$
as $m \rightarrow \infty$. One concludes that
$$\lim_{m \rightarrow \infty} \|w_{myy}^{(1)} \|_{L^2(|x|>T_m)} = 0,$$
and the complementary limit
$$\lim_{m \rightarrow \infty} \|w_{myy}^{(2)} \|_{L^2(|x|<T_m)} = 0$$
is obtained in a similar fashion.\qed

\begin{corollary} \label{Third splitting result}
The estimate
$$
\lim_{m \rightarrow \infty} \|\nabla u(\eta_m^{(1)}+\eta_m^{(2)})-\nabla u(\eta_m^{(1)})-\nabla u(\eta_m^{(2)})\|_1= 0
$$
holds under the hypotheses of Lemmata \ref{First splitting result} and \ref{Second splitting result}.
\end{corollary}

The proof of Theorem \ref{Main splitting theorem} is completed by applying the next lemma to the formula
for ${\mathcal T}^\prime$ given in Corollary \ref{Gradient of T}.

\begin{lemma}
\quad
\begin{list}{(\roman{count})}{\usecounter{count}}
\item
The estimates
\begin{eqnarray*}
& & \lim_{m \rightarrow \infty} \|f_1(\eta_m^{(1)}+\eta_m^{(2)})K(\eta_m^{(1)}+\eta_m^{(2)}) f_2(\eta_m^{(1)}+\eta_m^{(2)}) \\
& & \qquad\mbox{}
-f_1(\eta_m^{(1)})K(\eta_m^{(1)})f_2(\eta_m^{(1)})
-f_1(\eta_m^{(2)})K(\eta_m^{(2)})f_2(\eta_m^{(2)})\|_0 = 0
\end{eqnarray*}
and
\begin{eqnarray*}
& & \lim_{m \rightarrow \infty} \|f_1(\eta_m^{(1)}+\eta_m^{(2)})K(\eta_m^{(1)}+\eta_m^{(2)}) f_2(\eta_m^{(1)}+\eta_m^{(2)}) \\
& & \qquad\mbox{}
-f_1(\eta_m^{(1)})K(\eta_m^{(1)})f_2(\eta_m^{(1)})
-f_1(\eta_m^{(2)})K(\eta_m^{(2)})f_2(\eta_m^{(2)})\|_{L^1(\R)} = 0.
\end{eqnarray*}
hold for all real polynomials $f_1$, $f_2$.
\item
The estimate
\begin{eqnarray*}
& & \lim_{m \rightarrow \infty} \|{\mathcal H}^\prime(\eta_m^{(1)}+\eta_m^{(2)})(f_1(\eta_m^{(1)}+\eta_m^{(2)}),f_2(\eta_m^{(1)}+\eta_m^{(2)}) \\
& & \qquad\mbox{}
-{\mathcal H}^\prime(\eta_m^{(1)})(f_1(\eta_m^{(1)}),f_2(\eta_m^{(1)}))
-{\mathcal H}^\prime(\eta_m^{(2)})(f_1(\eta_m^{(2)}),f_2(\eta_m^{(2)}) )\|_0 =0
\end{eqnarray*}
holds for all real polynomials $f_1$, $f_2$.
\item
The estimate
$$\lim_{m \rightarrow \infty} \langle {\mathcal H}^\prime(\eta_m^{(1)})(f_1(\eta_m^{(1)}),f_2(\eta_m^{(1)})), \eta_m^{(2)} \rangle_0 = 0$$
holds for all real polynomials $f_1$, $f_2$.
\end{list}
\end{lemma}
{\bf Proof.} (i) Observe that
\begin{eqnarray*}
& & f_1(\eta_m^{(1)}+\eta_m^{(2)})K(\eta_m^{(1)}+\eta_m^{(2)}) f_2(\eta_m^{(1)}+\eta_m^{(2)}) \\
& & \qquad\mbox{}
-f_1(\eta_m^{(1)})K(\eta_m^{(1)})f_2(\eta_m^{(1)})
-f_1(\eta_m^{(2)})K(\eta_m^{(2)})f_2(\eta_m^{(2)}) \\
& = & f_1(\eta_m^{(1)})(u_x(\eta_m^{(1)}+\eta_m^{(2)})-u_x(\eta_m^{(1)}))
+ f_2(\eta_m^{(2)})(u_x(\eta_m^{(1)}+\eta_m^{(2)})-u_x(\eta_m^{(2)})).
\end{eqnarray*}
The $L^1(\R)$- and $L^2(\R)$-norms of this quantity can both be estimated by
\begin{eqnarray*}
& & \|f_1(\eta_m^{(1)})\|_1
\|u_x(\eta_m^{(1)}+\eta_m^{(2)})-u_x(\eta_m^{(1)}) |_{y=1} \|_{L^2(|x|<R_m)} \\
& & \quad\mbox{}+\|f_2(\eta_m^{(1)})\|_1
\|u_x(\eta_m^{(1)}+\eta_m^{(2)})-u_x(\eta_m^{(2)} |_{y=1} )\|_{L^2(|x|>S_m)} \hspace{0.75in} \\
& \leq & \underbrace{\|f_1(\eta_m^{(1)})\|_1}_{\textstyle O(1)}
\underbrace{\|\nabla u(\eta_m^{(1)}+\eta_m^{(2)})-\nabla u(\eta_m^{(1)})  \|_{H^1(|x|<T_m)}}_{\textstyle o(1)} \\
& & \quad\mbox{}+\underbrace{\|f_2(\eta_m^{(1)})\|_1}_{\textstyle O(1)}
\underbrace{\|\nabla u(\eta_m^{(1)}+\eta_m^{(2)})-\nabla u(\eta_m^{(2)}  )\|_{H^1(|x|>T_m)}}_{\textstyle o(1)} \hspace{0.75in} \\
& = & o(1)
\end{eqnarray*}
(use the Cauchy-Schwarz inequality or the maximum norm for the polynomials).

(ii) Observe that
\begin{eqnarray*}
& &{\mathcal H}^\prime(\eta_m^{(1)}+\eta_m^{(2)})(f_1(\eta_m^{(1)}+\eta_m^{(2)}),f_2(\eta_m^{(1)}+\eta_m^{(2)}) \\
& & \qquad\mbox{}
-{\mathcal H}^\prime(\eta_m^{(1)})(f_1(\eta_m^{(1)}),f_2(\eta_m^{(1)}))
-{\mathcal H}^\prime(\eta_m^{(2)})(f_1(\eta_m^{(2)}),f_2(\eta_m^{(2)}) )\\[1mm]
& & =
-u_x(\eta_m^{(1)}+\eta_m^{(2)})v_x(\eta_m^{(1)}+\eta_m^{(2)})
+u_x(\eta_m^{(1)})v_x(\eta_m^{(1)}) + u_x(\eta_m^{(2)})v_x(\eta_m^{(2)}) \\
& & \qquad\mbox{}+u_y(\eta_m^{(1)}+\eta_m^{(2)})v_y(\eta_m^{(1)}+\eta_m^{(2)})
-u_y(\eta_m^{(1)})v_y(\eta_m^{(1)}) - u_y(\eta_m^{(2)})v_y(\eta_m^{(2)}) \\
& & \qquad\mbox{}+h(\eta_m^{(1)}+\eta_m^{(2)})u_y(\eta_m^{(1)}+\eta_m^{(2)})v_y(\eta_m^{(1)}+\eta_m^{(2)}) \\
& & \qquad\mbox{}-h(\eta_m^{(1)})u_y(\eta_m^{(1)})v_y(\eta_m^{(1)})-h(\eta_m^{(2)})u_y(\eta_m^{(2)})v_y(\eta_m^{(2)})\Big|_{y=1},
\end{eqnarray*}
where
$$h(\eta)=\frac{\eta^{\prime 2}-\eta^2-2\eta}{(1+\eta)^2}$$
and $u(\eta)$, $v(\eta)$ are the solutions  to \eqn{BC for u 1}--\eqn{BC for u 3} with respectively $\xi = \partial_x f_1(\eta)$ and $\xi=\partial_x f_2(\eta)$, $\eta \in U$.

The estimates
\begin{eqnarray*}
\lefteqn{\|u_x(\eta_m^{(1)} + \eta_m^{(2)})v_x(\eta_m^{(1)} + \eta_m^{(2)})
- (u_x(\eta_m^{(1)}) + u_x(\eta_m^{(2)}))(v_x(\eta_m^{(1)}) + v_x(\eta_m^{(2)}))  |_{y=1} \|_0} \\
& & \leq \underbrace{\|v_x(\eta_m^{(1)}+\eta_m^{(2)}) |_{y=1}\|_1}_{\displaystyle = O(1)}\underbrace{\|u_x(\eta_m^{(1)} + \eta_m^{(2)}) - u_x(\eta_m^{(1)}) - u_x(\eta_m^{(2)}) |_{y=1}\|_0}_{\displaystyle = o(1)} \\
& & \qquad\mbox{}+\underbrace{\|u_x(\eta_m^{(1)})+u_x(\eta_m^{(2)}) |_{y=1}\|_1}_{\displaystyle = O(1)}\underbrace{\|v_x(\eta_m^{(1)} +\eta_m^{(2)}) - v_x(\eta_m^{(1)}) - v_x(\eta_m^{(2)}) |_{y=1}\|_0}_{\displaystyle = o(1)} \\
& & = o(1)
\end{eqnarray*}
and
\begin{eqnarray*}
& & \|(u_x(\eta_m^{(1)}) + u_x(\eta_m^{(2)}))(v_x(\eta_m^{(1)}) + v_x(\eta_m^{(2)}))
- u_x(\eta_m^{(1)})v_x(\eta_m^{(1)}) - u_x(\eta_m^{(2)})v_x(\eta_m^{(2)})   |_{y=1} \|_0 \\
& \leq & \|u_x(\eta_m^{(1)})v_x(\eta_m^{(2)}) |_{y=1} \|_0 + \|u_x(\eta_m^{(2)})v_x(\eta_m^{(1)}) |_{y=1} \|_0 \\
& \leq & c \big(\underbrace{\|u_x(\eta_n^{(1)}) |_{y=1} \|_{L^2(|x|>T_m)}}_{\displaystyle = o(1)}
\underbrace{\|v_x(\eta_n^{(2)}) |_{y=1} \|_1}_{\displaystyle = O(1)} +
\underbrace{\|u_x(\eta_n^{(1)}) |_{y=1} \|_1}_{\displaystyle = O(1)}\underbrace{\|v_x(\eta_n^{(2)}) |_{y=1} \|_{L^2(|x| <T_m)}}_{\displaystyle = o(1)}\\
& & \quad\mbox{}+\underbrace{\|u_x(\eta_n^{(2)}) |_{y=1} \|_{L^2(|x| <T_m)}}_{\displaystyle = o(1)}
\underbrace{\|v_x(\eta_n^{(1)}) |_{y=1} \|_1}_{\displaystyle = O(1)} +
\underbrace{\|u_x(\eta_n^{(2)}) |_{y=1} \|_1}_{\displaystyle = O(1)}\underbrace{\|v_x(\eta_n^{(1)}) |_{y=1} \|_{L^2(|x| >T_m)}}_{\displaystyle = o(1)}\big) \\
& = & o(1)
\end{eqnarray*}
imply that
$$
\|(u_x(\eta_m^{(1)}) + u_x(\eta_m^{(2)}))(v_x(\eta_m^{(1)}) + v_x(\eta_m^{(2)}))
- u_x(\eta_m^{(1)})v_x(\eta_m^{(1)}) - u_x(\eta_m^{(2)})v_x(\eta_m^{(2)})   |_{y=1} \|_0 = o(1)
$$
as $m \rightarrow \infty$; here we have used the estimate
$$\|u_x(\eta)|_{y=1}\|_1\ \leq\ c\|\nabla u\|_{3/2,1}\ \leq c\|f_1(\eta)\|_2, \qquad \eta \in U$$
and its counterpart for $v$. The same argument shows that
$$
\|(u_y(\eta_m^{(1)}) + u_y(\eta_m^{(2)}))(v_y(\eta_m^{(1)}) + v_y(\eta_m^{(2)}))
- u_y(\eta_m^{(1)})v_y(\eta_m^{(1)}) - u_y(\eta_m^{(2)})v_y(\eta_m^{(2)})   |_{y=1} \|_0 = o(1)
$$
as $m \rightarrow \infty$.

Because $h(\eta_m^{(1)}+\eta_m^{(2)})=h(\eta_m^{(1)})+h(\eta_m^{(2)})$ and
$$\|(u_y(\eta_m^{(1)}) + u_y(\eta_m^{(2)}))(v_y(\eta_m^{(1)}) + v_y(\eta_m^{(2)}))
- u_y(\eta_m^{(1)})v_y(\eta_m^{(1)}) - u_y(\eta_m^{(2)})v_y(\eta_m^{(2)})   |_{y=1} \|_0=o(1)$$
as $m \rightarrow \infty$ (see above), repeating the proof of part (i) above yields the estimate
\begin{eqnarray*}
& & \qquad\mbox{}\|h(\eta_m^{(1)}+\eta_m^{(2)})u_y(\eta_m^{(1)}+\eta_m^{(2)})v_y(\eta_m^{(1)}+\eta_m^{(2)}) \\
& & \qquad\qquad\mbox{}-h(\eta_m^{(1)})u_y(\eta_m^{(1)})v_y(\eta_m^{(1)})-h(\eta_m^{(2)})u_y(\eta_m^{(2)})v_y(\eta_m^{(2)}) |_{y=1}\|_0=o(1)
\end{eqnarray*}
as $m \rightarrow \infty$.

(iii) The methods used in part (ii) show that
$$
\|{\mathcal H}^\prime(\eta_m^{(1)})(f_1(\eta_m^{(1)}),f_2(\eta_m^{(1)}))\|_{L^2(|x|>T_m)} = o(1),
$$
so that
\begin{eqnarray*}
|\langle {\mathcal H}^\prime(\eta_m^{(1)})(f_1(\eta_m^{(1)}),f_2(\eta_m^{(1)})), \eta_m^{(2)} \rangle_0|
& \leq & \underbrace{\|{\mathcal H}^\prime(\eta_m^{(1)})(f_1(\eta_m^{(1)}),f_2(\eta_m^{(1)}))\|_{L^2(|x|>S_m)}}_{\displaystyle = O(1)}
\underbrace{\|\eta_m^{(2)}\|_0}_{\displaystyle = o(1)} \\
& \rightarrow & 0
\end{eqnarray*}
as $m \rightarrow \infty$.\qed

\section{Minimising sequences} \label{Minimising sequences}

The goal of this section is the proof of the following theorem, the existence of the sequence advertised in
which is a key ingredient
in the proof that the infimum of $\JJ_\mu$ over $U \sm \{0\}$ is a strictly sub-additive function of
$\mu$. The subadditivity property of $c_\mu$ is in turn used later to establish the convergence (up to subsequences
and translations) of \emph{any} minimising sequence for $\JJ_\mu$ over $U \sm \{0\}$ which does not approach the boundary of $U$.

\begin{theorem} \label{Special MS theorem}
There exists a minimising sequence $\{\tilde{\eta}_m\}$ for $\JJ_\mu$ over $U\sm\{0\}$
with the properties that $\|\tilde{\eta}_m\|_2^2 \leq c \mu$ for each $m \in {\mathbb N}$ and
$\lim_{m \rightarrow \infty}\|\JJ_\mu^\prime(\tilde{\eta}_m)\|_0=0$.
\end{theorem}

\subsection{The penalised minimisation problem} \label{Penalised minimisation}
We begin by studying the functional 
$\JJ_{\rho,\mu}: H^2(\R) \rightarrow \R \cup \{\infty\}$ defined by
$$\JJ_{\rho,\mu}(\eta) = \left\{\begin{array}{lll}\displaystyle \KK(\eta)+\frac{(\mu+\GG(\eta))^2}{\LL(\eta)}
+ \rho(\|\eta\|_2^2), & & \eta \in U\sm\{0\}, \\
\\
\infty, & & \eta \not\in U\sm\{0\},
\end{array}\right.$$
in which $\rho: [0,M^2) \rightarrow \R$ is a smooth, increasing `penalisation' function such that
$\rho(t)=0$ for $0 \leq t \leq \tilde{M}^2$ and $\rho(t) \rightarrow \infty$ as $t \uparrow M^2$. We allow
negative values of the small parameter, so that $0< |\mu| < \mu_0$ (see the comments below
Lemma \ref{Vanishing and concentration}) and the number $\tilde{M} \in (0,M)$ is chosen so that
$$\tilde{M}^2 > (c^\star+D\nu_0+D\nu_0^-)|\mu|;$$
the following analysis is valid for every such choice of $\tilde{M}$, which
in particular may be chosen arbitrarily close to $M$.
In this inequality $\nu_0$ and $\nu_0^-$ are the speeds of linear waves with frequency $k_0$ riding shear
flows with vorticities $\omega$ and $-\omega$ and $c^\star$, $D$ are constants identified
in Lemmata \ref{Test function}(i) and \ref{Critical point estimate} below.
In Section \ref{MSP}
we give a detailed description of the qualitative properties of an arbitrary minimising sequence
$\{\eta_m\}$ for $\JJ_{\rho,\mu}$; the penalisation function ensures that $\{\eta_m\}$ does not
approach the boundary of the set $U \sm \{0\}$ in which $\JJ_\mu$ is defined.

We first give some useful \emph{a priori} estimates. Lemma \ref{Test function}(i)
shows in particular that
$$c_{\rho,\mu}:=\inf \JJ_{\rho,\mu} < 2\nu_0^\mu|\mu| - c|\mu|^{r^\star}, \qquad
c_\mu := \inf_{\eta \in U\sm\{0\}} \JJ_\mu(\eta) < 2\nu_0^\mu|\mu| - c|\mu|^{r^\star},$$
where $\nu_0^\mu$ is the speed of linear waves with frequency $k_0$ riding a shear flow with vorticity
$(\sgn \mu) \omega$ (which depends only upon the sign of $\mu$),
while Lemma \ref{Critical point estimate}, whose proof is a straightforward modification of the argument
presented by Buffoni \emph{et al.} \cite[Propositions 2.34 and 3.2]{BuffoniGrovesSunWahlen13},
gives estimates on the size of critical points of $\JJ_\mu$ and a class of related functionals.

\begin{lemma} \label{Test function}
\quad
\begin{list}{(\roman{count})}{\usecounter{count}}
\item
There exists $\eta_{\mu}^\star \in U\sm\{0\}$ with compact support and a positive constant $c^\star$
such that
$\|\eta_\mu^\star\|_2 \leq c^\star |\mu|^\frac{1}{2}$, $\rho(\|\eta_\mu^\star\|_2^2)=0$ and
$$\JJ_{\rho,\mu}(\eta_\mu^\star)\ =\ \JJ_\mu(\eta_\mu^\star)\ <\ 2\nu_0^\mu|\mu| - c|\mu|^{r^\star},
\qquad
r^\star = \left\{\begin{array}{ll} \frac{5}{3},\quad\mbox{} & \beta>\beta_\mathrm{c}, \\[1.5mm]
3, & \beta<\beta_\mathrm{c}. \end{array}\right.$$
\item
The inequality
$$\KK_2(\eta) + \frac{(\mu+\GG_2(\eta))^2}{\LL_2(\eta)} \geq 2 \nu_0^\mu|\mu|$$
holds for each $\eta \in H^2(\R)\sm\{0\}$.
\end{list}
\end{lemma}
{\bf Proof.} First suppose that $\mu>0$.
The proof of part (i) is recorded in Appendix A, while part (ii) follows from the calculation
\begin{eqnarray*}
\lefteqn{\KK_2(\eta) + \frac{(\mu+\GG_2(\eta))^2}{\LL_2(\eta)}} \qquad \\
& = & \KK_2(\eta) + 2\nu_0 \GG_2(\eta) - \nu_0^2 \LL_2(\eta) + \frac{(\mu+\GG_2(\eta)-\nu_0\LL_2(\eta))^2}{\LL_2(\eta)}
+2\nu_0 \mu \\
& = & \frac{1}{2}\int_{-\infty}^\infty g(k) |\hat{\eta}|^2 + \frac{(\mu+\GG_2(\eta)-\nu_0\LL_2(\eta))^2}{\LL_2(\eta)}
+2\nu_0 \mu \\[1mm]
& \geq & 2\nu_0 \mu.
\end{eqnarray*}

For $\mu<0$ we observe that $\JJ_\mu(\eta)$, $\JJ_{\rho,\mu}(\eta)$ and
$\KK_2(\eta) + (\mu + \GG_2(\eta))^2/\LL_2(\eta)$ are invariant under the transformation $(\mu,\omega) \mapsto (-\mu,-\omega)$.\qed

\begin{lemma} \label{Critical point estimate}
Suppose that $\gamma_1$ and $\gamma_2$ belong to a bounded set of real numbers.
Any critical point $\eta$ of the functional $\tilde{\JJ}_\gamma: U \rightarrow \R$
defined by the formula
$$\tilde{\JJ}_\gamma(\eta)=\KK(\eta) - \gamma_1 \GG(\eta) - \gamma_2 \LL(\eta)+\gamma_3\|\eta\|_2^2, \qquad \gamma_3 \geq 0$$
satisfies the estimate
$$\|\eta\|_2^2 \leq D\KK(\eta),$$
where $D$ is a positive constant which does not depend upon $\gamma_1$, $\gamma_2$ or $\gamma_3$.
\end{lemma}
\begin{corollary} \label{Critical point estimate for weighted J}
Any critical point $\eta$ of $\JJ_{\rho,\mu}$ with $\JJ_{\rho,\mu}(\eta)<2\nu_0^\mu|\mu|$ satisfies the estimates
$$\|\eta\|_2^2 \leq 2D\nu_0^\mu|\mu|, \qquad  \rho(\|\eta\|_2^2)=0.$$
\end{corollary}
{\bf Proof.} Notice that any critical point $\eta$ of $\JJ_{\rho,\mu}$ is also a critical point of the functional
$\tilde{\JJ}_\gamma$, where
$$\gamma_1 = - \frac{2(\mu+\GG(\eta))}{\LL(\eta)}, \qquad
\gamma_2 = \frac{(\mu+\GG(\eta))^2}{\LL(\eta)^2}, \qquad
\gamma_3 = 2\rho^\prime(\|\eta\|_2^2).$$
Furthermore, any function $\eta \in U$ such that
$$
\frac{(\mu+\GG(\eta))^2}{\LL(\eta)} \leq 2\nu_0^\mu|\mu|
$$
satisfies
$$
\frac{\mu^2}{\LL(\eta)}
\ \leq\ 2\nu_0^\mu|\mu| - \frac{2\mu\GG(\eta)}{\LL(\eta)} - \frac{\GG(\eta)^2}{\LL(\eta)}
\ \leq\ 2\nu_0^\mu|\mu| + \frac{2|\mu||\GG(\eta)|}{\LL(\eta)} \\
\ \leq\ c |\mu|
$$
(see Proposition \ref{Quadratic estimates}), so that
\begin{equation}
\frac{|\mu|}{\mathcal L(\eta)}\leq c. \label{Bound on L from below}
\end{equation}
Observing that
$$
\frac{(\mu+\GG(\eta))^2}{\LL(\eta)}\ \leq\ \JJ_{\rho,\mu}(\eta)\ \leq\ 2\nu_0^\mu|\mu|,
$$
we find from Proposition \ref{Quadratic estimates} and inequality \eqn{Bound on L from below}
that $\gamma_1$ and $\gamma_2$ are bounded.
The previous lemma shows that $\|\eta\|_2^2 \leq D\KK(\eta) \leq D\JJ_{\rho,\mu}(\eta) < 2D\nu_0^\mu|\mu|$ and
hence $\rho(\|\eta\|_2^2)=0$ because of the choice of $\tilde{M}$.\qed

Finally, we establish some basic properties of a minimising sequence $\{\eta_m\}$ for $\JJ_{\rho,\mu}$.
Without loss of generality we may assume that
$$\sup_{m \in {\mathbb N}} \|\eta_m\|_2 < M$$
($\|\eta_m\|_2 \rightarrow M$ would imply that $\JJ_{\rho,\mu}(\eta_m) \rightarrow \infty$),
and it follows that $\{\eta_m\}$ admits a subsequence such that $\lim_{m \rightarrow \infty}
\|\eta_m\|_2$ exists and is positive ($\eta_m \rightarrow 0$ in $H^2(\R)$ would also imply
that $\JJ_{\rho,\mu}(\eta_m) \rightarrow \infty$). The following lemma records further
useful properties of $\{\eta_m\}$.

\begin{lemma} \label{General properties}
Every minimising sequence $\{\eta_m\}$ for $\JJ_{\rho,\mu}$ has the properties that
$$
\JJ_{\rho,\mu}(\eta_m) < 2\nu_0^\mu|\mu|-c|\mu|^{r^\star}, \qquad \LL(\eta_m) \geq c|\mu|,
\qquad \LL_2(\eta_m) \geq c|\mu|,$$
$$\MM_{\rho,\mu}(\eta_m) \leq - c|\mu|^{r^\star},
\qquad \|\eta_m\|_{1,\infty} \geq c |\mu|^{r^\star}
$$
for each $m \in {\mathbb N}$, where
$$\MM_{\rho,\mu}(\eta) = \JJ_{\rho,\mu}(\eta) - \KK_2(\eta) - \frac{(\mu+\GG_2(\eta))^2}{\LL_2(\eta)}.$$
\end{lemma}
{\bf Proof. }The first and second estimates are obtained from Lemma \ref{Test function}(i) and
the remark leading to \eqn{Bound on L from below}, while the third is a consequence of the
calculation
\begin{equation}
c\|\eta\|_{1/2}^2\ \leq\ \LL_2(\eta), \LL(\eta)\ \leq c\|\eta\|_{1/2}^2, \qquad \eta \in U.
\label{L and L2 equivalence}
\end{equation}
Turning to the fourth estimate, observe that
$$\MM_{\rho,\mu}(\eta_m)\ \leq\ \JJ_{\rho,\mu}(\eta_m) - 2\nu_0^\mu|\mu|\ \leq\ -c|\mu|^{r^\star}$$
because
$$\KK_2(\eta) + \frac{(\mu+\GG_2(\eta))^2}{\LL_2(\eta)} \geq 2\nu_0^\mu|\mu|$$
(see Lemma \ref{Test function}(ii)).

Finally, it follows from the calculation
\begin{eqnarray*}
\lefteqn{\MM_{\rho,\mu}(\eta_m)-\rho(\|\eta_m\|_2^2)}\quad \\
& & = \mathcal K_\mathrm{nl}(\eta_m)-\frac{\mu^2 \mathcal L_\mathrm{nl}(\eta_m)}{\mathcal L(\eta_m)\mathcal L_2(\eta_m)}
-\frac{2\mu\GG(\eta_m) \mathcal L_\mathrm{nl}(\eta_m)}{\mathcal L(\eta_m)\mathcal L_2(\eta_m)}
+\frac{2\mu\GG_\mathrm{nl}(\eta_m)}{\mathcal L(\eta_m)}\hspace{2cm} \\[1mm]
& &
\qquad\mbox{}-\frac{\GG_2(\eta_m) \mathcal L_\mathrm{nl}(\eta_m)}{\mathcal L(\eta_m)\mathcal L_2(\eta_m)}+\frac{(\GG(\eta_m)+\GG_2(\eta_m)) \GG_\mathrm{nl}(\eta_m)}{\mathcal L(\eta_m)},
\end{eqnarray*}
the inequalities
$$|\GG_2(\eta_m)|, |\GG(\eta_m)| \leq c\|\eta_m\|_{1/2}^2,$$
$$ |\GG_\mathrm{nl}(\eta_m)|,\ |\KK_\mathrm{nl}(\eta_m)|
\leq c\|\eta_m\|_{1,\infty}, \qquad
|\LL_\mathrm{nl}(\eta_m)| \leq c\|\eta_m\|_{1,\infty} \|\eta_m\|_{1/2}^2
$$
and \eqn{L and L2 equivalence} that
$$|\MM_{\rho,\mu}(\eta_m)-\rho(\|\eta_m\|_2^2)| \leq c\|\eta_m\|_{1,\infty}.$$
The fifth estimate is obtained from this result and the fact that
 $$\MM_{\rho,\mu}(\eta_m) - \rho(\|\eta_m\|_2^2) \leq -c|\mu|^{r^\star}.\eqno{\Box}$$
 
\begin{remark} \label{General properties without rho}
Replacing $\JJ_{\rho,\mu}(\eta)$ by $\JJ_\mu(\eta)$ and
$\MM_{\rho,\mu}(\eta)$ by
$$\MM_\mu(\eta) := \JJ_\mu(\eta) - \KK_2(\eta) - \frac{(\mu+\GG_2(\eta))^2}{\LL_2(\eta)}$$
in its statement, one finds that the above lemma is also valid for a minimising sequence $\{\eta_m\}$ for $\JJ_\mu$ over $U\sm\{0\}$.
\end{remark}

\subsection{Minimising sequences for the penalised problem} \label{MSP}

\subsubsection{Application of the concentration-compactness principle} \label{Application of cc}

The next step is to perform a more detailed analysis of the behaviour of a minimising sequence $\{\eta_m\}$
for $\JJ_{\rho,\mu}$ by applying the concentration-compactness principle (Lions \cite{Lions84a,Lions84b}); Theorem \ref{Concentration-compactness theorem} below states this result in a form suitable for the present situation.

\begin{theorem} \label{Concentration-compactness theorem}
Any sequence $\{u_m\} \subset L^1(\R)$ of non-negative functions with the
property that
$$\lim_{m \rightarrow \infty} \int_{-\infty}^\infty  u_m(x) \dx = \ell > 0$$
admits a subsequence for which precisely one of the following phenomena occurs.\\
\\
\underline{Vanishing}: For each $r>0$ one has that
$$\lim_{m \rightarrow \infty}\Bigg(\sup_{\tilde{x} \in \R} \int_{\tilde{x}-r}^{\tilde{x}+r}\!\!\! u_m(x) \dx\Bigg) = 0.
$$
\\
\underline{Concentration}: There is a sequence $\{x_m\} \subset \R$ with the property that for each
$\varepsilon>0$ there exists a positive real number $R$ with
$$\int_{-R}^R u_m(x+x_m) \dx \geq \ell-\varepsilon$$
for each $m \in {\mathbb N}$.\\
\\
\underline{Dichotomy}: There are sequences $\{x_m\} \subset \R$, $\{M_m^{(1)}\}, \{M_m^{(2)}\} \subset \R$ and a real number
$\kappa \in (0,\ell)$ with the properties that $M_m^{(1)}$, $M_m^{(2)} \rightarrow \infty$,
$M_m^{(1)}/M_m^{(2)} \rightarrow 0$,
$$\int_{-M_m^{(1)}}^{M_m^{(1)}} u_m(x+x_m) \dx \rightarrow \kappa, \qquad
\int_{-M_m^{(2)}}^{M_m^{(2)}} u_m(x+x_m) \dx \rightarrow \kappa$$
as $m \rightarrow \infty$. Furthermore
$$\lim_{m \rightarrow \infty} \Bigg(
\sup_{\tilde{x} \in \R} \int_{\tilde{x}-r}^{\tilde{x}+r} \!\!\! u_m(x) \dx\Bigg) \leq \kappa$$
for each $r>0$, and for each $\varepsilon>0$ there is a positive,
real number $R$ such that
$$\int_{-R}^R u_m(x+x_m) \dx \geq \kappa-\varepsilon$$
for each $m \in {\mathbb N}$.
\end{theorem}

Standard interpolation inequalities show that the norms
$\|\!\cdot\!\|_r$ are metrically equivalent on $U$ for $r \in [0,2)$; we therefore study
the convergence properties of $\{\eta_m\}$ in $H^r(\R)$ for $r \in [0,2)$, by focussing
on the concrete case $r=1$. One
may assume that $\|\eta_m\|_1 \rightarrow \ell$ as $m \rightarrow \infty$, where $\ell>0$
because $\eta_m \rightarrow 0$ in $H^r(\R)$ for $r>\frac{3}{2}$ would imply that
$\JJ_{\rho,\mu}(\eta_m) \rightarrow \infty$. This observation suggests applying 
Theorem \ref{Concentration-compactness theorem} to the sequence $\{u_m\}$ defined by 
$$
u_m = \eta_m^{\prime 2} + \eta_m^2,
$$
so that $\|u_m\|_{L^1(\R)} = \|\eta_m\|_1^2$.The following result
deals with `vanishing' and `concentration' (see Buffoni \emph{et al.}
\cite[Lemmata 3.7 and 3.9]{BuffoniGrovesSunWahlen13}.

\begin{lemma} \label{Vanishing and concentration}
\quad
\begin{list}{(\roman{count})}{\usecounter{count}}
\item
The sequence $\{u_m\}$ does not have the `vanishing' property.
\item
Suppose that $\{u_m\}$ has the `concentration' property. The sequence
$\{\eta_m(\cdot+x_m)\}$ admits a subsequence, with a slight abuse of
notation abbreviated to $\{\eta_m\}$, which satisfies
$$\lim_{m \rightarrow \infty} \|\eta_m\|_2 \leq \tilde{M}$$
and converges in $H^r(\R)$ for $r \in [0,2)$, to
$\eta^{(1)}$. The function $\eta^{(1)}$ satisfies the estimate
$$\|\eta^{(1)}\|_2^2\ \leq\ D\KK(\eta^{(1)})\ <\ 2D\nu_0^\mu|\mu|,$$
minimises $\JJ_{\rho,\mu}$ and minimises $\JJ_\mu$ over
$\tilde{U}\sm\{0\}$, where $\tilde{U}=\{\eta \in H^2(\R): \|\eta\|_2< \tilde{M}\}$.
\end{list}
\end{lemma}

We now present the more involved discussion of the remaining case (`dichotomy'),
again abbreviating the subsequence of $\{\eta_m(\cdot+x_m)\}$ identified by
Theorem \ref{Concentration-compactness theorem} to $\{\eta_m\}$.
The analysis is similar to that given by Buffoni \emph{et al.} \cite{BuffoniGrovesSunWahlen13}
in their study of three-dimensional irrotational solitary waves, the main difference being
that negative values of $\mu$ are also considered, so that $\mu$ is replaced by $|\mu|$ in
estimates (this change is necessary since the numbers $\mu^{(1)}$ and $\mu^{(2)}$
appearing in part (iv) of the following lemma, which are later used iteratively, may be negative).
We therefore omit proofs which are straightfoward modifications of those given
by Buffoni \emph{et al.}; note however that references in that paper to Appendix D
(in particular Theorem D.6) for `pseudo-local' properties of operators should be
replaced by references to Section \ref{Pseudo-local properties} (in particular
Theorem \ref{Main splitting theorem}) here.

Define sequences $\{\eta_m^{(1)}\}$, $\{\eta_m^{(2)}\}$ by the formulae
$$\eta_m^{(1)}(x)=\eta_m(x)\chi\left(\frac{x}{M_m^{(1)}}\right), \qquad
\eta_m^{(2)}(x) = \eta_m(x) \left(1-\chi\left(\frac{x}{M_m^{(2)}}\right)\right),$$
so that
$$\supp \eta_m^{(1)} \subset [-2M_m^{(1)},2M_m^{(1)}], \qquad
\supp \eta_m^{(2)} \subset \R \sm (-M_m^{(2)},M_m^{(2)}).$$

\begin{lemma} \label{Splitting properties}
\quad
\begin{list}{(\roman{count})}{\usecounter{count}}
\item
The sequences $\{\eta_m\}$, $\{\eta_m^{(1)}\}$ and $\{\eta_m^{(2)}\}$ have the limiting behaviour
$$\|\eta_m^{(1)}\|_2^2 \rightarrow \kappa, \qquad
\|\eta_m^{(2)}\|_2^2 \rightarrow \ell-\kappa, \qquad
\|\eta_m-\eta_m^{(1)}-\eta_m^{(2)}\|_2 \rightarrow 0$$
as $m \rightarrow \infty$ and satisfy the bounds
$$\sup_{m \in {\mathbb N}} \|\eta_m^{(1)}\|_2<M, \quad \sup_{m \in {\mathbb N}} \|\eta_m^{(2)}\|_2<M, \quad
\sup_{m \in {\mathbb N}} \|\eta_m^{(1)}+\eta_m^{(2)}\|_2 < M.$$
\item
The limits $\lim_{m \rightarrow \infty} \LL(\eta_m^{(1)})$ and $\lim_{m \rightarrow \infty} \LL(\eta_m^{(2)})$
are positive.
\item
The functionals $\GG$, $\KK$ and $\LL$ satisfy
$$\begin{Bmatrix} \GG \\ \KK \\ \LL \end{Bmatrix}
(\eta_m)
-\begin{Bmatrix} \GG \\ \KK \\ \LL \end{Bmatrix}
(\eta_m^{(1)})
-\begin{Bmatrix} \GG \\ \KK \\ \LL \end{Bmatrix}
(\eta_m^{(2)}) \rightarrow 0,$$
$$\left\|\begin{Bmatrix} \GG^\prime \\ \KK^\prime \\ \LL^\prime \end{Bmatrix}
(\eta_m)
-\begin{Bmatrix} \GG^\prime \\ \KK^\prime \\ \LL^\prime \end{Bmatrix}
(\eta_m^{(1)})
-\begin{Bmatrix} \GG^\prime \\ \KK^\prime \\ \LL^\prime \end{Bmatrix}
(\eta_m^{(2)})\right\|_0 \rightarrow 0$$
as $m \rightarrow \infty$.
\item
The sequences $\{\eta_m\}$, $\{\eta_m^{(1)}\}$ and $\{\eta_m^{(2)}\}$ satisfy
$$\lim_{m \rightarrow \infty} \JJ_\mu(\eta_m) = \lim_{m \rightarrow \infty} \JJ_{\mu^{(1)}}(\eta_m^{(1)})+
\lim_{m \rightarrow \infty} \JJ_{\mu^{(2)}}(\eta_m^{(2)}),$$
$$\lim_{m \rightarrow \infty} \JJ_\mu^\prime(\eta_m) = \lim_{m \rightarrow \infty} \JJ_{\mu^{(1)}}^\prime(\eta_m^{(1)})+
\lim_{m \rightarrow \infty} \JJ_{\mu^{(2)}}^\prime(\eta_m^{(2)}),$$
where
$$
\mu^{(1)}=\alpha^{(1)}(\mu+\lim_{m\rightarrow\infty}\GG(\eta_m))-\lim_{m\rightarrow\infty} \GG(\eta_m^{(1)}), \qquad
\mu^{(2)}=\alpha^{(2)}(\mu+\lim_{m\rightarrow\infty}\GG(\eta_m))-\lim_{m\rightarrow\infty} \GG(\eta_m^{(2)})
$$
and the positive numbers $\alpha^{(1)}$, $\alpha^{(2)}$ are defined by
$$
\alpha^{(1)}=\frac{\displaystyle\lim_{m \rightarrow \infty} \LL(\eta_m^{(1)})}{\displaystyle\lim_{m \rightarrow \infty} \vphantom{\LL^2}\LL(\eta_m)}, \qquad
\alpha^{(2)}=\frac{\displaystyle\lim_{m \rightarrow \infty} \LL(\eta_m^{(2)})}{\displaystyle\lim_{m \rightarrow \infty} \vphantom{\LL^2}\LL(\eta_m)}.
$$
\item
The sequence $\{\eta_m^{(1)}\}$ converges weakly in $H^2(\R)$ and
strongly  in $H^r(\R)$ for $r \in [0,2)$, to a function $\eta^{(1)} \in H^2(\R)$ with
$\|\eta^{(1)}\|_2^2 \leq D\KK(\eta^{(1)})$ and $\|\eta^{(1)}\|_1 \geq c|\mu|^{2r^\star}$.
\item
The sequence $\{\eta_m^{(2)}\}$ is a minimising sequence for the
functional
$\JJ_{\rho_2,\mu^{(2)}}: H^2(\R) \rightarrow \R \cup \{\infty\}$ defined by
$$\JJ_{\rho_2,\mu^{(2)}}(\eta) = \left\{\begin{array}{lll}\displaystyle \KK(\eta)+\frac{(\mu^{(2)}+\GG(\eta))^2}{\LL(\eta)}
+ \rho_2(\|\eta\|_2^2), & & \eta \in U_2\sm\{0\}, \\
\\
\infty, & & \eta \not\in U_2\sm\{0\},
\end{array}\right.$$
where
$$U_2 = \{\eta \in H^2(\R): \|\eta\|_2^2 \leq M^2 - \|\eta^{(1)}\|_2^2\}, \qquad
\rho_2(\|\eta\|_2^2) =\rho(\|\eta^{(1)}\|_2^2 + \|\eta\|_2^2).$$
\item
The sequences $\{\eta_m\}$ and $\{\eta_m^{(2)}\}$ satisfy
$$
\lim_{m \rightarrow \infty} \rho(\|\eta_m\|_2^2) = \lim_{m \rightarrow \infty} \rho_2(\|\eta_m^{(2)}\|_2^2),
$$
$$
\lim_{m \rightarrow \infty} \JJ_{\rho,\mu}(\eta_m) = \JJ_{\mu^{(1)}}(\eta^{(1)})+
\lim_{m \rightarrow \infty} \JJ_{\rho_2,\mu^{(2)}}(\eta_m^{(2)})
$$
and
$$
\|\eta^{(1)}\|_2^2 + \lim_{m \rightarrow \infty} \|\eta_m^{(2)}\|_2^2 \leq \lim_{m \rightarrow \infty} \|\eta_m\|_2^2
$$
with equality if $\lim_{m \rightarrow \infty} \rho(\|\eta_m\|_2^2) > 0$.
\end{list}
\end{lemma}
{\bf Proof.} For part (i) see Buffoni \emph{et al.} \cite[Lemma 3.10(i), (ii)]{BuffoniGrovesSunWahlen13}.

Turning to part (ii), observe that $\LL(\eta_m^{(1)}) \rightarrow 0$ as $m \rightarrow \infty$ implies that
$\|\eta_m^{(1)}\|_{1/2} \rightarrow 0$ and hence $\|\eta_m^{(1)}\|_1 \rightarrow 0$ as $m \rightarrow \infty$,
which contradicts part (i). The same argument shows that 
$\LL(\eta_m^{(2)}) \not \rightarrow 0$ as $m\rightarrow \infty$.
Because the derivative of $\GG$ is bounded on $U$, we find that
$$|\GG(\eta_m)-\GG(\eta_m^{(1)}+\eta_m^{(2)})|\ \leq\ c\|\eta_m-\eta_m^{(1)}-\eta_m^{(2)}\|_2\ \rightarrow\ 0$$
(see part (i)) and therefore that
$$
\GG(\eta_m)-\GG(\eta_m^{(1)})-\GG(\eta_m^{(2)})\ =\ 
\underbrace{\GG(\eta_m) - \GG(\eta_m^{(1)}+\eta_m^{(2)})}_{\displaystyle =o(1)}
+ \underbrace{\GG(\eta_m^{(1)}+\eta_m^{(2)}) - \GG(\eta_m^{(1)}) - \GG(\eta_m^{(2)})}_{\displaystyle =o(1)}$$
as $m \rightarrow \infty$, in which Theorem \ref{Main splitting theorem} has been used. The same argument
applies to $\KK$ and $\LL$ and establishes part (iii).

Part (iv) follows from part (iii) by a direct calculation (cf.\ Buffoni \emph{et al.} \cite[Corollary 3.11]{BuffoniGrovesSunWahlen13});
for parts (v), (vi) and (vii) see Buffoni \emph{et al.} \cite[Lemmata 3.12, 3.15(i), 3.15(ii)]{BuffoniGrovesSunWahlen13}.\qed

\subsubsection{Iteration} \label{Iteration of cc}

The next step is to apply the concentration-compactness principle to the sequence $\{u_{2,m}\}$
given by
$$
u_{2,m} = \eta_{2,m}^{\prime 2} + \eta_{2,m}^2,
$$
where $\eta_{2,m}=\eta_m^{(2)}$,
and repeat the above analysis.
We proceed iteratively in this fashion, writing $\{\eta_m\}$, $\mu$ and $U$ in iterative
formulae as respectively $\{\eta_{1,m}\}$, $\mu_1$ and $U_1$. The following lemma describes
the result of one step in this procedure (see Buffoni \emph{et al.} \cite[\S3.3]{BuffoniGrovesSunWahlen13}).

\begin{lemma} \label{Results of iteration}
Suppose there exist functions $\eta^{(1)}$, \ldots, $\eta^{(k)} \in H^2(\R)$
and a sequence $\{\eta_{{k+1},m}\}
\subset H^2(\R)$ with the following properties.
\begin{list}{(\roman{count})}{\usecounter{count}}
\item
The sequence $\{\eta_{{k+1},m}\}$ is a minimising sequence for
$\JJ_{\rho_{k+1},\mu_{k+1}}: H^2(\R) \rightarrow \R \cup \{\infty\}$ defined by
$$\JJ_{\rho_{k+1},\mu_{k+1}}(\eta) = \left\{\begin{array}{lll}\displaystyle \KK(\eta)+\frac{(\mu_{k+1}+\GG(\eta))^2}{\LL(\eta)}
+ \rho_{k+1}(\|\eta\|_2^2), & & \eta \in U_{k+1}\sm\{0\}, \\
\\
\infty, & & \eta \not\in U_{k+1}\sm\{0\},
\end{array}\right.$$
where
$$U_{k+1} = \left\{\eta \in H^2(\R): \|\eta\|_2^2 \leq M^2 - \sum_{j=1}^k\|\eta^{(j)}\|_2^2\right\}$$
and
$$\rho_{k+1}(\|\eta\|_2^2)=\rho\left(\sum_{j=1}^k\|\eta^{(j)}\|_2^2 + \|\eta\|_2^2\right),$$
$$
\mu_{k+1} = \frac{\displaystyle \lim_{m \rightarrow \infty} \LL(\eta_{k+1,m})}{\displaystyle \lim_{m \rightarrow \infty} \LL(\eta_m)}
\Big(\mu+\lim_{m\rightarrow\infty}\GG(\eta_m)\Big)-\lim_{m\rightarrow\infty} \GG(\eta_{k+1,m}).
$$

\item
The functions $\eta^{(1)}$, \ldots, $\eta^{(k)}$ satisfy
$$
0<\|\eta^{(j)}\|_2^2 \leq D\KK(\eta^{(j)}), \qquad j=1,\ldots,k 
$$
and
$$
c_{\rho,\mu} = \sum_{j=1}^k \JJ_{\mu_j^{(1)}} (\eta^{(j)}) + c_{\rho_{k+1},\mu_{k+1}},
$$
where
$$
\mu_j^{(1)} = \frac{\displaystyle\LL(\eta^{(j)})}{\displaystyle \lim_{m \rightarrow \infty} \LL(\eta_m)}
\Big(\mu+\lim_{m\rightarrow\infty}\GG(\eta_m)\Big)-\lim_{m\rightarrow\infty} \GG(\eta^{(j)}), \qquad j=1,\ldots,k
$$
and $c_{\rho_{k+1},\mu_{k+1}} = \inf \JJ_{\rho_{k+1},\mu_{k+1}}$.
\item
The sequences $\{\eta_m\}$, $\{\eta_{k+1,m}\}$ and functions $\eta^{(1)}$,
\ldots, $\eta^{(k)}$ satisfy
$$
\sum_{j=1}^k 
\begin{Bmatrix} \GG \\ \KK \\ \LL \end{Bmatrix}
(\eta^{(j)}) + \lim_{m \rightarrow \infty}
\begin{Bmatrix}\GG \\ \KK \\ \LL \end{Bmatrix}
(\eta_{k+1,m}) = \lim_{m \rightarrow \infty}
\begin{Bmatrix} \GG \\ \KK \\ \LL \end{Bmatrix}
(\eta_m),
$$
$$
\lim_{m \rightarrow \infty} \rho(\|\eta_m\|_2^2) = \lim_{m \rightarrow \infty} \rho_{k+1}(\|\eta_{k+1,m}\|_2^2)
$$
and
$$
\sum_{j=1}^k \|\eta^{(j)}\|_2^2
+ \lim_{m \rightarrow \infty} \|\eta_{k+1,m}\|_2^2 \leq \lim_{m \rightarrow \infty} \|\eta_m\|_2^2
$$
with equality if $\lim_{m \rightarrow \infty} \rho(\|\eta_m\|_2^2) > 0$.
\end{list}
Precisely one of the following phenomena occurs.

\begin{enumerate}
\item
There exists a sequence
$\{x_{k+1,m}\} \subset \R$ and a subsequence of
$\{\eta_{k+1,m}(\cdot+x_{k+1,m})\}$ which satisfies
$$\lim_{m \rightarrow \infty} \|\eta_{k+1,m}(\cdot+x_{k+1,m})\|_2^2 \leq \tilde{M}^2 - \sum_{j=1}^k \|\eta^{(j)}\|_2^2$$
and converges in $H^r(\R)$ for $r \in [0,2)$. The limiting function $\eta^{(k+1)}$ satisfies
$$\sum_{j=1}^{k+1}
\begin{Bmatrix} \GG \\ \KK \\ \LL \end{Bmatrix}
(\eta^{(j)}) = \lim_{m \rightarrow \infty}
\begin{Bmatrix} \GG \\ \KK \\ \LL \end{Bmatrix}
(\eta_m),
$$
$$0<\|\eta^{(k+1)}\|_2^2 \leq D\KK(\eta^{(k+1)}), \qquad c_{\rho,\mu} = \sum_{j=1}^{k+1} \JJ_{\mu_j^{(1)}}(\eta^{(j)}),$$
with $\mu_{k+1}^{(1)} =\mu_{k+1}$,
minimises $\JJ_{\rho_{k+1},\mu_{k+1}}$ and minimises $\JJ_{\mu_{k+1}^{(1)}}$ over
$\tilde{U}_{k+1}\sm\{0\}$, where
$$\tilde{U}_{k+1} = \left\{\eta \in H^2(\R): \|\eta\|_2^2 \leq \tilde{M}^2 - \sum_{j=1}^k\|\eta^{(j)}\|_2^2\right\}.$$
The step concludes the iteration.
\item
There exist sequences $\{\eta_{k+1,m}^{(1)}\}$, $\{\eta_{k+1,m}^{(2)}\}$ with the following properties.
\begin{list}{(\roman{count})}{\usecounter{count}}
\item
The sequence $\{\eta_{k+1,m}^{(1)}\}$ converges in $H^r(\R^2)$ for $r \in [0,2)$, to a function $\eta^{(k+1)}$
which satisfies the estimates
$$0<\|\eta^{(k+1)}\|_2^2 \leq D\KK(\eta^{(k+1)}), \qquad \|\eta^{(k+1)}\|_2 \geq c|\mu|_{k+1}^{2r^\star}.$$
\item
The sequence $\{\eta_{k+1,m}^{(2)}\}$ is a minimising sequence for
$\JJ_{\rho_{k+2},\mu_{k+1}^{(2)}}: H^2(\R) \rightarrow \R \cup \{\infty\}$ defined by
$$\JJ_{\rho_{k+2},\mu_{k+1}^{(2)}}(\eta) = \left\{\begin{array}{lll}\displaystyle \KK(\eta)+\frac{(\mu_{k+1}^{(2)}+\GG(\eta))^2}{\LL(\eta)}
+ \rho_{k+2}(\|\eta\|_2^2), & & \eta \in U_{k+2}\sm\{0\}, \\
\\
\infty, & & \eta \not\in U_{k+2}\sm\{0\},
\end{array}\right.$$
where
$$U_{k+2} = \left\{\eta \in H^2(\R): \|\eta\|_2^2 \leq M^2 - \sum_{j=1}^{k+1}\|\eta^{(j)}\|_2^2\right\}$$
and
$$\rho_{k+2}(\|\eta\|_2^2)=\rho\left(\sum_{j=1}^{k+1}\|\eta^{(j)}\|_2^2 + \|\eta\|_2^2\right),
$$

$$
\mu_{k+1}^{(2)} = \frac{\displaystyle \lim_{m \rightarrow \infty} \LL(\eta_{k+1,m}^{(2)})}{\displaystyle \lim_{m \rightarrow \infty} \LL(\eta_m)}
\Big(\mu+\lim_{m\rightarrow\infty}\GG(\eta_m)\Big)-\lim_{m\rightarrow\infty} \GG(\eta_{k+1,m}^{(2)});
$$
furthermore
$$c_{\rho,\mu} = \sum_{j=1}^{k+1} \JJ_{\mu_j^{(1)}} (\eta^{(j)}) + c_{\rho_{k+2},\mu_{k+1}^{(2)}},$$
where
$$\mu_{k+1}^{(1)} = \mu \frac{\LL(\eta^{(k+1)})}{\displaystyle \lim_{m \rightarrow \infty} \LL(\eta_m)}, \qquad
c_{\rho_{k+2},\mu_{k+1}^{(2)}} = \inf \JJ_{\rho_{k+2},\mu_{k+1}^{(2)}}.$$
\item
The sequences $\{\eta_m\}$, $\{\eta_{k+1,m}^{(2)}\}$ and functions $\eta^{(1)}$,
\ldots, $\eta^{(k+1)}$ satisfy
$$
\sum_{j=1}^k
\begin{Bmatrix} \GG \\ \KK \\ \LL \end{Bmatrix}
(\eta^{(j+1)}) + \lim_{m \rightarrow \infty}
\begin{Bmatrix} \GG \\ \KK \\ \LL \end{Bmatrix}
(\eta_{k+1,m}^{(2)}) = \lim_{m \rightarrow \infty}
\begin{Bmatrix} \GG \\ \KK \\ \LL \end{Bmatrix}
(\eta_m),
$$
$$
\lim_{m \rightarrow \infty} \rho(\|\eta_m\|_2^2) = \lim_{m \rightarrow \infty} \rho_{k+2}(\|\eta_{k+1,m}^{(2)}\|_2^2)
$$
and
$$
\sum_{j=1}^{k+1} \|\eta^{(j)}\|_2^2
+ \lim_{m \rightarrow \infty} \|\eta_{k+1,m}^{(2)}\|_2^2 \leq \lim_{m \rightarrow \infty} \|\eta_m\|_2^2
$$
with equality if $\lim_{m \rightarrow \infty} \rho(\|\eta_m\|_2^2) > 0$.
\end{list}
\end{enumerate}
The iteration continues to the next step with $\eta_{k+2,m} = \eta_{k+1,m}^{(2)}$, $m \in {\mathbb N}$.
\end{lemma}

The above construction does not assume that the iteration terminates (that is `concentration'
occurs after a finite number of iterations). If it does not terminate
we let $k \rightarrow \infty$ in Lemma \ref{Results of iteration} and find that
$\|\eta^{(k)}\|_2 \rightarrow 0$ (because
$$\sum_{j=1}^k \|\eta^{(j)}\|_2^2\ \leq\ D\sum_{j=1}^k \KK(\eta^{(j)})\ \leq\ D\sum_{j=1}^k \JJ_{\mu_j^{(1)}}(\eta^{(j)})
\ <\ D c_{\rho,\mu}\ <\ 2D\nu_0^\mu|\mu|$$
for each $k \in {\mathbb N}$, so that the series $\sum_{j=1}^\infty \|\eta^{(j)}\|_2^2$ converges),
$\mu_k \rightarrow 0$ (because $\|\eta^{(k)}\|_2^2 \geq c|\mu_k|^{2r^\star}$), $c_{\rho_k,\mu_k} \rightarrow 0$
(because $c_{\rho_k,\mu_k} < 2\nu_0^{\mu_k}|\mu_k|$) and
$$c_{\rho,\mu} = \sum_{j=1}^\infty \JJ_{\mu_j^{(1)}}(\eta^{(j)}).$$
For completeness we record the following corollary
of Lemma \ref{Results of iteration} which is not used in the remainder of the paper (cf.\ Buffoni \emph{et al.}
\cite[Corollary 3.17]{BuffoniGrovesSunWahlen13}).
\begin{corollary}
Every minimising sequence $\{\eta_m\}$ for $\JJ_{\rho,\mu}$ satisfies
$\lim_{m \rightarrow \infty} \|\eta_m\|_2 \leq \tilde{M}$.
\end{corollary}

\subsection{Construction of the special minimising sequence} \label{Special minimising sequence}

The sequence $\{\tilde{\eta}_m\}$ advertised in Theorem \ref{Special MS theorem}
is constructed by gluing together
the functions $\eta^{(j)}$ identified in Section \ref{Iteration of cc} above with increasingly large distances between them
(the index $j$ is taken between $1$ and $k$, where $k=\infty$ if the iteration does not terminate).
The minimal distance between the functions is chosen so that the interaction between
the `tails' of the indiviual functions is negligable and  $\|\tilde{\eta}_m\|_2^2$ is approximately
$\sum_{j=1}^k \|\eta^{(j)}\|_2^2 = O(\mu)$
(we return to the original physical setting in which $\mu$ is positive).
The algorithm
is stated precisely in part (ii) of the following proposition (which follows immediately from part (i)); for the proof of
part (i) see Buffoni \emph{et al.} \cite[Proposition 3.20]{BuffoniGrovesSunWahlen13}. 
\begin{proposition}
\quad
\begin{list}{(\roman{count})}{\usecounter{count}}
\item
There exists a constant $C>0$ such that
$$\left\|\sum_{j=1}^k \tau_{S_j} \eta^{(j)}\right\|_2^2 \leq 2C^2D\nu_0\mu,$$
where $(\tau_X \eta^{(j)})(x) := \eta^{(j)}(x+X)$,
for all choices of $\{S_j\}_{j=1}^k$. Moreover, in the case $k=\infty$ the series converges uniformly over all such sequences. 
\item
The sequence $\{\tilde{\eta}_m\}$ defined by the following algorithm satisfies $\|\tilde{\eta}_m\|_3^2 \leq 2C^2D\nu_0\mu$.
\begin{enumerate}
\item
Choose $R_j>1$ large enough so that
$$\|\eta^{(j)}\|_{H^2(|x|>R_j)} < \frac{\mu}{2^j}.$$
\item
Write $S_1=0$ and choose $S_j>S_{j-1}+2R_j+2R_{j-1}$ for $j=2,\ldots,k$.
\item
Define
$$\tilde{\eta}_m = \sum_{j=1}^k \tau_{S_j+(j-1)m}\eta^{(j)}, \qquad m \in {\mathbb N}.$$
\end{enumerate}
\end{list}
\end{proposition}

Observe that a local, translation-invariant, analytic operator $\TT: U \rightarrow \R$ has the property that
$$\lim_{m \rightarrow \infty} \TT(\tilde{\eta}_m) = \sum_{j=1}^k \TT(\eta^{(j)}).$$
Part (i) of the next lemma states that the functionals $\GG$, $\KK$ and $\LL$ behave in the same fashion
(with corresponding estimates for their $L^2$-gradients); it is deduced from Theorem \ref{Main splitting theorem}
using the method given by Buffoni \emph{et al.} \cite[Lemma 3.22]{BuffoniGrovesSunWahlen13}.
Part (ii) follows from part (i) by a straightforward calculation which shows that
$$\lim_{m \rightarrow \infty} \JJ_\mu(\tilde{\eta}_m) = \sum_{j=1}^k \JJ_{\mu_j^{(1)}}(\eta^{(j)}), \qquad
\lim_{m \rightarrow \infty} \left\| \JJ^\prime_\mu(\tilde{\eta}_m) - \sum_{j=1}^k \JJ^\prime_{\mu_j^{(1)}}(\eta^{(j)}) \right\|_0 =0
$$
(cf.\ Buffoni \emph{et al.} \cite[Corollary 3.23]{BuffoniGrovesSunWahlen13}).

\begin{lemma}
\quad
\begin{list}{(\roman{count})}{\usecounter{count}}
\item
The sequence $\{\tilde{\eta}_m\}$ and functions $\{\eta^{(i)}\}_{i=1}^m$ satisfy
$$\lim_{m \rightarrow \infty} \begin{Bmatrix} \GG \\ \KK \\ \LL \end{Bmatrix} (\tilde{\eta}_m) = \sum_{i=1}^k \begin{Bmatrix} \GG \\ \KK \\ \LL \end{Bmatrix}(\eta^{(i)}),\qquad
\lim_{m \rightarrow \infty} \left\| \begin{Bmatrix} \GG^\prime \\ \KK^\prime \\ \LL^\prime \end{Bmatrix}(\tilde{\eta}_m)
- \sum_{i=1}^k \begin{Bmatrix} \GG^\prime \\ \KK^\prime \\ \LL^\prime \end{Bmatrix}(\eta^{(i)})\right\|_0=0.$$
\item
The sequence $\{\tilde{\eta}_m\}$ has the properties that
$$\lim_{m \rightarrow \infty} \JJ_\mu(\tilde{\eta}_m) = c_{\rho,\mu}, \qquad
\lim_{m \rightarrow \infty} \|\JJ_\mu^\prime(\tilde{\eta}_m)\|_0=0.$$
\end{list}
\end{lemma}

The proof of Theorem \ref{Special MS theorem} is completed by the following proposition.

\begin{proposition} \label{Special sequence minimises}
The sequence $\{\tilde{\eta}_m\}$ is a minimising sequence for $\JJ_\mu$ over $U \sm \{0\}$.
\end{proposition}
{\bf Proof.} Let us first note that $\{\tilde{\eta}_m\}$ is a minimising sequence for $\JJ_\mu$ over
$\tilde{U}\sm\{0\}$ since the existence of a minimising sequence $\{v_m\}$ for $\JJ_\mu$ over $\tilde{U}\sm\{0\}$
with $\lim_{m \rightarrow \infty} \JJ_\mu(v_m) < \lim_{m \rightarrow \infty} \JJ_\mu(\tilde{\eta}_m)$ would
lead to the contradiction
$$\lim_{m \rightarrow \infty} \JJ_{\rho,\mu}(v_m) = \lim_{m \rightarrow \infty} \JJ_\mu(v_m)
 < \lim_{m \rightarrow \infty} \JJ_\mu(\tilde{\eta}_m) = \lim_{m \rightarrow \infty} \JJ_{\rho,\mu}(\tilde{\eta}_m)
=c_{\rho,\mu}.$$
It follows from this fact and the estimate $\|\tilde{\eta}_m\|_2^2 \leq 2C^2D\nu_0\mu$ that
$$\inf \{\JJ_\mu(\eta): \|\eta\|_2 \in (0,\tilde{M})\} = \inf \{\JJ_\mu(\eta): \|\eta\|_2 \in (0,\sqrt{2C^2D\nu_0\mu})\}$$
for all $\tilde{M} \in (\sqrt{2C^2D\nu_0\mu},M)$. The right-hand side of this equation does not depend
upon $\tilde{M}$; letting $\tilde{M} \rightarrow M$ on the left-hand side, one therefore
finds that
\begin{eqnarray*}
\inf \{\JJ_\mu(\eta): \|\eta\|_2 \in (0,M)\} & = & \inf \{\JJ_\mu(\eta): \|\eta\|_2 \in (0,\sqrt{2C^2D\nu_0\mu})\} \\
& = & \lim_{m \rightarrow \infty} \JJ_\mu(\tilde{\eta}_m).
\end{eqnarray*}
\qed

\section{Strict sub-additivity} \label{SSA section}

The goal of this section is to establish that $c_\mu$
is \emph{strictly sub-additive}, that is
\begin{equation}
c_{\mu_1+\mu_2} < c_{\mu_1} + c_{\mu_2}, \qquad 0<|\mu_1|, |\mu_2|, \mu_1+\mu_2 < \mu_0, \label{Strict SA}
\end{equation}
where negative values of the small parameter are again allowed. This fact is deduced from the
facts that $c_\mu$ is an increasing, \emph{strictly sub-homogeneous} function of $\mu>0$, that is
\begin{equation}
c_{a\mu} < ac_\mu, \qquad a>1. \label{Strict SH}
\end{equation}
The strict sub-homogeneity property of $c_\mu$ is established
by considering a `near minimiser'
of $\JJ_\mu$ over $U\sm\{0\}$, that is a function in $U \sm \{0\}$ with
$$\|\tilde{\eta}\|_2^2 \leq c\mu, \quad \JJ_\mu(\tilde{\eta}) < 2\nu_0\mu-c\mu^{r^\star}, \quad \|\JJ_\mu^\prime(\tilde{\eta})\|_0 \leq \mu^N$$
and hence $\LL(\tilde{\eta}), \LL_2(\tilde{\eta})>c\mu$ (see 
the remark above \eqn{Bound on L from below} and inequality \eqn{L and L2 equivalence}),
and identifying the dominant term in the `nonlinear' part $\MM_\mu(\tilde{\eta})$ of $\JJ_\mu(\tilde{\eta})$.
In Sections \ref{Strong ST} and \ref{Weak ST} below we show that
\begin{equation}
0 > \MM_\mu(\tilde{\eta}) = \left\{\begin{array}{ll} 
\displaystyle c\!\! \int_{-\infty}^\infty \tilde{\eta}_1^3 \dx + o(\mu^\frac{5}{3}), & \beta>\beta_\mathrm{c}, \\
\\
\displaystyle -c\!\! \int_{-\infty}^\infty \tilde{\eta}_1^4 \dx + o(\mu^3),\quad\mbox{} & \beta<\beta_\mathrm{c}, \end{array}\right.
\label{Behaviour of MM}
\end{equation}
where $\eta_1$ is obtained from $\eta \in H^2(\R)$ by multiplying its Fourier transform by
the characteristic function of the 
set $S=[-k_0-\delta_0,-k_0+\delta_0] \cup [k_0-\delta_0,k_0+\delta_0]$ with $\delta_0>0$ if $\beta>\beta_\mathrm{c}$ and
$\delta_0 \in (0,k_0/3)$ if $\beta<\beta_\mathrm{c}$;
inequality \eqn{Strict SH} is readily verified by approximating $\MM(\tilde{\eta}_m)$ by the homogeneous
term identified in \eqn{Behaviour of MM}. The details of this procedure are given in Section \ref{SSA subsection} below.

Straightforward estimates of the kind
$$\GG_j(\tilde{\eta}_m),\ \KK_j(\tilde{\eta}_m),\ \LL_j(\tilde{\eta}_m) = O(\|\tilde{\eta}_m\|_2^j)= O(\mu^{j/2})$$
do not suffice to establish \eqn{Behaviour of MM}.
According to the calculations presented in Appendix A, the function
$\eta^\star_\mu$, which is constructed using the KdV scaling for $\beta > \beta_\mathrm{c}$ and the
nonlinear Schr\"{o}dinger scaling for $\beta<\beta_\mathrm{c}$, satisfies the estimate \eqn{Behaviour of MM}
(with $\tilde{\eta}$ replaced by $\eta^\star_\mu$). The choice of $\eta_\mu^\star$ is of course motivated
by the expectation that a minimiser, and hence any near minimiser,
should have the KdV or nonlinear Schr\"{o}dinger length scales. Our
strategy is therefore to show that $\tilde{\eta}_1$ is $O(\mu^\frac{1}{2})$ with respect to a weighted norm.
To this end we consider the norm
$$
\nn \eta \nn_\alpha^2 := 
\int _{-\infty}^\infty (1+\mu^{-4\alpha}(|k|-k_0)^4) |\hat{\eta}(k)|^2\dk
$$
and choose $\alpha>0$ as large as possible so that $\nn \tilde{\eta}_1 \nn_\alpha$
is $O(\mu^\frac{1}{2})$; this more detailed description of the the behaviour
of $\tilde{\eta}$ allows one to obtain better estimates for
$\GG_j(\tilde{\eta})$, $\KK_j(\tilde{\eta})$ and $\LL_j(\tilde{\eta})$
and thus establish \eqn{Behaviour of MM} (see Sections \ref{Strong ST}
and \ref{Weak ST} for respectively $\beta>\beta_\mathrm{c}$ and $\beta<\beta_\mathrm{c}$).

\subsection{Preliminaries} \label{Preliminaries}

In this section we establish some basic facts which are used Sections \ref{Strong ST}--\ref{SSA subsection}.

\subsubsection{Splitting of \boldmath$\eta$\unboldmath}

In view of the expected frequency distribution of $\tilde{\eta}$ we split each $\eta \in U$ into the sum of a function $\eta_1$ with
spectrum near $k=\pm k_0$ and a function $\eta_2$ whose spectrum is bounded away from these points.
To this end we write the equation
\begin{eqnarray*}
\JJ_\mu^\prime(\eta)
& = & \KK_2^\prime(\eta) + \KK_\mathrm{nl}^\prime(\eta)+2 \left(\frac{\mu+\GG(\eta)}{\LL(\eta)}\right) \GG_2^\prime(\eta)
+2\left(\frac{\mu+\GG(\eta)}{\LL(\eta)}\right) \GG_\mathrm{nl}^\prime(\eta)\\
& & \qquad \mbox{}-\left(\frac{\mu+\GG(\eta)}{\LL(\eta)}\right)^{\!\!2}\LL_2^\prime(\eta)
-\left(\frac{\mu+\GG(\eta)}{\LL(\eta)}\right)^{\!\!2} \LL_\mathrm{nl}^\prime(\eta)
\end{eqnarray*}
\begin{eqnarray*}
& = & \KK_2^\prime(\eta) + 2\nu_0 \GG_2^\prime(\eta) - \nu_0^2 \LL_2^\prime(\eta) \\
& & \qquad\mbox{} + \KK_\mathrm{nl}^\prime(\eta)+2 \left(\frac{\mu+\GG(\eta)}{\LL(\eta)} - \nu_0\right) \GG_2^\prime(\eta)
+2\left(\frac{\mu+\GG(\eta)}{\LL(\eta)}\right) \GG_\mathrm{nl}^\prime(\eta) \\
& & \qquad \mbox{}-\left(\frac{\mu+\GG(\eta)}{\LL(\eta)} + \nu_0\right)\!\!\left(\frac{\mu+\GG(\eta)}{\LL(\eta)} - \nu_0\right)\LL_2^\prime(\eta)
-\left(\frac{\mu+\GG(\eta)}{\LL(\eta)}\right)^{\!\!2} \LL_\mathrm{nl}^\prime(\eta)
\end{eqnarray*}
in the form
\begin{eqnarray*}
\lefteqn{g(k)\hat{\eta}
= \FF\left[ \JJ_\mu^\prime(\eta) - \KK_\mathrm{nl}^\prime(\eta)-2 \left(\frac{\mu+\GG(\eta)}{\LL(\eta)} - \nu_0\right) \GG_2^\prime(\eta)
- 2\left(\frac{\mu+\GG(\eta)}{\LL(\eta)}\right) \GG_\mathrm{nl}^\prime(\eta)\right.} \qquad\qquad \\
& & \left.\mbox{}+\left(\frac{\mu+\GG(\eta)}{\LL(\eta)} + \nu_0\right)\!\!\left(\frac{\mu+\GG(\eta)}{\LL(\eta)} - \nu_0\right)\LL_2^\prime(\eta)
+\left(\frac{\mu+\GG(\eta)}{\LL(\eta)}\right)^{\!\!2} \LL_\mathrm{nl}^\prime(\eta)\right]
\end{eqnarray*}
and decompose it into two coupled equations by defining 
$\eta_2 \in H^2(\R)$ by the formula
\begin{eqnarray*}
\lefteqn{\eta_2 = \FF^{-1}\!\!\left[\frac{1-\chi_S(k)}{g(k)} \FF\!\!\left[ \JJ^\prime(\eta) - \KK_\mathrm{nl}^\prime(\eta)-2 \left(\frac{\mu+\GG(\eta)}{\LL(\eta)} - \nu_0\right) \GG_2^\prime(\eta)
- 2\left(\frac{\mu+\GG(\eta)}{\LL(\eta)}\right) \GG_\mathrm{nl}^\prime(\eta)\right.\right.}\\
& & \hspace{1.25in} \left.\left.\mbox{}+\left(\frac{\mu+\GG(\eta)}{\LL(\eta)} + \nu_0\right)\!\!\left(\frac{\mu+\GG(\eta)}{\LL(\eta)} - \nu_0\right)\LL_2^\prime(\eta)
+\left(\frac{\mu+\GG(\eta)}{\LL(\eta)}\right)^{\!\!2} \LL_\mathrm{nl}^\prime(\eta)\right]\!\right]
\end{eqnarray*}
and $\eta_1 \in H^2(\R)$ by $\eta_1=\eta-\eta_2$, so that $\hat{\eta}_1$ has support in $S$; here we have used the fact that
$$f \mapsto \FF^{-1}\left[\frac{1-\chi_S(k)}{g(k)} \hat{f}(k)\right]$$
is a bounded linear operator $L^2(\R) \rightarrow H^2(\R)$.

\subsubsection{Estimates for $\nn \cdot \nn_\alpha$}

\begin{proposition} \label{Estimates for nn}
\quad
\begin{list}{(\roman{count})}{\usecounter{count}}
\item
The estimates
$\|\eta\|_{1,\infty} \leq c \mu^\frac{\alpha}{2}\nn \eta \nn_\alpha$, $\|K^0\eta\|_\infty \leq c \mu^\frac{\alpha}{2}\nn \eta \nn_\alpha$
hold for each $\eta \in H^2(\R)$.
\item
The estimates
$$\|\eta^{\prime\prime}+k_0^2\eta\|_0 \leq c \mu^{\alpha}\nn \eta\nn_\alpha, \qquad k_0 \neq 0,$$
and
$$\| (K^0\eta)^{(n)}\|_\infty \leq \mu^\frac{\alpha}{2}\nn \eta \nn_\alpha, \quad n=0,1,2,\ldots,$$
hold for each $\eta \in H^2(\R)$ with $\supp \hat{\eta} \subseteq S$.
\end{list}
\end{proposition}\pagebreak
{\bf Proof.} (i) Observe that
\begin{eqnarray}
\|\eta^{(j)}\|_\infty^2 & \leq & c\||k|^j\hat{\eta}\|_{L^1(\R)}, \qquad j=0,1,\label{max norm} \\
\nonumber \\
\|K^0 \eta\|_\infty & \leq & \|(K^0-1)\eta\|_\infty + \|\eta\|_\infty \nonumber \\
& \leq & c(\| (|k|\coth|k|-1)\hat{\eta}\|_{L^1(\R)} + \|\eta\|_\infty) \nonumber \\
& \leq & c(\| |k| \hat{\eta} \|_{L^1(\R)} + \|\hat{\eta}\|_{L^1(\R)}) \label{max norm of K0}
\end{eqnarray}
and
\begin{eqnarray*}
\||k|^j \hat{\eta}\|_{L^1(\R)}^2
& \leq & \left(\int_{-\infty}^\infty \frac{k^{2j}}{1+\mu^{-4\alpha}(k-k_0)^4} \dk\right)
\int_0^\infty (1+\mu^{-4\alpha}(k-k_0)^4)|\hat{\eta}(k)|^2 \dk \\
& & \mbox{}+\left(\int_{-\infty}^\infty \frac{k^{2j}}{1+\mu^{-4\alpha}(k+k_0)^4} \dk\right)
\int_{-\infty}^0 (1+\mu^{-4\alpha}(k+k_0)^4)|\hat{\eta}(k)|^2 \dk \\
& \leq & c\mu^\alpha \nn \eta \nn^2, \qquad j=0,1.
\end{eqnarray*}

(ii) The first result follows from the calculation
\begin{eqnarray*}
\lefteqn{\|\eta^{\prime\prime}+k_0^2\eta\|_0^2} \\ & = & \|(k^2-k_0^2) \hat{\eta}\|_0^2 \\
& \leq & c \left( \int_{k_0-\delta_0}^{k_0+\delta_0} |k-k_0|^2 |\hat{\eta}(k)|^2 \dk + \int_{-k_0-\delta_0}^{-k_0+\delta_0} |k+k_0|^2 |\hat{\eta}(k)|^2 \dk \right) \\
& \leq & c \left( \int_{k_0-\delta_0}^{k_0+\delta_0} (\mu^{2\alpha} + \mu^{-2\alpha} |k-k_0|^4) |\hat{\eta}(k)|^2 \dk
+ \int_{-k_0-\delta_0}^{-k_0+\delta_0} (\mu^{2\alpha} + \mu^{-2\alpha} |k+k_0|^4) |\hat{\eta}(k)|^2 \dk \right) \\
& \leq & c\mu^{2\alpha} \left( \int_{k_0-\delta_0}^{k_0+\delta_0} (1 + \mu^{-4\alpha} |k-k_0|^4) |\hat{\eta}(k)|^2 \dk
+ \int_{-k_0-\delta_0}^{-k_0+\delta_0} (1 + \mu^{-4\alpha} |k+k_0|^4) |\hat{\eta}(k)|^2 \dk \right) \\
& = & c \mu^{2\alpha} \nn \eta \nn_\alpha^2,
\end{eqnarray*}
while the second is established by repeating the proof of the second inequality in part (i) and
estimating $|k| \leq k_0+\delta_0$.\qed

\subsubsection{Estimates for the wave speed}

The following proposition is used in particular to bound the deviation of the quantity\linebreak
$(\mu + \GG(\tilde{\eta}))/\LL(\tilde{\eta})$
(the speed of the corresponding travelling wave when $\tilde{\eta}$ is a minimiser of $\JJ_\mu$ over $U\sm\{0\}$)
from the linear wave speed $\nu_0$.

\begin{proposition} \label{Speed estimate}
The function $\tilde{\eta}$ satisfies the inequalities
$$
\RR_1(\tilde{\eta}) \leq \frac{\mu+\GG(\tilde{\eta})}{\LL(\tilde{\eta})} -\nu_0 \leq \RR_2(\tilde{\eta}),
$$
and
$$
\RR_1(\tilde{\eta}) -\tilde{\MM}_\mu(\tilde{\eta}) \leq \frac{\mu+\GG_2(\tilde{\eta})}{\LL_2(\tilde{\eta})} -\nu_0 \leq \RR_2(\tilde{\eta}) -\tilde{\MM}_\mu(\tilde{\eta})
$$
where
\begin{eqnarray*}
\RR_1(\tilde{\eta}) & = & -\frac{\langle \JJ_\mu^\prime(\tilde{\eta}),\tilde{\eta}\rangle}{4\mu}
+ \frac{1}{4\mu}\big(\langle \MM_\mu^\prime(\tilde{\eta}),\tilde{\eta}\rangle+4\mu\tilde{\MM}_\mu(\tilde{\eta})\big), \nonumber \\
\RR_2(\tilde{\eta}) & = & -\frac{\langle \JJ_\mu^\prime(\tilde{\eta}),\tilde{\eta}\rangle}{4\mu} 
+ \frac{1}{4\mu}\big(\langle \MM_\mu^\prime(\tilde{\eta}),\tilde{\eta}\rangle+4\mu\tilde{\MM}_\mu(\tilde{\eta})\big)
-\frac{\MM_\mu(\tilde{\eta})}{2\mu},
\end{eqnarray*}
and
$$
\tilde{\MM}_\mu(\tilde{\eta}) = \frac{\mu+\GG(\tilde{\eta})}{\LL(\tilde{\eta})}-\frac{\mu+\GG_2(\tilde{\eta})}{\LL_2(\tilde{\eta})}.
$$
\end{proposition}
{\bf Proof.} Taking the scalar product of the equation
$$\JJ_\mu^\prime(\tilde{\eta})=\KK_2^\prime(\tilde{\eta}) - \left(\frac{\mu+\GG_2(\tilde{\eta})}{\LL_2(\tilde{\eta})}\right)^{\!\!2}\LL_2^\prime(\tilde{\eta})
+2\left(\frac{\mu+\GG_2(\tilde{\eta})}{\LL_2(\tilde{\eta})}\right)\GG_2^\prime(\tilde{\eta})
+\MM_\mu^\prime(\tilde{\eta})$$
with $\tilde{\eta}$ yields the identity
$$\frac{\mu+\GG(\tilde{\eta})}{\LL(\tilde{\eta})} = - \frac{\langle \JJ_\mu^\prime(\tilde{\eta}),\tilde{\eta} \rangle}{4\mu}
+ \frac{1}{2\mu}\left(\KK_2(\tilde{\eta}) + \frac{(\mu+\GG_2(\tilde{\eta}))^2}{\LL_2(\tilde{\eta})}\right)
+ \frac{1}{4\mu}\big(\langle \MM_\mu^\prime(\tilde{\eta}),\tilde{\eta} \rangle + 4\mu \tilde{\MM}_\mu(\tilde{\eta}\big).$$
The first inequality is derived by estimating the quantity in brackets from above and below by means of the
estimate
$$2\nu_0\mu \leq \KK_2(\tilde{\eta}) + \frac{(\mu+\GG_2(\tilde{\eta}))^2}{\LL_2(\tilde{\eta})} = \JJ_\mu(\tilde{\eta}) - \MM_\mu(\tilde{\eta}) < 2\nu_0\mu - \MM_\mu(\tilde{\eta})$$
and the second inequality follows directly from the first.\qed

\subsubsection{Estimates for the functionals $\GG$, $\KK$ and $\LL$}

Turning to the functionals $\GG$, $\KK$ and $\LL: U \rightarrow {\mathbb R}$, denote
their non-quadratic parts by $\GG_\mathrm{nl}$, $\KK_\mathrm{nl}$, $\LL_\mathrm{nl}$
and write
$$\GG_\mathrm{nl}(\eta)=\sum_{k=3}^4 \GG_k(\eta)+\GG_\mathrm{r}(\eta),\quad
\KK_\mathrm{nl}(\eta)=\sum_{k=3}^4 \KK_k(\eta)+\KK_\mathrm{r}(\eta),\quad
\LL_\mathrm{nl}(\eta) = \sum_{k=3}^4 \LL_k(\eta) + \LL_\mathrm{r}(\eta),$$
so that
\begin{eqnarray}
\GG_\mathrm{r}(\eta) & = & \frac{\omega}{4} \int_{-\infty}^\infty \eta^2 (K(\eta)-K^0-K^1(\eta))\eta \dx, \label{Formula for GGr} \\
\KK_\mathrm{r}(\eta) & = & \beta \int_{-\infty}^\infty \big(\sqrt{1+\eta^{\prime 2}}-1-{\textstyle\frac{1}{2}}\eta^{\prime 2} + {\textstyle\frac{1}{8}}\eta^{\prime 4}\big)\dx
 - \frac{\omega^2}{2} \int_{-\infty}^\infty \frac{\eta^2}{2} (K(\eta)-K^0)\frac{\eta^2}{2} \dx,\quad \label{Formula for KKr} \\
\LL_\mathrm{r}(\eta) & = & \frac{1}{2}\int_{-\infty}^\infty \eta (K(\eta)-K^0-K^1(\eta)-K^2(\eta))\eta \dx. \label{Formula for LLr}
\end{eqnarray}
We now record useful explicit formulae for the cubic and quartic parts of the functionals in terms of the Fourier-multiplier
operator $K^0$ and give order-of-magnitude estimates for their cubic, quartic and higher-order parts.

\begin{proposition} \label{Formulae for the threes and fours}
The formulae
$$
\GG_3(\eta) = \frac{\omega}{4}\int_{-\infty}^\infty \eta^2 K^0 \eta \dx, \quad
\KK_3(\eta) = \frac{\omega^2}{6} \int_{-\infty}^\infty \eta^3 \dx, \quad
\LL_3(\eta) = \frac{1}{2}\int_{-\infty}^\infty \big(-(K^0 \eta)^2 \eta + \eta^{\prime 2}\eta\big)\dx
$$
and
\begin{eqnarray*}
\GG_4(\eta) & = & \frac{\omega}{2} \eta^2 \eta^{\prime 2} \dx- \frac{\omega}{4} \int_{-\infty}^\infty \eta^2 K^0 (\eta K^0 \eta)\dx , \\
\KK_4(\eta) & = & - \frac{\beta}{8}\int_{-\infty}^\infty \eta^{\prime 4} \dx- \frac{\omega^2}{8} \int_{-\infty}^\infty \eta^2 K^0 \eta^2 \dx, \\
\LL_4(\eta) & = & \frac{1}{2}\int_{-\infty}^\infty \big( K^0(\eta K^0 \eta)\eta K^0 \eta + (K^0 \eta) \eta^2 \eta^{\prime\prime}\big)\dx
\end{eqnarray*}
hold for each $\eta \in U$.
\end{proposition}
{\bf Proof.} The formulae for $\GG_3$ and $\KK_3$, $\KK_4$ follow directly from equations \eqn{Definition of GG} and \eqn{Definition of KK}.

Equations \eqn{Definition of LL} and \eqn{Definition of Lprime} imply that
$$\LL_3(\eta) = \frac{1}{2}\int_{-\infty}^\infty \eta K_1(\eta) \eta \dx, \qquad
\LL_3^\prime(\eta) = \frac{1}{2} \HH_1^\prime(\eta)(\eta,\eta) + K_1(\eta)\eta,$$
while Lemma \ref{Expression for HHprime} shows that
$$\HH_1^\prime(\eta)(\zeta_1,\zeta_2) = - u_{1x}^0u_{2x}^0 + u_{1y}^0u_{2y}^0 \Big|_{y=1} = -(K^0 \zeta_1)(K^0 \zeta_2) + \zeta_1^\prime \zeta_2^\prime,$$
where $u_j$ is the weak solution of \eqn{BC for u 1}--\eqn{BC for u 3}
with $\xi=\zeta_j^\prime$, $j=1,2$, so that
\begin{equation}
\LL_3^\prime(\eta) = - \frac{1}{2}(K^0 \eta)^2 + \frac{1}{2}\eta^{\prime 2} + K_1(\eta)\eta. \label{L3prime v1}
\end{equation}
Taking the inner product of this equation with $\eta$, we therefore find that
$$3\LL_3(\eta) = \frac{1}{2}\int_{-\infty}^\infty (-(K^0 \eta)^2 \eta + \eta^{\prime 2}\eta)\dx + 2\LL_3(\eta),$$
which yields the given formula for $\LL_3(\eta).$

Similarly, equations \eqn{Definition of GG} and \eqn{Definition of Gprime} imply that
$$\GG_4(\eta) = \frac{\omega}{4} \int_{-\infty}^\infty \eta^2 K_1(\eta)\eta \dx$$
and
\begin{eqnarray*}
\GG_4^\prime(\eta) & = & \frac{\omega}{4}\HH_1^\prime(\eta)(\eta^2,\eta) + \frac{\omega}{4} K_1(\eta)\eta^2 + \frac{\omega}{2}\eta K_1(\eta) \eta \\
& = & -\frac{\omega}{4}(K^0 \eta^2)K^0\eta + \frac{\omega}{4}(\eta^2)^\prime \eta^\prime + \frac{\omega}{4} K_1(\eta)\eta^2 + \frac{\omega}{2}\eta K_1(\eta) \eta.
\end{eqnarray*}
The formula for $\GG_4(\eta)$ follows by taking the inner product of the latter equation with $\eta$.

Finally, equations \eqn{Definition of LL} and \eqn{Definition of Lprime} imply that
$$\LL_4(\eta) = \frac{1}{2}\int_{-\infty}^\infty \eta K_2(\eta) \eta \dx, \qquad
\LL_4^\prime(\eta) = \frac{1}{2}\HH_2^\prime(\eta)(\eta,\eta) + K_2(\eta)\eta$$
and Lemma \ref{Expression for HHprime} shows that
$$\HH_2^\prime(\eta)(\zeta_1,\zeta_2) = -u_{1x}^0 u_{2x}^1 - u_{2x}^0 u_{1x}^1 + u_{1y}^0 u_{2y}^1 + u_{2y}^0 u_{1y}^1 
-2\eta u_{1y}^0 u_{2y}^0 \Big|_{y=1}.$$
Using equation \eqn{BC for un 2}, we find that
$$u_y^1|_{y=1}\ =\ G^1\cdot(0,1)\Big|_{y=1}\ =\ -(Q^1 \nabla u^0)\cdot(0,1)\Big|_{y=1} =
\eta u_y^0 + \eta^\prime u_x^0 \Big|_{y=1}\ =\ \eta\zeta^\prime - \eta^\prime K^0\zeta,$$
where $u$ is the weak solution of \eqn{BC for u 1}--\eqn{BC for u 3}
with $\xi=\zeta^\prime$, so that
$$\HH_2^\prime(\eta)(\eta,\eta) = -2 \eta^{\prime 2} K^0 \eta - 2K^0 \eta K^1(\eta) \eta.$$
Equating the expressions \eqn{L3prime v1} and
$$
\LL_3^\prime(\eta) = -K^0(\eta K^0 \eta) - \frac{1}{2}(K^0\eta)^2 - \frac{1}{2}\eta^{\prime 2} - \eta^{\prime\prime}\eta,
$$
which follows from the formula for $\LL_3(\eta)$, we find that
$$K^1(\eta)\eta = - K^0(\eta K^0 \eta) - (\eta^\prime\eta)^\prime,$$
so that
$$\LL_4^\prime(\eta) = -\eta^{\prime 2} K^0 \eta  + K^0\eta K^0(\eta K^0 \eta) + K^0\eta (\eta^\prime\eta)^\prime + K_2(\eta)\eta.$$
The formula for $\LL_4(\eta)$ is obtained by taking the inner product of the this expression with $\eta$.\qed
\begin{proposition} \label{Size of the functionals}
The estimates
\begin{eqnarray*}
\begin{Bmatrix} |\GG_3(\eta)| \\ |\KK_3(\eta)| \\ |\LL_3(\eta)| \end{Bmatrix}
& \leq & c\|\eta\|_2^2 (\|\eta\|_{1,\infty} + \|\eta^{\prime\prime}+k_0^2\eta\|_0), \\
\begin{Bmatrix} |\GG_4(\eta)| \\ |\KK_4(\eta)| \\ |\LL_4(\eta)| \end{Bmatrix}
& \leq & c \|\eta\|_2^2 (\|\eta\|_{1,\infty} + \|\eta^{\prime\prime}+k_0^2\eta\|_0)^2, \\
\begin{Bmatrix} |\GG_\mathrm{r}(\eta)| \\ |\KK_\mathrm{r}(\eta)| \\ |\LL_\mathrm{r}(\eta)| \end{Bmatrix}
& \leq & c\|\eta\|_2^3 (\|\eta\|_{1,\infty} + \|\eta^{\prime\prime}+k_0^2\eta\|_0)^2
\end{eqnarray*}
hold for each $\eta \in U$.
\end{proposition}
{\bf Proof.} These results are obtained by estimating the right-hand sides of the formulae given in Propositions
\ref{Formulae for the threes and fours}
and equations \eqn{Formula for GGr}--\eqn{Formula for LLr} using Proposition \ref{Analyticity of K 3}.\qed
\begin{proposition} \label{Size of the gradients}
The estimates
\begin{eqnarray*}
\begin{Bmatrix} \|\GG_3^\prime(\eta)\|_0 \\ \|\KK_3^\prime(\eta)\|_0 \\ \|\LL_3^\prime(\eta)\|_0 \end{Bmatrix}
& \leq & c\|\eta\|_2 (\|\eta\|_{1,\infty} + \|\eta^{\prime\prime}+k_0^2\eta\|_0 + \|K^0 \eta\|_\infty), \\
\begin{Bmatrix} \|\GG_4^\prime(\eta)\|_0 \\ \|\KK_4^\prime(\eta)\|_0 \\ \|\LL_4^\prime(\eta)\|_0 \end{Bmatrix}
& \leq & c\|\eta\|_2 (\|\eta\|_{1,\infty} + \|\eta^{\prime\prime}+k_0^2\eta\|_0 + \|K^0 \eta\|_\infty)^2, \\
\begin{Bmatrix} \|\GG_\mathrm{r}^\prime(\eta)\|_0 \\ \|\KK_\mathrm{r}^\prime(\eta)\|_0 \\ \|\LL_\mathrm{r}^\prime(\eta)\|_0 \end{Bmatrix}
& \leq & c\|\eta\|_2^2  (\|\eta\|_{1,\infty} + \|\eta^{\prime\prime}+k_0^2\eta\|_0)^2
\end{eqnarray*}
hold for each $\eta \in U$.
\end{proposition}
{\bf Proof.}
We estimate the right-hand sides of the formulae
\begin{equation}
\GG_3^\prime(\eta) = \frac{\omega}{4}K^0 \eta^2 + \frac{\omega}{2} \eta K^0 \eta, \quad
\KK_3^\prime(\eta) = \frac{\omega^2}{2} \eta^2, \quad
\LL_3^\prime(\eta) = -K^0(\eta K^0 \eta) - \frac{1}{2}(K^0\eta)^2 - \frac{1}{2}\eta^{\prime 2} - \eta^{\prime\prime}\eta,
\label{Formulae for the three gradients}
\end{equation}
\begin{eqnarray*}
\GG_4^\prime(\eta) & = & - \frac{\omega}{4}(K^0 \eta^2) K^0 \eta - \frac{\omega}{4}K^0(\eta K^0 \eta^2) - \omega \eta \eta^{\prime 2}
- \omega \eta^2 \eta^{\prime\prime} - \frac{\omega}{2} \eta K^0(\eta K^0 \eta), \\
\KK_4^\prime(\eta) & = & \frac{3\beta}{2} \eta^{\prime 2} \eta^{\prime\prime} - \frac{\omega^2}{4}\eta^2 K^0 \eta^2, \\
\LL_4^\prime(\eta) & = & -2\eta^{\prime 2} K^0 \eta - 2K^0 \eta K^1(\eta)\eta + K_2(\eta)\eta
\end{eqnarray*}
and
\begin{eqnarray*}
\GG_\mathrm{r}^\prime(\eta) & = & \frac{\omega}{4}(\HH^\prime(\eta)-\HH_1^\prime(\eta))(\eta^2,\eta)
+ \frac{\omega}{4}\big(K(\eta)\!-\!K^0\!-\!K^1(\eta)\big)\eta^2 + \frac{\omega}{2}\eta\big(K(\eta)\!-\!K^0\!-\!K^1(\eta)\big)\eta, \\
\KK_\mathrm{r}^\prime(\eta) & = & \beta \left(1-\frac{3}{2}\eta^{\prime 2}-\frac{1}{(1+\eta^{\prime 2})^\frac{3}{2}}\right)
\eta^{\prime\prime} - \frac{\omega^2}{8}\HH^\prime(\eta)(\eta^2,\eta^2)
-\frac{\omega^2}{2}\eta^2 \big(K(\eta)-K^0\big)\eta, \\
\LL_\mathrm{r}^\prime(\eta) & = & \frac{1}{2}\big(\HH^\prime(\eta)\!-\!\HH_1^\prime(\eta)\!-\!\HH_2^\prime(\eta)\big)(\eta,\eta)
+\big(K(\eta)\!-\!K^0\!-\!K^1(\eta)\!-\!K^2(\eta)\big)\eta
\end{eqnarray*}
using Proposition \ref{Analyticity of K 3} and the estimate
$$\|\HH_{j+1}^\prime(\eta)(\zeta_1,\zeta_2)\|_0 \leq C B^j (\|\eta\|_{1,\infty} + \|\eta^{\prime\prime} + k_0^2 \eta\|_0)^j
\|\zeta_1\|_{3/2} \|\zeta_2\|_{3/2}, \qquad j=0,1,2,\ldots.\eqno{\Box}$$

It is also helpful to write
$$
\KK_3^\prime(\eta) = m_1(\eta,\eta), \quad \GG_3^\prime(\eta) = m_2(\eta,\eta), \quad \LL_3^\prime(\eta) = m_3(\eta,\eta),
$$
where $m_j \in {\mathcal L}_\mathrm{s}^2(H^2(\R), L^2(\R))$, $j=1,2,3$, are defined by
\begin{eqnarray*}
m_1(u_1,u_2) & = & \frac{\omega^2}{2}u_1u_2, \\
m_2(u_1,u_2) & = & \frac{\omega}{4}K^0(u_1u_2) + \frac{\omega}{4}u_1 K^0 u_2 + \frac{\omega}{4}u_2 K^0 u_1, \\
m_3(u_1,u_2) & = & -\frac{1}{2}K^0(u_1 K^0 u_2) - \frac{1}{2}K^0(u_2 K^0 u_1)  \\
& & \qquad\mbox{} - \frac{1}{2}K^0 u_1 K^0 u_2 -\frac{1}{2}u_{1x}u_{2x} - \frac{1}{2}u_{1xx}u_2 - \frac{1}{2}u_1u_{2xx},
\end{eqnarray*}
and similarly
$$\KK_3(\eta) = n_1(\eta,\eta,\eta), \quad \GG_3(\eta) = n_2(\eta,\eta,\eta), \quad \LL_3(\eta) = n_3(\eta,\eta,\eta),$$
where $n_j \in {\mathcal L}_\mathrm{s}^3(H^2(\R), \R)$, $j=1,2,3$, are defined by
\begin{eqnarray*}
n_1(u_1,u_2,u_3) & = & \frac{\omega^2}{6}\int_{-\infty}^\infty u_1 u_2 u_3 \dx, \\
n_2(u_1,u_2,u_3) & = & \frac{\omega}{12}\int_{-\infty}^\infty \PP[u_1u_2K^0u_3] \dx, \\
n_3(u_1,u_2,u_3) & = & \frac{1}{6}\int_{-\infty}^\infty \PP[u_1^\prime u_2^\prime u_3] \dx
- \frac{1}{6}\int_{-\infty}^\infty \PP[(K^0u_1)(K^0u_2)u_3] \dx
\end{eqnarray*}
and the symbol $\PP[\cdot]$ denotes the sum of all distinct expressions resulting from permutations of the variables 
appearing in its argument.
\begin{proposition} \label{Estimates for mj and nj}
The estimates
$$
\|m_j(\eta_1,u_2)\|_0 \leq c(\|\eta_1\|_{1,\infty} + \|\eta_1^{\prime\prime}+k_0^2\eta_1\|_0 + \|K^0 \eta_1\|_{1,\infty})
\|u_2\|_2, \qquad j=1,2,3,
$$
and
$$|n_j(\eta_1,u_2,u_3)| \leq c(\|\eta_1\|_{1,\infty} + \|\eta_1^{\prime\prime}+k_0^2\eta_1\|_0 + \|K^0 \eta_1\|_{1,\infty}) \|u_2\|_2 \|u_3\|_2, \qquad j=1,2,3,$$
hold for each $\eta \in U$ and $u_2$, $u_3 \in H^2(\R)$.
\end{proposition}

\subsubsection{Formulae for the functionals $\MM_\mu$ and $\tilde{\MM}_\mu$}

\begin{lemma} \label{Formulae for M}
The estimates
\begin{eqnarray*}
 \MM_\mu(\eta) & = & \KK_3(\eta) + 2 \nu_0 \GG_3(\eta) - \nu_0^2 \LL_3(\eta)
+\KK_4(\eta) + 2 \nu_0 \GG_4(\eta) - \nu_0^2 \LL_4(\eta) \\
& & \mbox{}+2\left(\frac{\mu+\GG_2(\eta)}{\LL_2(\eta)}-\nu_0\right)(\GG_3(\eta)+\GG_4(\eta)) \\
& & \mbox{}-\left(\frac{\mu+\GG_2(\eta)}{\LL_2(\eta)}-\nu_0\right)\!\!\left(\frac{\mu+\GG_2(\eta)}{\LL_2(\eta)}+\nu_0\right)(\LL_3(\eta)+\LL_4(\eta)) \\
& & \mbox{}+\frac{1}{\LL_2(\eta)}\left(\GG_3(\eta)-\left(\frac{\mu+\GG_2(\eta)}{\LL_2(\eta)}\right)\LL_3(\eta)\right)^{\!\!2}
+ O(\mu^\frac{3}{2}(\|\eta\|_{1,\infty} + \|\eta^{\prime\prime}+k_0^2 \eta\|_0)^2),
\end{eqnarray*}
\begin{eqnarray*}
\lefteqn{\langle \MM_\mu(\eta),\eta \rangle+4 \mu \tilde{\MM}_\mu(\eta)} \qquad\\
& = & 3\big(\KK_3(\eta) + 2 \nu_0 \GG_3(\eta) - \nu_0^2 \LL_3(\eta)\big)
+4\big(\KK_3(\eta) + 2 \nu_0 \GG_3(\eta) - \nu_0^2 \LL_3(\eta)\big) \\
& & \mbox{}+2\left(\frac{\mu+\GG_2(\eta)}{\LL_2(\eta)}-\nu_0\right)(3\GG_3(\eta)+4\GG_4(\eta)) \\
& & \mbox{}-\left(\frac{\mu+\GG_2(\eta)}{\LL_2(\eta)}-\nu_0\right)\!\!\left(\frac{\mu+\GG_2(\eta)}{\LL_2(\eta)}+\nu_0\right)(3\LL_3(\eta)+4\LL_4(\eta)) \\
& & \mbox{}+\frac{4}{\LL_2(\eta)}\left(\GG_3(\eta)-\left(\frac{\mu+\GG_2(\eta)}{\LL_2(\eta)}\right)\LL_3(\eta)\right)^{\!\!2} 
+ O(\mu^\frac{3}{2}(\|\eta\|_{1,\infty} + \|\eta^{\prime\prime}+k_0^2 \eta\|_0)^2)
\end{eqnarray*}
and
$$\tilde{\MM}_\mu(\eta) = \mu^{-1} (\GG_3(\eta)+\GG_4(\eta))
+\mu^{-1}\!\left(\frac{\mu+\GG_2(\eta)}{\LL_2(\eta)}\right)\!(\LL_3(\eta)+\LL_4(\eta)) + O(\mu^\frac{1}{2}(\|\eta\|_{1,\infty} + \|\eta^{\prime\prime}+k_0^2\eta\|_0)^2)$$
hold for each $\eta \in U$ with $\|\eta\|_2 \leq c\mu^\frac{1}{2}$ and $\LL_2(\eta)>c\mu$.
\end{lemma}
{\bf Proof.} Using the formulae
$$
\MM_\mu(\eta) = \KK_\mathrm{nl}(\eta) + \frac{(\mu+\GG(\eta))^2}{\LL(\eta)} - \frac{(\mu+\GG_2(\eta))^2}{\LL_2(\eta)}$$
and
$$\frac{1}{\LL(\eta)} = \frac{1}{\LL_2(\eta)}\left(1-\frac{\LL_\mathrm{nl}(\eta)}{\LL(\eta)}\right),$$
one finds that
\begin{eqnarray*}
\MM_\mu(\eta) & = & \KK_\mathrm{nl}(\eta) + 2\left(\frac{\mu+\GG_2(\eta)}{\LL_2(\eta)}\right)\GG_\mathrm{nl}(\eta)
- \left(\frac{\mu+\GG_2(\eta)}{\LL_2(\eta)}\right)^2 \LL_\mathrm{nl}(\eta) \\
& & \mbox{} + \frac{\GG_\mathrm{nl}(\eta)^2}{\LL(\eta)} 
- 2\left(\frac{\mu+\GG_2(\eta)}{\LL_2(\eta)}\right)\frac{\GG_\mathrm{nl}(\eta)\LL_\mathrm{nl}(\eta)}{\LL(\eta)}
+ \left(\frac{\mu+\GG_2(\eta)}{\LL_2(\eta)}\right)^2 \frac{\LL_\mathrm{nl}(\eta)^2}{\LL(\eta)}.
\end{eqnarray*}
We estimate the first line by substituting
$$\begin{Bmatrix} \GG_\mathrm{nl}(\eta) \\ \KK_\mathrm{nl}(\eta) \\ \LL_\mathrm{nl}(\eta) \end{Bmatrix}
=
\begin{Bmatrix}
\GG_3(\eta)+\GG_4(\eta) \\
\KK_3(\eta)+\KK_4(\eta) \\
\LL_3(\eta)+\LL_4(\eta) \\
\end{Bmatrix}
+O(\mu^\frac{3}{2}(\|\eta\|_{1,\infty} + \|\eta^{\prime\prime}+k_0^2 \eta\|_0)^2)
$$
(see Proposition \ref{Size of the functionals}) and
$$\frac{\mu+\GG_2(\eta)}{\LL_2(\eta)}=O(1).$$
Writing
$$\GG_\mathrm{nl}(\eta) = \GG_3(\eta) + O(\mu(\|\eta\|_{1,\infty} + \|\eta^{\prime\prime}+k_0^2 \eta\|_0)^2)$$
(see Proposition \ref{Size of the functionals}) and estimating
$$\GG_3(\eta)\ =\ O(\|\eta\|_\infty \|\eta\|_2^2)\ =\ O(\mu\|\eta\|_\infty)$$
(using the formula for $\GG_3(\eta)$ given in Proposition \ref{Formulae for the threes and fours}) yields
$$\GG_\mathrm{nl}(\eta)^2 = \GG_3(\eta)^2 + O(\mu^2(\|\eta\|_{1,\infty} + \|\eta^{\prime\prime}+k_0^2 \eta\|_0)^3)$$
and
$$\frac{\LL_\mathrm{nl}(\eta)\GG_3(\eta)^2}{\LL_2(\eta)\LL(\eta)} = O(\mu^2(\|\eta\|_{1,\infty} + \|\eta^{\prime\prime}+k_0^2 \eta\|_0)^3)$$
(recall that $\LL(\eta) \geq c\LL_2(\eta)$ for $\eta \in U$), so that
$$\frac{\GG_\mathrm{nl}(\eta)^2}{\LL(\eta)} = \frac{\GG_3(\eta)^2}{\LL_2(\eta)} + O(\mu^\frac{3}{2}(\|\eta\|_{1,\infty} + \|\eta^{\prime\prime}+k_0^2 \eta\|_0)^2);$$
the remaining terms on the second line are estimated in the same fashion.

Altogether we find that
\begin{eqnarray*}
\MM_\mu(\eta) & = & \KK_3(\eta) + 2 \left(\frac{\mu+\GG_2(\eta)}{\LL_2(\eta)}\right)\GG_3(\eta) 
- \left(\frac{\mu+\GG_2(\eta)}{\LL_2(\eta)}\right)^2\LL_3(\eta) \\
& & \mbox{}+\KK_4(\eta) + 2 \left(\frac{\mu+\GG_2(\eta)}{\LL_2(\eta)}\right)\GG_4(\eta) 
- \left(\frac{\mu+\GG_2(\eta)}{\LL_2(\eta)}\right)^2\LL_4(\eta) \\
& & \mbox{} + \frac{1}{\LL_2(\eta)}\left(\GG_3(\eta) - \LL_3(\eta)\left(\frac{\mu+\GG_2(\eta)}{\LL_2(\eta)}\right)\right)^2
+ O(\mu^\frac{3}{2}(\|\eta\|_{1,\infty} + \|\eta^{\prime\prime}+k_0^2 \eta\|_0)^2),
\end{eqnarray*}
from which the stated formula for $\MM_\mu(\eta)$ follows by an algebraic manipulation.

The other estimates are derived by similar calculations.
\qed

\subsection{The case $\beta>\beta_\mathrm{c}$} \label{Strong ST}

We begin by estimating the wave speed.

\begin{proposition} \label{Strong ST speed estimate}
The function $\tilde{\eta}$ satisfies the estimates
$$
\begin{Bmatrix}
\displaystyle \left|\frac{\mu+\GG(\tilde{\eta})}{\LL(\tilde{\eta})} - \nu_0\right| \\
\\
\displaystyle \left|\frac{\mu+\GG_2(\tilde{\eta})}{\LL_2(\tilde{\eta})} - \nu_0\right|
\end{Bmatrix}
\leq c(\|\tilde{\eta}\|_{1,\infty}+\|\tilde{\eta}^{\prime\prime}\|_0 + \mu^{N-\frac{1}{2}}).
$$
\end{proposition}
{\bf Proof.} Proposition \ref{Size of the functionals} implies that
$$
\begin{Bmatrix}
|\GG_j(\tilde{\eta})| \\
|\KK_j(\tilde{\eta})| \\
|\LL_j(\tilde{\eta})| \\
\end{Bmatrix}
 \leq c \mu (\|\tilde{\eta}\|_{1,\infty} + \|\tilde{\eta}^{\prime\prime}\|_0), \qquad j=3,4,
$$
and Lemma \ref{Formulae for M} shows that
$$|\MM_\mu(\tilde{\eta})|,\ |\langle \MM_\mu^\prime(\tilde{\eta}),\tilde{\eta} \rangle+4 \mu \tilde{\MM}_\mu(\tilde{\eta})| \leq c\mu (\|\tilde{\eta}\|_{1,\infty} + \|\tilde{\eta}^{\prime\prime}\|_0), \qquad
|\tilde{\MM}_\mu(\tilde{\eta})| \leq c(\|\tilde{\eta}\|_{1,\infty} + \|\tilde{\eta}^{\prime\prime}\|_0).
$$
The results are obtained by combining these estimates with Proposition \ref{Speed estimate}.\qed

\begin{corollary}
The quantity
\begin{eqnarray*}
\lefteqn{\SS(\tilde{\eta})
= \JJ_\mu^\prime(\tilde{\eta}) - \KK_\mathrm{nl}^\prime(\tilde{\eta})-2 \left(\frac{\mu+\GG(\tilde{\eta})}{\LL(\tilde{\eta})} - \nu_0\right) \GG_2^\prime(\tilde{\eta})
- 2\left(\frac{\mu+\GG(\tilde{\eta})}{\LL(\tilde{\eta})}\right) \GG_\mathrm{nl}^\prime(\tilde{\eta})} \qquad\qquad \\
& & \mbox{}+\left(\frac{\mu+\GG(\tilde{\eta})}{\LL(\tilde{\eta})} + \nu_0\right)\!\!\left(\frac{\mu+\GG(\tilde{\eta})}{\LL(\tilde{\eta})} - \nu_0\right)\LL_2^\prime(\tilde{\eta})
+\left(\frac{\mu+\GG(\tilde{\eta})}{\LL(\tilde{\eta})}\right)^{\!\!2} \LL_\mathrm{nl}^\prime(\tilde{\eta})
\end{eqnarray*}
satisfies
$$\|\SS(\tilde{\eta})\|_0 \leq c (\mu^\frac{1}{2}(\|\tilde{\eta}\|_{1,\infty} + \|\tilde{\eta}^{\prime\prime}\|_0+\|K^0\tilde{\eta}\|_\infty)+\mu^N).$$
\end{corollary}

The next step is an estimate for $\nn \tilde{\eta}_1 \nn_\alpha$ and $\| \tilde{\eta}_2\|_2$.

\begin{lemma} \label{Strong ST weighted norm estimate}
The function $\tilde{\eta}$ satisfies $\nn \tilde{\eta}_1 \nn_\alpha^2 \leq c\mu$ and $\| \tilde{\eta}_2 \|_2^2 \leq c\mu^{2+\alpha}$
 for $\alpha<\frac{1}{3}$.
\end{lemma}
{\bf Proof.} Using the equations
$$g(k)\tilde{\eta}_1 = \FF[\SS(\eta)], \qquad \tilde{\eta}_2 = \FF^{-1}\left[\frac{1-\chi_S(k)}{g(k)}\FF[\SS(\tilde{\eta})]\right],$$
we find from the previous corollary that
$$
\|\tilde{\eta}_2\|_2\ \leq\ c (\mu^\frac{1}{2}(\|\tilde{\eta}_1\|_{1,\infty} + \|\tilde{\eta}_1^{\prime\prime}\|_0
+\|K^0\tilde{\eta}_1\|_\infty)+\mu^\frac{1}{2}\|\tilde{\eta}_2\|_2+\mu^N)
$$
and therefore
\begin{equation}
\|\tilde{\eta}_2\|_2\ \leq\ c (\mu^\frac{1}{2}(\|\tilde{\eta}_1\|_{1,\infty} + \|\tilde{\eta}_1^{\prime\prime}\|_0
+\|K^0\tilde{\eta}_1\|_\infty)+\mu^N),
\label{Strong ST Gives result for eta2}
\end{equation}
and
$$
\int_{-\infty}^\infty g(k)^2 |\tilde{\eta}_1(k)|^2\dk
\leq c (\mu(\|\tilde{\eta}_1\|_{1,\infty} + \|\tilde{\eta}_1^{\prime\prime}\|_0+\|K^0\tilde{\eta}_1\|_\infty)^2+\mu\|\tilde{\eta}_2\|_2^2 +\mu^{2N})
$$
(see Proposition \ref{Estimates for nn}).
Multiplying the above inequality by $\mu^{-4\alpha}$, using \eqn{Strong ST Gives result for eta2}
and adding $\|\tilde{\eta}_1\|_0^2 \leq \|\tilde{\eta}\|_0^2 \leq c\mu$, one finds that
\begin{eqnarray}
\nn \tilde{\eta} \nn_\alpha^2
& \leq & c (\mu^{1-4\alpha}(\|\tilde{\eta}_1\|_{1,\infty} + \|\tilde{\eta}_1^{\prime\prime}\|_0+\|K^0\tilde{\eta}_1\|_\infty)^2+\mu) 
\label{Strong ST Gives result for eta1} \\
& \leq &
c(\mu^{1-3\alpha} \nn \tilde{\eta} \nn_\alpha^2+\mu), \nonumber
\end{eqnarray}
so that $\nn \tilde{\eta} \nn_\alpha^2 \leq c\mu$ for $\alpha<\frac{1}{3}$. The estimate for $\tilde{\eta}_2$ follows from
inequality \eqn{Strong ST Gives result for eta2}.\qed

It remains to identify the dominant terms in the formulae for $\MM_\mu(\tilde{\eta})$ and\linebreak
$\langle \MM_\mu^\prime(\tilde{\eta}), \tilde{\eta} \rangle + 4\mu \tilde{\MM}_\mu(\tilde{\eta})$
given in Lemma \ref{Formulae for M}; this task is accomplished by combining the estimates in
Propositions \ref{Strong ST threes estimate 1},  \ref{Strong ST threes estimate 2} and
Lemma \ref{Strong ST MM estimate} below.

\begin{proposition} \label{Strong ST threes estimate 1}
The function $\tilde{\eta}$ satisfies the estimate
$$
\begin{Bmatrix}
\GG_3(\tilde{\eta}) \\
\KK_3(\tilde{\eta}) \\
\LL_3(\tilde{\eta})
\end{Bmatrix}
=
\begin{Bmatrix}
\GG_3(\tilde{\eta}_1) \\
\KK_3(\tilde{\eta}_1) \\
\LL_3(\tilde{\eta}_1)
\end{Bmatrix}
+o(\mu^{\frac{5}{3}}).
$$
\end{proposition}
{\bf Proof.} Using Proposition \ref{Estimates for mj and nj}, we find that
\begin{eqnarray*}
\left|n_j \left (
\tilde{\eta}_1,
\begin{Bmatrix}
\tilde{\eta}_1\\
\tilde{\eta}_2
\end{Bmatrix},
\tilde{\eta}_2
\right)\right|
& \leq & c \mu^\frac{\alpha}{2} \nn \tilde{\eta}_1 \nn_\alpha
\begin{Bmatrix}
\|\tilde{\eta}_1\|_2\\
\|\tilde{\eta}_2\|_2
\end{Bmatrix}
\|\tilde{\eta}_2\|_2 \\
& \leq & c\mu^{2+\alpha} \\
&  =& o(\mu^\frac{5}{3}),
\end{eqnarray*}
while
$$|n_j(\tilde{\eta}_2, \tilde{\eta}_2, \tilde{\eta}_2)|\ \leq\ c\|\tilde{\eta}_2\|_2^3\ \leq\ c\mu^{3+\frac{3\alpha}{2}}\ =\ o(\mu^\frac{5}{3});$$
it follows that
$$n_j(\tilde{\eta}_1+\tilde{\eta}_2,\tilde{\eta}_1+\tilde{\eta}_2,\tilde{\eta}_1+\tilde{\eta}_2) - n_j(\tilde{\eta}_1,\tilde{\eta}_1,\tilde{\eta}_1) = o(\mu^\frac{5}{3})$$
for $j=1,2,3$.\qed

\begin{proposition} \label{Strong ST threes estimate 2}
The function $\tilde{\eta}$ satisfies the estimate
$$\KK_3(\tilde{\eta}_1) + 2 \nu_0 \GG_3(\tilde{\eta}_1) - \nu_0^2 \LL_3(\tilde{\eta}_1) = \frac{1}{2}\left(\frac{\omega^2}{3}+1\right)\int_{-\infty}^\infty \tilde{\eta}_1^3 \dx+ o(\mu^\frac{5}{3}).$$
\end{proposition}
{\bf Proof.} Note that
\begin{eqnarray*}
\GG_3(\tilde{\eta}_1) & = & \frac{\omega}{4}\int_{-\infty}^\infty \tilde{\eta}_1^3  \dx+ \frac{\omega}{4}\int_{-\infty}^\infty \tilde{\eta}_1^2 (K^0\tilde{\eta}_1-\tilde{\eta}_1)\dx, \\
\\
\KK_3(\tilde{\eta}_1) & = & \frac{\omega^2}{6}\int_{-\infty}^\infty \tilde{\eta}_1^3 \dx,
\end{eqnarray*}
 \begin{eqnarray*}
\LL_3(\tilde{\eta}_1) & = & -\frac{1}{2}\int_{-\infty}^\infty \tilde{\eta}_1^3\dx - \int_{-\infty}^\infty (K^0 \tilde{\eta}_1 - \tilde{\eta}_1)\tilde{\eta}_1^2 \dx  \\
& & \qquad\quad\mbox{}- \frac{1}{2}\int_{-\infty}^\infty (K^0 \tilde{\eta}_1-\tilde{\eta}_1)^2 \tilde{\eta}_1\dx
+\frac{1}{2}\int_{-\infty}^\infty \tilde{\eta}_1^{\prime 2} \tilde{\eta}_1\dx
\end{eqnarray*}
(see Proposition \ref{Formulae for the threes and fours}) and estimate
\begin{eqnarray*}
& &\hspace{-7mm}
\left| \int_{-\infty}^\infty \tilde{\eta}_1^{\prime 2} \tilde{\eta}_1 \dx \right|
 \leq \|\tilde{\eta}_1\|_\infty \|\tilde{\eta}_1^\prime\|_0^2
 \leq c\mu^\frac{5\alpha}{2}\nn \tilde{\eta}_1 \nn_\alpha^3
 \leq c\mu^{\frac{3}{2}+\frac{5\alpha}{2}}
=o(\mu^\frac{5}{3}),\\
& &\hspace{-7mm}
\left| \int_{-\infty}^\infty \tilde{\eta}_1^2(K^0\tilde{\eta}_1-\tilde{\eta}_1)\dx \right|
 \leq \|\tilde{\eta}_1\|_\infty \|\tilde{\eta}_1\|_0 \|K^0 \tilde{\eta}_1 - \tilde{\eta}_1\|_0
 \leq c\mu^{\frac{1}{2}+\frac{5\alpha}{2}}\nn \tilde{\eta}_1 \nn_\alpha^2
 \leq c\mu^{\frac{3}{2}+\frac{5\alpha}{2}}
=o(\mu^\frac{5}{3}),\\
& &\hspace{-7mm}
\left| \int_{-\infty}^\infty \tilde{\eta}_1(K^0\tilde{\eta}_1-\tilde{\eta}_1)^2 \dx \right|
 \leq \|\tilde{\eta}_1\|_\infty \|K^0 \tilde{\eta}_1 - \tilde{\eta}_1\|_0^2
 \leq c\mu^\frac{9\alpha}{2} \nn \tilde{\eta}_1 \nn_\alpha^3
 \leq c\mu^{\frac{3}{2}+\frac{9\alpha}{2}}
=o(\mu^\frac{5}{3}),
\end{eqnarray*}
in which the calculation
$$\|K^0 \eta - \eta\|_0^2\ =\ \int_{-\infty}^\infty (|k|\coth |k|-1)^2|\hat{\eta}(k)|^2 \dk\ \leq\ c \int_{-\infty}^\infty k^4 |\hat{\eta}(k)|^2 \dk
\ =\ c\|\eta^{\prime\prime}\|_0^2\ \leq\ c\mu^{4\alpha}\nn \eta\nn_\alpha^2$$
for $\eta \in H^2(\R)$ has been used. One concludes that
$$\KK_3(\tilde{\eta}_1) + 2 \nu_0 \GG_3(\tilde{\eta}_1) - \nu_0^2 \LL_3(\tilde{\eta}_1) = \frac{1}{2}\Bigg(\frac{\omega^2}{3}+\underbrace{\omega \nu_0 + \nu_0^2}_{\displaystyle =1}\Bigg)\int_{-\infty}^\infty \tilde{\eta}_1^3 + o(\mu^\frac{5}{3}).\eqno{\Box}$$
\begin{lemma} \label{Strong ST MM estimate}
The estimates
\begin{eqnarray*}
& & \MM_{a^2\mu}(a\tilde{\eta}) = a^3\big(\KK_3(\tilde{\eta}) + 2 \nu_0 \GG_3(\tilde{\eta}) - \nu_0^2 \LL_3(\tilde{\eta})\big) + a^3 o(\mu^\frac{5}{3}),\\
& & \langle \MM_{a^2\mu}^\prime(a\tilde{\eta}), a\tilde{\eta} \rangle + 4a^2\mu \tilde{\MM}_{a^2\mu}(a\tilde{\eta})
=3a^3\big(\KK_3(\tilde{\eta}) + 2\nu_0\GG_3(\tilde{\eta})-\nu_0^2 \LL_3(\tilde{\eta})\big)  + a^3o(\mu^\frac{5}{3})
\end{eqnarray*}
hold uniformly over $a \in [1,2]$.
\end{lemma}
{\bf Proof.} Using Lemma \ref{Formulae for M}, the estimates given in Proposition \ref{Size of the functionals} and
$$\frac{\mu+\GG_2(\eta)}{\LL_2(\eta)}=O(1),$$
we find that
\begin{eqnarray*}
\MM_{a^2\mu}(a\tilde{\eta}) & = & a^3\Bigg[\KK_3(\tilde{\eta}) +2\nu_0\GG_3(\tilde{\eta}) -\nu_0^2\LL_3(\tilde{\eta})
+ 2\left(\frac{\mu+\GG_2(\tilde{\eta})}{\LL_2(\tilde{\eta})}-\nu_0\right)\GG_3(\tilde{\eta}) \\
& & \qquad\mbox{}-\left(\frac{\mu+\GG_2(\tilde{\eta})}{\LL_2(\tilde{\eta})}-\nu_0\right)\!\!
\left(\frac{\mu+\GG_2(\tilde{\eta})}{\LL_2(\tilde{\eta})}+\nu_0\right)\LL_3(\tilde{\eta})\Bigg] \\
& & \mbox{}+O(a^4\mu^\frac{3}{2}(\|\tilde{\eta}\|_{1,\infty} + \|\tilde{\eta}^{\prime\prime}\|_0))
\end{eqnarray*}
uniformly over $a \in [1,2]$. The first result follows by estimating
$$
\|\tilde{\eta}\|_{1,\infty} + \|\tilde{\eta}^{\prime\prime}\|_0\ \leq\ c(\mu^\frac{\alpha}{2}\nn \tilde{\eta}_1 \nn_\alpha + \|\tilde{\eta}_2\|_2)
\ \leq\ c\mu^{\frac{1}{2}+\frac{\alpha}{2}},
$$
$$
\frac{\mu+\GG_2(\tilde{\eta})}{\LL_2(\tilde{\eta})} - \nu_0 = O(\mu^{\frac{1}{2}+\frac{\alpha}{2}}),
\qquad
\begin{Bmatrix}
\GG_3(\tilde{\eta}) \\ \LL_3(\tilde{\eta})
\end{Bmatrix}
= O(\mu^\frac{3}{2})
$$
and $a^4 \leq 2a^3$. The second result is derived in a similar fashion.\qed

\begin{corollary} \label{Strong ST final MM estimates}
The estimates
\begin{eqnarray*}
& & \MM_{a^2\mu}(a\tilde{\eta}) = \frac{1}{2}a^3\left(\frac{\omega^2}{3}+1\right) \int_{-\infty}^\infty \tilde{\eta}_1^3 \dx+ a^3 o(\mu^\frac{5}{3}),\\
& & \langle \MM_{a^2\mu}^\prime(a\tilde{\eta}), a\tilde{\eta} \rangle + 4a^2\mu \tilde{\MM}_{a^2\mu}(a\tilde{\eta})
=\frac{3}{2}a^3\left(\frac{\omega^2}{3}+1\right) \int_{-\infty}^\infty \tilde{\eta}_1^3 \dx+ a^3 o(\mu^\frac{5}{3})
\end{eqnarray*}
hold uniformly over $a \in [1,2]$ and
$$\int_{-\infty}^\infty \tilde{\eta}_1^3 \dx \leq -c\mu^\frac{5}{3}.$$
\end{corollary}
{\bf Proof.} The estimates follow by combining Propositions \ref{Strong ST threes estimate 1}
and \ref{Strong ST threes estimate 2} with Lemma \ref{Strong ST MM estimate}, while the inequality for $\tilde{\eta}$ is a consequence of
the first estimate (with $a=1$) and the fact that $\MM_\mu(\tilde{\eta}) \leq -c\mu^\frac{5}{3}$.\qed

\subsection{The case $\beta<\beta_\mathrm{c}$} \label{Weak ST}

\subsubsection{Estimates for near minimisers}

We begin with an observation which shows that the equation for $\eta_1$ may be written as
\begin{equation}
g(k)\hat{\eta}_1 = \chi_S(k) \FF[ \SS(\eta)], \label{eta1 equation}
\end{equation}
where
\begin{eqnarray*}
\lefteqn{\SS(\eta) = \JJ_\mu^\prime(\eta) - \KK_\mathrm{nl}^\prime(\eta)+\KK_3(\eta_1)-2 \left(\frac{\mu+\GG(\eta)}{\LL(\eta)} - \nu_0\right) \GG_2^\prime(\eta)}\\
& & \qquad\quad\left.-\mbox{} 2\left(\frac{\mu+\GG(\eta)}{\LL(\eta)}\right) (\GG_\mathrm{nl}^\prime(\eta)-\GG_3^\prime(\eta))
+\left(\frac{\mu+\GG(\eta)}{\LL(\eta)} + \nu_0\right)\!\!\left(\frac{\mu+\GG(\eta)}{\LL(\eta)} - \nu_0\right)\LL_2^\prime(\eta)\right.\\
& & \qquad\quad\mbox{}
+\left(\frac{\mu+\GG(\eta)}{\LL(\eta)}\right)^{\!\!2} (\LL_\mathrm{nl}^\prime(\eta)-\LL_3^\prime(\eta)).
\end{eqnarray*}

\begin{proposition}
The identity
$$\chi_S \FF\left[\begin{Bmatrix} \GG_3^\prime(\eta_1) \\ \KK_3^\prime(\eta_1) \\ \LL_3^\prime(\eta_1) \end{Bmatrix} \right] =0$$
holds for each $\eta \in U$.
\end{proposition}
{\bf Proof.} Using \eqn{Formulae for the three gradients}, we find that the supports of $\GG_3^\prime(\eta_1)$, $\KK_3^\prime(\eta_1)$
and $\LL_3^\prime(\eta_1)$ lie in the set\linebreak
$[-2k_0-2\delta_0, -2k_0 + 2\delta_0] \cup [-2\delta_0, 2\delta_0] \cup [2k_0-2\delta_0, 2k_0+2\delta_0]$.\qed

In keeping with equation \eqn{eta1 equation} we write the equation for $\eta_2$ in the form
\begin{equation}
\underbrace{\eta_2+H(\eta)}_{\displaystyle :=\eta_3} = \FF^{-1}\left[\frac{1-\chi_S(k)}{g(k)} \FF[\SS(\eta)]\right],
\label{eta3 equation}
\end{equation}
where
\begin{equation}
H(\eta) = \FF^{-1} \left[\frac{1}{g(k)} \FF \left[ \KK_3^\prime(\eta_1) + 2 \left(\frac{\mu+\GG(\eta)}{\LL(\eta)}\right) \GG_3^\prime(\eta_1)-\left(\frac{\mu+\GG(\eta)}{\LL(\eta)}\right)^{\!\!2} \LL_3^\prime(\eta_1)\right]\right];
\label{Definition of H}
\end{equation}
the decomposition $\eta=\eta_1-H(\eta)+\eta_3$ forms the basis of the calculations presented below.
An estimate on the size of $H(\eta)$ is obtained from \eqn{Definition of H} and Proposition \ref{Estimates for mj and nj}.

\begin{proposition} \label{Estimate for H}
The estimate
$$\|H(\eta)\|_2 \leq c (\|\eta_1\|_{1,\infty} + \|\eta_1^{\prime\prime}+k_0^2\eta_1\|_0 + \|K^0 \eta_1\|_{1,\infty} + \|\eta_3\|_2) \|\eta_1\|_2$$
holds for each $\eta \in U$.
\end{proposition}

The above results may be used to derive estimates for the gradients of the cubic parts of the functionals
which are used in the analysis below.

\begin{proposition} \label{Three gradients eta to eta1}
The function $\tilde{\eta}$ satisfies the estimates
$$
\begin{Bmatrix}
\|\GG_3^\prime(\tilde{\eta})-\GG_3^\prime(\tilde{\eta}_1)\|_0 \\
\|\KK_3^\prime(\tilde{\eta})-\KK_3^\prime(\tilde{\eta}_1)\|_0 \\
\|\LL_3^\prime(\tilde{\eta})-\LL_3^\prime(\tilde{\eta}_1)\|_0
\end{Bmatrix}
\leq c \mu^\frac{1}{2}((\|\tilde{\eta}_1\|_{1,\infty} + \|\tilde{\eta}_1^{\prime\prime}+k_0^2\tilde{\eta}_1\|_0 + \|K^0 \tilde{\eta}_1\|_{1,\infty} )^2+ \|\tilde{\eta}_3\|_2).
$$
\end{proposition}
{\bf Proof.} Observe that
$$\GG_3^\prime(\eta)-\GG_3^\prime(\eta_1)=m_2(H(\eta),H(\eta))+m_2(\eta_3,\eta_3) - 2m_2(\eta_1,H(\eta))
-2m_2(\eta_3,H(\eta))+2m_2(\eta_1,\eta_3)$$
and estimate the right-hand side of this equation using Propositions \ref{Estimates for mj and nj}
and \ref{Estimate for H}. The same method yields the results for $\KK_3^\prime$ and $\LL_3^\prime$.\qed

Estimates for $\GG_3(\tilde{\eta})$, $\KK_3(\tilde{\eta})$ and $\LL_3(\tilde{\eta})$ are obtained in a similar fashion.

\begin{proposition} \label{Threes}
The function $\tilde{\eta}$ satisfies the estimates
$$
\begin{Bmatrix}
|\GG_3(\tilde{\eta})| \\
|\KK_3(\tilde{\eta})| \\
|\LL_3(\tilde{\eta})| \\
\end{Bmatrix}
\leq c\big(\mu (\|\tilde{\eta}_1\|_{1,\infty} + \|\tilde{\eta}_1^{\prime\prime}+k_0^2\tilde{\eta}_1\|_0 + \|K^0 \tilde{\eta}_1\|_{1,\infty}) + \mu\|\tilde{\eta}_3\|_2\big) .
$$
\end{proposition}
{\bf Proof.} Observe that
$$
\GG_3(\eta_1)=\frac{1}{3}\langle \GG_3^\prime(\eta_1), \eta_1\rangle =\frac{1}{3}\int_{-\infty}^\infty \FF[\GG_3^\prime(\eta_1)]\overline{\hat{\eta}_1} \dk
=\frac{1}{3}\int_{-\infty}^\infty \underbrace{\chi_S(k)\FF[\GG_3^\prime(\eta_1)]}_{\displaystyle = 0}\overline{\hat{\eta}_1}\dk=0,
$$
(since $\hat{\eta}_1 = \chi_S(k) \hat{\eta}_1$), so that
\begin{eqnarray*}
\GG_3(\eta) & = & \GG_3(\eta)-\GG_3(\eta_1) \\
& = & -n_2(H(\eta),H(\eta),H(\eta)) + n_2(\eta_3,\eta_3,\eta_3) - 6n_2(\eta_1,H(\eta),\eta_3) - 3n_2(\eta_1,\eta_1,H(\eta)) \\
& & \mbox{} + 3n_2(\eta_1,\eta_1,\eta_3) + 3n_2(H(\eta),H(\eta),\eta_3) + 3n_2(H(\eta),H(\eta),\eta_1)
+ 3n_2(\eta_3,\eta_3,\eta_1) \\
& & \mbox{}-3n_2(\eta_3,\eta_3,H(\eta))
\end{eqnarray*}
and estimate the right-hand side of this equation using Propositions
\ref{Estimates for mj and nj} and \ref{Estimate for H}.
The same method yields the results for $\KK_3$ and $\LL_3$.\qed

Estimating the right-hand sides of the inequalities
\begin{eqnarray*}
\|\GG_\mathrm{nl}^\prime(\tilde{\eta})-\GG_3^\prime(\tilde{\eta}_1)\|_0 & \leq &
\|\GG_\mathrm{r}^\prime(\tilde{\eta})\|_0 + \|\GG_4^\prime(\tilde{\eta})\|_0 + \|\GG_3^\prime(\tilde{\eta})-\GG_3^\prime(\tilde{\eta}_1)\|_0, \\
|\GG_\mathrm{nl}(\tilde{\eta})| & \leq & |\GG_\mathrm{r}(\tilde{\eta})|+|\GG_4(\tilde{\eta})| + |\GG_3(\tilde{\eta})|
\end{eqnarray*}
(together with the corresponding inequalities for $\KK$ and $\LL$)
using Propositions \ref{Size of the functionals} and \ref{Size of the gradients}, the calculation
\begin{eqnarray}
\lefteqn{\|\eta\|_{1,\infty} + \|\eta^{\prime\prime}+k_0^2\eta\|_0 + \|K^0 \eta\|_\infty} \quad \nonumber \\
& & \leq c(\|\eta_1\|_{1,\infty} + \|\eta_1^{\prime\prime}+k_0^2\eta_1\|_0 + \|K^0 \eta_1\|_\infty + \|H(\eta)\|_2 + \|\eta_3\|_2) \nonumber \\
& & \leq c(\|\eta_1\|_{1,\infty} + \|\eta_1^{\prime\prime}+k_0^2\eta_1\|_0 + \|K^0 \eta_1\|_{1,\infty} + \|\eta_3\|_2). \label{General estimate pt 1}
\end{eqnarray}
and Propositions \ref{Three gradients eta to eta1} and \ref{Threes}
yields the following estimates for the `nonlinear' parts of the functionals.

\begin{lemma} \label{Estimates for nl}
The function $\tilde{\eta}$ satisfies the estimates
\begin{eqnarray*}
& & \begin{Bmatrix}
\|\GG_\mathrm{nl}^\prime(\tilde{\eta})-\GG_3^\prime(\tilde{\eta}_1)\|_0 \\
\|\KK_\mathrm{nl}^\prime(\tilde{\eta})-\KK_3^\prime(\tilde{\eta}_1)\|_0 \\
\|\LL_\mathrm{nl}^\prime(\tilde{\eta})-\LL_3^\prime(\tilde{\eta}_1)\|_0
\end{Bmatrix}
\leq c\big(\mu^\frac{1}{2}(\|\tilde{\eta}_1\|_{1,\infty} + \|\tilde{\eta}_1^{\prime\prime}+k_0^2\tilde{\eta}_1\|_0)^2 + \|K^0 \tilde{\eta}_1\|_{1,\infty} )^2
+ \mu^\frac{1}{2}\|\tilde{\eta}_3\|_2\big), \\
& & \begin{Bmatrix}
|\GG_\mathrm{nl}(\tilde{\eta})| \\
|\KK_\mathrm{nl}(\tilde{\eta})| \\
|\LL_\mathrm{nl}(\tilde{\eta})|
\end{Bmatrix}
\leq c\big(\mu(\|\tilde{\eta}_1\|_{1,\infty} + \|\tilde{\eta}_1^{\prime\prime}+k_0^2\tilde{\eta}_1\|_0)^2 + \|K^0 \tilde{\eta}_1\|_{1,\infty} )^2
+ \mu\|\tilde{\eta}_3\|_2\big).
\end{eqnarray*}
\end{lemma}

We now have all the ingredients necessary to estimate the wave speed and the
quantity $\nn \tilde{\eta}_1 \nn_\alpha$.

\begin{proposition} \label{Weak ST speed estimate}
The function $\tilde{\eta}$ satisfies the estimates
$$
\begin{Bmatrix}
\displaystyle \left|\frac{\mu+\GG(\tilde{\eta})}{\LL(\tilde{\eta})} - \nu_0\right| \\
\\
\displaystyle \left|\frac{\mu+\GG_2(\tilde{\eta})}{\LL_2(\tilde{\eta})} - \nu_0\right|
\end{Bmatrix}
\leq c\big((\|\tilde{\eta}_1\|_{1,\infty} + \|\tilde{\eta}_1^{\prime\prime}+k_0^2\tilde{\eta}_1\|_0 + \|K^0 \tilde{\eta}_1\|_{1,\infty} )^2+ \|\tilde{\eta}_3\|_2
 + \mu^{N-\frac{1}{2}}\big).
$$
\end{proposition}
{\bf Proof.} Combining Lemma \ref{Formulae for M}, inequality \eqn{General estimate pt 1}
and Lemma \ref{Estimates for nl}, one finds that
$$|\MM(\tilde{\eta})|,\ |\langle \MM^\prime(\tilde{\eta}),\tilde{\eta}\rangle+4\mu\tilde{\MM}_\mu(\tilde{\eta})| \leq
c\big(\mu(\|\tilde{\eta}_1\|_{1,\infty} + \|\tilde{\eta}_1^{\prime\prime}+k_0^2\tilde{\eta}_1\|_0)^2 + \|K^0 \tilde{\eta}_1\|_{1,\infty} )^2
+ \mu\|\tilde{\eta}_3\|_2\big),$$
$$|\tilde{\MM}_\mu(\tilde{\eta})| \leq c\big((\|\tilde{\eta}_1\|_{1,\infty} + \|\tilde{\eta}_1^{\prime\prime}+k_0^2\tilde{\eta}_1\|_0)^2 + \|K^0 \tilde{\eta}_1\|_{1,\infty} )^2
+ \|\tilde{\eta}_3\|_2\big),$$
from which the given estimates follow by Proposition \ref{Speed estimate}.\qed

\begin{lemma} \label{Weak ST nn estimate}
The function $\tilde{\eta}$ satisfies
$\nn \tilde{\eta}_1 \nn_\alpha^2 \leq c\mu$,
$\|\tilde{\eta}_3\|_2^2 \leq c\mu^{3+2\alpha}$ and $\|H(\tilde{\eta})\|_2^2 \leq c\mu^{2+\alpha}$ for $\alpha<1$.
\end{lemma}
{\bf Proof.} Lemma \ref{Estimates for nl} and Proposition \ref{Weak ST speed estimate} assert that
$$\|\SS(\tilde{\eta})\|_0 \leq c\big( \mu^\frac{1}{2} (\|\tilde{\eta}_1\|_{1,\infty} + \|\tilde{\eta}_1^{\prime\prime}+k_0^2\tilde{\eta}_1\|_0 + \|K^0 \tilde{\eta}_1\|_{1,\infty} )^2+ \mu^\frac{1}{2}\|\tilde{\eta}_3\|_2
+ \mu^N\big),$$
which inequality shows that
$$\|\tilde{\eta}_3\|_2 \leq c\big( \mu^\frac{1}{2} (\|\tilde{\eta}_1\|_{1,\infty} + \|\tilde{\eta}_1^{\prime\prime}+k_0^2\tilde{\eta}_1\|_0 + \|K^0 \tilde{\eta}_1\|_{1,\infty} )^2+ \mu^\frac{1}{2}\|\tilde{\eta}_3\|_2
+ \mu^N\big)$$
and therefore
\begin{equation}
\|\tilde{\eta}_3\|_2 \leq c\big( \mu^\frac{1}{2} (\|\tilde{\eta}_1\|_{1,\infty} + \|\tilde{\eta}_1^{\prime\prime}+k_0^2\tilde{\eta}_1\|_0 + \|K^0 \tilde{\eta}_1\|_{1,\infty} )^2 + \mu^N\big),
\label{Weak ST Gives result for eta3}
\end{equation}
and
\begin{eqnarray*}
\int_{-\infty}^\infty g(k)^2 |\tilde{\eta}_1|^2 \dk
& \leq & c\big( \mu(\|\tilde{\eta}_1\|_{1,\infty} + \|\tilde{\eta}_1^{\prime\prime}+k_0^2\tilde{\eta}_1\|_0 + \|K^0 \tilde{\eta}_1\|_{1,\infty} )^4+ \mu\|\tilde{\eta}_3\|_2^2
+ \mu^{2N}\big) \\
& \leq & c\big( \mu(\|\tilde{\eta}_1\|_{1,\infty} + \|\tilde{\eta}_1^{\prime\prime}+k_0^2\tilde{\eta}_1\|_0 + \|K^0 \tilde{\eta}_1\|_{1,\infty} )^4
+ \mu^{2N}\big).
\end{eqnarray*}
Multiplying the above inequality by $\mu^{-4\alpha}$ and adding $\|\tilde{\eta}_1\|_0^2 \leq \|\tilde{\eta}\|_0^2 \leq c \mu$, one finds that
\begin{eqnarray}
\nn \tilde{\eta}_1 \nn_\alpha^2 & \leq & c\big( \mu^{1-4\alpha}(\|\tilde{\eta}_1\|_{1,\infty} + \|\tilde{\eta}_1^{\prime\prime}+k_0^2\tilde{\eta}_1\|_0 + \|K^0 \tilde{\eta}_1\|_{1,\infty} )^4
+ \mu\big) \label{Weak ST Gives result for eta1}\\
& \leq & c(\mu^{1-2\alpha}\nn \tilde{\eta}_1 \nn_\alpha ^4 +\mu) \nonumber
\end{eqnarray}
where Proposition \ref{Estimates for nn} and the fact that $g(k) \geq c(|k|-k_0)^2$ for $k \in S$ have also been used.

The estimate for $\tilde{\eta}_1$ follows from the previous inequality
using the argument given by Groves \& Wahl\'{e}n \cite[p.\ 401]{GrovesWahlen10}, while those for
$\tilde{\eta}_3$ and $H(\tilde{\eta})$ are derived by estimating $\nn \tilde{\eta}_1 \nn_\alpha^2 \leq c\mu$
in equation \eqn{Weak ST Gives result for eta3} and Proposition \ref{Estimate for H}.\qed

\subsubsection{Estimates for the variational functional}

The next step is to identify the dominant terms in the formulae for $\MM_\mu(\tilde{\eta})$ and\linebreak
$\langle \MM_\mu^\prime(\tilde{\eta}), \tilde{\eta} \rangle + 4\mu \tilde{\MM}_\mu(\tilde{\eta})$
given in Lemma \ref{Formulae for M}. We begin by examining the quantities
$\GG_4(\tilde{\eta})$, $\KK_4(\tilde{\eta})$ and $\LL_4(\tilde{\eta})$.

\begin{proposition}
The function $\tilde{\eta}$ satisfies the estimates
$$
\begin{Bmatrix}
\GG_4(\tilde{\eta}) \\
\KK_4(\tilde{\eta}) \\
\LL_4(\tilde{\eta})
\end{Bmatrix}
=
\begin{Bmatrix}
\GG_4(\tilde{\eta}_1) \\
\KK_4(\tilde{\eta}_1) \\
\LL_4(\tilde{\eta}_1)
\end{Bmatrix}
+o(\mu^3).
$$
\end{proposition}
{\bf Proof.} Write
$$\KK_4(\eta) = p_1(\eta,\eta,\eta,\eta), \qquad
\GG_4(\eta)=p_2(\eta,\eta,\eta,\eta), \qquad
\LL_4(\eta)=p_3(\eta,\eta,\eta,\eta),$$
where $p_j \in {\mathcal L}_\mathrm{s}^4(H^2(\R), \R)$, $j=1,2,3$, are defined by
\begin{eqnarray*}
p_1(u_1,u_2,u_3,u_4)
& = & - \frac{1}{8}\int_{-\infty}^\infty u_1^\prime u_2^\prime u_3^\prime u_4^\prime  \dx
-\frac{\omega^2}{48}\int_{-\infty}^\infty \PP[u_1u_2 K^0(u_3u_4)] \dx , \\
\\
p_2(u_1,u_2,u_3,u_4)
& = & \frac{\omega}{12} \int_{-\infty}^\infty \PP[u_1u_2 u_3^\prime u_4^\prime]\dx
-\frac{\omega}{48}\int_{-\infty}^\infty \PP[u_1u_2K^0(u_3K^0u_4)] \dx, \\
\\
p_3(u_1,u_2,u_3,u_4)
& = & \frac{1}{24}\int_{-\infty}^\infty \PP[u_1u_2 (K^0u_3) u_4^{\prime\prime}] \dx
+\frac{1}{48} \int_{-\infty}^\infty \PP[ K^0(u_1K^0u_2)u_3K^0u_4] \dx,
\end{eqnarray*}
and estimate each term in the expansion of
$$
p_j(\tilde{\eta}_1-H(\tilde{\eta})+\tilde{\eta}_3,\tilde{\eta}_1-H(\tilde{\eta})+\tilde{\eta}_3,\tilde{\eta}_1-H(\tilde{\eta})+\tilde{\eta}_3,\tilde{\eta}_1-H(\tilde{\eta})+\tilde{\eta}_3) - p_j(\tilde{\eta}_1,\tilde{\eta}_1,\tilde{\eta}_1,\tilde{\eta}_1)$$
for $j=1,2,3$.
Terms with zero, one or two occurrences of $\tilde{\eta}_1$ are estimated by
$$
\left|p_j \left (
\begin{Bmatrix}
\tilde{\eta}_1 \\
H(\tilde{\eta}) \\
\tilde{\eta}_3
\end{Bmatrix}^{\!\!(2)}, 
\begin{Bmatrix}
H(\tilde{\eta}) \\
\tilde{\eta}_3
\end{Bmatrix}^{(2)}\right)\right|
\ \leq \ c
\begin{Bmatrix}
\|\tilde{\eta}_1\|_2 \\
\|H(\tilde{\eta})\|_2 \\
\|\tilde{\eta}_3\|_2
\end{Bmatrix}^2
\begin{Bmatrix}
\|H(\tilde{\eta})\|_2 \\
\|\tilde{\eta}_3\|_2
\end{Bmatrix}^{2} \\
\ \leq \  c\mu\mu^{2+\alpha}\\
\ =\ o(\mu^3),
$$
while terms with three occurrences of $\tilde{\eta}_1$ are estimated by
\begin{eqnarray*}
\left|p_j \left (
\{\tilde{\eta}_1\}^{(3)},
\begin{Bmatrix}
H(\tilde{\eta}) \\
\tilde{\eta}_3
\end{Bmatrix}
\right)\right|
& \leq & c
\begin{Bmatrix}
\|\tilde{\eta}_1\|_\infty \\
\|K^0 \tilde{\eta}_1 \|_{1,\infty} \\
\|\tilde{\eta}_1^{\prime\prime}\|_0
\end{Bmatrix}
\|\tilde{\eta}_1\|_2^2
\begin{Bmatrix}
\|H(\tilde{\eta})\|_2 \\
\|\tilde{\eta}_3\|_2
\end{Bmatrix} \\
& \leq & c\mu^\frac{\alpha}{2}\nn \tilde{\eta}_1 \nn_\alpha \mu \mu^{1+\frac{\alpha}{2}} \\
& \leq & c\mu^{\frac{5}{2}+\alpha} \\
&  =& o(\mu^3).
\end{eqnarray*}

To identify the dominant terms in $\GG_4(\tilde{\eta}_1)$, $\KK_4(\tilde{\eta}_1)$ and
$\LL_4(\tilde{\eta}_1)$ we use the following result, which shows how Fourier-mutliplier operators
acting upon the function $\eta_1$, whose spectrum is concentrated near $k=\pm k_0$, may be
approximated by multiplication by constants.

\begin{lemma} \label{Approximate operators}
For each $\eta \in H^2(\R)$ with $\|\eta\|_2 \leq c\mu^\frac{1}{2}$ the quantities $\eta_1^+:=\FF^{-1}[\chi_{[0,\infty)}\hat \eta_1]$ and
$\eta_1^-:=\FF^{-1}[\chi_{(-\infty,0]}\hat \eta_1]$ (that is $\eta_1^-=\overline{\eta_1^+}$)
satisfy the estimates
\begin{list}{(\roman{count})}{\usecounter{count}}
\item
$\eta_1^{\pm\prime} = \pm \i  k_0\eta_1^\pm + \underline{O}(\mu^{\frac{1}{2}+\alpha})$,
\item
$K^0(\eta_1^\pm) = f(k_0)\eta_1^\pm + \underline{O}(\mu^{\frac{1}{2}+\alpha})$,
\item
$((\eta_1^\pm)^2)^\prime= \pm 2k_0\i (\eta_1^\pm)^2
 + \underline{O}(\mu^{1+\frac{3\alpha}{2}})$,
 \item
$(\eta_1^+\eta_1^-)^\prime = \underline{O}(\mu^{1+\frac{3\alpha}{2}})$,
\item
$K^0((\eta_1^\pm)^2) = f(2k_0)(\eta_1^\pm)^2 + \underline{O}(\mu^{1+\frac{3\alpha}{2}})$,
\item
$K^0(\eta_1^+\eta_1^-) = \eta_1^+\eta_1^-+\underline{O}(\mu^{1+\frac{3\alpha}{2}})$,
 \item
 $\FF^{-1}\left[g(k)^{-1}\FF[(\eta_1^\pm)^2]\right] = g(2k_0)(\eta_1^\pm)^2
 + \underline{O}(\mu^{1+\frac{3\alpha}{2}})$,
 \item
$\FF^{-1}\left[g(k)^{-1}\FF[\eta_1^+\eta_1^-]\right] = g(0)^{-1}\eta_1^+\eta_1^-
+ \underline{O}(\mu^{1+\frac{3\alpha}{2}})$.
\end{list}
Here the symbol $\underline{O}(\mu^\gamma)$ denotes a quantity whose Fourier transform has
compact support and whose $L^2(\R)$-norm (and hence $H^s(\R)$-norm for $s \geq 0$)
is $O(\mu^\gamma)$.
\end{lemma}
{\bf Proof.} Estimates (i) and (ii) follow from the calculations
$$\|(\i k \mp \i k_0)\hat{\eta}_1^\pm\|_0^2 = \|(|k|-k_0)\hat \eta_1\|_0^2, \qquad
\|(K^0-f(k_0))(\eta_1^\pm)\|_0^2 \leq c \|(|k|-k_0)\hat \eta_1\|_0^2$$
(because $f(k)=f(k_0) + O(|k|-k_0)$ for $k \in S$) and
$$\|(|k|-k_0)\hat \eta_1\|_0^2
 \leq \frac{1}{2}\int_{-\infty}^\infty (\mu^{2\alpha} + \mu^{-2\alpha}(|k|-k_0)^4)|\hat{\eta}_1|^2\dk
 \leq c\mu^{2\alpha} 
\nn \eta_1\nn_\alpha^2\le c\mu^{1+2\alpha},$$
while (iii) and (iv) are obtained from the observations
\begin{eqnarray*}
\|(\partial_x \mp 2\i k_0) (\eta_1^\pm)^2\|_0
& = & \|2((\partial_x\mp k_0\i )\eta_1^\pm)\eta_1^\pm\|_0\\
& \leq & 2\|(\partial_x\mp \i k_0 )\eta_1^
\pm \|_0 \|\eta_1^\pm\|_\infty \\
& \leq & c\mu^{\frac{1}{2}+\frac{3\alpha}{2}}\nn \eta_1^\pm\nn_\alpha \\
& \leq & c\mu^{1+\frac{3\alpha}{2}},
\end{eqnarray*}
and
\begin{eqnarray*}
\|(\eta_1^+\eta_1^-)^\prime\|_0 & = & \|((\partial_x-\i k_0)\eta_1^+)\eta_1^- + \eta_1^+((\partial_x+\i k_0)\eta_1^-)\|_0 \\
& \leq & \|(\partial_x-\i k_0)\eta_1^+\|_0\|\eta_1^-\|_\infty+\|\eta_1^+\|_\infty \|(\partial_x+\i k_0)\eta_1^-\|_0\\
& \leq & c\mu^{1+\frac{3\alpha}{2}},
\end{eqnarray*}
in which Proposition \ref{Estimates for nn} has been used. Estimates (v) and (vi) are deduced
from respectively (iii) and (iv) by means of the inequalities
$$\|(K^0-f(2k_0))(\eta_1^\pm)^2\|_0^2 \leq c \|(|k|-2k_0)\FF[(\eta_1^\pm)^2]\|_0^2=\|(\i k \mp \i k_0)\FF[(\eta_1^\pm)^2]\|_0^2$$
(because $f(k)=f(2k_0) + O(|k|-2k_0)$ for $k \in 2S$) and
$$\|(K^0-\underbrace{f(0)}_{\displaystyle = 1})\eta_1^+\eta_1^-\|_0^2 \leq c \||k|\FF[\eta_1^+\eta_1^-]\|_0^2=\|\i k \FF[\eta_1^+\eta_1^-]\|_0^2$$
(because $f(k)=f(0) + O(|k|)$ for $k \in [-2\delta_0,2\delta_0]$),
and (vii) and (viii) are deduced from (iii) and (iv) in the same fashion.\qed

\begin{proposition}
The function $\tilde{\eta}_1$ satisfies the estimates
\begin{eqnarray*}
&& \KK_4(\tilde{\eta}_1) = A_4^1 \int_{-\infty}^\infty \tilde{\eta}_1^4 \dx+ o(\mu^3), \qquad A_4^1 = -\frac{\beta\omega k_0^4}{8} - \frac{\omega^2}{24}(f(2k_0)+2), \\
&& \GG_4(\tilde{\eta}_1) = A_4^2 \int_{-\infty}^\infty \tilde{\eta}_1^4 \dx+ o(\mu^3), \qquad A_4^2 = \frac{\omega k_0^2}{6} - \frac{\omega}{12}f(k_0)(f(2k_0)+2), \\
&& \LL_4(\tilde{\eta}_1) = A_4^3 \int_{-\infty}^\infty \tilde{\eta}_1^4 \dx+ o(\mu^3), \qquad A_4^3 = \frac{1}{6}f(k_0)^2(f(2k_0)+2) - \frac{k_0^2f(k_0)}{2}.
\end{eqnarray*}
\end{proposition}
{\bf Proof.} Using the formulae given in Lemma \ref{Approximate operators}, we find that
\begin{eqnarray*}
\int_{-\infty}^\infty \tilde{\eta}_1^2 \tilde{\eta}_1^{\prime 2} \dx
& = & \int_{-\infty}^\infty \big((\tilde{\eta}_1^+)^2((\tilde{\eta}_1^-)^\prime)^2 + (\tilde{\eta}_1^-)^2((\tilde{\eta}_1^+)^\prime)^2 
+4\tilde{\eta}_1^+\tilde{\eta}_1^-(\tilde{\eta}_1^+)^\prime(\tilde{\eta}_1^-)^\prime  \big) \dx\\
& = & 2k_0^2 \int_{-\infty}^\infty (\tilde{\eta}_1^+)^2(\tilde{\eta}_1^-)^2 \dx+ o(\mu^3),
\end{eqnarray*}
and similarly
\begin{eqnarray*}
\int_{-\infty}^\infty K^0(\tilde{\eta}_1^2)\tilde{\eta}_1K^0\tilde{\eta}_1 \dx
& = & (2f(2k_0)f(k_0)+4f(k_0))\int_{-\infty}^\infty (\tilde{\eta}_1^+)^2(\tilde{\eta}_1^-)^2 \dx+ o(\mu^3),\\
\int_{-\infty}^\infty (\tilde{\eta}_1^\prime)^4 \dx & = & 6k_0^4  \int_{-\infty}^\infty (\tilde{\eta}_1^+)^2 (\tilde{\eta}_1^-)^2 \dx+ o(\mu^3), \\
\int_{-\infty}^\infty \tilde{\eta}_1^2 K^0 (\tilde{\eta}_1^2) \dx
& = & (2f(2k_0)+4)\int_{-\infty}^\infty (\tilde{\eta}_1^+)^2(\tilde{\eta}_1^-)^2 \dx + o(\mu^3), \\
\int_{-\infty}^\infty K^0(\tilde{\eta}_1 K^0 \tilde{\eta}_1)\tilde{\eta}_1K^0\tilde{\eta}_1 \dx
& = & (2f(2k_0)f(k_0)^2 + 4f(k_0)^2) \int_{-\infty}^\infty (\tilde{\eta}_1^+)^2(\tilde{\eta}_1^-)^2 \dx + o(\mu^3), \\
\int_{-\infty}^\infty (K^0 \tilde{\eta}_1)\tilde{\eta}_1^2 \tilde{\eta}_1^{\prime\prime} \dx
& = & -6k_0^2f(k_0) \int_{-\infty}^\infty (\tilde{\eta}_1^+)^2(\tilde{\eta}_1^-)^2 \dx + o(\mu^3).
\end{eqnarray*}
The result is obtained by substituting the above expressions into the explicit formulae for $\KK_4$, $\GG_4$ and $\LL_4$ given in Proposition
\ref{Formulae for the threes and fours}.\qed

\begin{corollary}  \label{lot in fours}
The function $\tilde{\eta}$ satisfies the estimate
$$\KK_4(\tilde{\eta}) + 2\nu_0 \GG_4(\tilde{\eta}) - \nu_0^2 \LL_4(\tilde{\eta})
= A_4 \int_{-\infty}^\infty \tilde{\eta}_1^4 \dx + o(\mu^3),$$
where
$$A_4=A_4^1+2\nu_0A_4^2-\nu_0^2A_4^3.$$
\end{corollary}

We now turn to the corresponding result for $\GG_3(\tilde{\eta})$, 
$\KK_3(\tilde{\eta})$ and $\LL_3(\tilde{\eta})$.

\begin{proposition} \label{lot in threes step 1}
The function $\tilde{\eta}$ satisfies the estimate
$$
\begin{Bmatrix}\GG_3(\tilde{\eta}) \\ \KK_3(\tilde{\eta}) \\ \LL_3(\tilde{\eta}) \end{Bmatrix}
= - \int_{-\infty}^\infty 
\begin{Bmatrix} \GG_3^\prime(\tilde{\eta}_1) \\ \KK_3^\prime(\tilde{\eta}_1) \\ \LL_3^\prime(\tilde{\eta}_1) \end{Bmatrix}
H(\tilde{\eta}) \dx+ o(\mu^3).
$$
\end{proposition}
{\bf Proof.} Each term in the expansion of
$$
n_2(\tilde{\eta}_1-H(\tilde{\eta})+\tilde{\eta}_3,\tilde{\eta}_1-H(\tilde{\eta})+\tilde{\eta}_3,\tilde{\eta}_1-H(\tilde{\eta})+\tilde{\eta}_3)
$$
with zero or one occurrence of $\tilde{\eta}_1$ can be estimated by
$$
\left| n_2 \left (
\begin{Bmatrix}
\tilde{\eta}_1 \\
H(\tilde{\eta}) \\
\tilde{\eta}_3
\end{Bmatrix}, 
\begin{Bmatrix}
H(\tilde{\eta}) \\
\tilde{\eta}_3
\end{Bmatrix}^{(2)}\right)\right|
\ \leq \ c
\begin{Bmatrix}
\|\tilde{\eta}_1\|_2 \\
\|H(\tilde{\eta})\|_2 \\
\|\tilde{\eta}_3\|_2
\end{Bmatrix}
\begin{Bmatrix}
\|H(\tilde{\eta})\|_2 \\
\|\tilde{\eta}_3\|_2
\end{Bmatrix}^{2} \\
\ \leq \  c\mu^\frac{1}{2}\mu^{2+\alpha}\\
\ =\ o(\mu^3),
$$
while
$$|n_2(\tilde{\eta}_1,\tilde{\eta}_1,\tilde{\eta}_3)|\ \leq\ c \|\tilde{\eta}\|_2^2\|\tilde{\eta}_3\|_2\ \leq\ c\mu\mu^{\frac{3}{2}+\alpha}\ =\ o(\mu^3)$$
and
$$n_2(\tilde{\eta}_1,\tilde{\eta}_1,\tilde{\eta}_1) = \GG_3(\tilde{\eta}_1)=0.$$
It follows that
\begin{eqnarray*}
\GG_3(\tilde{\eta}) & = & -3n_2(\tilde{\eta}_1,\tilde{\eta}_1,H(\tilde{\eta})) + o(\mu^3) \\
& = & -\mathrm{d}\GG_3[\tilde{\eta}_1](H(\tilde{\eta})) + o(\mu^3) \\
& = & - \int_{-\infty}^\infty \GG_3^\prime(\tilde{\eta}_1) H(\tilde{\eta}) \dx + o(\mu^3).
\end{eqnarray*}
The same argument yields the results for $\KK_3(\tilde{\eta})$ and $\LL_3(\tilde{\eta})$.\qed

\begin{proposition} \label{lot in threes step 2}
The function $\tilde{\eta}$ satisfies the estimate
$$H(\tilde{\eta}) = \FF^{-1} \left[\frac{1}{g(k)} \FF[\KK_3^\prime(\tilde{\eta}_1) + 2\nu_0 \GG_3^\prime(\tilde{\eta}_1) - \nu_0^2 \LL_3^\prime(\tilde{\eta}_1)]\right]
+\underline{o}(\mu^3).$$
\end{proposition}
{\bf Proof.} Noting that
$$\left| \frac{\mu+\GG(\tilde{\eta})}{\LL(\eta)} - \nu_0 \right|\ \leq\ c(\mu^\alpha \nn \tilde{\eta}_1 \nn_\alpha^2 + \|\tilde{\eta}_3\|_2 + \mu^{N-\frac{1}{2}})\ =\ O(\mu^{1+\alpha})$$
(see Corollary \ref{Weak ST speed estimate}) and
$$
\begin{Bmatrix} \|\GG_3^\prime(\tilde{\eta}_1)\|_0 \\ \|\KK_3^\prime(\tilde{\eta}_1)\|_0 \\ \|\LL_3^\prime(\tilde{\eta}_1)\|_0 \end{Bmatrix}
\ \leq\ c\mu^\frac{\alpha}{2}\nn \tilde{\eta}_1 \nn_\alpha \|\tilde{\eta}_1\|_2\ =\ O(\mu^{1+\frac{\alpha}{2}}),
$$
(see Proposition \ref{Size of the gradients}) one finds that
$$H(\tilde{\eta}) = \FF^{-1} \left[\frac{1}{g(k)} \FF[\KK_3^\prime(\tilde{\eta}_1) + 2\nu_0 \GG_3^\prime(\tilde{\eta}_1) - \nu_0^2 \LL_3^\prime(\tilde{\eta}_1)]\right]
+\underbrace{O(\mu^{1+\alpha})\underline{O}(\mu^{1+\frac{\alpha}{2}})}_{\displaystyle = \underline{o}(\mu^3)}.\eqno{\Box}$$

Combining Propositions \ref{lot in threes step 1} and \ref{lot in threes step 2}, one finds that
\begin{eqnarray}
\lefteqn{\KK_3(\tilde{\eta}) + 2\nu_0 \GG_3(\tilde{\eta}) - \nu_0^2\LL_3(\tilde{\eta})} \nonumber \\
& & \hspace{-7mm}= -\int_{-\infty}^\infty (\KK_3^\prime(\tilde{\eta}_1) + 2\nu_0 \GG_3^\prime(\tilde{\eta}_1) - \nu_0^2 \LL_3^\prime(\tilde{\eta}_1))
\FF^{-1} \left[\frac{1}{g(k)} \FF[\KK_3^\prime(\tilde{\eta}_1) + 2\nu_0 \GG_3^\prime(\tilde{\eta}_1) - \nu_0^2 \LL_3^\prime(\tilde{\eta}_1)]\right] \dx\nonumber \\
& & \hspace{-3.5mm}\mbox{} +o(\mu^3), \label{lcomb threes preliminary}
\end{eqnarray}
which we write as
\begin{eqnarray}
\lefteqn{\KK_3(\tilde{\eta}) + 2\nu_0 \GG_3(\tilde{\eta}) - \nu_0^2\LL_3(\tilde{\eta})} \nonumber \\
& & \hspace{-7mm}= -2\!\!\int_{-\infty}^\infty M(\tilde{\eta}_1^+,\tilde{\eta}_1^+)\FF^{-1}[g(k)^{-1}M(\tilde{\eta}_1^-,\tilde{\eta}_1^-)]\dx
-4\int_{-\infty}^\infty M(\tilde{\eta}_1^+,\tilde{\eta}_1^-)\FF^{-1}[g(k)^{-1}M(\tilde{\eta}_1^+,\tilde{\eta}_1^-)] \dx\nonumber \\
& & \hspace{-3.5mm}\mbox{} +o(\mu^3), \qquad \label{lcomb threes}
\end{eqnarray}
where
$$M=m_1+2\nu_0m_2-\nu_0^2m_3,$$
in order to determine the dominant term on its right-hand side.

\begin{proposition} \label{lot in threes}
The function $\tilde{\eta}$ satisfies
$$\KK_3(\tilde{\eta}) + 2\nu_0 \GG_3(\tilde{\eta}) - \nu_0^2\LL_3(\tilde{\eta}) = A_3 \int_{-\infty}^\infty \tilde{\eta}_1^4 \dx + o(\mu^3),$$
where
\begin{eqnarray*}
A_3 & = & - \frac{g(2k_0)^{-1}}{3}(A_3^1)^2 - \frac{2g(0)^{-1}}{3}(A_3^2)^2, \\
A_3^1 & = & \frac{\omega \nu_0}{2}f(2k_0)+\omega \nu_0f(k_0)+\frac{\omega^2}{2}+\nu_0^2 f(2k_0)f(k_0)+\frac{\nu_0^2}{2}f(k_0)^2-\frac{3k_0^2\nu_0^2}{2}, \\
A_3^2 & = & \frac{\omega \nu_0}{2}+\omega \nu_0f(k_0)+\frac{\omega^2}{2}+\nu_0^2f(k_0)+\frac{\nu_0^2}{2}f(k_0)^2-\frac{\nu_0^2k_0^2}{2}.
\end{eqnarray*}
\end{proposition}
{\bf Proof.} Lemma \ref{Approximate operators} implies that
$$
M(\tilde{\eta}_1^+,\tilde{\eta}_1^+) = A_3^1 (\tilde{\eta}_1^+)^2 + \underline{O}(\mu^{1+\alpha}),
$$
so that
$$\FF^{-1}[g(k)^{-1}M(\tilde{\eta}_1^-,\tilde{\eta}_1^-)]\ =\ \FF^{-1}[g(k)^{-1}\overline{M(\tilde{\eta}_1^+,\tilde{\eta}_1^+)}]
\ =\ g(2k_0)^{-1}A_3^1(\tilde{\eta}_1^-)^2 + \underline{O}(\mu^{1+\alpha}),$$
and
$$
M(\tilde{\eta}_1^+,\tilde{\eta}_1^-) = A_3^2 \tilde{\eta}_1^+\tilde{\eta}_1^- + \underline{O}(\mu^{1+\alpha}),
$$
so that
$$\FF^{-1}[g(k)^{-1}M(\tilde{\eta}_1^+,\tilde{\eta}_1^-)] = g(0)^{-1}A_3^2 \tilde{\eta}_1^+ \tilde{\eta}_1^- + \underline{O}(\mu^{1+\alpha});$$
the result follows from these calculations and equation \eqn{lcomb threes}.\qed

The requisite estimates for $\MM_\mu(\tilde{\eta})$ and
$\langle \MM_\mu^\prime(\tilde{\eta}), \tilde{\eta} \rangle + 4\mu \tilde{\MM}_\mu(\tilde{\eta})$
may now be derived from Corollary \ref{lot in fours} and Proposition \ref{lot in threes}.

\begin{lemma} \label{Weak ST MM estimate}
The estimates
\begin{eqnarray*}
& & \hspace{-6mm}\MM_{a^2\mu}(a\tilde{\eta}) = a^3\big(\KK_3(\tilde{\eta}) + 2 \nu_0 \GG_3(\tilde{\eta}) - \nu_0^2 \LL_3(\tilde{\eta})\big)
+a^4\big(\KK_4(\tilde{\eta}) + 2 \nu_0 \GG_4(\tilde{\eta}) - \nu_0^2 \LL_4(\tilde{\eta})\big)
+ a^3 o(\mu^3),\\
& &  \hspace{-6mm}\langle \MM_{a^2\mu}^\prime(a\tilde{\eta}), a\tilde{\eta} \rangle + 4a^2\mu \tilde{\MM}_{a^2\mu}(a\tilde{\eta})
=3a^3\big(\KK_3(\tilde{\eta}) + 2\nu_0\GG_3(\tilde{\eta})-\nu_0^2 \LL_3(\tilde{\eta})\big) \\
& & \hspace{2.7in}\mbox{}+4a^4\big(\KK_4(\tilde{\eta}) + 2\nu_0 \GG_4(\tilde{\eta}) - \nu_0^2 \LL_4(\tilde{\eta}) \big) + a^3o(\mu^3)
\end{eqnarray*}
hold uniformly over $a \in [1,2]$.
\end{lemma}
{\bf Proof.} Lemma \ref{Formulae for M} asserts that
\begin{eqnarray*}
\MM_{a^2\mu}(a\tilde{\eta}) & = & a^3\big(\KK_3(\tilde{\eta}) + 2 \nu_0 \GG_3(\tilde{\eta}) - \nu_0^2 \LL_3(\tilde{\eta}))
+a^4\big(\KK_4(\tilde{\eta}) + 2 \nu_0 \GG_4(\tilde{\eta}) - \nu_0^2 \LL_4(\tilde{\eta})\big) \\
& & \mbox{}+2\left(\frac{\mu+\GG_2(\tilde{\eta})}{\LL_2(\tilde{\eta})}-\nu_0\right)(a^3\GG_3(\tilde{\eta})+a^4\GG_4(\tilde{\eta})) \\
& & \mbox{}-\left(\frac{\mu+\GG_2(\tilde{\eta})}{\LL_2(\tilde{\eta})}-\nu_0\right)\!\!\left(\frac{\mu+\GG_2(\tilde{\eta})}{\LL_2(\tilde{\eta})}+\nu_0\right)(a^3\LL_3(\tilde{\eta})+a^4\LL_4(\tilde{\eta})) \\
& & \mbox{}+\frac{a^4}{\LL_2(\tilde{\eta})}\left(\GG_3(\tilde{\eta})-\left(\frac{\mu+\GG_2(\tilde{\eta})}{\LL_2(\tilde{\eta})}\right)\LL_3(\tilde{\eta})\right)^{\!\!2}
+ O(a^5\mu^\frac{3}{2}(\|\tilde{\eta}\|_{1,\infty} + \|\tilde{\eta}^{\prime\prime}+k_0^2 \tilde{\eta}\|_0)^2)
\end{eqnarray*}
uniformly over $a \in [1,2]$.

The first result follows by estimating
$$
\begin{Bmatrix}
\GG_3(\tilde{\eta}) \\ \LL_3(\tilde{\eta})
\end{Bmatrix}
= O(\mu^\frac{3}{2}),
\quad
\begin{Bmatrix}
\GG_4(\tilde{\eta}) \\ \LL_4(\tilde{\eta})
\end{Bmatrix}
= O(\mu^2),
$$
$$
\|\tilde{\eta}\|_{1,\infty} + \|\tilde{\eta}^{\prime\prime}+k_0^2 \tilde{\eta}\|_0
\ \leq\ c(\mu^\frac{\alpha}{2}\nn \tilde{\eta} \nn_\alpha+\|\tilde{\eta}_3\|_2)
\ \leq\ c\mu^{\frac{1}{2}+\frac{\alpha}{2}}
$$
(see equation \eqn{General estimate pt 1}),
$$\left|\frac{\mu+\GG_2(\tilde{\eta})}{\LL_2(\tilde{\eta})} - \nu_0\right| \leq c( \mu^\alpha \nn \tilde{\eta}_1 \nn_\alpha^2 + \|\eta_3\|_2 + \mu^{N-\frac{1}{2}}) \leq c \mu^{1+\alpha}$$
and noting that
\begin{eqnarray*}
\lefteqn{\GG_3(\tilde{\eta})-\left(\frac{\mu+\GG_2(\tilde{\eta})}{\LL_2(\tilde{\eta})}\right)\LL_3(\tilde{\eta})} \\
& = & \GG_3(\tilde{\eta}) - \nu_0 \LL_3(\tilde{\eta}) + o(\mu^3) \\
& = & -\int_{-\infty}^\infty (\GG_3^\prime(\tilde{\eta}_1) - \nu_0 \LL_3^\prime(\tilde{\eta}_1))
\FF^{-1} \left[\frac{1}{g(k)} \FF[\KK_3^\prime(\tilde{\eta}_1) + 2\nu_0 \GG_3^\prime(\tilde{\eta}_1) - \nu_0^2 \LL_3^\prime(\tilde{\eta}_1)]\right] \dx
+o(\mu^3) \\
& = & -\int_{-\infty}^\infty \big(\tilde{M}(\tilde{\eta}_1^+,\tilde{\eta}_1^+)\FF^{-1}[g(k)^{-1}M(\tilde{\eta}_1^-,\tilde{\eta}_1^-)]
+ \tilde{M}(\tilde{\eta}_1^-,\tilde{\eta}_1^-)\FF^{-1}[g(k)^{-1}M(\tilde{\eta}_1^+,\tilde{\eta}_1^+)] \big)\dx \\
& & \mbox{}-4\int_{-\infty}^\infty \tilde{M}(\tilde{\eta}_1^+,\tilde{\eta}_1^-)\FF^{-1}[g(k)^{-1}M(\tilde{\eta}_1^+,\tilde{\eta}_1^-)] \dx+o(\mu^3)\\
& = & \gamma \int_{-\infty}^\infty \tilde{\eta}_1^4 \dx + o(\mu^3) \\
& = & O(\mu^{2+\alpha}) + o(\mu^3),
\end{eqnarray*}
where $\tilde{M}=m_2-\nu_0m_3$ and $\gamma$ is a (possibly negative) constant. Here the
third line follows from the second by Propositions \ref{lot in threes step 1} and
\ref{lot in threes step 2} and the fifth from the fourth by repeating the proof of Proposition
\ref{lot in threes}.

The second result is derived in a similar fashion.\qed

\begin{corollary} \label{Weak ST final MM estimates}
The estimates
\begin{eqnarray*}
& & \MM_{a^2\mu}(a\tilde{\eta}) = (a^3A_3+a^4A_4)\int_{-\infty}^\infty \tilde{\eta}_1^4 \dx+ a^3 o(\mu^3),\\
& & \langle \MM_{a^2\mu}^\prime(a\tilde{\eta}), a\tilde{\eta} \rangle + 4a^2\mu \tilde{\MM}_{a^2\mu}(a\tilde{\eta})
=(3a^3A_3+4a^4A_4)\int_{-\infty}^\infty \tilde{\eta}_1^4 \dx+ a^3 o(\mu^3),\\
\end{eqnarray*}
hold uniformly over $a \in [1,2]$ and
$$\int_{-\infty}^\infty \tilde{\eta}_1^4 \dx \geq c\mu^3.$$
\end{corollary}
{\bf Proof.} The estimates follow by combining
Corollary \ref{lot in fours}, Proposition \ref{lot in threes} and Lemma \ref{Weak ST MM estimate},
while the inequality for $\tilde{\eta}_1$ is a consequence of
the first estimate (with $a=1$) and the fact that $\MM_\mu(\tilde{\eta}) \leq -c\mu^3$.\qed

\subsection{Derivation of the strict sub-additivity property} \label{SSA subsection}

In this section we derive the strict sub-additivity property \eqn{Strict SA}. We begin with by showing that
$c_\mu$ is a strictly sub-homogeneous, increasing function of $\mu>0$. The first of these properties is
a corollary of the next proposition.

\begin{proposition}
There exists $a_0 \in (1,2]$ and $q>2$ with the property that the function
$$a \mapsto a^{-q}\MM_{a^2\mu}(a\tilde{\eta}), \qquad a \in [1,a_0]$$
is decreasing and strictly negative.
\end{proposition}
{\bf Proof.} This result follows from the calculations
\begin{eqnarray*}
\frac{\mathrm{d}}{\mathrm{d}a}\left(a^{-\frac{5}{2}}\MM_{a^2\mu}(a\tilde{\eta})\right)
& = & a^{-\frac{7}{2}}\left(
{\textstyle -\frac{5}{2}}\MM_{a^2\mu}(a\tilde{\eta})+\langle \MM^\prime_{a^2\mu}(a\tilde{\eta}),a\tilde{\eta}\rangle_0
+ 4a^2\mu\tilde{\MM}_{a^2\mu}(a\tilde{\eta})\right) \\
& = & \frac{1}{4}a^{-\frac{7}{2}}\left(\frac{a^3}{4}\left(\frac{\omega^3}{3}+1\right)\int_{-\infty}^\infty \tilde{\eta}_1^3 \dx+ a^3 o(\mu^\frac{5}{3})\right) \\
& = & a^{-\frac{1}{2}}\left(\frac{1}{4}\left(\frac{\omega^3}{3}+1\right)\int_{-\infty}^\infty \tilde{\eta}_1^3 \dx+ o(\mu^\frac{5}{3})\right) \\
& \leq & -c\mu^\frac{5}{3} \\
& <  & 0, \hspace{2in}a \in (1,2)
\end{eqnarray*}
for $\beta>\beta_\mathrm{c}$ (see Corollary \ref{Strong ST final MM estimates}) and
\begin{eqnarray*}
\lefteqn{\frac{\mathrm{d}}{\mathrm{d}a}\left(a^{-q}\MM_{a^2\mu}(a\tilde{\eta})\right)}\qquad \\
& = & a^{-(q+1)}\left(
-q\MM_{a^2\mu}(a\tilde{\eta})+\langle \tilde{\MM}^\prime_{a^2\mu}(a\tilde{\eta}),a\tilde{\eta}\rangle_0
+ 4a^2\mu\tilde{\MM}_{a^2\mu}(a\tilde{\eta})\right) \\
& = & a^{-(q+1)}\left(\big(-q(a^3A_3+a^4A_4)+3a^3A_3+4a^4A_4\big)\int_{-\infty}^\infty \tilde{\eta}_1^4  \dx+ a^3o(\mu^3)\right) \\
& = & a^{2-q}\left(\big((3-q)A_3+a(4-q)A_4\big)\int_{-\infty}^\infty \tilde{\eta}_1^4 \dx + o(\mu^3)\right) \\
& \leq & -c\mu^3 \\
& <  & 0, \hspace{2.5in}a \in (1,a_0),\ q \in (2,q_0)
\end{eqnarray*}
for $\beta<\beta_\mathrm{c}$ (see Corollary \ref{Weak ST final MM estimates}); here $a_0>1$ and $q_0>2$ are chosen
so that $(3-q)A_3+a(4-q)A_4$, which is negative for $a=1$ and $q=2$ (see Appendix B), is also negative for $a \in (1,a_0]$ and $q \in (2,q_0]$.\qed
\begin{corollary}
The number $c_\mu$ is a strictly sub-homogeneous function of $\mu>0$.
\end{corollary}
{\bf Proof.} 
The previous lemma implies that
\[
\mathcal M_{a\mu}(a^\frac{1}{2} \tilde \eta_m)\le a^\frac{q}{2}\mathcal M_\mu (\tilde \eta_m)<0, 
\quad a \in [1,a_0^2],
\]
from which it follows that
\begin{eqnarray*}
c_{a\mu} & \leq & \JJ_{a\mu}(a^\frac{1}{2}\tilde{\eta}_m) \\
& = & \KK_2(a^{\frac{1}{2}}\tilde{\eta}_m) + \frac{(a\mu+\GG_2(a^{\frac{1}{2}}\tilde{\eta}_m))^2}{\LL_2(a^{\frac{1}{2}}\tilde{\eta}_m)}
+\MM(a^{\frac{1}{2}}\tilde{\eta}_m) \\
& \leq & a\left(\KK_2(\tilde{\eta}_m)+\frac{(\mu+\GG(\tilde{\eta}_m))^2}{\LL(\tilde{\eta}_m)}\right) + a^\frac{q}{2}\MM_\mu(\tilde{\eta}_m) \\
& = & a\left(\KK_2(\tilde{\eta}_m)+\frac{\mu^2}{\LL(\tilde{\eta}_m)}+\MM_\mu(\tilde{\eta}_m)\right)
+ (a^\frac{q}{2}-a)\MM_\mu(\tilde{\eta}_m) \\
& \leq & a\JJ(\tilde{\eta}_m) - c(a^\frac{q}{2}-a)\mu^{r^\star}
\end{eqnarray*}
for $a \in [1,a_0^2]$. In the limit $n \rightarrow \infty$ the above inequality yields
$$c_{a\mu}\ \leq\ ac_\mu - c(a^\frac{q}{2}-a)\mu^{r^\star}\ < ac_\mu.\eqno{\Box}$$
for $a\in (1,a_0^2]$.

For $a>a_0^2$ we choose $p \geq 2$ such that $a \in (1,a_0^{2p}]$ (and hence $a^\frac{1}{p} \in (1,a_0^2]$) and observe that
$$c_{a\mu}<a^\frac{1}{p} c_{a^{(p-1)/p} \mu}<a^\frac{2}{p} c_{a^{(p-2)/p} \mu}
< \cdots < ac_\mu.\eqno{\Box}$$

\begin{lemma}
The number $c_\mu$ is an increasing function of $\mu>0$.
\end{lemma}
{\bf Proof.}
Using Proposition \ref{Strong ST speed estimate} for $\beta>\beta_\mathrm{c}$ and Proposition
\ref{Weak ST speed estimate} for $\beta<\beta_\mathrm{c}$, one finds that
$$\mu + \GG(\tilde{\eta}_m) =  \nu_0\LL(\tilde{\eta}_m) + O(\mu^\frac{3}{2}) \geq c\mu + O(\mu^\frac{3}{2})$$
so that
$$\mu + \GG(\tilde{\eta}_m) \geq c_\star \mu$$
for some $c_\star \in (0,1)$. Let $d_\star = 1- c_\star$, so that $d_\star \in (0,1)$.

First suppose that $\mu_1 \in [d_\star\mu_2, \mu_2]$. Let $\{\tilde{\eta}_m^2\}$ be the special minimising
sequence constructed in Theorem \ref{Special MS theorem} for $\mu=\mu_2$ and note that
$$\mu_1 + \GG(\tilde{\eta}_m^2) = \mu_2 + \GG(\tilde{\eta}_m^2) - (\mu_2-\mu_1) \geq \mu_1-d_\star \mu_2 \geq 0,$$
so that $\JJ_{\mu_1}(\tilde{\eta}_m^2) \leq \JJ_{\mu_2}(\tilde{\eta}_m^2)$. It follows that
$$c_{\mu_1} \leq \JJ_{\mu_1}(\tilde{\eta}_m^2) \leq \JJ_{\mu_2}(\tilde{\eta}_m^2) \rightarrow c_{\mu_2}$$
as $n \rightarrow \infty$,
that is
$$c_{\mu_1} \leq c_{\mu_2}.$$

For $\mu_1 < d_\star\mu_2$ we choose $p \geq 2$ such that $\mu_1 \in [d_\star^p \mu_2, \mu_2]$
(and hence $\mu_1 \in [d_\star d_\star^{p-1} \mu_2, d_\star^{p-1} \mu_2]$
and obviously $d_\star^{q+1}\mu_2 \in [d_\star d_\star^q \mu_2, d_\star^q\mu_2]$, $q=0,\ldots, p-2$) and
observe that
$$c_{\mu_1} \leq c_{d_\star^{p-1}\mu_2} \leq c_{d_\star^{p-2}\mu_2} \leq \cdots \leq c_{\mu_2}.\eqno{\Box}$$

Our final result is stated in the following theorem.

\begin{theorem} \label{SSA theorem}
The number $c_\mu$
has the strict sub-additivity property
$$
c_{\mu_1+\mu_2} < c_{\mu_1} + c_{\mu_2}, \qquad 0<|\mu_1|, |\mu_2|, \mu_1+\mu_2 < \mu_0.
$$
\end{theorem}
{\bf Proof.} Using the strict sub-homogeneity of $c(\mu)$ for $\mu>0$, we find that
$$c_{\mu_1+\mu_2} < \frac{\mu_1+\mu_2}{\mu_1}c_{\mu_1} = c_{\mu_1} + \frac{\mu_2}{\mu_1}c_{\mu_1} \leq c_{\mu_1}+c_{\mu_2}$$
for $0<\mu_1\leq \mu_2$,
and for $\mu_1<0$, $\mu_2>0$ with $\mu_1+\mu_2 > 0$ its monotonicity for $\mu>0$ shows that
$$c_{\mu_1+\mu_2} \leq c_{\mu_2} < c_{\mu_1} + c_{\mu_2}.\eqno{\Box}$$

\section{Existence theory and consequences}

\subsection{Minimisation} \label{Minimisation}

The following theorem, which is proved using the results of Sections \ref{Minimising sequences}
and \ref{SSA section}, is our final result concerning the set of minimisers
of $\JJ_\mu$ over $U\sm\{0\}$. 

\begin{theorem} \label{Key minimisation theorem}
\quad
\begin{list}{(\roman{count})}{\usecounter{count}}
\item
The set $B_\mu$ of minimisers of $\JJ_\mu$ over $U \sm \{0\}$ is non-empty.
\item
Suppose that $\{\eta_m\}$ is a minimising sequence for $\JJ_\mu$ on $U\sm\{0\}$ which satisfies
$$
\sup_{m\in{\mathbb N}} \|\eta_m\|_2 < M.
$$
There exists a sequence $\{x_m\} \subset \R$ with the property that
a subsequence of $\{\eta_m(x_m+\cdot)\}$ converges
in $H^r(\R)$, $r \in [0,2)$, to a function $\eta \in B_\mu$.
\end{list}
\end{theorem}
{\bf Proof.} It suffices to prove part (ii), since an application of this result to the sequence
$\{\tilde{\eta}_m\}$ constructed in Theorem \ref{Special MS theorem} yields part (i).

In order to establish part (ii) we choose $\tilde{M} \in (\sup_{m \in {\mathbb N}}  \|\eta_m\|_2, M)$, so that
$\{\eta_m\}$ is also a minimising sequence for the functional $\JJ_{\rho,\mu}$
introduced in Section \ref{Penalised minimisation} (the existence of a minimising sequence $\{v_m\}$
for $\JJ_{\rho,\mu}$ with $\lim_{m \rightarrow \infty} \JJ_{\rho,\mu}(v_m) <
\lim_{m \rightarrow \infty} \JJ_{\rho,\mu}(\eta_m)$ would lead to the contradiction
$$\lim_{m \rightarrow \infty} \JJ_\mu(v_m)\ \leq\ \lim_{m \rightarrow \infty} \JJ_{\rho,\mu}(v_m)
\ <\  \lim_{m \rightarrow \infty} \JJ_{\rho,\mu}(\eta_m)
\ =\ \lim_{m \rightarrow \infty} \JJ_\mu(\eta_m)\ =\ c_\mu\ ).$$
We may therefore study $\{\eta_m\}$ using
the theory given in Section \ref{MSP}, noting that  the sequence $\{u_m\}$
with $u_m = (\eta_m^\prime)^2 + \eta_m^2$
does not have the `dichotomy' property:
the existence of sequences $\{\eta_m^{(1)}\}$, $\{\eta_m^{(2)}\}$ with the features
listed in Lemma \ref{Splitting properties} is
incompatible with the  strict sub-additivity property of $c_\mu$ (Theorem \ref{SSA theorem}).
Recall that the numbers $\mu^{(1)}$, $\mu^{(2)}$ sum to $\mu$;  this fact leads to the contradiction
\begin{eqnarray*}
c_\mu & < & c_{\mu^{(1)}} + c_{\mu^{(2)}} \\
& \leq & \lim_{m \rightarrow \infty}\JJ_{\mu^{(1)}}(\eta_m^{(1)}) + \lim_{m \rightarrow \infty}\JJ_{\mu^{(2)}}(\eta_m^{(2)}) \\
& = & \lim_{m \rightarrow \infty}\JJ_\mu(\eta_m) \\
& = & c_\mu.
\end{eqnarray*}
We conclude that $\{u_m\}$ has the `concentration' property and hence
$\eta_m(\cdot+x_m) \rightarrow \eta^{(1)}$ as $n \rightarrow \infty$ in
$H^r(\R)$ for every $r \in [0,2)$, (see Lemma \ref{Vanishing and concentration}(ii)), whereby
$\JJ_\mu(\eta)=\lim_{m \rightarrow \infty} \JJ_\mu(\eta_m(\cdot+x_m)) = c_\mu$, so that $\eta^{(1)}$
is a minimiser of $\JJ_\mu$ over $U\sm\{0\}$.\qed

The next step is to relate the above result to our original problem of
finding minimisers of $\HH(\eta,\Phi)$ subject to the constraint $\II(\eta,\Phi) = 2\mu$,
where $\HH$ and $\II$ are defined in equations \eqn{Definition of HH} and
\eqn{Definition of II}.

\begin{theorem} \label{Result for constrained minimisation}
\quad
\begin{list}{(\roman{count})}{\usecounter{count}}
\item
The set $D_\mu$ of minimisers of $\HH$ on the set
$$S_\mu = \{(\eta,\xi) \in U \times H_\star^{1/2}(\R): \II(\eta,\Phi) = 2\mu\}$$
is non-empty.
\item
Suppose that $\{(\eta_m,\xi_m)\} \subset S_\mu$ is a minimising sequence for $\HH$ with the
property that $\sup_{m \in \N} \|\eta_m\|_2 < M$. There exists a sequence $\{x_m\} \subset \R$
with the property that a subsequence of $\{(\eta_m(x_m + \cdot), \xi_m(x_m+\cdot)\}$ converges
in $H^r(\R) \times H_\star^{1/2}(\R)$, $r \in [0,2)$, to a function in $D_\mu$.
\end{list}
\end{theorem}
{\bf Proof.} (i) We consider the minimisation problem in two steps.
\begin{enumerate}
\item
\emph{Fix $\eta \in U\sm\{0\}$ and 
minimise $\HH(\eta,\cdot)$ over $T_\mu = \{\xi \in H_\star^{1/2}(\R): \II(\eta,\xi)=2\mu\}$.}
Notice that 
$\HH(\eta,\cdot)$ is weakly lower semicontinuous 
on $H_\star^{1/2}(\R)$ (since $\xi \mapsto \langle G(\eta)\xi,\xi\rangle_0^\frac{1}{2}$ 
is equivalent to its usual norm), while
$\II (\eta,\cdot)$ is weakly continuous on $H_\star^{1/2}(\R)$; furthermore 
$\HH(\eta,\cdot)$ is convex and coercive. A familiar argument shows that $\HH(\eta,\cdot)$
has a unique minimiser $\xi_\eta$ over $T_\mu$.

\emph{2. Minimise $\HH(\eta,\xi_\eta)$ over $U\sm\{0\}$.}
Because $\Phi_\eta$ minimises $\HH(\eta,\cdot)$ over $T_\mu$ there exists a
Lagrange multiplier $\nu_\eta$ such that
$$G(\eta)\xi_\eta+\omega\eta\eta^\prime=\nu_\eta \eta^\prime,$$
and a straightforward calculation shows that
\begin{equation}
\xi_\eta = G(\eta)^{-1}(\nu_\eta \eta^\prime -\omega \eta\eta^\prime), \quad \nu_\eta =\frac{\mu+\GG(\eta)}{\mathcal L(\eta)}. \label{Formula for minimiser and LM}
\end{equation}
According to Theorem \ref{Key minimisation theorem}(i)
the set $B_\mu$ of minimisers of $\JJ_\mu(\eta):=\HH(\eta,\xi_\eta)$ over $U \sm \{0\}$ is not empty; it follows
that $D_\mu$ is also not empty.
\end{enumerate}

(ii)
Let $\{(\eta_m,\xi_m)\}\subset U\times H_\star^{1/2}(\R)$ be a minimising 
sequence for $\HH$ over $S_\mu$ with $\sup_{m \in {\mathbb N}} \| \eta_m \|_2 < M$. The inequality
\[
\HH(\eta_m,\xi_{\eta_m})\le \HH(\eta_m,\xi_m)
\]
shows that $\{(\eta_k,\xi_{\eta_k})\}\subset U\times H_\star^{1/2}(\R)$ 
is also a minimising sequence; it follows that $\{\eta_m\}\subset U\sm\{0\}$ 
is a minimising sequence for $\JJ_\mu$ which therefore converges (up to translations
and subsequences) in $H^r(\R)$, $r \in [0,2)$, to a minimiser $\eta$ of $\JJ_\mu$ over $U\sm\{0\}$.

The relations \eqn{Formula for minimiser and LM} show that
$\xi_{\eta_m}\rightarrow \xi_\eta$ in $H_\star^{1/2}(\R)$, and using this result and the calculation
\begin{eqnarray*}
c\|\xi_m-\xi_{\eta_m}\|_{*,1/2}^2 & \le & \frac{1}{2}\langle G(\eta_m)(\xi_m-\xi_{\eta_m}),(\xi_m-\xi_{\eta_m})\rangle \\
&=  & 2\HH(\eta_m,\xi_m)+2\HH(\eta_m,\xi_{\eta_m})-4 \HH(\eta_m, 
{\textstyle \frac{1}{2}}(\xi_m+\xi_{\eta_m}))\\
&\le & 2\HH(\eta_m,\xi_m)+2\HH(\eta_m,\xi_{\eta_m})-4c_\mu \\
&\rightarrow & 2c_\mu+2c_\mu-4 c_\mu \\
& = & 0
\end{eqnarray*}
as $n \rightarrow \infty$ (recall that $\HH(\eta_m, \xi) \geq \HH(\eta_m,\xi_{\eta_m}) = \JJ(\eta_m) \geq c_\mu$
for all $\xi \in H_\star^{1/2}(\R)$),
one finds that $\xi_m\rightarrow \xi_\eta$ in $H_\star^{1/2}(\R)$
as $m \rightarrow \infty$.\qed

\subsection{Convergence to solitary-wave solutions of model equations} \label{Convergence}

\subsubsection{The case \boldmath$\beta>\beta_\mathrm{c}$\unboldmath}

Suppose that $\eta$ is a minimiser of $\JJ$ over $U\sm\{0\}$, write $\eta=\eta_1+\eta_2$
according to the decomposition introduced in Section \ref{Preliminaries},
and define $\phi_\eta \in H^2(\R)$ by the formula
$$\eta_1(x)=\mu^\frac{2}{3}\phi_\eta(\mu^\frac{1}{3} x).$$
In this section we prove that $\dist(\phi_\eta, D_\mathrm{KdV}) \rightarrow 0$
as $\mu \downarrow 0$, uniformly over $\eta \in B_\mu$, where $D_\mathrm{KdV}$ is the
set of solitary-wave solutions to the Korteweg-deVries equation and
`$\dist$' denotes the distance in $H^1(\R)$.

\begin{remark}
Observe that
$$
\begin{Bmatrix} \KK_2(\eta) \\ \GG_2(\eta) \\ \LL_2(\eta) \end{Bmatrix}
= \begin{Bmatrix} \KK_2(\eta_1) \\ \GG_2(\eta_1) \\ \LL_2(\eta_1) \end{Bmatrix}
+\underbrace{\begin{Bmatrix} \KK_2(\eta_2) \\ \GG_2(\eta_2) \\ \LL_2(\eta_2) \end{Bmatrix}}_{\begin{array}{c}
=O(\|\eta\|_2^2) \\[1mm]= O(\mu^{2+\alpha})\end{array}}
$$
because $\hat{\eta}_1$ and $\hat{\eta}_2$ have disjoint supports, and
$$\GG_2(\eta_1) = -\frac{\mu\omega}{4}\int_{-\infty}^\infty \phi_\eta^2 \dx, \qquad
\KK_2(\eta_1) = \frac{\mu}{2}\int_{-\infty}^\infty \phi_\eta^2 \dx,$$
while the estimates
$$\int_{-\infty}^\infty (|k|\coth |k|-1) |\hat{\eta}_1|^2 \dk \leq c \int_{-\infty}^\infty k^2 |\hat{\eta}_1|^2 \dk
 = c\|\eta^\prime\|_0^2 \leq c\mu^{2\alpha}\nn \eta\nn_\alpha^2 \leq c\mu^{1+2\alpha},$$
$$\int_{-\infty}^\infty (|k|\coth |k|-1-{\textstyle\frac{1}{3}}k^2)|\hat{\eta}_1|^2 \dk \leq c \int_{-\infty}^\infty k^4 |\hat{\eta}_1|^2 \dk
 = c\|\eta^{\prime\prime}\|_0^2 \leq c\mu^{4\alpha}\nn \eta\nn_\alpha^2=c\mu^{1+4\alpha}$$
 show that
 $$\LL_2(\eta_1) = \frac{\mu}{2}\int_{-\infty}^\infty \phi_\eta^2 \dx + O(\mu^{1+2\alpha})$$
 and
$$\LL_2(\eta_1) = \frac{\mu}{2}\int_{-\infty}^\infty \phi_\eta^2 \dx - \frac{\beta}{3} \mu^\frac{5}{3}\int_{-\infty}^\infty (\phi_\eta^\prime)^2 \dx
+ O(\mu^{1+4\alpha}).$$
Furthermore, Corollary \ref{Strong ST final MM estimates} implies that
$$\MM_\mu(\eta) = \frac{1}{2}\left(\frac{\omega^2}{3}+1\right)\mu^\frac{5}{3} \int_{-\infty}^\infty \phi_\eta^3\dx + o(\mu^\frac{5}{3}).$$
\end{remark}

Our first result concerns the convergence of the $L^2(\R)$-norm of minimisers of $\JJ_\mu$ over $U\sm\{0\}$.

\begin{proposition} \label{Strong ST convergence in L2}
The estimate $\|\phi_\eta\|_0^2 = 4(\omega^2+4)^{-\frac{1}{2}} + O(\mu^{2\alpha})$ holds
for each $\eta \in B_\mu$.
\end{proposition}
{\bf Proof.} It follows from
$$\left|\frac{\mu + \GG_2(\eta)}{\LL_2(\eta)}-\nu_0\right| \leq c \mu^{\frac{\alpha}{2}+\frac{1}{2}}, \qquad \LL(\eta) \leq c \mu$$
that
$$\nu_0 \LL_2(\eta)-\GG_2(\eta) = \mu + O(\mu^{\frac{\alpha}{2}+\frac{3}{2}}),$$
and the result is obtained by combining this estimate with 
$$\nu_0 \LL_2(\eta)-\GG_2(\eta) = \frac{1}{4}\hspace{-1mm}\underbrace{(2\nu_0+\omega)}_{\displaystyle = \sqrt{\omega^2+4}}\hspace{-1mm}\mu \int_{-\infty}^\infty \phi_\eta^2 \dx + O(\mu^{1+2\alpha}).\eqno{\Box}$$

The next step is to show that the Korteweg-deVries energy $\EE_\mathrm{KdV}(\phi_\eta)$
corresponding to a minimiser $\eta$ of $\JJ_\mu$ over $U\sm\{0\}$ approaches
$c_\mathrm{KdV}$ in the limit $\mu \downarrow 0$.

\begin{theorem} \label{Convergence to KdV 1}
\quad
\begin{list}{(\roman{count})}{\usecounter{count}}
\item
$c_\mu = 2\nu_0\mu+c_\mathrm{KdV} \mu^\frac{5}{3}+ o(\mu^\frac{5}{3})$;
\item Each $\eta \in B_\mu$ satisfies $\EE_\mathrm{KdV}(\phi_\eta) \rightarrow c_\mathrm{KdV}$ as $\mu \downarrow 0$.
\end{list}
\end{theorem}
{\bf Proof.} Notice that
\begin{eqnarray}
c_\mu & = &  \JJ_\mu(\eta) \nonumber \\
& = & \KK_2(\eta) + \frac{(\mu+\GG_2(\eta))^2}{\LL_2(\eta)} +\MM_\mu(\eta) \nonumber \\
& = & 2\nu_0\mu + \KK_2(\eta) + 2\nu_0 \GG_2(\eta) - \nu_0^2 \LL_2(\eta) + \left(\frac{\mu+\GG_2(\eta)}{\sqrt{\LL_2(\eta)}} - \nu_0\sqrt{\LL_2(\eta)}\right)^{\!\!2} + \MM_\mu(\eta) \nonumber \\
& \geq & 2\nu_0\mu + \KK_2(\eta) + 2\nu_0 \GG_2(\eta) - \nu_0^2 \LL_2(\eta) + \MM_\mu(\eta) \nonumber \\
& = & 2\nu_0\mu + \frac{1}{2}\mu^\frac{5}{3}\int_{-\infty}^\infty \left(\left(\beta-\frac{\nu_0^2}{3}\right)(\phi_\eta^\prime)^2+\left(\frac{\omega^2}{3}+1\right)\phi_\eta^3
\right)\dx  + o(\mu^\frac{5}{3}) \nonumber \\
& = & 2\nu_0\mu + \mu^\frac{5}{3}\EE_\mathrm{KdV}(\phi_\eta) + o(\mu^\frac{5}{3}), \label{Strong ST minimum upper bound}
\end{eqnarray}
and combining this estimate with Lemma \ref{Strong ST test function} yields
$$
\EE_\mathrm{KdV}(\phi_\eta) \leq c_\mathrm{KdV} + o(1).
$$

A straightforward scaling argument shows that
$$\inf \{\EE_\mathrm{KdV}(\phi): \phi \in H^1(\R), \|\phi\|_0^2 = 4(\omega^2+4)^{-\frac{1}{2}}a\} = a^\frac{5}{3} c_\mathrm{KdV},$$
whence
$$\EE_\mathrm{KdV}(\phi_\eta)\ \geq\ (1+ O(\mu^{2\alpha}))^\frac{5}{3}c_\mathrm{KdV}\ =\ c_\mathrm{KdV} + o(1)$$
because $\|\phi_\eta\|_0^2 = 4(\omega^2+4)^{-\frac{1}{2}} + O(\mu^{2\alpha})$ (see
Proposition \ref{Strong ST convergence in L2}),
and it follows from inequality \eqn{Strong ST minimum upper bound} that
$$c_\mu \geq 2 \nu_0\mu + \mu^\frac{5}{3}c_\mathrm{KdV} + o(\mu^\frac{5}{3}).$$
The complementary estimate
$$c_\mu \leq 2 \nu_0\mu + \mu^\frac{5}{3}c_\mathrm{KdV} + o(\mu^\frac{5}{3}).$$
is a consequence of Lemma \ref{Strong ST test function}.\qed

We now present our main convergence result.

\begin{theorem} \label{Convergence to KdV 2}
The set $B_\mu$ of minimisers of $\JJ_\mu$ over $U\sm\{0\}$ satisfies
$$
\sup_{\eta \in B_\mu} \inf_{x \in \R} \|\phi_\eta - \phi_\mathrm{KdV}(\cdot+x)\|_1 \rightarrow 0
$$
as $\mu \downarrow 0$.
\end{theorem}
{\bf Proof.} Suppose that the limit is positive, so that there exists $\varepsilon>0$ and a sequence $\{\mu_m\}$ with $\mu_m \downarrow 0$ such that
$$\sup_{\eta \in C_{\mu_m}}  \inf_{x \in \R} \|\phi_\eta - \phi_\mathrm{KdV}(\cdot+x)\|_1\geq \varepsilon, \qquad m \in \N$$
and hence a further sequence $\{\eta_m\} \subset U \sm\{0\}$ with $\eta_m \in C_{\mu_m}$ and
$$
\dist(\phi_{\eta_m},D_\mathrm{KdV})\ =\  \inf_{x \in \R} \|\phi_\eta - \phi_\mathrm{KdV}(\cdot+x)\|_1\ \geq\ \frac{\varepsilon}{2}, \qquad m \in \N.
$$
On the other hand $\EE_\mathrm{KdV}(\phi_{\eta_m})\to c_\mathrm{KdV}$
and $\|\phi_{\eta_m}\|_0^2\to 4(\omega^2+4)^{-\frac{1}{2}}$ as $n\to \infty$ (see Proposition \ref{Strong ST convergence in L2} and
Theorem \ref{Convergence to KdV 1}(ii)); combining Lemma \ref{Variational KdV}(ii)
with a straightforward scaling argument, we arrive at the contradiction of the existence
of a sequence $\{x_m\} \subset \R$ such that a subsequence of $\{\phi_{\eta_m}(x_m+\cdot)\}$
converges in $H^1(\R)$ to an element of $D_\mathrm{KdV}$.\qed

\begin{remark}
The previous theorem implies that $\{\|\phi_\eta\|_1: \eta \in B_\mu\}$ is bounded, so that
\begin{eqnarray*}
\|\hat \eta_1\|_{L^1(\R)}^2 & \leq & \left(\int_{-\infty}^\infty\frac{1}{1+\mu^{-\frac{2}{3}}k^2} \dk\right)
\!\!\!\left(\int_{-\infty}^\infty (1+\mu^{-\frac{2}{3}}k^2)|\hat \eta_1(k)|^2\dk\right)  \\
& = & \mu^\frac{2}{3}\left(\int_{-\infty}^\infty\frac{1}{1+\mu^{-\frac{2}{3}}k^2} \dk\right)
\!\!\!\left(\int_{-\infty}^\infty (1+\mu^{-\frac{2}{3}}k^2)\left|\hat \phi_\eta\left(\frac{k}{\mu^\frac{1}{3}}\right)\right|^2\dk\right)\\
& = & 2\pi\mu^\frac{4}{3}\|\phi_\eta\|_1^2 \\
& \leq & c\mu^\frac{4}{3}
\end{eqnarray*}
and hence $\|\eta_1\|_{1,\infty}$, $\| K^0 \eta_1\|_\infty \leq c\mu^\frac{2}{3}$ (see equations \eqn{max norm},
\eqn{max norm of K0}), and it follows from inequalities \eqn{Strong ST Gives result for eta2},
\eqn{Strong ST Gives result for eta1}
that $\nn \eta_1 \nn_{1/3}^2 \leq c\mu$, $\| \eta_2 \|_2^2 \leq \mu^\frac{7}{3}$.
For $\eta \in B_\mu$ Lemma \ref{Strong ST weighted norm estimate} therefore also holds with $\alpha=\frac{1}{3}$
(the result predicted in the Korteweg-deVries scaling limit).
\end{remark}

Our final result shows that the speed $\nu_\mu$ of a solitary wave corresponding to $\eta \in B_\mu$, which is given by the formula
$$\nu_\mu = \frac{\mu+\GG(\eta)}{\LL(\eta)},$$
satisfies
$$\nu_\mu = \nu_0 + 2(\omega^2+4)^{-\frac{1}{2}}\nu_\mathrm{KdV}\mu^\frac{2}{3} + o(\mu^\frac{2}{3})$$
uniformly over $\eta \in B_\mu$.

\begin{theorem}
The set $B_\mu$ of minimisers of $\JJ_\mu$ over $U\sm\{0\}$ satisfies
$$\sup_{\eta \in B_\mu} \left| \frac{\mu + \GG(\eta)}{\LL(\eta)} - (\nu_0 + 2(\omega^2+4)^{-\frac{1}{2}}\nu_\mathrm{KdV}\mu^\frac{2}{3})\right| = o(\mu^\frac{2}{3}).$$
\end{theorem}
{\bf Proof.} Using the identity
$$
\frac{\mu + \GG(\eta)}{\LL(\eta)} = \frac{1}{2\mu}(c_\mu-\MM_\mu(\eta)) + \frac{1}{4\mu}(\langle \MM_\mu^\prime(\eta),\eta \rangle + 4\mu\tilde{\MM}_\mu(\eta))
$$
(see the proof of Proposition \ref{Speed estimate}), we find that
\begin{eqnarray*}
\frac{\mu + \GG(\eta)}{\LL(\eta)} & = & \nu_0 + \frac{1}{2}c_\mathrm{KdV}\mu^\frac{2}{3} + \frac{1}{8\mu}\left(\frac{\omega^2}{3}+1\right)\int_{-\infty}^\infty \eta_1^3 \dx
 +o(\mu^\frac{2}{3}) \\
 & = & \nu_0 + \frac{1}{2}c_\mathrm{KdV}\mu^\frac{2}{3} + \frac{1}{8}\left(\frac{\omega^2}{3}+1\right)\mu^\frac{2}{3} \int_{-\infty}^\infty \phi_\eta^3\dx
 +o(\mu^\frac{2}{3}) \\
  & = & \nu_0 + \frac{1}{2}\EE_\mathrm{KdV}(\phi_\mathrm{KdV})\mu^\frac{2}{3} + \frac{1}{8}\left(\frac{\omega^2}{3}+1\right)\mu^\frac{2}{3} \int_{-\infty}^\infty \phi_\mathrm{KdV}^3 \dx
 +o(\mu^\frac{2}{3}) \\
  & = & \nu_0 + \frac{1}{4}\mu^\frac{2}{3}  \underbrace{\int_{-\infty}^\infty \left(\left(\beta-\frac{\nu_0^2}{3}\right) (\phi_\mathrm{KdV}^\prime)^2
  + \frac{3}{2}\left(\frac{\omega^2}{3}+1\right)\phi_\mathrm{KdV}^3 \right)\dx}_{\displaystyle = 8(\omega^2+4)^{-\frac{1}{2}}\nu_\mathrm{KdV}}
+ o(\mu^\frac{2}{3}) \\
 & = & \nu_0 + 2(\omega^2+4)^{-\frac{1}{2}}\nu_\mathrm{KdV}\mu^\frac{2}{3} + o(\mu^\frac{2}{3}),
 \end{eqnarray*}
in which Theorem \ref{Convergence to KdV 1}(i),
Corollary \ref{Strong ST final MM estimates} and Theorem \ref{Convergence to KdV 2}
have been used.\qed.

\subsubsection{The case \boldmath$\beta<\beta_\mathrm{c}$\unboldmath}

Suppose that $\eta$ is a minimiser of $\JJ_\mu$ over $U\sm\{0\}$, write
$\eta=\eta_1-H(\eta_1)+\eta_3$ and $\eta_1=\eta_1^++\eta_1^-$ according to the decompositions
introduced in Section \ref{Weak ST},
and define $\phi_\eta \in H^2(\R)$ by the formula
$$\eta_1^+(x) = \frac{1}{2}\mu \phi_\eta(\mu x)\e^{\i k_0 x}.$$
In this section we prove that $\dist(\phi_\eta, D_\mathrm{NLS}) \rightarrow 0$
as $\mu \downarrow 0$, uniformly over $\eta \in B_\mu$, where $D_\mathrm{NLS}$ is the
set of solitary-wave solutions to the nonlinear Schr\"{o}dinger equation and
`$\dist$' denotes the distance in $H^1(\R)$.

\begin{remark}
Note that
\begin{equation}
\begin{Bmatrix} \KK_2(\eta) \\ \GG_2(\eta) \\ \LL_2(\eta) \end{Bmatrix}
= \begin{Bmatrix} \KK_2(\eta_1) \\ \GG_2(\eta_1) \\ \LL_2(\eta_1) \end{Bmatrix}
+\begin{Bmatrix} \KK_2(-H(\eta)+\eta_3) \\ \GG_2(-H(\eta)+\eta_3) \\ \LL_2(-H(\eta)+\eta_3) \end{Bmatrix}
\label{Splitting of the twos}
\end{equation}
because $\hat{\eta}_1$ and $\FF[-H(\eta)+\eta_3]$ have disjoint supports.
\end{remark}

Our first result concerns the convergence of the $L^2(\R)$-norm of minimisers of $\JJ_\mu$
over $U_2\sm\{0\}$.

\begin{proposition} \label{Weak ST convergence in L2}
The estimate
$\|\phi_\eta\|_0^2 = \left(\frac{1}{4}\nu_0f(k_0)+\frac{\omega}{8}\right)^{-1}
+ O(\mu^\alpha)$
holds for each $\eta \in B_\mu$.
\end{proposition}
{\bf Proof.} It follows from
$$\left|\frac{\mu + \GG_2(\eta)}{\LL_2(\eta)}-\nu_0\right| \leq c \mu^{1+\alpha}, \qquad \LL_2(\eta) \leq c \mu$$
that
\begin{equation}
\nu_0 \LL_2(\eta)-\GG_2(\eta) = \mu + O(\mu^{2+\alpha}). \label{Weak ST convergence in L2 part 1}
\end{equation}
On the other hand
\begin{eqnarray*}
\nu_0\LL_2(\eta) - \GG_2(\eta) & = & \nu_0 \LL_2(\eta_1)-\GG_2(\eta_1) + O(\|H(\eta)\|_2^2 + \|\eta_3\|_2^2) \\
& = & \nu_0 \LL_2(\eta_1)-\GG_2(\eta_1) + O(\mu^{2+\alpha}) \\
& = & \nu_0\int_{-\infty}^\infty \eta_1^+ K^0 \eta_1^- \dx + \frac{\omega}{2}\int_{-\infty}^\infty \eta_1^+\eta_1^-\dx
+ O(\mu^{2+\alpha}) \\
& = & \left(\nu_0f(k_0)+\frac{\omega}{2}\right) \int_{-\infty}^\infty \eta_1^+\eta_1^- \dx + O(\mu^{1+\alpha}) \\
& = & \left(\frac{1}{4}\nu_0f(k_0)+\frac{\omega}{8}\right)\mu \int_{-\infty}^\infty |\phi_\eta|^2 \dx+ O(\mu^{1+\alpha}),
\end{eqnarray*}
and the result is obtained by combining this estimate with \eqn{Weak ST convergence in L2 part 1}.\qed

The next step is to show that the nonlinear Schr\"{o}dinger energy $\EE_\mathrm{NLS}(\phi_\eta)$
corresponding to a minimiser $\eta$ of $\JJ_\mu$ over $U\sm\{0\}$ approaches $c_\mathrm{NLS}$
in the limit $\mu \downarrow 0$.

\begin{theorem} \label{Convergence to NLS 1}
\quad
\begin{list}{(\roman{count})}{\usecounter{count}}
\item
$c_\mu = 2\nu_0\mu+c_\mathrm{NLS} \mu^3+ o(\mu^3)$;
\item Each $\eta \in B_\mu$ satisfies $\EE_\mathrm{NLS}(\phi_\eta) \rightarrow c_\mathrm{NLS}$ as $\mu \downarrow 0$.
\end{list}
\end{theorem}
{\bf Proof.} Notice that
\begin{eqnarray}
c_\mu & = &  \JJ_\mu(\eta) \nonumber \\
& = & \KK_2(\eta) + \frac{(\mu+\GG_2(\eta))^2}{\LL_2(\eta)} +\MM_\mu(\eta) \nonumber \\
& = & 2\nu_0\mu + \KK_2(\eta) + 2\nu_0 \GG_2(\eta) - \nu_0^2 \LL_2(\eta) + \left(\frac{\mu+\GG_2(\eta)}{\sqrt{\LL_2(\eta)}} - \nu_0\sqrt{\LL_2(\eta)}\right)^{\!\!2} + \MM_\mu(\eta) \nonumber \\
& \geq & 2\nu_0\mu + \KK_2(\eta) + 2\nu_0 \GG_2(\eta) - \nu_0^2 \LL_2(\eta) + \MM_\mu(\eta),
\label{Conv est 0}
\end{eqnarray}
where
\begin{eqnarray}
\lefteqn{\KK_2(\eta) + 2\nu_0 \GG_2(\eta) - \nu_0^2 \LL_2(\eta)}\qquad \nonumber \\
& = & (\KK_2 + 2\nu_0 \GG_2 - \nu_0^2 \LL_2)(\eta_1) + (\KK_2 + 2\nu_0 \GG_2 - \nu_0^2 \LL_2)(-H(\eta)+\eta_3).
\label{Conv est 1}
\end{eqnarray}

The second term on the right-hand side of \eqn{Conv est 1} is estimated using the calculation
\begin{eqnarray*}
\lefteqn{(\KK_2 + 2\nu_0 \GG_2 - \nu_0^2 \LL_2)(-H(\eta)+\eta_3)} \qquad \\
& = & (\KK_2 + 2\nu_0 \GG_2 - \nu_0^2 \LL_2)(H(\eta)) + O(\|H(\eta)\|_2\|\eta_3\|_2) + O(\|\eta_3\|_2^2) \\
& = & \frac{1}{2}\int_{-\infty}^\infty g(k) |\FF[H(\eta)]|^2 \dk+ o(\mu^3) \\
& = & \frac{1}{2}\int_{-\infty}^\infty g(k)^{-1}|\FF[\KK_3(\eta_1) + 2\nu_0 \GG_3(\eta_1) - \nu_0^2 \LL_3(\eta_1)]|^2 \dk + o(\mu^3) \\
&= & -\frac{1}{2}\big(\KK_3(\eta) + 2\nu_0 \GG_3(\eta) - \nu_0^2 \LL_3(\eta) \big)+ o(\mu^3) \\
& = & - \frac{A_3}{2} \int_{-\infty}^\infty \eta_1^4 \dx + o(\mu^3) \\
& = & -\frac{3A_3}{16}\mu^3 \int_{-\infty}^\infty |\phi_\eta|^4 \dx + o(\mu^3),
\end{eqnarray*}
where we have used
Proposition \ref{lot in threes step 2}, equation \eqn{lcomb threes preliminary} and Proposition \ref{lot in threes}.
Turning to the first term on the right-hand side of \eqn{Conv est 1}, write
$$
(\KK_2 + 2\nu_0 \GG_2 - \nu_0^2 \LL_2)(\eta_1) = \frac{1}{2}\int_{-\infty}^\infty g(k) |\hat{\eta}_1|^2 \dk =
 \int_{-\infty}^\infty g(k)|\hat{\eta}_1^+(k)|^2 \dk.
$$
and note that
$$g(k) = \frac{1}{2}g^{\prime\prime}(k_0)(k-k_0)^2 + O(|k-k_0|^3), \qquad k \in [k_0-\delta_0,k_0+\delta_0].$$
One finds that
$$
\int_{-\infty}^\infty (k-k_0)^2 |\hat{\eta}_1^+(k)|^2 \dk=\!\!\int_{-\infty}^\infty k^2 |\hat{\eta}_1^+ (k+k_0)|^2 \dk
\!=\!\frac{\mu^2}{4}\! \int_{-\infty}^\infty \left|\frac{\mathrm{d}}{\mathrm{d}x} \phi_\eta(\mu x)\right|^2 \dx
\!=\!\frac{\mu^3}{4}\int_{-\infty}^\infty |\phi_\eta^\prime|^2 \dx
$$
(because $\hat{\eta}_1^+(k+k_0) = \frac{\mu}{2}\FF[\phi_\eta(\mu x)]$)
and
$$
\int_{-\infty}^\infty (k-k_0)^3 |\hat{\eta}_1^+(k)|^2 \dk \leq c\mu^{3\alpha} \nn \eta_1 \nn_\alpha^2 = O(\mu^{1+3\alpha}),
$$
so that
$$\int_{-\infty}^\infty \big(g(k)-{\textstyle\frac{1}{2}}(k-k_0)^2\big) |\hat{\eta}_1^+(k)|^2 \dk= o(\mu^3).$$
Altogether these calculations show that
\begin{eqnarray}
\lefteqn{(\KK_2 + 2\nu_0 \GG_2 - \nu_0^2 \LL_2)(\eta_1)} \hspace{1cm} \nonumber \\
& =& \frac{1}{8}g^{\prime\prime}(k_0)\mu^3\int_{-\infty}^\infty |\phi_\eta^\prime|^2\dx
-\frac{3A_3}{16}\mu^3 \int_{-\infty}^\infty |\phi_\eta|^4 \dx+ o(\mu^3). \label{Conv est 2}
\end{eqnarray}

Substituting \eqn{Conv est 2} and
$$
\MM_\mu(\eta) 
\ =\ (A_3+A_4)\int_{-\infty}^\infty \eta_1^4 \dx+ o(\mu^3) \\
\ =\ \frac{3}{8}(A_3+A_4)\mu^3 \int_{-\infty}^\infty |\phi_\eta|^4 \dx + o(\mu^3)
$$
(see Corollary \ref{Weak ST final MM estimates}) into inequality \eqn{Conv est 0} yields
\begin{eqnarray}
c_\mu & \geq & 2\nu_0\mu + \frac{1}{8}g^{\prime\prime}(k_0)\mu^3\int_{-\infty}^\infty |\phi_\eta^\prime|^2 \dx
+ \frac{3}{8}\left(\frac{A_3}{2}+A_4\right)\mu^3 \int_{-\infty}^\infty |\phi_\eta|^4 \dx + o(\mu^3) \nonumber \\
& = & 2\nu_0\mu + \mu^3 \EE_\mathrm{NLS}(\phi_\eta) + o(\mu^3), \label{Weak ST minimum upper bound}
\end{eqnarray}
and combining this estimate with Lemma \ref{Weak ST test function} yields
$$
\EE_\mathrm{NLS}(\phi_\eta) \leq c_\mathrm{NLS} + o(1).
$$

A straightforward scaling argument shows that
$$\inf \{\EE_\mathrm{NLS}(\phi): \phi \in H^1(\R), \|\phi\|_0^2 = \left(\textstyle\frac{1}{4}\nu_0f(k_0)+\frac{\omega}{8}\right)^{-1}a\} = a^3 c_\mathrm{NLS},$$
whence
$$\EE_\mathrm{NLS}(\phi_\eta)\ \geq\ (1+ O(\mu^\alpha))^3c_\mathrm{NLS}\ =\ c_\mathrm{NLS} + o(1)$$
because $\|\phi_\eta\|_0^2 = \left(\frac{1}{4}\nu_0f(k_0)+\frac{\omega}{8}\right)^{-1}
+ O(\mu^\alpha)$ (see
Proposition \ref{Weak ST convergence in L2}),
and it follows from inequality \eqn{Weak ST minimum upper bound} that
$$c_\mu \geq 2 \nu_0\mu + \mu^3c_\mathrm{NLS} + o(\mu^3).$$
The complementary estimate
$$c_\mu \leq 2 \nu_0\mu + \mu^3c_\mathrm{NLS} + o(\mu^3).$$
is a consequence of Lemma \ref{Weak ST test function}.\qed

Our main convergence result is derived from Theorem \ref{Convergence to NLS 1} in the same
way as the corresponding result for $\beta>\beta_\mathrm{c}$ (see Appendix A.1).

\begin{theorem} \label{Convergence to NLS 2}
The set $B_\mu$ of minimisers of $\JJ_\mu$ over $U\sm\{0\}$ satisfies
$$\sup_{\eta \in B_\mu} \inf_{\substack{\omega \in [0,2\pi],\\ x \in \R}} \|\phi_\eta-e^{\mathrm{i}\omega}\phi_\mathrm{NLS}(\cdot+x)\|_1 \rightarrow 0$$
as $\mu \downarrow 0$.
\end{theorem}
\begin{remark}
The previous theorem implies that $\{\|\phi_\eta\|_1: \eta \in B_\mu\}$ is bounded, so that
\begin{eqnarray*}
\|\hat \eta_1\|_{L^1(\R)}^2 & \leq & 2\left(\int_{k_0-\delta_0}^{k_0+\delta_0}\frac{1}{1+\mu^{-2}(k-k_0)^2} \dk\right)
\!\!\!\left(\int_{k_0-\delta_0}^{k_0+\delta_0} (1+\mu^{-2}(k-k_0)^2)|\hat \eta_1(k)|^2\dk\right)  \\
& \leq & 2\left(\int_{-\infty}^\infty\frac{1}{1+\mu^{-2}(k-k_0)^2} \dk\right)
\!\!\!\left(\int_{-\infty}^\infty (1+\mu^{-2}(k-k_0)^2)\left|\hat \phi_\eta\left(\frac{k-k_0}{\mu}\right)\right|^2\dk\right)\\
& = & 2\pi\mu^2\|\phi_\eta\|_1^2 \\
& \leq & c\mu^2
\end{eqnarray*}
and hence $\|\eta_1\|_{1,\infty}$, $\| K^0 \eta_1\|_{1,\infty} \leq c\mu$ (see equations \eqn{max norm} and
\eqn{max norm of K0}), and it follows from Proposition \ref{Estimate for H} and inequalities
\eqn{Weak ST Gives result for eta3}, \eqn{Weak ST Gives result for eta1}
that
$$\nn \eta_1 \nn_1^2 \leq c\mu, \qquad
\|H(\eta_1)\|_2^2 \leq c\mu^3, \qquad
\|u_3\|_2^2 \leq c\mu^5.$$
For $\eta \in B_\mu$ Lemma \ref{Weak ST nn estimate} therefore also holds with $\alpha=1$
(the result predicted in the nonlinear Schr\"{o}dinger scaling limit).
\end{remark}

Our final result shows that the speed $\nu_\mu$ of a solitary wave corresponding to $\eta \in B_\mu$, which is given by the formula
$$\nu_\mu+\frac{\mu+\GG(\eta)}{\LL(\eta)},$$
satisfies
$$\nu_\mu = \nu_0 + 4(\omega+2\nu_0f(k_0))^{-1}\nu_\mathrm{NLS}\mu^2 + o(\mu^2)$$
uniformly over $\eta \in B_\mu$.

\begin{theorem}
The set $B_\mu$ of minimisers of $\JJ_\mu$ over $U\sm\{0\}$ satisfies
$$\sup_{\eta \in B_\mu} \left| \frac{\mu + \GG(\eta)}{\LL(\eta)} - 
(\nu_0 + 4(\omega+2\nu_0f(k_0))^{-1}\nu_\mathrm{NLS}\mu^2)\right| = o(\mu^2).$$
\end{theorem}
{\bf Proof.} Using the identity
$$
\frac{\mu + \GG(\eta)}{\LL(\eta)} = \frac{1}{2\mu}(c_\mu-\MM_\mu(\eta)) + \frac{1}{4\mu}(\langle \MM_\mu^\prime(\eta),\eta \rangle + 4\mu\tilde{\MM}_\mu(\eta))
$$
(see the proof of Proposition \ref{Speed estimate}), we find that
\begin{eqnarray*}
\frac{\mu + \GG(\eta)}{\LL(\eta)} & = & \nu_0 + \frac{1}{2}c_\mathrm{NLS}\mu^2 + 
\frac{1}{2\mu}\left(\frac{1}{2}A_3+A_4\right)\int_{-\infty}^\infty \eta_1^4 \dx + o(\mu^2) \\
& = & \nu_0 + \frac{1}{2}c_\mathrm{NLS}\mu^2 + 
\frac{3}{16}\left(\frac{1}{2}A_3+A_4\right)\mu^2\int_{-\infty}^\infty |\phi_\eta|^4 \dx + o(\mu^2) \\
& = & \nu_0 + \frac{1}{2}\EE_\mathrm{NLS}(\phi_\mathrm{NLS})\mu^2 + 
\frac{3}{16}\left(\frac{1}{2}A_3+A_4\right)\mu^2\int_{-\infty}^\infty |\phi_\mathrm{NLS}|^4 \dx + o(\mu^2) \\
& = & \nu_0 + \frac{1}{4} \mu^2\underbrace{\int_{-\infty}^\infty\left(
\frac{1}{4}g^{\prime\prime}(k_0) |\phi^\prime_\mathrm{NLS}|^2
+\frac{3}{2}\left(\frac{1}{2}A_3+A_4\right) |\phi_\mathrm{NLS}|^4\right) \dx}_{\textstyle
= 2\left(\frac{1}{4}\nu_0f(k_0)+\frac{\omega}{8}\right)^{\!\!-1}\nu_\mathrm{NLS}}
+ o(\mu^2) \\
& = & \nu_0 + 4(\omega+2\nu_0f(k_0))^{-1}\nu_\mathrm{NLS}\mu^2 + o(\mu^2),
\end{eqnarray*}
in which Theorem \ref{Convergence to NLS 1}(i), Corollary \ref{Weak ST final MM estimates}
and Theorem \ref{Convergence to NLS 2} have been used.\qed

\appendix
\section*{Appendix A: Proof of Lemma \ref{Test function}(i)}
\setcounter{section}{1}

\subsection*{A.1\quad The case $\beta>\beta_\mathrm{c}$}

\begin{lemma} \label{Strong ST test function}
Suppose that $\mu>0$.
There exists a continuous, invertible mapping $\mu \rightarrow \alpha(\mu)$ such
that
$$\JJ_\mu(\eta^\star) = 2\nu_0\mu + c_\mathrm{KdV}\mu^\frac{5}{3} + o(\mu^\frac{5}{3}),$$
where
$$\eta^\star(x) = \alpha^2 \phi_\mathrm{KdV}(\alpha x).$$
\end{lemma}
{\bf Proof.} Let us first note that
$$K^0 \eta^\star - \eta^\star + {\textstyle\frac{1}{3}}(\eta^\star)^{\prime\prime}=
\FF^{-1}[\underbrace{(|k|\coth|k|-1-{\textstyle\frac{1}{3}}|k|^2)}_{\displaystyle \leq c|k|^4}\hat{\eta}^\star] = \underline{O}(\alpha^\frac{11}{2})$$
and hence
$$K^0 \eta^\star - \eta^\star = \FF^{-1}[(|k|\coth|k|-1)\hat{\eta}^\star] = \underline{O}(\alpha^\frac{7}{2}).$$
Using these estimates and $\|\eta^\star\|_0 = O(\alpha^\frac{3}{2})$, one finds that
$$
\KK_2(\eta^\star) = \frac{\alpha^3}{2}\int_{-\infty}^\infty \phi_\mathrm{KdV}^2 \dx+ \frac{\alpha^5}{2}\beta\int_{-\infty}^\infty \phi_\mathrm{KdV}^{\prime 2} \dx, \qquad
\GG_2(\eta^\star) = -\frac{\alpha^3}{4}\omega \int_{-\infty}^\infty \phi_\mathrm{KdV}^2 \dx,
$$
$$\LL_2(\eta^\star)\ =\ \frac{1}{2}\int_{-\infty}^\infty \eta^\star K^0 \eta^\star \dx \ =\ \frac{\alpha^3}{2}\int_{-\infty}^\infty \phi_\mathrm{KdV}^2 \dx+
\frac{\alpha^5}{6}\int_{-\infty}^\infty \phi_\mathrm{KdV}^{\prime 2} \dx + O(\alpha^7),$$
and
\begin{eqnarray*}
\KK_3(\eta^\star) & = & \frac{\alpha^5}{6}\omega^2\int_{-\infty}^\infty \phi_\mathrm{KdV}^3 \dx, \\
\\
\GG_3(\eta^\star) & = & \frac{\omega}{4}\int_{-\infty}^\infty (\eta^\star)^2 K^0 \eta^\star \dx\\
& = & \frac{\omega}{4}\int_{-\infty}^\infty (\eta^\star)^3 \dx+ \frac{\omega}{4}\int_{-\infty}^\infty  (\eta^\star)^2 (K^0\eta^\star-\eta^\star)\dx \\
& = &  \frac{\alpha^5}{4}\omega \int_{-\infty}^\infty \phi_\mathrm{KdV}^3\dx + O(\alpha^7),\\
\\
\LL_3(\eta^\star) & = & \frac{1}{2}\int_{-\infty}^\infty \big(-(K^0 \eta^\star)^2 \eta^\star + (\eta^{\star\prime})^2\eta^\star\big) \dx \\
& = & -\frac{1}{2}\int_{-\infty}^\infty (\eta^\star)^3 \dx + \frac{1}{2} \int_{-\infty}^\infty \!\!\big(\! -2(K^0\eta^\star\!-\!\eta^\star)(\eta^\star)^2 
- (K^0\eta^\star\!-\!\eta^\star)^2\eta^\star + (\eta^{\star\prime})^2\eta^\star \big) \dx\\
& = & -\frac{\alpha^5}{2} \int_{-\infty}^\infty \phi_\mathrm{KdV}^3 \dx + O(\alpha^7),
\end{eqnarray*}
in which the further estimate $\|\eta^\star\|_\infty = O(\alpha^2)$ has been used (see Proposition
\ref{Formulae for the threes and fours} for the formulae for $\GG_3$, $\KK_3$ and $\LL_3$).
Finally, Proposition \ref{Size of the functionals} shows that $\GG_4(\eta^\star)$, $\KK_4(\eta^\star)$, $\LL_4(\eta^\star)$
and $\GG_\mathrm{r}(\eta^\star)$, $\KK_\mathrm{r}(\eta^\star)$, $\LL_\mathrm{r}(\eta^\star)$ are all $O(\alpha^7)$.

The above calculations show that
\begin{eqnarray*}
\lefteqn{\KK(\eta^\star) + 2 \nu_0 \GG(\eta^\star) - \nu_0^2 \LL(\eta^\star)} \qquad \\
& = & \frac{\alpha^3}{2}\underbrace{(1-\omega \nu_0 -\nu_0^2)}_{\displaystyle = 0}\int_{-\infty}^\infty \phi_\mathrm{KdV}^2 \dx
+ \frac{1}{2}\left(\beta - \frac{\nu_0^2}{3}\right)\alpha^5 \int_{-\infty}^\infty \phi_\mathrm{KdV}^{\prime 2}\dx \\
& & \qquad \mbox{}+\frac{1}{2}\Bigg(\frac{\omega^2}{3}+\underbrace{\omega \nu_0 + \nu_0^2}_{\displaystyle =1}\Bigg)\alpha^5\int_{-\infty}^\infty \phi_\mathrm{KdV}^3 \dx+ O(\alpha^7) \\
& = & \alpha^5 \EE_\mathrm{KdV}(\phi_\mathrm{KdV}) + O(\alpha^7) \\
& = & c_\mathrm{KdV}\alpha^5 + O(\alpha^7).
\end{eqnarray*}
The mapping
\begin{eqnarray*}
\alpha & \mapsto & \nu_0 \LL(\eta^\star) - \GG(\eta^\star) \\
& & = \alpha^3\left(\frac{\nu_0}{2}+\frac{\omega}{4}\right) \int_{-\infty}^\infty \phi_\mathrm{KdV}^2 \dx+ O(\alpha^5) \\
& & = \frac{\alpha^3}{4}\sqrt{\omega^2+4}  \int_{-\infty}^\infty \phi_\mathrm{KdV}^2 \dx + O(\alpha^5)
\end{eqnarray*}
is continuous and strictly increasing and therefore has a continuous inverse $\mu \mapsto \alpha(\mu)$; furthermore
$\alpha(\mu) = \mu^\frac{1}{3} + o(\mu^\frac{1}{3})$ and
$$\JJ_\mu(\eta^\star) - 2\nu_0 \mu\ =\ \KK(\eta^\star) + 2 \nu_0 \GG(\eta^\star) - \nu_0^2 \LL(\eta^\star)\ =\ c_\mathrm{KdV}\mu^\frac{5}{3} + o(\mu^\frac{5}{3}).\eqno{\Box}$$

\subsection*{A.2\quad The case $\beta<\beta_\mathrm{c}$}

\begin{lemma} \label{Weak ST test function}
Suppose that $\mu>0$.
There exists a continuous, invertible mapping $\mu \rightarrow \alpha(\mu)$ such
that
$$\JJ_\mu(\eta^\star) = 2\nu_0\mu + c_\mathrm{NLS}\mu^3 + o(\mu^3),$$
where
$$\eta^\star(x) = \alpha \phi_\mathrm{NLS}(\alpha x) \cos k_0 x - \frac{\alpha^2}{2}g(2k_0)^{-1}A_3^1 \phi_\mathrm{NLS}(\alpha x)^2 \cos 2k_0 x - \frac{\alpha^2}{2}g(0)^{-1}A_3^2 \phi_\mathrm{NLS} (\alpha x)^2.$$
\end{lemma}
{\bf Proof.} We seek a test function $\eta^\star$ of the form
$$\eta^\star(x) = \alpha \phi_\mathrm{NLS}(\alpha x) \cos k_0 x + \alpha^2 \psi(\alpha x) \cos 2k_0 x + \alpha^2 \xi (\alpha x)$$
with $\psi$, $\xi \in \SS(\R)$.

Choose $n \in {\mathbb N}$ and $\chi \in C_0^\infty(\R)$.
Straightforward calculations yield the formulae
$$K^0 (\chi(\alpha x)) = \chi(\alpha x) + S_2(x),$$
where
$$S_2(x) = \frac{1}{\alpha} \FF^{-1}\left[(|k|\coth|k|-1)\hat{\chi}\left(\frac{k}{\alpha}\right)\right],$$
and
\begin{eqnarray*}
\lefteqn{K^0(\chi(\alpha x)\cos n k_0 x) =} \\
& & f(nk_0)\chi(\alpha x) \cos n k_0 x + \alpha f^\prime(nk_0)\chi^\prime(\alpha x)\sin nk_0 x
- \frac{\alpha^2}{2}f^{\prime\prime}(nk_0) \chi^{\prime\prime}(\alpha x) \cos nk_0 x + S_1(x),
\end{eqnarray*}
where
$$S_1(x) = \frac{1}{2}\FF^{-1}\left[R_{nk_0}(k)(k-nk_0)^3\hat{\chi}\left(\frac{k-nk_0}{\alpha}\right)\right]
+\frac{1}{2}\FF^{-1}\left[R_{-nk_0}(k)(k+nk_0)^3\hat{\chi}\left(\frac{k+nk_0}{\alpha}\right)\right]
$$
and $R_\omega(k) = \frac{1}{6}f^{\prime\prime\prime}(k_\omega)$ for some $k_\omega$ between $k$ and $\omega$;
the remainder terms $S_1$ and $S_2$ satisfy the estimates
$\|S_1\|_\infty = O(\alpha^3)$, $\|S_1\|_1 = O(\alpha^\frac{7}{2})$ and
$\|S_2\|_m = O(\alpha^{n+\frac{3}{2}})$.
Furthermore, repeated integration by parts shows that
$$\int_{-\infty}^\infty \chi(\alpha x) \left\{\begin{array}{c}\ \sin \\ \cos \end{array}\right\} (mx) \dx
= O(\alpha^n)$$
for each $m \in \N$, so that
$$\int_{-\infty}^\infty \chi(\alpha x) \left\{\begin{array}{c}\ \sin \\ \cos \end{array}\right\} (m_1 x)
\cdots
\left\{\begin{array}{c}\ \sin \\ \cos \end{array}\right\} (m_\ell x)
\dx
= O(\alpha^n)$$
for all $m_1,\ldots,m_\ell \in \N$ with $m_1 \pm \ldots \pm m_\ell \neq 0$.

Estimating using the above rules, one finds that
\begin{eqnarray*}
\KK_2(\eta^\star) & \!\!\! = & \!\!\!\frac{\alpha}{4}(1+\beta k_0^2)\int_{-\infty}^\infty \phi_\mathrm{NLS}^2 \dx
+ \frac{\alpha^3}{4}\beta \int_{-\infty}^\infty \phi_\mathrm{NLS}^{\prime 2} \dx \\
& & \qquad\mbox{}
+\frac{\alpha^3}{4}(1+4\beta k_0^2) \int_{-\infty}^\infty \psi^2\dx
+\frac{\alpha^3}{2}\int_{-\infty}^\infty \xi^2 \dx+ O(\alpha^4), \\
\GG_2(\eta^\star) & \!\!\! = & \!\!\!-\frac{\alpha}{8}\omega\int_{-\infty}^\infty \phi_\mathrm{NLS}^2 \dx
-\frac{\alpha^3}{8} \omega \int_{-\infty}^\infty \psi^2\dx
- \frac{\alpha^3}{4}\omega \int_{-\infty}^\infty \xi^2\dx+ O(\alpha^4), \\
\LL_2(\eta^\star) & \!\!\! = & \!\!\!\frac{\alpha}{4}f(k_0)\int_{-\infty}^\infty \phi_\mathrm{NLS}^2 \dx
+\frac{\alpha^3}{8}f^{\prime\prime}(k_0) \int_{-\infty}^\infty \phi_\mathrm{NLS}^{\prime 2} \dx\\
& & \qquad\mbox{}
+\frac{\alpha^3}{4}f(2k_0)\int_{-\infty}^\infty \psi^2\dx + \frac{\alpha^3}{2} \int_{-\infty}^\infty \xi^2 \dx + O(\alpha^4), \\
\KK_3(\eta^\star) & \!\!\! = & \!\!\! \frac{\alpha^3}{8}\omega^2 \int_{-\infty}^\infty \phi_\mathrm{NLS}^2 \psi \dx
+ \frac{\alpha^3}{4}\omega^2 \int_{-\infty}^\infty \phi_\mathrm{NLS}^2 \xi \dx + O(\alpha^4), \\
\GG_3(\eta^\star) & \!\!\! = & \!\!\! \frac{\alpha^3}{8}\left(f(k_0)+{\textstyle\frac{1}{2}}f(2k_0)\right)\omega
\int_{-\infty}^\infty \phi_\mathrm{NLS}^2 \psi \dx+ \frac{\alpha^3}{4}\left( f(k_0) +{\textstyle\frac{1}{2}} \right) \omega \int_{-\infty}^\infty \phi_\mathrm{NLS}^2 \xi \dx+ O(\alpha^4), \\
\LL_3(\eta^\star) & \!\!\! = & \!\!\! \frac{\alpha^3}{4}\left(-f(k_0)f(2k_0)-{\textstyle\frac{1}{2}}f(k_0)^2 +{\textstyle\frac{3}{2}} k_0^2\right) \int_{-\infty}^\infty \phi_\mathrm{NLS}^2\psi \dx \\
& & \qquad\mbox{}
+\frac{\alpha^3}{4}(-2f(k_0)-f(k_0)^2+k_0^2)\int_{-\infty}^\infty \phi_\mathrm{NLS}^2 \xi \dx + O(\alpha^4)
\end{eqnarray*}
and
\begin{eqnarray*}
\KK_4(\eta^\star) & \!\!\! = & \!\!\! -\frac{\alpha^3}{64}\big(3\beta k_0^4 +\omega^2(f(2k_0)+2)\big) \int_{-\infty}^\infty \phi_\mathrm{NLS}^4 \dx + O(\alpha^4), \\
\GG_4(\eta^\star) & \!\!\! = & \!\!\! \frac{\alpha^3}{16}\left(k_0^2 - {\textstyle\frac{1}{2}}f(k_0)(f(2k_0)+2)\right)\omega \int_{-\infty}^\infty \phi_\mathrm{NLS}^4 \dx + O(\alpha^4), \\
\LL_4(\eta^\star) & \!\!\! = & \!\!\! \frac{\alpha^3}{16}(f(k_0)^2(f(2k_0)+2)-3k_0^2f(k_0))\int_{-\infty}^\infty \phi_\mathrm{NLS}^4 \dx + O(\alpha^4)
\end{eqnarray*}
(see Proposition \ref{Formulae for the threes and fours}
for the formulae for $\KK_3$, $\GG_3$, $\LL_3$ and $\KK_4$, $\GG_4$, $\LL_4$). Finally, observe that
\begin{eqnarray*}
\eta^{\star \prime\prime}(x)+k_0^2\eta^\star(x)
& = & \alpha^3 \phi_\mathrm{NLS}^{\prime\prime}(\alpha x) \cos k_0 x -2 \alpha^2k_0\phi_\mathrm{NLS}^\prime(\alpha x)\sin k_0 x
+\alpha^4 \psi^{\prime\prime}(\alpha x) \cos 2 k_0 x\\
& & \qquad\mbox{} - 4 \alpha^3 k_0 \psi^\prime(\alpha x) \sin 2k_0 x
-3 k_0^2 \alpha^2 \psi(\alpha x) \cos 2k_0 x + \alpha^4 \xi^{\prime\prime}(\alpha x),
\end{eqnarray*}
so that $\|\eta^{\star\prime\prime}+k_0^2\eta^\star\|_0 = O(\alpha^\frac{3}{2})$,
and using the further estimates
$\|\eta^\star\|_2 = O(\alpha^\frac{1}{2})$ and\linebreak
$\|\eta^\star\|_{1,\infty} = O(\alpha)$, one finds
from Proposition \ref{Size of the functionals} that
$\KK_\mathrm{r}(\eta^\star)$, $\GG_\mathrm{r}(\eta^\star)$, $\LL_\mathrm{r}(\eta^\star)$ are all $O(\alpha^\frac{7}{2})$.

The above calculations show that
\begin{eqnarray*}
\lefteqn{\KK(\eta^\star) + 2 \nu_0 \GG(\eta^\star) - \nu_0^2 \LL(\eta^\star)} \\
& = & \frac{\alpha^3}{8}(2\beta-\nu_0^2f^{\prime\prime}(k_0))\int_{-\infty}^\infty \phi_\mathrm{NLS}^{\prime 2} \dx
+\frac{\alpha^3}{4}\int_{-\infty}^\infty \big(g(2k_0)\psi^2 + A_3^1 \phi_\mathrm{NLS}^2 \psi \big)\dx \\
& & \qquad\mbox{}
+\frac{\alpha^3}{2}\int_{-\infty}^\infty \big(g(0)\xi^2 + A_3^2 \phi_\mathrm{NLS}^2 \xi\big) \dx+ \frac{3\alpha^3}{8}A_4 \int_{-\infty}^\infty \phi_\mathrm{NLS}^4\dx + O(\alpha^\frac{7}{2}) \\
& = & \frac{\alpha^3}{8}(2\beta-\nu_0^2f^{\prime\prime}(k_0))\int_{-\infty}^\infty \phi_\mathrm{NLS}^{\prime 2} \dx
+\frac{\alpha^3}{4}g(2k_0)\int_{-\infty}^\infty \left( \psi + \frac{g(2k_0)^{-1}}{2}A_3^1 \phi_\mathrm{NLS}^2 \right)^{\!\!2}\dx \\
& & \qquad\mbox{}+\frac{\alpha^3}{4}g(0)\int_{-\infty}^\infty \left( \xi + \frac{g(0)^{-1}}{2}A_3^2 \phi_\mathrm{NLS}^2 \right)^{\!\!2} \dx\\
& & \qquad\mbox{}+\alpha^3\left(\frac{3}{8}A_4 - \frac{g(2k_0)^{-1}}{16}(A_3^1)^2 - \frac{g(0)^{-1}}{8}(A_3^2)^2\right)
\int_{-\infty}^\infty \phi_\mathrm{NLS}^4 \dx+ O(\alpha^\frac{7}{2}),
\end{eqnarray*}
in which the second line follows from the first by the definitions of $A_3^1$, $A_3^2$, $A_4$ and the
third from the second by completing the square. The choice
$$\psi = - \frac{g(2k_0)^{-1}}{2}A_3^1\phi_\mathrm{NLS}^2, \qquad \xi = - \frac{g(0)^{-1}}{2}A_3^1\phi_\mathrm{NLS}^2$$
therefore minimises the value of $\KK(\eta^\star) + 2 \nu_0 \GG(\eta^\star) - \nu_0^2 \LL(\eta^\star)$
up to $O(\alpha^\frac{7}{2})$, whereby
\begin{eqnarray*}
\KK(\eta^\star) + 2 \nu_0 \GG(\eta^\star) - \nu_0^2 \LL(\eta^\star)
& = & \alpha^3\EE_\mathrm{NLS}(\phi_\mathrm{NLS})+O(\alpha^\frac{7}{2}) \\
& = & c_\mathrm{NLS}\alpha^3 + O(\alpha^\frac{7}{2}).
\end{eqnarray*}

The mapping
\begin{eqnarray*}
\alpha & \mapsto & \nu_0 \LL(\eta^\star) - \GG(\eta^\star) \\
& & = \alpha\left(\frac{\nu_0}{4}f(k_0)+\frac{\omega}{8}\right) \int_{-\infty}^\infty \phi_\mathrm{NLS}^2\dx + O(\alpha^2)
\end{eqnarray*}
is continuous and strictly increasing and therefore has a continuous inverse $\mu \mapsto \alpha(\mu)$; furthermore
$\alpha(\mu) = \mu + o(\mu)$ and
$$\JJ_\mu(\eta^\star) - 2\nu_0 \mu\ =\ \KK(\eta^\star) + 2 \nu_0 \GG(\eta^\star) - \nu_0^2 \LL(\eta^\star)\ =\ c_\mathrm{NLS}\mu^3 + o(\mu^3).\eqno{\Box}$$

\section*{Appendix B: The sign of \boldmath$A_3+2A_4$\unboldmath}

The quantities $\beta$, $\omega$, $k_0$ and $\nu_0$ are related by the fact that $g(k) \geq 0$ with equality
precisely when $k=\pm k_0$. It follows from the simultaneous equations $g(k_0)=0$, $g^\prime(k_0)=0$ that
$$\beta = \frac{\nu_0^2f^\prime(k_0)}{2k_0}, \qquad \omega = \frac{1+\beta k_0^2 - \nu_0^2f(k_0)}{\nu_0},$$
and inserting these expressions for $\beta$ and $\omega$ into the formulae for
$A_3$ and $A_4$ (Corollary \ref{lot in fours} and Proposition \ref{lot in threes}), one finds that
\begin{equation}
\nu_0^6 (A_3 + 2A_4) = a_8 \nu_0^8 + a_6 \nu_0^ 6 + a_4 \nu_0^4 + a_2 \nu_0^2 + a_0, \label{Speed poly}
\end{equation}
in which
\begin{eqnarray*}
a_0 & = & \textstyle- \frac{1}{12}h_2(k_0)^{-1}(1+2h_1(k_0)), \\[2mm]
a_2 & = & \textstyle-\frac{1}{3}h_2(k_0)^{-1}
\Big(\frac{1}{2}f(2k_0)+\frac{1}{2}k_0 f^\prime(k_0)+ 2h_1(k_0)\left(\frac{1}{2}+\frac{1}{2}k_0f^\prime(k_0)\right)\Big) , \\[2mm]
a_4 & = & \textstyle-\frac{1}{3}h_2(k_0)^{-1}\left(\left(\frac{1}{2}f(2k_0)+\frac{1}{2}k_0 f^\prime(k_0)\right)^2 + 2h_1(k_0)\left(\frac{1}{2}+\frac{1}{2}k_0f^\prime(k_0)\right)^2\right) \\
& & \textstyle \quad\mbox{}
- 2\left(\frac{1}{12}+\frac{1}{24}f(2k_0)\right), \\[2mm]
a_6 & = & \textstyle - \frac{2}{3}h_2(k_0)^{-1}\!\!
\left(\frac{1}{2}f(k_0)f(2k_0)-\frac{3}{2}k_0^2+\frac{1}{4}k_0f^\prime(k_0)f(2k_0)\!+\!\frac{1}{8}f^\prime(k_0)^2\right)\!\!\left(\frac{1}{2}f(2k_0)+\frac{1}{2}k_0f^\prime(k_0)\right) \\
& & \textstyle\quad\mbox{}-\frac{4}{3}h_2(k_0)^{-1}h_1(k_0)\left(\frac{1}{4}k_0f^\prime(k_0)+\frac{1}{2}f(k_0)-\frac{1}{2}k_0^2 + \frac{1}{8}k_0^2 f^\prime(k_0)^2\right)\!\!\left(\frac{1}{2}+\frac{1}{2}k_0f^\prime(k_0)\right) \\
& & \textstyle\quad\mbox{}+2\left(-\frac{1}{24}k_0f^\prime(k_0)f(2k_0)+\frac{1}{3}k_0^2
-\frac{1}{12}k_0f^\prime(k_0)-\frac{1}{6}f(k_0)-\frac{1}{12}f(k_0)f(2k_0)\right), \\[1mm]
a_8 & = & \textstyle- \frac{1}{3}h_2(k_0)^{-1}
\left(\frac{1}{2}f(k_0)f(2k_0)-\frac{3}{2}k_0^2+\frac{1}{4}f^\prime(k_0)f(2k_0)+\frac{1}{8}f^\prime(k_0)^2\right)^2 \\
& & \textstyle \quad\mbox{}-\frac{2}{3}h_2(k_0)^{-1}h_1(k_0)
\left(\frac{1}{4}k_0f^\prime(k_0)+\frac{1}{2}f(k_0)-\frac{1}{2}k_0^2+\frac{1}{8}k_0^2f^\prime(k_0)^2\right)^2 \\
& & \textstyle\quad\mbox{}-2\bigg(\frac{1}{16}k_0^3f^\prime(k_0)+\frac{1}{6}f(k_0)^2(f(k_0)+2)
- \frac{1}{2}k_0^2f(k_0) \\
& & \qquad\qquad\mbox{}\textstyle-2\left(\frac{1}{2}k_0f^\prime(k_0)-f(k_0)\right)\!\!\left(\frac{1}{6}k_0^2-\frac{1}{12}f(k_0)(f(2k_0)+2)\right) \\
& & \qquad\qquad\mbox{}\textstyle
+\frac{1}{24}\left(\frac{1}{2}k_0f^\prime(k_0)-f(k_0)\right)^2(f(2k_0)+2)\bigg)
\end{eqnarray*}
and
$$h_1(k_0)=\frac{-2f(2k_0)+2f(k_0)+3k_0f^\prime(k_0)}{-2-k_0f^\prime(k_0)+2f(k_0)},
\qquad
h_2(k_0)=\frac{3}{2}k_0f^\prime(k_0)+f(k_0)-f(2k_0).$$
The right-hand side of \eqn{Speed poly} defines a polynomial function of $\nu_0$ with coefficients
which depend upon $k_0$, and the following argument shows that it is negative for all positive values of
$\nu_0$.

First note that $a_0$, $a_2$ and $a_4$ are negative because
$$h_1(k_0)=g(0)^{-1}g(2k_0)^{-1}>0, \qquad h_2(k_0)=\frac{g(2k_0)}{\nu_0^2} > 0.$$
A lengthy calculation shows that
$$a_8 = - \frac{k_0^3}{\sinh^6 k_0}
\left( \sum_{j=0}^\infty \frac{a_{8,2j+1}}{(2j+1)!}k_0^{2j+1}\right)^{\!\!-1}\sum_{j=0}^\infty \frac{a_{8,2j}}{(2j)!}k_0^{2j},$$
in which explicit formulae for the coefficients $a_{8,j}$ are computed from the above expression for $a_8$.
Elementary estimates are used to establish that $a_{8,j}>0$,
so that $a_8$ is also negative. The argument is completed by
demonstrating that $4a_4a_8-a_6^2$ is positive. For this purpose we use the calculation
$$4a_4a_8-a_6^2 = \frac{k_0^4}{\sinh^8 k_0}\left(\sum_{j=0}^\infty \frac{b_j}{(2j)!}k_0^{2j}\right)^{\!\!-1}\sum_{j=0}^\infty \frac{c_j}{(2j)!}k_0^{2j}$$
with explicit formulae for the coefficients $b_j$ and $c_j$, which are also found to be positive.

\noindent\\
{\bf Acknowledgement.} E. Wahl\'{e}n was supported by an Alexander von Humboldt Research Fellowship, the Royal Physiographic Society in Lund, 
and the Swedish Research Council (grant no. 621-2012-3753). We would like to thank Boris Buffoni (EPFL Lausanne),
Per-Anders Ivert (Lund) and David Lannes (\'{E}cole Normale Sup\'{e}rieure) for many helpful discussions during the preparation of this article.

\bibliographystyle{standard}
\bibliography{mdg}
\end{document}